\documentclass[12pt,reqno]{amsart}
\usepackage{amssymb}

\input epsf.tex
\newdimen\xsize
\newdimen\oldbaselineskip
\newdimen\oldlineskiplimit
\def\nolineskip{\oldbaselineskip=\baselineskip\baselineskip=0pt%
\oldlineskiplimit=\lineskiplimit\lineskiplimit=0pt}
\def\restorelineskip{\baselineskip=\oldbaselineskip%
\lineskiplimit=\oldlineskiplimit}
\def\putm[#1][#2]#3{
\hbox{\vbox to 0pt{\parindent=0pt%
\vskip#2\xsize\hbox to0pt{\hskip#1\xsize $#3$\hss}\vss}}}%
\def\putt[#1][#2]#3{
\vbox to 0pt{\noindent\hskip#1\xsize\lower#2\xsize%
\vtop{\restorelineskip#3}\vss}}

\DeclareFontFamily{U}{rsf}{\skewchar\font'177}%
\DeclareFontShape{U}{rsf}{m}{n}{<-6>rsfs5<6-8>rsfs7<8->rsfs10}{}%
\DeclareFontShape{U}{rsf}{b}{n}{<-6>rsfs5<6-8>rsfs7<8->rsfs10}{}%
\DeclareMathAlphabet\RSFS{U}{rsf}{m}{n}
\SetMathAlphabet\RSFS{bold}{U}{rsf}{b}{n}
\let\scr=\rfs

\def\mib#1{\boldsymbol{#1}}

\def\sf#1{{\mathsf{#1}}}

\def\msmall#1{\mathchoice{\hbox{\small$\displaystyle {#1}$}}{#1}{#1}{#1}}        

\let\3=\ss
\def\cc{{\mathbb C}}

\def\rr{{\mathbb R}}
\def\nn{{\mathbb N}}
\def\pp{{\mathbb P}}
  \def\cp{\cc\pp}
\def\zz{{\mathbb Z}}

\def\ttt{{\mathbb T}}
\def\bfg{{\mathbf G}}

\def\aut{{\sf{Aut}}}
\def\loc{{\sf{loc}}}
\def\st{_{\mathsf{st}}}
\def\area{\sf{area}}

\def\cosh{\sf{cosh}}
\def\cos{\sf{cos}}
\def\coker{\sf{Coker}\,}

\def\deg{\sf{deg}\,}
\def\diam{\sf{diam}\,}
\def\dim{\sf{dim}\,}
\def\dimc{\dim_\cc}
\def\dimr{\dim_\rr}
\def\dfrm{\sf{dfrm}}
\def\dist{\sf{dist}\,}
\def\endo{\sf{End}}

\def\ev{\sf{ev}}

\def\ind{\sf{ind}}

\def\sind{\sf{S\hbox{-}ind}\,}

\def\gcd{\sf{gcd}}
\def\sfh{\sf{H}}
\def\hom{\sf{Hom}\,}
\def\homr{\sf{Hom}\vph_\rr}
\def\id{\sf{Id}}
\def\im{\sf{Im}\,}
\def\bint{{-}\mathchoice{\mkern-18mu}{\mkern-15mu}{\mkern-15mu}{\mkern-15mu}\int}
\def\sfl{\sf{L}}
\def\re{\sf{Re}\,}
\def\ker{\sf{Ker}\,}
\def\lim{\mathop{\sf{lim}}}

\def\log{\sf{log}\,}
\def\max{\sf{max}}
\def\min{\sf{min}}
\def\nod{{\sf{nod}}}
\def\ord{\sf{ord}}
\def\osc{\sf{osc}}
\def\op{\sf{op}}
\def\pr{\sf{pr}}

\def\rank{\sf{rank}}
\def\res{\sf{Res}}
\def\sinh{\sf{sinh}}
\def\sin{\sf{sin}}
\def\supp{\sf{supp}\,}
\def\sup{\sf{sup}\,}

\def\bfone{\boldsymbol{1}}
\def\mbfe{{\mib{e}}}
\def\mbfk{{\mib{k}}}
\def\mbfl{{\mib{l}}}
\def\mbfp{{\mib{p}}}

\def\mbfw{{\mib{w}}}
\def\mbfx{{\mib{x}}}
\def\mbfz{{\mib{z}}}

\def\vpp{\,\,\vec{\vphantom{p}}\!\!p'{}}
\def\vdp{\,\,\vec{\vphantom{d}}\!\!d'{}}

\def\epsi{\varepsilon}
\let\vkappa=\varkappa
\def\<{\langle}\let\la=\<
\def\>{\rangle}\let\ra=\>
 \let\bs=\bss
\def\comp{\Subset}
\def\d{\partial}
\def\dbar{{\barr\partial}}

\def\ddef{\mathrel{{=}\raise0.3pt\hbox{:}}}
\def\deff{\mathrel{\raise0.3pt\hbox{\rm:}{=}}}
\def\inv{^{-1}}
\def\hook{\hookrightarrow}
\def\fraction#1/#2{\mathchoice{{\msmall{ #1\over#2}}}%
{{ #1\over #2 }}{{#1/#2}}{{#1/#2}}}
\def\norm#1{\Vert #1 \Vert}
\def\half{{\fraction1/2}}
\def\le{\leqslant}
\def\vph{^{\mathstrut}}
\def\lrar{\longrightarrow}
\def\emptyset{\varnothing}
\def\scirc{\mathop{\mathchoice{\hbox{\small$\circ$}}{\hbox{\small$\circ$}}%
{{\scriptscriptstyle\circ}}{{\scriptscriptstyle\circ}}}}
\def\longpoints{\leaders\hbox to 0.5em{\hss.\hss}\hfill \hskip0pt}
\def\stateskip{\smallskip}
\def\state#1. {\stateskip\noindent{\bf#1. }} 
\def\statep#1. {\stateskip\noindent{\bf#1 }} 
\def\step#1{{\sl Step #1)}}
\def\proof{\state Proof. \.}
\def\Chi{\raise 2pt\hbox{$\chi$}}
\def\eg{\hskip1pt plus1pt{\sl{e.g.\/\ \hskip1pt plus1pt}}}
\def\ie{\hskip1pt plus1pt{\sl i.e.\/\ \hskip1pt plus1pt}}
\def\iff{if and only if }
\def\wrt{with respect to }
\def\isl{{\mathrm i}}
\def\sli{{\sl i)} } \def\slip{{\sl i)}}
\def\slii{{\sl i$\!$i)} } \def\sliip{{\sl i$\!$i)}}
\def\sliii{{\sl i$\!$i$\!$i)} }\def\sliiip{{\sl i$\!$i$\!$i)}}
\def\sliv{{\sl i$\!$v)} } 
\def\slv{{\sl v)} } 

\def\reg{^\sf{reg}}
\def\sing{^\sf{sing}}

\def\diff{{\scr D}\mskip -2mu i\mskip -5mu f\mskip-6.3mu f}

\def\barr#1{\mskip1mu\overline{\mskip-1mu{#1}\mskip-1mu}\mskip1mu}
\def\Chi{\raise 2pt\hbox{$\chi$}}
\def\yps{\Upsilon}
\let\phI=\phi\let\phi=\varphi\let\varphi=\phI
\def\bfnu{{\boldsymbol\nu}}%
\def\bfxi{{\boldsymbol\xi}}%
\def\bftau{{\boldsymbol\tau}}%
\def\scra{\scr{A}}
\def\scrb{\scr{B}}

\def\scrh{\scr{H}}
\def\scre{\scr{E}}

\def\scrj{\scr{J}}

\def\scrm{\scr{M}}
  \def\whcalm{\wh{\scr M}}
  \def\barm{{\mskip5mu\barr{\mskip-5mu\scr M}}{}}
  
\def\scrn{{\scr N}}
\def\scro{{\scr O}}
\def\scrp{{\scr P}}
\def\scrpp{{\scr{P\!\!P}}}

\def\scrs{{\scr S}}

\def\scru{{\scr U}}
\def\scrv{{\scr V}}
\def\scrw{{\scr W}}
\def\scrx{{\scr X}}
\def\scry{{\scr Y}}
\def\scrz{{\scr Z}}

\def\trans{\pitchfork}

\def\epsi{\varepsilon}
\def\bs{\backslash}
\def\ogran{{\hskip0.7pt\vrule height6pt depth2.5pt\hskip0.7pt}}
\def\comp{\Subset}
\def\d{\partial}
\def\dbar{{\barr\partial}}
\def\1{{1\mkern-5mu{\rom l}}}

\def\ge{\geqslant}
\def\inv{^{-1}}
\let\wh=\widehat
\let\wt=\widetilde
\def\fraction#1/#2{\mathchoice{{\msmall{ #1\over#2}}}%
{{ #1\over #2 }}{{#1/#2}}{{#1/#2}}}
\def\half{{\fraction1/2}}
\def\le{\leqslant}
\def\vph{^{\mathstrut}}
\def\emptyset{\varnothing}
\def\.{\thinspace}

\def\ti#1{{\tilde{#1}}}
\def\term#1#2{\underline{#2}_{[{#1}]}}

\def\qed{\ \ \hfill\hbox to .1pt{}\hfill\hbox to .1pt{}\hfill $\square$\par}

\def\comment#1\endcomment{}

 \belowdisplayskip=\abovedisplayskip

\def\lineeqqno(#1){\hfill\llap{\vbox to 10pt%
{\vss\begin{align} \eqqno(#1)\end{align}\vss}}\vskip1pt}

\advance\headheight 1.2pt

\numberwithin{equation}{subsection}

\def\newsection[#1]#2{\section{#2}\label{sec:#1}}
   \def\refsection#1{{\sl Section \ref{sec:#1}$\,$}}
\def\newsubsection[#1]#2{\subsection{#2}\label{sec:#1}\showlabel{\tt #1}}
   \def\refsubsection#1{{\sl Paragraph \ref{sec:#1}$\,$}}

\def\inDex#1{\index{#1}}

\newtheorem{thm}{Theorem}[subsection]
   \def\newthm#1{\begin{thm} \showlabel{\tt #1}\label{#1}%
      \index{{\sl Theorem \ref{#1}},  \ label {\tt #1}|Zpage}} 
   \def\refthm#1{{\sl Theorem \ref{#1}$\,$}}
\newtheorem{lem}[thm]{Lemma}
   \def\newlemma#1{\begin{lem} \showlabel{\tt #1}\label{#1}%
      \inDex{{\sl Lemma \ref{#1}},  \ label {\tt #1}|Zpage}} 
   \def\lemma#1{{\sl Lemma \ref{#1}$\,$}}
\newtheorem{prop}[thm]{Proposition}
   \def\newprop#1{\begin{prop} \showlabel{\tt #1}\label{#1}%
      \inDex{{\sl Proposition \ref{#1}},  \ label {\tt #1}|Zpage}}
   \def\propo#1{{\sl Proposition \ref{#1}$\,$}}
\newtheorem{corol}[thm]{Corollary}
   \def\newcorol#1{\begin{corol} \showlabel{\tt #1}\label{#1}%
      \inDex{{\sl Corollary \ref{#1}},  \ label {\tt #1}|Zpage}}
   \def\refcorol#1{{\sl Corollary \ref{#1}$\,$}}
\newtheorem{defi}{Definition}[subsection]
   \def\newdefi#1{\begin{defi} \showlabel{\tt #1}\label{#1}\rm %
      \inDex{{\sl Definition \ref{#1}},  \ label {\tt #1}|Zpage}}
   \def\refdefi#1{{\sl Definition \ref{#1}$\,$}}

\def\eqqno(#1){\label{eq#1}%
      \index{{\sl  Eqtn \ref{eq#1}},  \ label {\tt(#1)}|Zpage}}
\def\eqqref(#1){\eqref{eq#1}}
\def\showlabel#1{}

\makeindex{}

\begin{document}
\baselineskip=14.0pt plus 2pt

\title[Pseudoholomorphic curves and the symplectic isotopy problem]%
{Pseudoholomorphic curves 
\\[5pt]
and 
the symplectic isotopy problem%
}
\author[V.~Shevchishin]{Vsevolod V.~Shevchishin}
\address{Fakult\"at f\"ur Mathematik\\
Ruhr-Universit\"at Bochum\\ 
Universit\"atsstrasse 150\\
44780 Bochum\\
Germany}
\email{sewa@@cplx.ruhr-uni-bochum.de}
\dedicatory{}
\begin{abstract}
The deformation problem for pseudoholomorphic curves and related
geometrical properties of the total moduli space of pseudoholomorphic curves
are studied. A sufficient condition for the saddle point property
of the total moduli space is established. The local symplectic isotopy problem
is formulated and solved for the case of imbedded pseudoholomorphic curves. 
It is shown that any two symplectically imbedded surfaces $\Sigma_0, \Sigma_1
\subset \cp^2$ of the same degree $d\le 6$ are symplectically isotopic.
\end{abstract}
\maketitle
\setcounter{tocdepth}{2}

\setcounter{section}{-1}

\newsection[intro]{Introduction}

\baselineskip =14.0pt plus .5pt
The symplectic isotopy problem can be formulated as follows:

{\it For a given  symplectic 4-dimensional manifold $(X,\omega)$ and symplectically
imbedded compact connected oriented surfaces $\Sigma_0, \Sigma_1 \subset X$ 
in the same homology class $[A] \in \sfh_2(X,\zz)$, does there exist an isotopy 
$\{\Sigma_t\}_{t\in  [0,1]}$ connecting $\Sigma_0$ with $\Sigma_1$ such that all\/ 
$\Sigma_t$ are also symplectically imbedded?} 

In this case $\Sigma_0$ and $\Sigma_1$ are called {\sl symplectically isotopic}.
Note that by the {\sl genus formula} for pseudoholomorphic curves (see 
\refsection{1}) $\Sigma_0$ and $\Sigma_1$ have the same genus. 
The example of Fintushel and Stern \cite{Fi-St} (see \refsection{6} for details)
shows that in general the answer is negative. On the other hand, Sikorav 
\cite{Sk-3} gives an affirmative answer in the case of surfaces of positive degree
$d\le 3$ in $\cp^2$. So it is natural to ask under which conditions on $(X,\omega)$
and $[A] \in \sfh_2(X,\zz)$ the symplectic isotopy does exist.

In this paper techniques involving pseudoholomorphic curves are developed in
directions needed for solving the symplectic isotopy problem. The main result 
is the following 

\state Theorem 1. {\it Two symplectically imbedded compact connected oriented 
surfaces $\Sigma_0, \Sigma_1 \subset \cp^2$ of the same positive degree $d\le 6$ 
are symplectically isotopic}.

\smallskip
The proof is based on the solution of two problems which are closely related to 
the symplectic isotopy problem and concern the geometry of the total moduli space
of pseudoholomorphic curves. The first result is a sufficient condition for the 
{\it saddle point property} of the total moduli space of pseudoholomorphic 
curves. This removes one of the obstacles to the existence of a symplectic isotopy 
and is applied to solve the {\it local symplectic isotopy problem} for imbedded 
pseudoholomorphic curves. The latter appears as a necessary part of the global 
problem.

\smallskip
The obtained progress gives hope that the symplectic isotopy problem has
an affirmative solution in the case $c_1(X,\omega)[A] >0$. Note that compact 
symplectic 4-manifolds with this property are classified: Up to the case when 
$[A]$ is represented by an exceptional sphere, these are symplectic blow-ups of 
a rational or ruled symplectic manifold, see \cite{McD-Sa-3}, {\sl Corollary 1.5},
and also \refsection{6}.

\newsubsection[0.1]{Overview of main results} There is essentially only one known
method of constructing symplectic isotopies. Having its origin in Gromov's 
celebrated article \cite{Gro}, it utilizes moduli spaces of pseudoholomorphic 
curves. One fixes a homotopy $h(t) \deff J_t$, $t\in [0,1]$, of 
$\omega$-tame almost complex structures and considers the relative moduli space,
$$
\scrm_h \deff \bigr\{ (C,t) : 
C \text{ is an imbedded $J_t$-holomorphic 
curve in the homology class }[A] 
\bigl\}, 
$$
equipped with a natural topology and with the projection $\pi_h: (C,t) \in \scrm_h 
\mapsto t \in [0,1]$. It follows essentially from \cite{Gro} that for a generic 
homotopy $h(t) = J_t$ the space $\scrm_h$ is a smooth manifold of expected 
dimension $\dimr\scrm_h = 2c_1(X,\omega)[A] + 2g-1$ and the projection $\pi_h$ is 
also smooth. Moreover, if this dimension is positive, then there exists a generic 
path $h(t)=J_t$ such that the original surface $\Sigma_0$ (resp.\ $\Sigma_1$) is 
$J_0$-holomorphic (resp.\ $J_1$-holomorphic) so that $(\Sigma_i, i) \in \scrm_h$ 
for $i=0,1$. Using a trivial but crucial observation that for every $(C,t) \in 
\scrm_h$ the curve $C$ is an $\omega$-symplectic real surface, one tries to find 
the desired isotopy $\Sigma_t$ by constructing a continuous section $\sigma: t 
\in [0,1] \mapsto (\Sigma_t, t) \in \scrm_h$ connecting $(\Sigma_0,0)$ with 
$(\Sigma_1,1)$.

One can easily see possible obstacles to the existence of such a section $\sigma:
 [0,1]\to \scrm_h$. The first one is that the projection $\pi_h: \scrm_h \to [0,1]$,
considered as a function, can have local maxima and minima. Indeed, if some 
$(C^*, t^*)\in \scrm_h$ appears as a local maximum of $\pi_h$, then for all $t>
t^*$ there exists no $J_t$-holomorphic curves sufficiently close to $C^*$. Observe 
that the mere fact of existence of a $J_t$-holomorphic curve $C_t$ 
for $t>t^*$ does not help much, because this does not imply that such a curve $C_t$ 
is symplectically isotopic to $\Sigma_0$ or $\Sigma_1$. Note that exactly the 
existence of $J$-holomorphic curves is the main technical tool in the Gromov's 
article \cite{Gro}. 

On the other hand, this obstacle does not appear in the case when the projection 
$\pi_h$ has the following {\sl saddle point property}: the Hesse matrix 
of $\pi_h$ at any critical point has at least one positive and at least one negative
eigenvalues. In \refsection{3} we prove 

\state Theorem 2. {\it Assume that $c_1(X,\omega)[A]>0$. Then for a generic 
homotopy $h(t)$ all critical points of the projection $\pi_h: \scrm_h \to 
 [0,1]$ are saddle}. 

\smallskip
On the other hand, we also show that in the case $c_1(X,\omega)[A]\le 0$ local 
maxima and minima of the projection $\pi_h: \scrm_h \to [0,1]$ do exist in the case
of an apropriate generic choice of $h$. Note that in the example of Fintushel 
and Stern \cite{Fi-St} one has $c_1(X,\omega)[A]= 0$.

\medskip
The next obstacle to existence of the desired section $\sigma: [0,1]\to \scrm_h$
comes from the fact that in general $\scrm_h$ is not compact and the projection
$\pi_h$ is not proper. This means that while attempting to construct the section 
$\sigma: [0,1]\to \scrm_h$ we go to ``infinity'' in the space $\scrm_h$. The 
Gromov compactness theorem provides control on the limiting behavior of curves 
$C_t$ in this case. It says that some subsequence, say $C_{t_n}$, converges in 
a certain weak sense to a pseudoholomorphic curve $C^*$ in the same homology class
$[A]$. The curve $C^*$ can have several irreducible components, some of them 
multiple, and also several singular points. Thus we come to a problem of 
describing  symplectic isotopy classes of imbedded pseudoholomorphic
curves close to a given singular curve $C^*$.

Here we have essentially two different difficulties. The first one comes from
multiple components. At the moment, we have no remedy for this problem. So we 
appempt to avoid the appearence of multiple components. This is done
by imposing the following additional constraint. We consider the curves $C_t$ 
which pass through fixed points $\mbfx=(x_1,\ldots,x_k)$ on $X$. This means, 
however, that now we must consider a new moduli space $\scrm_{h,\mbfx}$ with a 
new projection $\pi_{h,\mbfx}: \scrm_{h,\mbfx} \to [0,1]$. In \refsection{3}
we also show that if the number $k$ of the fixed points is strictly less than
$c_1(X,\omega)[A]$, then the saddle point property for $\pi_{h,\mbfx}$ remains
valid. An easy calculation shows that fixing $k=3d-1$ generic points on $\cp^2$, 
we can avoid the appearance of multiple components in pseudoholomorphic curves 
of degree $d\le6$. This explains the restriction in the main theorem.

\smallskip
In the case when $C^*$ has no multiple components we still have to describe the
possible symplectic isotopy classes of imbedded pseudoholomorphic curves close to 
$C^*$. Recall that by the result of Micallef and White (see \refsubsection{1.2}) 
every singular point of $C^*$ is isolated and topologically equivalent to a 
singular point of a usual holomorphic curve. In this way we come to the {\sl local
symplectic isotopy problem} which asks about possible symplectic isotopy types of 
imbedded pseudoholomorphic curves in a neighborhood of a given isolated 
singularity, see \refsubsection{6.2} for details. 

The solution of the local symplectic isotopy problem is based on the simple
observation that for {\sl holomorphic curves} this problem has a trivial solution.
Namely, a generic holomorphic deformation of a given (local) holomorphic curve 
$C^*$ gives a non-singular curve $C$, and the set of such curves is open and 
connected. Using this fact, we expoit essentially the same method as in the case 
of the global problem and prove that it is possible to deform isotopically an 
imbedded pseudoholomorphic curve $C$ sufficiently close to a given singular curve 
$C^*$ into a genuine holomorphic curve. In this way we prove

\state Theorem 3. {\it There exists a unique symplectic isotopy class of 
non-singular pseudoholomorphic curves which are close to a given pseudoholomorphic 
curve $C^*$ without multiple components}.

\smallskip
It should be noted that in the proof the saddle point property from {\sl Theorem 
2\/} is used in an essential way. {\sl Theorem 1\/} follows now by 
the procedure of avoiding multiple components as it is explained above. 

\medskip
Further technical results of the paper are as follows. \refsection{4} is devoted 
to the deformation problem of pseudoholomorphic maps with prescribed singularities.
It is shown that the subspace of such maps is an immerced submanifold of expected 
codimension in the total moduli space of pseudoholomorphic maps. In \refsection{3}
the second variation of the $\dbar$-equation is computed. The result establishes 
the relationship between the geometry of a pseudoholomorphic curve $C^*$ 
corresponding to a critical point of the projection $\pi_h: \scrm_h \to [0,1]$
and the eigenvalues of the Hesse matix $d^2\pi_h$ at this point. Combined with
transversality results, this yields the proof of {\sl Theorem 2}. Finally,
in \refsection{5} the problem of smoothing of nodal points on pseudoholomorphic 
curves is studied.

\bigskip
{\it Acknowledgements.}
The author would like to express his gratitude
to  A.~Huckle\-berry, S.~Ivash\-ko\-vich,
St.~Ne\-mi\-rov\-ski, St.~Orev\-kov, B.~Sie\-bert and J.-C.~Siko\-rav for 
numerous useful conversations, suggestions and remarks.


\newsection[1]{Deformation and the normal sheaf of pseudoholomorphic curves}

In this section we give a brief description of pseudoholomorphic 
curves and the related deformation theory.

\newsubsection[1.1]{Pseudoholomorphic curves}
First we collect some facts from the Gromov's theory. 
Since there are several books devoted to or
treating this theme (see \eg \cite{McD-Sa-1} or \cite{McD-Sa-2})
we only mention the basic definitions and results we shall use later.

\newdefi{def1.1.0} 
An {\sl almost complex structure} on a manifold $X$ is an
endomorphism  $J \in TM$ of the tangent bundle such that $J^2 = - \id$.
The pair $(X,J)$ is called an {\sl almost complex manifold}.
One of the most important classes of such manifolds appears in symplectic
geometry. An almost complex structure on a symplectic manifold $(X,\omega)$ is
called {\sl$\omega$-tame} if $\omega(v,Jv)>0$ for any non-zero tangent vector
$v$. It is well-known that the set $\scrj_\omega$ of $\omega$-tame almost complex
structures is non-empty and contractible, (see \eg \cite{Gro},
\cite{McD-Sa-1}, or \cite{McD-Sa-2}). In
particular, any two $\omega$-tame almost complex structures $J_0$ and $J_1$ can be
connected by a homotopy (path) $J_t, t\in [0,1],$ inside $\scrj_\omega$.
\end{defi}

\newdefi{def1.1.1} A {\sl parameterized $J$-holomorphic curve} in an almost 
complex manifold $(X,J)$ is given by a (connected) Riemann surface $S$ with a
complex structure $J_S$ on $S$ and a non-constant $C^1$-map $u: S \to X$ 
satisfying the Cauchy-Riemann equation
\begin{equation}
 du + J \scirc du \scirc J_S =0.
\eqqno(1.1.1)
\end{equation}
In this case we call $u$ a $(J_S, J)$-holomorphic map, or simply
$J$-holomorphic map. Here we use the fact that if $u$ is not constant, then
such a structure $J_S$ is unique. We shall also use the notion {\sl $J$-curve}
which, depending on the context, will mean a map $u:S \to X$, \ie a
parameterized curve in $X$, or an image $u(S)$ of $J$-holomorphic map, taken
with appropriate multiplicity, \ie a {\sl non-parameterized $J$-curve}.
\end{defi}

The equation \eqqref(1.1.1) is elliptic with the Cauchy-Riemann symbol. This 
provides regularity properties for $u$. In particular, $u$ is H\"older 
$C^{l+1,\alpha}$-smooth, $u\in C^{l+1,\alpha}(S,X)$, provided $J\in C^{l,\alpha}$
with integer $l\ge 1$ and $0< \alpha <1$. To simplify the notations we set $\ell
\deff l+ \alpha$ and write $C^\ell$ to indicate $C^{l,\alpha}$-smoothness.
In what follows we shall assume that almost complex structures $J$ on $X$ are
$C^\ell$-smooth for some fixed sufficiently big non-integer $\ell$.

An easy consequence of the tameness condition is that any $J$-holomorphic
{\sl imbedding} $u: S \to X$ with $J \in \scrj_\omega$ is {\sl symplectic}
\ie the pull-back $u^*\omega$ is non-degenerate on $S$. The converse is also
true: Any $C^{\ell+1}$-smooth symplectic imbedding $u:S \to X$ with $\ell >1$
is $J$-holomorphic for some $C^\ell$-smooth $\omega$-tame structure
$J$. For {\sl immersions} the situation is more complicated. We state a result 
in a setting which will be relevant later on 
(see, \eg \cite{Gro} or \cite{McD-Sa-2} for details). 

\newlemma{lem1.1.1}
Let $(X,\omega)$ be a symplectic manifold with $\dim_\rr X =4$, and $u: S \to X$ 
an $\omega$-symplectic $C^1$-map such that $u(S)$ has only simple transversal
{\sl positive} self-intersection points.

Then there exist an $\omega$-tame almost complex structure $J$ on $X$ and a complex
structure $J_S$ on $S$ and making $u$ a $J$-holomorphic map.
\end{lem}

It is worth to make the following remark. If $x \in X$ is a self-intersection 
point of $u$, $x= u(z_1)= u(z_2)$ with $z_1 \not= z_2$, such that the tangent 
planes $du(T_{z_i}) \subset T_xX$ are transversal and {\sl complex} \wrt some 
structure $J_x$ in $T_xX$, then the intersection index of planes $du(T_{z_i})$ 
in $x$ is {\sl positive}. However, it is possible that two {\sl symplectic} planes
$L_i$ in $(\rr^4, \omega)$ have {\sl negative} intersection index.

\smallskip
More detailed considerations lead to 
the {\sl genus formula} (also called {\sl adjunction formula}) for immersed 
symplectic surfaces in symplectic four-folds. For this let $(X,\omega)$ be 
a symplectic 
manifold of dimension 4, $S\deff \bigsqcup_{j=1}^dS_j$ a compact oriented surface 
and $u:  S \to X$ an immersion with only transversal self-intersection points.
Denote by $g_j$ the genus of $S_j$, by $[C]$ the homology class of the image
$C \deff u(S)$, $[C]^2$ the homological self-intersection number of $[C]$, and
by $c_1(X)$ the first Chern class of $(X, \omega)$. Define the {\sl geometric
self-intersection number} $\delta$ of $M=u(S)$ as the algebraic number of pairs
$z' \not= z'' \in S$ with $u(z') = u(z'')$, taken with the sign corresponding
to the intersection index.

\newlemma{lem1.1.2}
Suppose that $u: S \to X$ is a symplectic immersion which is 
compatible with the orientation on each component $S_j$ of $S$. Then
\begin{equation}
\sum_{j=1}^d g_j = {[C]^2 - c_1(X)[C]\over2} + d - \delta .
\eqqno(1.1.3)
\end{equation}
\end{lem}

\smallskip
An elementary proof uses the fact that for a {\sl symplectic} immersion
$u:S \to (X,\omega)$ one has $c_1(X)[C]= \chi(S) + \chi(N)$, where
$N$ is the normal bundle and $\chi$ denotes the Euler characteristic. Finally,
one observes that $\chi(N)= [C]^2 -2\delta$. For details, see \cite{Iv-Sh-1}.

\newsubsection[1.2]{Local structure of pseudoholomorphic curves}
For the most results of this paragraph we refer to \cite{Mi-Wh} where a very 
precise description of the local structure of pseudoholomorphic  curves is 
given. As a rough summary, one can say that the 
local behavior of (non-parameterized) pseudoholomorphic  curves is essentially 
the same as for usual holomorphic curves.

\newlemma{lem1.2.1} {\rm (\cite{Mi-Wh})} Let $(X,J)$ be an almost complex
manifold of $\dimc X 
\allowbreak
=n$, $u: S \to X$  a~$J$-holomorphic map, and $x\in X$ a 
point. Suppose that $J\in C^2$ and that for any $x'\in X$ sufficiently close to
$x$ the pre-image $u\inv(x)$ is finite. Then there exist neighborhoods $U
\subset X$ of $x$, $U' \subset \cc^n$ of\/ $0\in \cc^n$ and a 
$C^1$-diffeomorphism $\phi: U \to U'$ such that $C' \deff \phi\bigl( u(S) 
\cap U\bigr)$ is a proper analytic curve in $U'$ and such that $\phi_*(J_x)=J\st$,
where $J\st$ is the standard complex structure in $\cc^n$.
\end{lem}

In particular, the notion of a (local) irreducible component of a
$J$-holomorphic curve $C=u(S)$ is well-defined. Further, in the special  case
when $(X,J)$ is an {\sl almost complex surface} one can correctly define

\noindent
\sli  {\sl the intersection index $\delta_{ij}(x)\in \nn$} of two local
components $C_i$ and $C_j$ at $x\in X$, and

\noindent
\slii {\sl the nodal number $\delta_i(x) \in \nn$} of a local component
$C_i$ at $x\in X$ (see \cite{Mil}, \S\.10 and \refdefi{def6.2.1}).

The main properties of these local invariants are summarized in

\newlemma{lem1.2.2}
\sli If $x\in C_i \cap C_j$, then $\delta_{ij}(x)
\ge 1$. The equality holds \iff $C_i$ and $C_j$ are smooth and
intersect transversally in $x$;

\slii The set $\{ z \in S\;:\; \delta_i(u(z)) >0 \}$ is discrete in $S$;

\sliii Suppose additionally that $S= \sqcup_{j=1}^d S_j$ is a closed surface
and $u: S \to X$ is an imbedding almost everywhere on $S$. Set $C\deff u(S)$.
Denote by $\delta$ the sum of all local intersection indices $\delta_{ij}(x)$ and
all local nodal numbers $\delta_i(x)$, the homology class of $C$ by $[C]$, and
the genera of particular components $S_j$ by $g_j$. Then
\begin{equation}
{\textstyle \sum_{j=1}^d} g_j = {[C]^2 - c_1(X)[C]\over2} + d -\delta.
\eqqno(1.2.1)
\end{equation}
\end{lem}

The formula \eqqref(1.2.1) is the {\sl genus formula for pseudoholomorphic
curves}. We shall also apply a local version of this result. Here we say that a 
pseudoholomorphic curve $C$ in a symplectic manifold $X$ is {\sl parameterized by 
a real surface $S$} if there exists a map $u: S \to X$ which is an imbedding 
outside a discrete subset in $S$. Such a surface $S$, possibly not connected, can
be constructed as the normalization of $C \subset X$.

\newlemma{lem1.2.2a} Let $B \subset \rr^4$ be the unit ball equipped with the 
standard symplectic structure $\omega\st$, and $C_1$, $C_2$ pseudoholomorphic 
curves in $B$.
Assume that the boundaries of the curves $\d C_i$ are imbedded in the boundary of
the ball $\d B$, are sufficiently close to each other, and that every $C_i$ meets
transversally $\d B$. Denote by $\Chi_i$ the Euler characteristic of the surface 
$S_i$ parameterizing $C_i$ and by $\delta_i$ the sum of the nodal number of 
singular points of $C_i$. Then
\begin{equation}\eqqno(1.2.1a)
 \Chi_1 - 2\delta_1 = \Chi_2 - 2\delta_2
\end{equation}
\end{lem}

\proof It is shown in \cite{Iv-Sh-1} that every $C_i$ can be perturbed 
to a nearby pseudoholomorphic curve $C'_i$ in such a way that every singular point 
$x \in C_i$
with the nodal number $\delta_x(C_i)\ge 2$ ``splits'' into $\delta_x(C_i)$ {\sl 
nodal points} on the perturbed curve $C'_i$, \ie the points where exactly two 
branches of $C'_i$ meet transversally. By this procedure the topology of each $S_i$ 
and the whole nodal number of every $C_i$ remain unchanged. After this, one can 
replace a sufficiently small neighborhood of every nodal point $x \in C'_i$ by a 
symplectically imbedded handle. This ``symplectic surgery of $C'_i$'' produce {\sl 
imbedded} pseudoholomorphic curves $C''_i$ with $\Chi(C''_i)= \Chi_i - 2\delta_i$.

Moreover, all this can be carried out with the boundaries $\d C_i$ unchanged. 
Further, the hypothesis of the lemma implies that the boundaries $\d C_i$ are
transversal to the standard symplectic structure on $\d B=S^3$ and are isotopic
as {\sl transversal links}, see \cite{Iv-Sh-1} and \cite{Eli}. Now one applies 
the theorem of Bennequin \cite{Bn}, see also \cite{Eli}, which claims that, up
to sign convention, $\Chi(C''_i)$ is the {\sl Bennequin index of $\d C_i$} and
depends only on the {\sl transversal isotopy} class of $\d C_i$. The lemma follows.
\qed

\medskip
The result of Micallef and White, \lemma{lem1.2.1}, is not sufficient for
our purpose, because it does not allow us to control local structure of
pseudoholomorphic  curves under deformation. A necessary tool is provided
by the following statement proven in \cite{Iv-Sh-1}. Here and thereafter 
$\Delta$ denotes the unit disc in $\cc$ equipped with the standard complex 
structure.

\newlemma{lem1.2.3} 
Suppose that a $f\in L^{1,2}_\loc(\Delta, \cc^n)$ is not identically $0$ and 
satisfies {\sl a.e.} 
the~inequality
\begin{equation}
| \dbar f(z)| \le h(z)\cdot |z|^k \cdot |f(z)| 
\eqqno(1.2.2)
\end{equation}
for some $k\in \nn$ and nonnegative $h\in L^p_\loc(\Delta)$ with $2<p<\infty$.
Then
\begin{equation}
f(z)=z^\mu\bigl(P^{(k)}(z) + z^k g(z)\bigr),
\eqqno(1.2.3)
\end{equation}
where $\mu\in\nn$, $P^{(k)}$ is a polynomial in $z$ of degree $\le k$ with
$P^{(k)}(0) \not=0$, and $g\in L^{1, p}_\loc(\Delta,\cc^n)\hookrightarrow C^{0,
\alpha}$, $\alpha=1-{2\over p}$, with $g(0)=0$.
\end{lem}

\medskip
Using this result one can obtain the following  description of the local
behavior of a pseudoholomorphic  map. Note that on a given almost complex
manifold $(X,J)$ in a neighborhood of a given point $x_0\in X$ there exist
an (integrable) complex structure $J^*$ with $J^*(x_0)= J(x_0)$ and
$J^*$-holomorphic coordinates $w_1, \ldots, w_n$, $n=\dimc X$.

\newlemma{lem1.2.4} \sli Assume that $J$ is $C^1$-smooth and $u: \Delta \to
X$ is a non-constant $J$-holomorphic map with $u(0)=x_0$. Then in coordinates
$w_1, \ldots, w_n$ chosen as above in a neighborhood of $x_0\in X$ the map $u$
has the form
\begin{equation}
u(z) = z^\mu\cdot P^{(\mu-1)}(z) + z^{2\mu-1} \cdot v(z),
\eqqno(1.2.4)
\end{equation}
where $\mu\in \nn$, $P^{(\mu-1)}(z)$ is a complex $\cc^n$-valued polynomial
of degree $\le \mu-1$ with $P^{(\mu-1)}(0)\not=0$, and $v(z) \in
L^{1,p}(\Delta, \cc^n)$ with $v(0)=0$. 

\slii  Assume that $J$ is $C^1$-smooth and let $u_1, u_2: \Delta \to X$
be $J$-holomorphic maps such that $u_1$ is an immersion and $u_2(0) \in 
u_1(\Delta)$. Then there exists $r\in ]0,1[$ such that either $u_2(\Delta(r))
\subset u_1(\Delta)$ or $u_2(\Delta(r)) \cap u_1(\Delta)= u_2(0)$.

\sliii Let $J$ be a $C^\ell$-smooth almost complex structure
on the ball $B \subset \cc^n$ with $J(0) = J\st(0)$, and let $u_1, u_2 : 
\Delta \to B$ be $J$-holomorphic maps with $u_1(0) =u_2(0) =0 \in B$,
such that $u_1 \not = u_2$.

Then there exists a uniquely defined $\nu \in \nn$ and $w(z) \in C^1(
\Delta, \cc^n)$ such that 
\begin{equation}
u_1(z) - u_2(z) = z^\nu w(z).
\eqqno(1.2.5)
\end{equation}
\end{lem}

\proof The first and second parts of the lemma are proven in \cite{Iv-Sh-1}. 
For the third see {\sl Remark 1.6} and {\sl Section 6} of \cite{Mi-Wh}.

\smallskip
\newdefi{def1.2.1} If for an appropriate local complex coordinate $z$ on $S$ 
a $J$-holomorphic map $u:S \to X$ has the form \eqqref(1.2.4), then we call 
the (uniquely defined) $\mu$ the {\sl multiplicity of $u$ at the point $z=0$}.
\end{defi}

\newdefi{def1.2.2} A $J$-holomorphic map $u: S \to X$ is {\sl multiple}
if there exists a non-empty $U \subset S$ such that the restriction
$u\ogran_U$ can be represented as a composition $u\ogran_U = u' \scirc \phi$
where $u': \Delta \to X$ is a $J$-holomorphic map and $\phi: U \to \Delta$ is
a (branched) covering of degree $m\ge 2$. In other words, $u$ is multiple
if some part of the image $u(S)$ is multiply covered by $u$. 
\end{defi}

\newsubsection[1.3]{Deformation of pseudoholomorphic maps and the Gromov operator
$D_{u, J}$}
Roughly speaking, the main idea of the Gromov's theory is to construct and
study $J$-holomorphic curves in a symplectic manifold $(X, \omega)$ for some
special (\eg {\sl integrable}) $J$. Often, one can show the existence of a
$J_0$-holomorphic curve with some other almost complex structure $J_0$, see \eg
{\sl Lemma \ref{lem1.1.1}}. If both $J$ and $J_0$ are $\omega$-tame, then 
there exists a homotopy $\{J_t\}_{t\in [0,1]}$ from $J_0$ to $J_1=J$. Hence 
one could try to deform the constructed $J_0$-holomorphic map $u_0:S \to X$ 
into a $J_1$-holomorphic one using the continuity principle. The first step in 
this direction is to study the linearization of (\ie the first variation) the 
equation \eqqref(1.1.1). This means that we are interested in the first
differential of the section $\sigma_\dbar$.

\smallskip
Fix a compact surface $S$ of genus $g$. Denote by $\scrj_S$ the Banach
manifold of $C^{1,\alpha}$-smooth complex structures on $S$ with some fixed
$\alpha \in\; ]0,1[$. Thus
\begin{equation}
\scrj_S = \{ J_S\in C^{1,\alpha}(S, \endo(TS)):
J_S^2=-\id \}
\end{equation}
and the tangent space to $\scrj_S$ at $J_S$ is
\begin{equation}
T_{J_S}\scrj_S = \{ I\in C^{1,\alpha}(S, \endo(TS)) :
J_S I + I J_S = 0\} \equiv C^{1,\alpha}(S, \Lambda^{0,1}S\otimes TS),
\end{equation}
where $\Lambda^{0,1}S$ denotes the line bundle of $(0,1)$-form on $S$.

Let $\scrj$ be an open {\sl connected} subset in the Banach manifold of all
$C^\ell$-smooth almost complex structures on $X$ for some fixed non-integer $\ell>
2$. In our context the most interesting example is the set $\scrj_\omega$ of
$C^\ell$-smooth $\omega$-tame almost complex structures on $X$. The tangent space
to $\scrj$ at $J$ consists of $C^\ell$-smooth $J$-antilinear endomorphisms of
$TX$,
\begin{equation}
T_J\scrj = \{ I\in C^\ell(X, \endo(TX)) : JI + IJ = 0\}
\equiv C^\ell(X, \Lambda^{0,1}X \otimes TX),
\end{equation}
where $\Lambda^{0,1}X$ denotes the complex bundle of $(0,1)$-forms on $X$.

Fix $p$ with $2<p<\infty$. Then the set $L^{1,p}(S, X)$ of all Sobolev
$L^{1,p}$-smooth maps from $S$ to $X$ is a Banach manifold. For $u\in L^{1,p}(
S, X)$ we denote
\begin{equation}
E_u \deff u^*TX.
\eqqno(1.3.1)
\end{equation}
In this notation, the tangent space at $u \in L^{1,p}(S, X)$ is $T_uL^{1,p}(S,
X) = L^{1,p}(S, E_u)$, the space of $L^{1,p}$-smooth sections
of the pulled-back tangent bundle of $X$.

Fix a homology class $[C] \in \sfh_2(X, \zz)$ and consider the set
\begin{equation}
\scrs = \{ u\in L^{1,p}(S, X):u(S)\in [C]\}
\end{equation}
of maps $u$ representing the class $[C]$. Then $\scrs$ is open
in $L^{1,p}(S, X)$ and has the same tangent space, $T_u\scrs = L^{1,p}(S,
E_u)$. Since $\scrj$ is connected, the first Chern class $c_1(X,J)$ is
constant on $\scrj$. We shall denote it simply by $c_1(X)$. Set
\begin{equation}
\mu \deff \la c_1(X), [C] \ra.
\eqqno(1.3.1a)
\end{equation}

\smallskip
Consider the subset $\scrp\subset \scrs\times \scrj_S \times \scrj$
consisting of all triples $(u,J_S,J)$ with $u$ being $(J_S,J)$-holomorphic,
\begin{equation}
\scrp = \{(u, J_S, J)\in \scrs\times \scrj_S\times \scrj: du +
J\scirc du\scirc  J_S = 0 \}.
\eqqno(1.3.2)
\end{equation}

Let $\nabla$ be some symmetric connection on $TX$. Covariant differentiation
of \eqqref(1.1.1) gives the equation for the tangent space to $\scrp$.
Namely, a vector $(v,\dot J_S,\dot J)$ is tangent to $\scrp$ at the point
$(u,J_S,J)$ if it satisfies the equation
\begin{equation}
\nabla v + J\scirc \nabla v\scirc J_S + (\nabla_vJ)\scirc (du\scirc J_S) +
J\scirc du\scirc \dot J_S + \dot J\scirc du\scirc J_S = 0. \eqqno(1.3.3)
\end{equation}

\newdefi{def1.3.1}\.{\bf a)} For a complex bundle $E$ over $S$ let 
\begin{equation}
L^p_{(0,1)}(S, E)\deff L^p(S, E\otimes \Lambda^{(0,1)}S)
\end{equation}
denote the Banach space of $L^p$-integrable $E$-valued (0,1)-forms on $S$.

\statep b). Let $u$ be a $J$-holomorphic curve in $X$. Define the operator
$D_{u, J}: L^{1,p}(S, E_u) \to L^p_{(0,1)}(S, E_u)$ as
\begin{equation}
D_{u, J}(v) \deff \nabla v + J\scirc\nabla v\scirc J_S
+ (\nabla_vJ) \scirc du\scirc J_S
\eqqno(1.3.4)
\end{equation}

\statep c). Define complex Banach bundles $\scre$ and $\scre'$ over $\scrs
\times \scrj_S \times \scrj$ by 
\begin{equation}
\scre_{(u,J_S,J)} \deff L^{1,p}(S, E_u)
\quad\hbox{and}\quad
\scre'_{(u,J_S,J)} \deff L^p_{(0,1)}(S, E_u).
\eqqno(1.3.5)
\end{equation}
These bundles are $C^\ell$-smooth and the formula \eqqref(1.3.4) defines a
$\rr$-linear homomorphism $D=D_{u,J_S,J}: \scre \to \scre'$ which is $C^{\ell-1}
$-smooth. The bundle $\scre$ is essentially the tangent bundle to $\scrs$, whereas 
$\scre'$ appears as the space where the equation \eqqref(1.1.1) ``lives''. 
More precisely, \eqqref(1.1.1) defines a section $\sigma_\dbar$ of
$\scre'$,
\begin{equation}
\sigma_\dbar: (u,J_S,J) \in \scrs \times \scrj_S \times \scrj
\;\mapsto \;
 (du + J \scirc du \scirc J_S) \in \scre'_{(u,J_S,J)},
\eqqno(1.3.6)
\end{equation}
such that the equation \eqqref(1.1.1) reads $\sigma_\dbar(u, J_S, J)=0$. 
The space $\scr P$ appears then as the zero set of the section $\sigma_\dbar$.
\end{defi}

\state Remark. Here and thereafter we use the normalization ${\d \over \d\bar z} 
= {\d\over \d x} + \isl {\d\over \d y}$, deviating from the usual convention 
${\d \over \d\bar z} = \mathbf{\half} \cdot ({\d\over \d x} +\isl {\d\over\d y})$.
The same normalization is used for all operators with Cauchy-Riemann symbol.

\newlemma{lem1.3.1}
Let $\scrx$ be a Banach manifold, $\scre\to \scrx$
and $\scre'\to\scrx$ $C^1$-smooth Banach bundles over $\scrx$, $\nabla$ and
$\nabla'$ linear connections in $\scre$ and $\scre'$ respectively, $\sigma$
a $C^1$-smooth section of $\scre$, and $D: \scre \to \scre'$ a $C^1$-smooth
bundle homomorphism.

\sli If $\sigma(x)=0$ for some $x\in \scrx$, then the map $\nabla\sigma_x:
T_x\scrx \to \scre_x$ is independent of the choice of the connection $\nabla$
in $\scre$;

\slii Set $K_x\deff \ker(D_x:\scre_x\to \scre'_x)$ and $Q_x\deff \coker(D_x:
\scre_x\to \scre'_x)$ with the corresponding imbedding $i_x:K_x \to \scre_x$
and projection $p_x:\scre'_x \to Q_x$. Let $\nabla^\hom$ be the connection in
$\hom(\scre,\scre')$ induced by the connections $\nabla$ and $\nabla'$. Then
the map
\begin{equation}
p_x \circ(\nabla^\hom\!\! D_x)\scirc i_x:
T_x\scrx\to \hom(K_x,Q_x)
\end{equation}
is independent of the choice of connections
$\nabla$ and $\nabla'$.
\end{lem}

\state Remark. Taking this lemma into account, we shall use the following
notation. For $\sigma\in \Gamma(\scrx, \scre)$, $D\in \Gamma(\scrx,
\hom(\scre, \scre'))$ and $x\in \scrx$ as in the hypothesis of the lemma we
shall denote by $\nabla\sigma_x:T_x\scrx \to \scre_x$ and ${\nabla\!D}:
T_x\scrx\times \ker D_x \to \coker D_x$ the corresponding operators without
pointing out which connections were used to define them.

\proof \sli Let $\wt\nabla$ be another connection in $\scre$. Then
$\wt\nabla$ has the form $\wt\nabla=\nabla +A$ for some
$A\in\Gamma(\scrx,\hom(T\scrx,\endo(\scre)))$. So for $\xi\in T_x\scrx$
we obtain $\wt\nabla_\xi\sigma -\nabla_\xi\sigma = A(\xi,\sigma(x))=0$.

\slii Similarly, let $\wt\nabla'$ be another connection in $\scre'$, and let
$\wt\nabla^\hom$ be the connection in $\hom(\scre,\scre')$ induced by
$\wt\nabla$ and $\wt\nabla'$. Then $\wt\nabla'$ also has the form $\wt\nabla=
\nabla +A'$ for some $A'\in\Gamma(\scrx,\hom(T\scrx, \endo(\scre')))$.
So for $\xi\in T_x\scrx$ we obtain $\wt\nabla_\xi^\hom D -\nabla_\xi^\hom D=
A'(\xi)\scirc D_x - D_x \scirc A(\xi)$. The statement of the lemma now follows
from the identities $p_x\scirc D_x=0$ and $D_x\scirc i_x=0$.
\qed

\smallskip
\state Remark. The operator $D_{u, J}$ is the linearization of the equation
\eqqref(1.1.1). Thus \lemma{lem1.3.1} shows that the definition of $D_{u, J}$
is independent of the choice of $\nabla$. In particular, one can also use
non-symmetric connections, \eg those compatible with $J$, as it is done in
\cite{Gro}. However, it is convenient to have a fixed
connection considering varying almost complex structures $J$ on $X$. But this 
is impossible for $\nabla$ compatible with $J$. On the other hand, with a
symmetric connection computations become simpler.

\smallskip
The operator $D_{u, J}$, as well as the equation \eqqref(1.1.1) itself, is 
elliptic
of order 1 with the Cauchy-Riemann symbol. This implies standard regularity 
properties for $D_{u, J}$. In particular, the kernel and the cokernel are of 
finite dimension. The Riemann-Roch formula gives the index of $D_{u, J}$:
\begin{equation}
\dim_\rr \ker D_{u, J} - \dim_\rr \coker D_{u, J} = 2\cdot
\bigl(\mu  + n (1-g) \bigr),
\eqqno(1.3.7)
\end{equation}
where $\mu \deff c_1(X) \cdot [u(S)]$, $g$ is the genus of $S$, and $n$ the
{\sl complex} dimension of $X$, \ie $n \deff \half \dim_\rr X$. The factor
2 appears because we compute {\sl real, not complex} dimensions of the
(co)kernel.

\newsubsection[1.4]{Holomorphic structure on the induced bundle}
Now we want to understand the structure of the operator $D_{u,J}$ in more
detail. Note that the pulled-back bundle $E_u = u^*TX$ carries a complex
structure, namely $J$ itself, or more accurately $u^*J$. However, the
operator $D\deff D_{u,J}$ is only $\rr$-linear. So we decompose it into
$J$-linear and $J$-antilinear parts. Namely, for $v\in
L^{1,p}(S, E)$ we write $D v = \half \bigl( D v - JD(Jv)\bigr) +
\half (D v + JD (Jv)\bigr) = \dbar_{u, J}v+ R(v)$.

\newdefi{def1.4.1} The $J$-linear part $\dbar_{u, J}$ of the operator $D_{u, J}$
is called the {\sl $\dbar$-operator associated with a $J$-holomorphic map $u$}.
\end{defi}

\smallskip
By the definition, the operator $\dbar_{u, J} : L^{1,p}(S, E_u) \to L^p
_{(0,1)}(S, E_u)$ is $J$-linear. The following  statement is well
known in the smooth case.

\newlemma{lem1.4.1A} Let $S$ be a Riemann surface with a complex
structure $J_S$ and $E$ a $L^{1,p}$-smooth complex vector bundle of rank $r$
over $S$. Let also $\dbar_E : L^{1,p}(S, E) \to L^p_{(0,1)}(S,E)$ be a complex
linear differential operator satisfying the condition
\begin{equation}
\dbar_E(f\xi) = \dbar_S f \otimes \xi
+ f\cdot \dbar_E\xi,  \eqqno(1.3.6a)
\end{equation}
where $\dbar_S$ is the Cauchy-Riemann operator on $S$ associated to $J_S$.
Then the sheaf
\begin{equation}
U\subset S \mapsto \scro(E)(U) := \{\,\xi \in  L^{1,p}(U, E) \, :\,
\dbar_E\xi=0 \,\}   \eqqno(1.3.7a)
\end{equation}
is coherent and locally free of rank $r$. This defines a holomorphic
structure on $E$ for which $\dbar_E$ is the associated Cauchy-Riemann
operator.
\end{lem}

\state Remark. The condition $\eqqref(1.3.6a)$ means that $\dbar_E$ is of 
order 1 and has the Cauchy-Riemann symbol. For the proof we refer to 
\cite{Iv-Sh-1} and \cite{Iv-Sh-2} for the general case, or to \cite{H-L-S} 
for  the case of line bundles.

\smallskip
Thus, according to \lemma{lem1.4.1A}, the operator $\dbar_{u, J}$ defines
a holomorphic structure on the bundle $E_u$. We shall denote by $\scro(E_u)$
the sheaf of holomorphic sections of $E_u$. The tangent bundle $TS$ to our
Riemann surface also carries a natural holomorphic structure. We shall
denote by $\scro(TS)$ the~corresponding coherent sheaf.

\medskip
Denote by $N_J(v, w)$ the Nijenhuis torsion tensor of the almost complex
structure $J$, (see \eg \cite{Ko-No}, Vol.II., p.123.)\.\footnote{\.
Note that in \cite{Ko-No} another normalization constant is used. However,
this is not essential for our purpose.}

\newlemma{lem1.4.1}
\sli The $J$-antilinear part $R$ of $D_{u,J}$ is related to $u$ and $J$ by
the formula
\begin{equation}
R(v)(\xi ) = N_J(v, du(\xi ))\qquad \xi \in TS.
\eqqno(1.4.1)
\end{equation}
Thus $R$ is a continuous $J$-antilinear operator from $E$ to $\Lambda^{0,1}S
\otimes E_u$ of order zero which satisfies $R\scirc du \equiv 0 $, \ie
$R( du(\eta), \xi )=0$ for all $\eta ,\xi \in TS$.

\smallskip
\slii If $u$ is non-constant, then $du$ defines an injective analytic
morphism of coherent sheaves
\begin{equation}
0\lrar \scro(TS) \buildrel du \over \lrar \scro(E_u). \eqqno(1.4.2)
\end{equation}
\end{lem}

\proof \sli Formula \eqqref(1.4.1) can be found in \cite{McD-2}. The rest of
part \sli follows from the well-known fact that $N_J(v, w)$ is skew-symmetric
and $J$-antilinear in both arguments.

\noindent \slii The fact that $du: TS \to E_u$ defines a morphism between
coherent sheaves $\scro(TS)$ and $\scro(E_u)$ means that $du$ is
a {\sl holomorphic} section of $T^*S \otimes E_u$. This is equivalent to
relation
\begin{equation}
du \scirc \dbar_S = \dbar_{u,J} \scirc du.
\end{equation}
For the proof of this fact we refer to \cite{Iv-Sh-1} and \cite{Iv-Sh-2}.

Injectivity of the sheaf homomorphism $du$ is equivalent to its
nondegeneracy which is the case in our context. \qed

\smallskip
The zeros of the analytic morphism $du : \scro(TS) \to \scro (E_u)$ are
isolated. So we obtain

\newcorol{cor1.4.2} {\rm(\cite{Mi-Wh}, \cite{Sk-1})}
The set of critical points of a $J$-holo\-mor\-phic map is discrete, 
provided $J$ is of class $C^1$.
\end{corol}

\newdefi{def1.4.2} By the {\sl order of zero $\ord_p du$} of the
differential $du$ at a point $p\in S$ we shall understand the order of
vanishing at $p$ of the holomorphic morphism $du : \scro(TS)\to \scro(E_u)$.
\end{defi}

It follows from \lemma{lem1.4.1} that $\ord_p du$ is a well-defined
non-negative integer.

\newsubsection[1.5]{The normal sheaf of a pseudoholomorphic curve}
From \eqqref(1.4.2) we obtain the following short exact sequence of coherent 
sheaves
\begin{equation}
0\lrar \scro(TS) \buildrel du \over\lrar \scro(E_u)
\lrar \scrn_u\lrar 0,
\eqqno(1.5.1)
\end{equation}
where $\scrn_u\deff \scro(E)/du(\scro(TS))$ is the quotient sheaf. It follows from
\lemma{lem1.4.1} \slii that there is a decomposition $\scrn_u= \scro(N_u)\oplus 
\scrn_u\sing$ where $N_u$ is a holomorphic
vector bundle and $\scrn_u\sing = \bigoplus_{z\in S} \cc_z^{\ord_z du}$ is a
discrete sheaf with support in the set of critical points $a_i$ of $u$
with the stalk $\cc^{n_i}$ of dimension $n_i \deff \ord_{a_i}du$ at every
such point $a_i$. 

\newdefi{def1.5.1} The quotient sheaf $\scrn_u \deff \scro(E)/du(TS)$ is called 
the {\sl normal sheaf of a $J$-curve $u:S \to X$}, $N_u$ the {\sl normal 
bundle to the $J$-curve $u:S \to X$}, and $[A] \deff \sum n_i[a_i]$ 
the {\sl branching divisor of the $J$-curve $u:S \to X$}.
\end{defi}

Denote by $\scro([A])$ the sheaf of meromorphic functions on $S$ having poles
of order at most $n_i$ at $a_i$. Then \eqqref(1.5.1) gives rise to the exact 
sequence of coherent sheaves
\begin{equation}
0\lrar \scro(TS)\otimes \scro([A])\buildrel{du}\over
{\lrar} \scro(E_u)
\lrar \scro(N_u)\lrar 0.
\eqqno(1.5.2)
\end{equation}

The holomorphic structure in $N_u$ defines the Cauchy-Riemann operator
$\dbar_N :L^{1,p}(S, N_u) \allowbreak \lrar L^p_{(0,1)}(S, N_u)$.
\lemma{lem1.4.1} implies that
the homomorphism $R: E_u \to E_u \otimes \Lambda^{(0,1)}S$
induces a $J$-antilinear bundle homomorphism $R_N : N_u \to N_u \otimes
\Lambda^{(0,1)}S$. Define the operator
\begin{equation}
D_{u, J}^N : L^{1,p}(S, N_u) \allowbreak \lrar L^p_{(0,1)}(S, N_u)
\quad\hbox{by}\quad
D_{u, J}^N \deff\dbar_N + R_N.
\eqqno(1.5.3)
\end{equation}

\smallskip
\newdefi{def1.5.2} Let $E$ be a holomorphic vector bundle over a compact Riemann 
surface $S$ of genus $g$ and let $D:L^{1,p}(S, E)\to L^p(S, \Lambda^{0,1}S \otimes 
E)$ be an operator of the~form $D=\dbar + R$ where $R\in L^p\bigl(S,\,\hom_\rr(E,
\,\Lambda^{0, 1}S\otimes E) \bigr)$ with $2<p<\infty$. Define $\sfh^0_D(S, E)\deff 
\ker D$ and $\sfh^1_D(S, E) \deff \coker D$.  The groups $\sfh^i_D(S, E)$ are 
referred to as {\sl $D$-cohomology groups of $E$}.
\end{defi}

\medskip
The Riemann-Roch formula gives the {\sl index} of $D$,
\begin{equation}
\ind_\rr D \deff \dimr\sfh^0_D(S, E) - \dimr\sfh^1_D(S, E)=
2\bigl( c_1(E) + \rank(E) (1-g)\bigr).
\eqqno(1.5.?4)
\end{equation}

\state Remark. Taking into account the elliptic regularity of the
Cauchy-Riemann operator, for given $S$, $E$ and $R\in L^p$, $2<p<\infty$, one
can define $\sfh^i_D(S, E)$ as the (co)kernel of the operator $\dbar +R:
L^{1,q}(S, E) \to L^q(S,\,\Lambda^{0, 1}S\otimes E)$ for any $q\in\; ]1, p]$.
So the definition is independent of the choice of the functional spaces. Note also
that the $\sfh^i_D(S, E)$ are of finite dimension provided that $S$ is closed.
For details, see \eg \cite{Iv-Sh-1}.

\smallskip
The following lemmas contain main properties of $D$-cohomologies which will
be used later. For complete proofs we refer to \cite{Iv-Sh-1} and \cite{Iv-Sh-2}.

\newlemma{lem1.5.1} {\sl (Serre duality for $D$-cohomologies.)}
Let $E$ be a holomorphic vector bundle over a~compact Riemann
surface $S$ and let $D:L^{1,p}(S, E)\to L^p_{(0,1)}(S, E)$ be
an operator of the~form $D=\dbar + R$, where $R\in L^p\bigl(S,\,
\hom_\rr(E,\,\Lambda^{0, 1}S\otimes E) \bigr)$ with $2<p <\infty$. Let
$K_S\deff \Lambda^{1, 0}S$ be the~canonical holomorphic line
bundle of $S$. Then there exists the~naturally defined operator
\begin{equation}
D^*=\dbar- R^*: L^{1,p}(S, E^* \otimes K_S)
\to L^p_{(0,1)}(S, E^* \otimes K_S)
\end{equation}
with $R^* \in L^p\bigl(S,\,\homr(E^*\otimes K_S,\,
\Lambda^{0, 1}S \otimes E^*\otimes K_S) \bigr)$ and the~natural isomorphisms
\begin{equation}
\matrix
\sfh^0_D(S,\, E)^* &\cong& \sfh^1_{D^*}(S,\, E^*\otimes K_S),
\cr
\sfh^1_D(S,\, E)^* &\cong& \sfh^0_{D^*}(S,\, E^*\otimes K_S),
\endmatrix
\eqqno(1.5.4)
\end{equation}
induced by the pairings
\begin{equation}
\begin{matrix}
\phi \in \sfh^0_D(S,\, E), \quad
\psi \in L^p_{(0,1)}(S,\, E^*\otimes K_S)&
\mapsto \<\phi,\psi\> &
\deff \re \int_S \psi \scirc \phi
\cr
\psi \in \sfh^0_D(S,\, E^*\otimes K_S), \quad
\phi \in L^p_{(0,1)}(S,\, E)&
\mapsto \<\phi,\psi\> &
\deff \re \int_S \psi \scirc \phi
\end{matrix}
\eqqno(1.5.5)
\end{equation}
If, in addition, $R$ is $\cc$-antilinear, then $R^*$ is also
$\cc$-antilinear.
\end{lem}

\state Remark. The lemma expresses the well-known relation $\ker(D^*) =(\im D)
^\perp$ between a linear operator $D$ and its adjoint. It is worth observing
that the spaces themselves and the duality are defined only over the real numbers
$\rr$ and not over $\cc$.

\smallskip
\newlemma{lem1.5.2}
{\rm(\cite{H-L-S}, {\sl Vanishing theorem for $D$-cohomologies.})}
Let $S$ be a closed Riemann surface of genus $g$ and $L$ 
a~holomorphic {\sl line} bundle over $S$, equipped with a~differential
operator $D=\dbar + R$ with $R\in L^p\bigr(S, \homr(L,\, \Lambda^{0,
1}S\otimes L) \bigl)$, $p>2$. If $c_1(L)<0$, then $\sfh^0_D (S,\, L)=0$. If
$c_1(L)>2g-2$, then $\sfh^1_D (S,\, L)=0$.
\end{lem}

\medskip
The importance of the operator $D=\dbar + R$ lies in the fact that we can
associate with the short exact sequence \eqqref(1.5.1) the long exact 
sequence of $D$-cohomo\-logies. Note, that due to \lemma{lem1.4.1} we obtain
the short exact sequence of complexes
\begin{equation}
\eqqno(1.5.6)
\def\normalbaselines{\baselineskip20pt\lineskip3pt \lineskiplimit3pt }
\setbox1=\hbox{$\lrar$}
\def\mapright#1{\!\!\!\smash{\mathop{{\hbox to
\wd1{\hss\hbox{$\displaystyle\longrightarrow$}\hss}}}\limits^{#1}}\!\!\!}
\def\mapdown#1{\Big\downarrow\rlap{$\vcenter{\hbox{$\scriptstyle#1$}}$}}
\begin{matrix}
\llap{$0\lrar$}  L^{1,p}(S, TS)&
\mapright{du}&
L^{1,p}(S, E)&
\mapright{\barr\pr}&
L^{1,p}(S, E)\bigm/ du(L^{1,p}(S, 
\rlap{$TS))\lrar0$}
\\
\mapdown{\barr \partial_S}\qquad&& \mapdown{D}\qquad&&
\mapdown{\barr D}\qquad\qquad\qquad\qquad
\\
\llap{$0\lrar$} L^p_{(0,1)}(S, TS)&
\mapright{du}&
L^p_{(0,1)}(S, E)&
\mapright{\barr\pr}&
L^p_{(0,1)}(S, E)\bigm/
du(L^p_{(0,1)}(S, 
\rlap{$TS)) \lrar0$}
\end{matrix}
\end{equation}
where $\barr D$ is induced by $D\equiv D_{u, J}$.

\newlemma{lem1.5.3} For $\barr D$ as in \eqqref(1.5.6), $\ker \barr D =
\sfh^0_D (S, N_u) \oplus \sfh^0(S, \scrn_u\sing)$ and $\coker \barr D
\allowbreak= \sfh^1_D (S, N_u)$.
\end{lem}

\proof For an open set $U \subset S$ let $\Gamma_D(U, E_u) \deff \{ v\in L^{1,p}
_\loc(U, E_u) \;:\; Dv = 0\}$. Use the analogous notation for $N_u$. Consider the 
sheaves $U \mapsto \Gamma\, (U, \scro(TS)\,)$, $U \mapsto \Gamma_D (U, E_u)$, and 
$U \mapsto \Gamma_D(U, N_u) \oplus \Gamma(U, \scrn_u \sing)$. It is easy to show 
that the first two columns of the diagram \eqqref(1.5.6)
define fine resolutions of the sheaves $\scro(TS)$ and $\Gamma_D (\, \cdot\,,
E_u)$. Moreover, $du$ defines injective homomorphisms between these sheaves
and between their resolutions. An explicit computation shows that $\Gamma_D(\,
\cdot \,, N_u) \oplus \Gamma(\,\cdot\,, \scrn_u\sing)$ is the
corresponding quotient sheaf and that the third column of \eqqref(1.5.6) is its
resolution. For details, see \cite{Iv-Sh-1}. \qed

\newcorol{cor1.5.4} The~short exact sequence $\eqqref(1.5.1)$ induces
the~long exact sequence of $D$-cohomologies
\begin{equation*}
\def\normalbaselines{\baselineskip20pt\lineskip3pt \lineskiplimit3pt }
\def\mapright#1{\;\smash{\mathop{\longrightarrow}\limits^{#1}}\;}
\def\mapdown{\Big\downarrow}
\matrix
0& \mapright{}& \sfh^0(S, TS) &\mapright{} & \sfh^0_D(S, E)
 & \mapright{}& \sfh^0_D(S, N_u)\oplus \sfh^0(S, \scrn_u\sing)
 &\mapright{\delta} &\vphantom{\mapdown}\\
 & \mapright{}& \sfh^1(S, TS) &\mapright{} & \sfh^1_D(S, E)
 & \mapright{}& \sfh^1_D(S, N_u)         &\mapright{} &0.
\endmatrix
\end{equation*}
\end{corol}

\medskip
\newsection[2]{The total moduli space of pseudoholomorphic curves}

\newsubsection[2.1]{Transversality}
Any deformation of a given $J$-holomorphic map $u: S \to X$ defines a path in
the space $\scrp$ of pseudoholomorphic  maps. Thus to construct such a deformation 
we want to equip the space $\scrp$ with a structure of a smooth Banach manifold.

Note that by definition the set $\scrp$ is the zero set of the
section $\sigma_\dbar$ of the bundle $\scre' $, \ie the intersection of
the images of $\sigma_\dbar$ and the zero-section $\sigma_0$.
Thus we are interested in which points these sections meet
transversally. The analysis of the problem leads to the following

\newdefi{def2.1.1} Let $\scrx$, $\scry$, and $\scrz$ be Banach manifolds
with $C^\ell$-smooth maps $f:\scry \to \scrx$ and $g:\scrz\to\scrx$, $\ell
\ge1$. Define the {\sl fiber product} $\scry\times_\scrx \scrz $ by setting
$\scry\times_\scrx \scrz \deff \{\,(y,z)\in \scry\times\scrz \;:\; f(y)= g(z)
\,\}$. The map $f$ is called {\sl transversal to $g$} at a point $(y,z)\in
\scry\times_\scrx \scrz$ with $x\deff f(y)=g(z)$, and $(y,z)$ is called a
{\sl transversality point}, if the map $df_y\oplus dg_z: T_y\scry \oplus
T_z\scrz \to T_x\scrx$ is {\sl surjective} and admits a {\sl closed
complement} to its kernel. The set of transversality points $(y,z)\in \scry
\times_\scrx \scrz$ will be denoted by $\scry\times^\trans_\scrx \scrz$, with
$\trans$ symbolizing the transversality condition.

We say that $f:\scry \to \scrx$ is {\sl transversal} to $g:\scrz\to\scrx$ if
every point $(y,z) \in \scry\times_\scrx \scrz$ is a transversality point. In
particular, if $\scry$ consists of a point $x \in \scrx$ and the imbedding $\{
x\} \hook \scrx$ is transversal to $g$, we call $x$ a {\sl regular value of
$g$}. Note that by this definition any $x\in \scrx\bs g(\scrz)$ is a regular
value of $g$.

In the special case when the map $g:\scrz \to \scrx$ is a closed imbedding,
the fiber product $\scry\times_\scrx \scrz$ is simply the preimage
$f\inv\scrz$ of $\scrz\subset\scrx$. In particular, every point
$(y,z)\in\scry\times_\scrx \scrz$ is completely defined by its component
$y\in\scry$, $z=f(y)\in\scrz\subset \scrx$. In this case we simply say that
$f:\scry \to \scrx$ is {\sl transversal to $\scrz$ at $y\in\scry$},
\iff $(y,f(y))$ is a transversal point of $\scry\times_\scrx \scrz\cong
f\inv\scrz$.
\end{defi}

\newlemma{lem2.1.1}
The set $\scry\times^\trans_\scrx \scrz$ is open in $\scry\times_\scrx \scrz$
and is a $C^\ell$-smooth Banach manifold with tangent space
\begin{equation}
T_{(y,z)}\scry\times^\trans_\scrx \scrz = \ker\bigl(
df_y\oplus d(-g_z): T_y\scry \oplus T_z\scrz \to T_x\scrx \bigr).
\end{equation}
\end{lem}

\proof Fix $w_0\deff(y_0,z_0)\in \scry\times^\trans_\scrx \scrz$
and set $K_0 \deff \ker( df_{y_0}\oplus dg_{z_0}:
T_{y_0}\scry \oplus T_{z_0}\scrz \to T_x\scrx \bigr)$. Let $Q_0$ be a closed
complement to $K_0$. Then the map $df_{y_0}\oplus dg_{z_0}:Q_0\to T_x\scrx$
is an isomorphism.

Due to the choice of $Q_0$, there exists a neighborhood $V\subset
\scry\times\scrz$ of $(y_0,z_0)$ and $C^\ell$-smooth maps $w':V\to K_0$ and
$w'':V\to Q_0$, such that $dw'_{w_0}$ (resp.\ $dw''_{w_0}$) is the projection
from $T_{y_0}\scry \oplus T_{z_0}\scrz$ onto $K_0$ (resp.\ onto $Q_0$),
so that $(w',w'')$ are coordinates in some smaller  neighborhood $V_1\subset
\scry\times\scrz$ of $w_0=(y_0,z_0)$. It remains to consider the equation
$f(y)=g(z)$ in new coordinates $(w',w'')$ and apply the implicit function
theorem.
\qed

\medskip
Due to \lemma{lem2.1.1}, the set $\scrp$ is a Banach manifold at those points
$(u,J_S,J)\in\scrp$ where $\sigma_\dbar$ is transversal to $\sigma_0$.
However, at any point $(u,J_S,J;0)$ on the zero section $\sigma_0$ of
$\scre'$ we have the natural decomposition
\begin{equation}
T_{(u,J_S,J;0)}\scre' = d\sigma_0\bigl( T_{(u,J_S,J)}
(\scrs\times \scrj_S\times \scrj) \bigr)
\oplus \scre'_{(u,J_S,J)},
\end{equation}
where the first component is the tangent space to the zero section of $\scre'$
and the second one is the tangent space to the fiber $\scre'_{(u,J_S,J)}$.

Let $p_2$ denote the projection on the second component. Then the
transversality $\sigma_\dbar$ and $\sigma_0$ is equivalent to the surjectivity
of the map $p_2\scirc d\sigma_\dbar: T_{(u,J_S,J)}(\scrs\times \scrj_S \times
\scrj)\to \scre'_{(u,J_S,J)}$, \ie to the surjectivity of the operator
\begin{align*}
\nabla\sigma_\dbar:\;&
 T_uL^{1,p}(S,X)\oplus T_{J_S}\ttt_g \oplus T_J\scrj \lrar
\scre'_{(u,J_S,J)}
\\
\nabla\sigma_\dbar:\;&(v,\dot J_S,\dot J) \longmapsto
D_{(u,J)}v + J\scirc du \scirc \dot J_S
+ \dot J\scirc du \scirc J_S.
\end{align*}
By \refdefi{def1.5.2}, the quotient of $\scre'_{(u,J_S,J)}$ with respect
to the image of $D_{u,J}$ is $\sfh^1_D(S,E_u)$. The induced map
\begin{equation}
T_{J_S}\scrj_S \ni \dot J_S \mapsto J\scirc du \scirc \dot J_S
\in \sfh^1_D(S,E_u)
\eqqno(2.1.1)
\end{equation}
is also easy to describe. Recall that for a given complex
structure $J_S$ on $S$ one has the Dolbeault isomorphism
\begin{equation}
\sfh^1(S,TS)= \coker \bigl(\,\dbar :
C^{2,\alpha}(S, TS) \lrar C^{1,\alpha}_{(0,1)}(S, TS) \,\bigr)
\end{equation}
with the operator $\dbar$ associated to $J_S$. Recall also that $C^{1,\alpha}
_{(0,1)}(S, TS)$ is the tangent space $T_{J_S}\scrj_S$. This shows that the map 
\eqqref(2.1.1) is the same as the map $J\scirc du: \sfh^1(S,TS)
\to \sfh^1_D(S,E_u)$ and, due to identity $J\scirc du = du\scirc J_S$ and 
\refcorol{cor1.5.4}, its cokernel is $\sfh^1_D(S,N_u)$.

It remains to study the image of $T_J\scrj$ in $\sfh^1_D(S,N_u)$.

\newdefi{def2.1.2} For $(u,J_S,J)\in \scrp$ we define  $\Psi=
\Psi_{(u,J)}:T_J\scrj \to \scre'_{(u,J_S,J)}$ by setting $\Psi_{(u,J
)}( \dot J) \deff \dot J\scirc du \scirc J_S$. Let $\barr\Psi=\barr\Psi_{(u,
J)}:T_J\scrj \to \sfh^1_D(S,N_u)$ be induced by $\Psi$. Finally, define
$\scrp^* \deff \{ (u,J_S,J) \in \scrp: u \hbox{ is injective in generic }z
\in S\}$.
\end{defi}

\state Remark. One can show that $\scrp\bs \scrp^*$ consists of {\sl multiple
curves} for which the map $u:(S,J_S) \to (X,J)$ admits a factorization
$u= u' \scirc g$ for some non-trivial holomorphic branched covering $g:(S,J_S)
\to (S',J'_S)$ and a $J$-holomorphic map  $u':(S',J'_S) \to X$. On the other
hand, for any $(u,J_S,J) \in \scrp^*$ the map $u$ is a smooth imbedding outside 
finitely many points in $S$. For details see \cite{Mi-Wh} or \cite{Iv-Sh-1}.

\newlemma{lem2.1.2}  {\sl(Infinitesimal transversality)}.
Let $(u,J_S,J)\in \scrp^*$. Then the operator $\barr\Psi:T_J\scrj \to
\sfh^1_D(S, N_u)$ is surjective.
\end{lem}

\proof Choose some nonempty open set $V\subset S$, such that
$u\ogran_V$ is an imbedding. Use Serre duality (\lemma{lem1.5.1}) and
find a basis $\psi_1,\ldots \psi_l\in  \sfh^0_D(S, N^*\otimes K_S)\cong
\sfh^1_D(S, N)^*$.

Note that $\psi_i$ satisfy the equation $D \psi_i =0$, where the operator
$D= D_{N^*\otimes K_S}$ is of the form $\dbar + R$. One can show
(see \eg \cite{Iv-Sh-1} or \cite{H-L-S}) that any solution $v$ of the equation
$(\dbar + R)v=0$ is $L^{1,p}$-smooth and furthermore such a $v$ is either
identically zero or has isolated zeros.

This implies that there exist $I_1,\ldots I_l\in  C^\ell(S, N \otimes
\Lambda^{(0,1)} )$ with supports $\supp(I_i)$ in $V$ such that the
matrix $\bigl(\re\int_S \psi_i \scirc I_j\bigr)_{i,j=1} ^l$ is non-degenerate.
Since $u\ogran_V$ is a $C^{\ell+1}$-smooth imbedding, any such $I_i$ can be
represented in the form $I_i = \dot J_i \scirc du \scirc J_S$ with some
$\dot J_i \in C^\ell(X, \endo(TX))$ with $ J \scirc \dot J_i + \dot J_i\scirc
J =0$. The latter relation means that $\dot J_i \in T_J\scrj$. \qed

\smallskip
\newcorol{cor2.1.3} $\scrp^*$ is a $C^\ell$-smooth Banach manifold
with the tangent space
\begin{equation}
T_{(u,J_S,J)} \scrp^* = \bigl\{ (v,\dot J_S, \dot J) \;:\;
D_{u,J}v + J \scirc du \scirc \dot J_S + \dot J \scirc du \scirc J_S
=0 \;\bigr\}.
\eqqno(2.1.2)
\end{equation}
\end{corol}

\state Remark. $\scrp^*$ is in general smaller than the set $\scrp^\trans
\deff \sigma_\dbar \times^\trans_{\scre'} \sigma_0$ of all transversality
points of $\scrp$. On the other hand, it is sufficient for applications to
consider only the space $\scrp^*$.

\newsubsection[2.2]{Moduli space of pseudoholomorphic curves}
The space $\scrp$ (see \eqqref(1.3.2)) of all pseudoholomorphic maps is too big.
Indeed, one has a natural right action of the group $\diff _+(S)$ of the
orientation preserving $C^{2,\alpha}$-smooth diffeomorphisms of $S$ on the
product $\scrs \times \scrj_S \times \scrj$ given by formula
\begin{equation}
(u,J_S,J) \in \scrs \times \scrj_S \times \scrj, g\in \diff_+(S)
\; \lrar \; (u,J_S,J)\cdot g \deff (u\scirc g, J_S \scirc g, J),
\end{equation}
such that $\scrp$ and $\scrp^*$ are invariant \wrt this action. It is natural
to consider $(u,J_S,J) \in \scrp$ and $(u,J_S,J) \cdot g$ as two
parameterization of the same a $J$-curve. In other words, we are interested
in the quotient space $\scrp/ \diff_+(S)$ rather than the space $\scrp$
itself. As in Yang-Mills theory one can treat the group
$\diff_+(S)$ as the gauge group of the problem and the quotient as
the corresponding moduli space. Again as in Yang-Mills theory it is useful to
know in which points the quotient spaces $\scrp/\diff_+(S)$ is
a Banach manifold.

\smallskip
First we consider the action of the group $\diff_+(S)$ on the space $\scrj_S$.
Note that if $J'_S, J''_S \in \scrj_S$ are related by $J''_S = J'_S \scirc g$
for some $C^1$-diffeomorphism $g:S \to S$, then $g$ is $(J'_S, J''_S)
$-holomorphic. Since $J'_S$ and $J''_S$ are $C^{1,\alpha}$-smooth, elliptic
regularity implies that $g$ is $C^{2,\alpha}$-smooth.

Further, it is known that the action of $\diff_+(S)$ on $\scrj_S$
admits a global finite-dimen\-sional slice. To describe this slice we recall
some standard facts from Teichm\"uller theory.

\smallskip
Denote by $\ttt_g$ the Teichm\"uller space of marked complex structures on $S$.
This is a complex manifold of dimension
\begin{equation}
\dim_\cc \ttt_g = \cases  0      &\text{ if $g=0$;}  \\
                       1      &\text{ if $g=1$;}  \\
                       3g-3   &\text{ if $g\ge2$;}
               \endcases
\end{equation}
which can be completely characterized in the following way.

\newprop{prop2.2.A1} The product
$S\times \ttt_g$ admits a (non-unique) complex (\ie holomorphic) structure
$J_{S\times \ttt}$ such that:

\sli The natural projection $\pi_\ttt: S\times \ttt_g \to \ttt_g$ is
holomorphic, so that for any $\tau\in \ttt_g$ the identification $S \cong S
\times \{\tau\}$ induces the complex structure $J_S(\tau)\deff J_{S\times
\ttt}\ogran_{S\times \{\tau\}}$ on $S$;

\slii For any complex structure $I_S$ on $S$ there exist a uniquely
defined $\tau\in\ttt_g$ and a diffeomorphism $f: S\to S$ homotopic to the
identity map $\id_S: S\to S$ such that $I_S = f^*J_S(\tau)$;

\sliii Moreover, for any finite-dimen\-sional manifold $Y$ and any smooth
map $H: Y \to \scrj_S$ there exist maps $F:Y \to \diff_+(S)$ and $h:Y \to
\ttt$ such that $H(y) = (F(y))^* h(y)$;

\sliv The group $\bfg$ of automorphisms of $S\times \ttt_g$ preserving the
projection onto $\ttt_g$ is
\begin{equation}
\bfg=\cases \mathbf{PGl}(2,\cc)       &\text{ for $g=0$,} \\
            \mathbf{Sp}(2,\;\;\zz) \ltimes T^2   &\text{ for $g=1$,} \\
            \hbox{\rm discrete}                  &\text{ for $g\ge2$;}
\endcases
\eqqno(2.2.1)
\end{equation}

\slv The tangent space to $\ttt_g$ at $\tau$ is canonically isomorphic to
$\sfh^1(S, TS)$ where $S$ is equipped with the structure $J_S(\tau)$. The
group $\sfh^0(S, TS)$ is canonically isomorphic to the Lie algebra of
$\bfg$.
\end{prop}

\smallskip
We shall assume that such a structure $J_{S \times \ttt}$ is fixed. Then we
obtain an imbedding $\ttt \hook \scrj_S$ given by $\tau \in \ttt \mapsto
J_S(\tau) \in \scrj_S$. Using it, we identify $\ttt$ with its image in
$\scrj_S$. For any $J_S \in \ttt_g$ this induces a monomorphism
$T_{J_S}\ttt_g \hook T_{J_S}\scrj_S=C^{1,\alpha}_{(0,1)}(S, TS)$. Now, the
isomorphism $T_{J_S}\ttt_g \cong \sfh^1(S, TS)$ mentioned in
\slv is obtained as the composition
\begin{equation}
T_{J_S}\ttt_g \hook T_{J_S}\scrj_S=C^{1,\alpha}_{(0,1)}(S, TS) \lrar
C^{1,\alpha}_{(0,1)}(S, TS)/\dbar\bigl(C^{2,\alpha}(S, TS)\bigr) =
\sfh^1(S, TS).
\eqqno(2.2.2)
\end{equation}

\smallskip
By our construction, any $\diff_+(S)$-orbit in $\scrj_S$ intersects
$\ttt$. This implies that instead of $\scrp/\diff_+(S)$ we can consider the
quotient $\scrp \cap (\scrs \times \ttt \times \scrj_S)$ by the action of
$\bfg$.

\newdefi{def2.2.1} Let $\whcalm \deff \scrp^* \cap (\scrs \times
\ttt \times \scrj_S)$ and use the same the notations to the restrictions of
the bundles $\scre$ and $\scre'$ onto $\whcalm$ and for the induced
operator $D: \scre\to \scre'$. The quotient $\scrm \deff \whcalm/ \bfg$ is
the {\sl total moduli space of parameterized pseudoholomorphic curves}. It
is equipped with the projection $\pr_{\!\!\scrj}:\scrm \to \scrj$. Elements of 
$\scrm$ are denoted by $[u,J]$. To indicate the surface $S$, the ambient manifold
$X$, and the homology class $[C]\in \sfh_2(X, \zz)$ involved in the definition
of $\scrm$ we shall also use the notation $\scrm(S, X, [C])$. The same meaning
have the notations $\whcalm(S, X, [C])$ and $\scrp(S, X, [C])$.
\end{defi}

\newlemma{lem2.2.1}.
\sli The projection $\wh\pr:\whcalm \lrar \scrm$ is a principal $\bfg$-bundle.

\slii The bundles\/ $\scre$ and $\scre'$ over $\whcalm$ admit a natural
$C^\ell$-smooth $\bfg$-action such that $D:\scre \to \scre'$ is
$\bfg$-invariant.

\sliii For any $J \in \scrj$ and any non-multiple $J$-holomorphic map $u: S
\to X$ there exists a diffeomorphism $\phi: S \to S$ such that $[u\scirc \phi,
J]$ lies in $\scrm$. Moreover, such element of $\scrm$ is unique.
\end{lem}

\state Remarks.~1. Part \slii of the lemma is equivalent to the existence
of $C^\ell$-smooth bundles $\scre_\scrm$ and $\scre'_\scrm$ over $\scrm$ and
a $C^{\ell-1}$-smooth bundle homomorphism $D_\scrm: \scre_\scrm \to \scre'_\scrm$
which lift to the corresponding objects over $\whcalm$. Later on we drop the
sub-index $_\scrm$, so that, \eg $\scre$ will denote also the corresponding 
bundle over $\scrm$.

\state 2. Our main interest is the space $\scrm$. However, in the proofs below
we shall mostly work with $\whcalm$. The reason is that an element $(u, J_S,
J) \in \whcalm$ fixes a parameterization of a pseudoholomorphic curve, whereas
$[u,J] \in \scrm$ defines only an appropriate equivalence class of
parameterizations.

\proof Part \slip.

{\sl Case $g\ge2$}. It is known that in this case $\bfg$ acts properly
discontinuously on $\ttt_g$. This implies that the same is true for the action
of $\bfg$ on $\whcalm$. Moreover, it is clear that $\bfg$ acts freely on $\whcalm$.
Consequently, the map $\whcalm\lrar \scrm = \whcalm/\bfg$ is simply an (unbranched)
covering.

\smallskip
{\sl Case $g=0$}. In this case $S=S^2$, $\ttt_0=\{J\st\}$, and the action of
$\bfg$ on $S$ is generated by holomorphic vector fields, \ie by the space
$\sfh^0(S,TS)$. One can show that the action of $\bfg$ on
$\whcalm$ is generated by vector fields
\begin{equation}
(u, J\st,J)\in \whcalm \mapsto
(du(v), 0,0) \in T_{(u, J\st, J)}\whcalm
\quad\hbox{with $v\in \sfh^0(S,TS)$ fixed}.
\end{equation}
In particular, the action is {\sl smooth} or, more precisely, $C^\ell$-smooth.

Consequently, for a given $(u^0,J\st,J^0)\in \whcalm$ we can find a closed
complementing space $\scrv \subset T_{(u^0,J\st,J^0)}\whcalm$ to $(u^0 (\sfh
^0(S,TS)),0,0)$. Represent it as the tangent space of a submanifold $\scrw
\subset \whcalm$ through $(u^0,J\st,J^0)$, $T_{(u^0,J\st,J^0)} \scrw=\scrv$.
If $\scrv$ is chosen sufficiently small, then it intersects every orbit $\bfg
\cdot (u,J\st,J)$ transversally in exactly one point. Moreover, we have a
$\bfg$-invariant diffeomorphism $\bfg\cdot \scrw \cong \bfg
\times \scrw$, so that $\scrw$ is a local slice of $\bfg$-action at $(u^0,
J\st,J^0)$. This equips the quotient $\whcalm/\bfg$ with a structure
of a smooth Banach manifold such that the projection $\scrw \to \scrm=
\whcalm/\bfg$ is a $C^\ell$-smooth chart.

\smallskip
{\sl Case $g=1$} is a combination of the above two cases. First we consider
the action of $T^2 \triangleleft \bfg$. The existence of a local $T^2$-slice
$\scrw$ through any given $(u^0,J_S^0, J^0)\in \whcalm$ can be
shown by copying the construction of {\sl Case $g=0$}. This implies that
$\whcalm \lrar \whcalm/T^2$ is a principle $T^2$-bundle. Then we repeat
the argument of {\sl Case $g\ge2$} and show that $\whcalm/T^2 \to
\whcalm/\bfg$ is an unbranched covering with the group $\mathbf{Sp}(2,\zz) =
\bfg/ T^2$. This completes the proof of part \slip.

\medskip {\sl Part} \sliip.
The action of $\bfg$ extends in a natural way to an action on
$\scrz\deff S \times \scrs \times \ttt \times \scrj$. The evaluation map
$\ev: \scrz \to X$, $\ev(z, u, J_S, J) \deff u(z)$ is $\bfg$-equivariant. 
Consequently, the bundle $E\deff \ev^* TX$ over $\scrz$ is equipped
with the natural $\bfg$-action. The action of $\bfg$ on $E$ induces the
actions on section spaces $\scre$ and $\scre'$. Since all constructions are
natural, $D: \scre \to \scre'$ is $\bfg$-invariant.

Finally, it remains to note that the action of $\bfg$ on the bundles $\scre$
and $\scre'$ over $\whcalm$ is $C^\ell$-smooth.

\smallskip\noindent
{\sl Part} \sliii of the lemma states the universality property of
$\scrm$ which easily follows from the definitions.
\qed

\newcorol{cor2.2.2} $\scrm$ is a $C^\ell$-smooth Banach
manifold and $\pi_{\!\!\scrj}:\scrm \to \scrj$ is a Fredholm map. For
$[u,J]\in\scrm$ there exist natural isomorphisms
\begin{align*}
\ker(d\pi_{\!\!\scrj}:T_{[u,J]}\scrm \to T_J\scrj)&
\;\cong\; \sfh^0_D(S, \scrn_u),
\\
\coker(d\pi_{\!\!\scrj}:T_{[u,J]} \scrm \to T_J\scrj)
&\;\cong\; \sfh^1_D(S, \scrn_u).
\end{align*}
In particular, the index of $\pi_{\!\!\scrj}$ is equal to
\begin{equation}
\ind_\rr \pi_{\!\!\scrj} = \chi_\rr(\scrn_u)= 2\bigl(\mu + (n-3)(1-g) \bigr),
\eqqno(2.2.3)
\end{equation}
where $\mu = \la c_1(X), [C] \ra$.
\end{corol}

\proof The $C^\ell$-smooth structure on $\scrm$ is is the quotient structure 
defined by the the $C^\ell$-smooth $\bfg$-action on $\whcalm$.

\smallskip
Using \refcorol{cor2.1.3} we see that the tangent space
to $\whcalm$ is
\begin{equation}
T_{(u,J_S,J)} \whcalm = \bigl\{ (v,\dot J_S, \dot J) \,:\,
\dot J_S \in T_{J_S}\ttt,\;\,
D_{u,J}v + J \scirc du \scirc \dot J_S + \dot J \scirc du \scirc J_S
=0 \,\bigr\}.
\eqqno(2.2.4)
\end{equation}
Consider the natural projection $\pi_\scrp: \scrp^* \to \scrj$,  $(u,J_S,J)
\mapsto J$ with the differential $d\pi_\scrp: T_{(u,J_S,J)}\scrp^* \to
T_J\scrj$ given by $(v,\dot J_S,\dot J)\in T_{(u,J_S,J)}\scrp^*\mapsto
\dot J\in  T_J\scrj$.

The kernel $\ker(d\pi_\scrp)$ consists of solutions $v\in \scre_{(u,J)}$
of the equation
\begin{equation}
D_{u,J}v +  J\scirc du \scirc \dot J_S =0
\end{equation}
with $\dot J_S\in T_{J_S}\ttt$. Since the map
$\wh\pi: \wh\scrm\to \scrm$ is a principle $\bfg$-bundle, the kernel
$\ker( d\pi_{\!\!\scrj}:T_{(M,J)}\scrm \to T_J\scrj)$ is obtained from $\ker(d\pi)$
by taking the quotient by the tangent space to the fiber $\bfg
\cdot (u,J_S,J)$ which is equal to $du(\sfh^0(S,TS))$. Using the relations
$\sfh^0(S,TS)= \ker(\dbar_{TS}:L^{1,p}(S,TS)\to L^p(S,TS\otimes
\Lambda^{(0,1)}S)$, $T_{J_S}\ttt_g\cong\sfh^1(S,TS)=\coker(\dbar_{TS})$, and
$du\scirc \dbar_{TS}= D_{(u,J)}\scirc du$, we conclude that the space
$\ker(d\pi_{\!\!\scrj})$ is isomorphic to the quotient
\begin{equation}
\{v\in L^{1,p}(S,E_u)\;:\; Dv=du(\phi) \hbox{ for some }\phi\in
L^p(S,TS\otimes\Lambda^{(0,1)}S)\}
\!\!\bigm/\!\!
du\bigl(L^{1,p}(S,TS)\bigr).
\end{equation}
Hence, by \lemma{lem1.5.3}, $\ker(d\pi_{\!\!\scrj}:T_{[u,J]}\scrm\to T_J\scrj)
\cong \sfh^0_D(S,\scrn_u)$. In particular, $\ker(d\pi_{\!\!\scrj})$ is finite 
dimensional.

\smallskip
Similarly, the image of $d\pi_{\!\!\scrj}$ consists of those $\dot J$ for which
the equation
\begin{equation}
D_{u,J}v + J\scirc du \scirc \dot J_S + \dot J\scirc
du \scirc J_S =0
\end{equation}
has a solution $(v,\dot J_S)$. Using \eqqref(1.5.6) and \lemma{lem1.5.3}
we obtain the relations $\im(d\pi_{\!\!\scrj})=\ker\barr\Psi$ and
$\coker(d\pi) \cong \sfh^1_D(S,N_u)$.

This implies the Fredholm property for the projection $\pi_{\!\!\scrj}:\scrm \to
\scrj$ and the formula $\ind(d\pi_{\!\!\scrj}) = \ind (\scrn_u)$.
\qed

\medskip
We conclude the paragraph with a description of the deformations of non-closed 
pseudoholomorphic curves. 

\newdefi{def2.2.2} Let $\barr S=S \cup \d S$ be a compact {\sl non-closed} 
oriented surface with the boundary $\d S$ consisting of finitely many circles.
Denote by $\scrj_S$ the (Banach) space of complex structure on $S$ which are 
compatible with the orientation of $S$ and $C^\ell$-smooth up to boundary 
$\d S$. As usual let 
\begin{equation}
\eqqno(2.2.11)
\scrp(S, X) \deff \{ (u,J_S,J) \in L^{1,p}(S, X) \times \scrj_S \times 
\scrj: \dbar_{J_S,J} u=0\},
\end{equation}
the space of pseudoholomorphic maps equipped it with the natural projections
$\pr_{\scrj_S}: \scrp(S, X) \to \scrj_S$ and $\pr_{\!\!\scrj}: \scrp(S, X) \to 
\scrj$. The fibers of the projections are denoted by $\scrp(S, X, J)=\pr_{\!\!\scrj}
\inv(J)$, $\scrp(S,J_S, X)=\pr_{\scrj_S}\inv(J_S)$, and $\scrp(S,J_S, X,J)= 
\scrp(S,J_S, X) \cap \scrp(S, X, J)$ respectively.
\end{defi}

\newlemma{lem2.2.4} \sli Let $S$ be a non-closed oriented surface. Then

\sli the space $\scrp(S, X)$ is a Banach submanifolds of $L^{1,p}(S, X) 
\times \scrj_S\times \scrj$;

\slii For any $(u,J_S, J)\in \scrp(S, X)$, the operators $d\pr_{\!\!\scrj}: T_{(u,
J_S,J)}\scrp(S, X) \to T_J\scrj$ and $d\pr_{\scrj_S}: T_{(u,J_S,J)} \scrp(S, 
X) \to T_{J_S}\scrj$ are surjective and split. In particular, $\scrp(S; X,J)$ 
and $\scrp(S,J_S; X)$ are are Banach submanifolds of $\scrp(S, X)$;

\sliii For any $(u,J_S, J)\in \scrp(S, X)$, the submanifolds $\scrp(S; X,J)$ 
and $\scrp(S,J_S; X)$ are transversal in $(u,J_S,J)$; in particular, $\scrp(S, J_S; 
X,J)= \scrp(S; X,J) \cap \scrp(S,J_S; X)$ is also a Banach submanifold;

\sliv The the tangent spaces are given by
\begin{align}
T_{(u,J_S,J)}\scrp(S, X) &= \bigl\{ (v,\dot J_S,\dot J) \in T_uL^{1,p}(S, X)
\times T_{J_S}\scrj_S \times T_J\scrj : 
\notag\\
&\qquad\quad D_{u,J}v + \dot J \scirc du \scirc J_S + 
J \scirc du \scirc \dot J_S =0\bigr\};
\eqqno(2.2.12)\\
T_{(u,J_S)}\scrp(S; X,J)& = \bigl\{ (v,\dot J_S,\dot J) \in 
T_{(u,J_S,J)}\scrp(S, X): \dot J=0 \bigr\};
\eqqno(2.2.13)\\
T_{(u,J)}\scrp(S,J_S; X)& = \bigl\{ (v,\dot J_S,\dot J) \in 
T_{(u,J_S,J)}\scrp(S, X): \dot J_S=0 \bigr\};
\eqqno(1.2.14)\\ 
T_u\scrp(S,J_S; X,J)& = \bigl\{ (v,\dot J_S,\dot J) \in 
T_{(u,J_S,J)}\scrp(S, X): \dot J_S=0 = \dot J \bigr\}.
\eqqno(2.2.15) 
\end{align}
\end{lem}

\proof The lemma is obtained by the transversality techniques of 
\refsubsection{2.1} using the following claim: {\sl For any $(u,J_S,J) \in 
\scrp(S,X)$ the operator $D_{u,J} : L^{1,p}(C, E_u) \to L^p_{(0,1)}(C, E_u)$ 
is surjective}. Since $D_{u,J}$ is elliptic, this is a standard fact following 
from compactness and non-closedness of $S$. \qed

\newsubsection[2.3]{Transversality I$\!$I}
Before stating further results, we introduce some new notation. Here $S$ is a 
{\sl closed} real surface.

\newdefi{def2.3.1}
Let $Y$ be a $C^\ell$-smooth finite-dimen\-sional manifold, possibly with
non-empty $C^\ell$-smooth boundary $\d Y$, and $h:Y \to \scrj$ a $C^\ell
$-smooth map. Define the {\sl relative Moduli space}
\begin{equation}
\scrm_h \deff Y\times_{\!\!\scrj} \scrm \cong
\{\,(u,J_S,y)\in \scrs\times \ttt_g\times Y\,:\, (u,J_S, h(y))\in \scrp^*
\,\}/\bfg
\end{equation}
with the natural projection $\pi_h:\scrm_h \to Y$. In the special case
$Y=\{J\}\hook \scrj$, we obtain the Moduli space of $J$-holomorphic curves
$\scrm_J\deff \pi_{\!\!\scrj}\inv(J)$. The projection $\pi_h:\scrm_h \to Y$
is a fibration with a fiber $\pi_h\inv(y)=\scrm_{h(y)}$. We shall denote
elements of $\scrm_h$ by $[u,y]$, where $u:S\to X$ is a $h(y)$-holomorphic map.
\end{defi}

\medskip
The next two lemmas follow from the transversality theory.

\newlemma{lem2.3.1}
Let $Y$ be a $C^\ell$-smooth finite-dimen\-sional
manifold, and $h:Y \to \scrj$ a $C^\ell$-smooth map. Then $\scrm_h$ is a
$C^\ell$-smooth manifold in some neighborhood of a point $[u,y] \in
\scrm_h$ with $J \deff h(y)$ \iff the map $\barr\Psi_{u, J} \scirc dh
: T_uY \to \sfh^1_D(S, N_u)$ is surjective.
In this case the tangent space to $\scrm_h$ is
\begin{equation}
T_{[u,y]}\scrm_h=\ker\bigl(D\;\oplus\; \Psi\scirc dh:
\scre_{u,h(y)}\oplus T_yY \lrar \scre'_{u,h(y)} \bigr)
\bigm/ du(\sfh^0(S, TS))
\eqqno(2.3.1)
\end{equation}
\end{lem}

\proof We reformulate the transversality condition and
use \lemma{lem2.1.1}.
\qed

\newlemma{lem2.3.2}
\sli There exists a Baire subset $\scrj\reg \subset
\scrj$ such that any $J \in \scrj$ is a regular value of $\pi_{\!\!\scrj} :\scrm
\to \scrj$.

\slii There exists a Baire subset $\scrv$ in the space $C^\ell([0,1], \;
\scrj)$, such that any map $h:[0,1]\to \scrj$ from $\scrv$ is transversal to
$\pi_{\!\!\scrj} :\scrm \to \scrj$ and both $h(0)$ and $h(1)$ are
regular values of $\pi_{\!\!\scrj}$.
\end{lem}

\state Remark. In general, for any finite-dimen\-sional manifold $Y$ with
boundary $\d Y$ there exists a Baire subset $\scrv \subset C^\ell(Y,\; \scrj)$
such that any $h \in \scrv$, as well as its restriction $h\ogran_{\d Y}$ are
transversal to $\pi_{\!\!\scrj}$. The proof uses the Sard lemma.

\newlemma{lem2.3.3}
Suppose that $S$ is the sphere $S^2$ and $\dim_\rr(X) =4$.
Then the exists a {\sl connected} Baire subset $\scrj\reg \subset
\scrj$ such that any $J \in \scrj$ is a regular value of $\pi_{\!\!\scrj} :\scrm
\to \scrj$. Moreover, any $J_0, J_1 \in \scrj\reg$ can be connected by
a smooth path $h: [0,1] \to \scrj\reg$.
\end{lem}

\proof By \lemma{lem2.3.2}, there exists a Baire subset $\scrj\reg
\subset\scrj$ such that any $J \in \scrj\reg$ is a regular value of
$\pi_{\!\!\scrj}$. Further, any $J_0, J_1 \in \scrj\reg$ can be connected by a
smooth path $h: [0,1] \to \scrj$, transversal to $\pi_{\!\!\scrj}$.

For any such path $h: [0,1] \to \scrj$ and any $[u,t] \in \scrm_h$ with
$h(t)=J$ the map $\barr\Psi_{u,J}\scirc dh: T_t[0,1]\cong\rr \; \lrar
\sfh^1(S, N_u)$ is surjective by \lemma{lem2.3.1}. Consequently,
for such $u$ and $J$ we have $\dimr \sfh^1_D(S, \scrn_u)= \dimr \sfh^1_D(S,
N_u) \le 1$.

Recall that the difference $\dimr\sfh^0_D(S, N_u) -\dimr\sfh^1_D(S, N_u)$ is
{\sl even}, see \eqqref(1.5.?4). Hence, if $\dimr \sfh^1(S, N_u)= 1$ then
$\sfh^1_D(S, N_u)$ should also be nontrivial. On the other hand, the condition
$\dimr(X)=4$ implies that $N_u$ is a {\sl line bundle}. But, in view of {\sl
Lemma 1.5.2}, on the sphere $S=S^2$ one of the spaces $\sfh^i_D(S, N_u)$ must
be trivial.

Thus, we see that for any such path $h$ and any $[u,t]\in \scrm_h$ one has
$\sfh^1_D(S, \scrn_u)=0$. This means that $h$ takes values in $\scrj\reg$.
\qed

\smallskip
For higher genus $g(S) \ge1$ we have a similar, but weaker result.

\newlemma{lem2.3.4}
Assume that $\dim_\rr X =4$ and $g\deff g(S) \ge1$
and set $\mu \deff c_1(X)[u(S)]$. Then
\begin{equation}
\mu  \le \dim_\cc \sfh^0(S, \scrn\sing_u) \le \mu + g-1
\eqqno(2.3.2)
\end{equation}
for any $[u,J] \in \scrm$ with $\sfh^1(S, N_u)\cong \rr$.
\end{lem}

\proof The condition $\dim_\rr X =4$ means that $N_u$ is a {\sl line} bundle. 
Thus, by \lemma{lem1.5.2},  $\sfh^1(S, N_u)\cong \rr$ implies $c_1(N_u) \le 2g 
-2$. Further, since $\dim_\rr \sfh^1(S, N_u)=1$ and $\ind_\rr D^N_{u,J} = \dimr
\sfh^0_D(S, N_u) -\dimr\sfh^1_D(S, N_u)$ is {\sl even}, see \eqqref(1.5.?4), we 
conclude that $\dimr\sfh^0_D(S, N_u) \ge 1$. Consequently, $\ind_\rr D^N_{u,J} 
\ge 0$. The formula \eqqref(1.5.?4) for $\ind_\rr D^N_{u,J}$ yields $c_1(N_u) 
\ge g-1$. Finally, the definition of $\scrn_u$ yields the relation $\mu= c_1(X)
 [u(S)] = c_1(E_u) = c_1(TS) + c_1(N_u) + \dimc \sfh^0(S, \scrn\sing_u)$. \qed

\newsubsection[2.4]{Pseudoholomorphic curves through fixed points} In this 
paragraph we consider the total moduli space of pseudoholomorphic curves passing 
through given fixed points $x_1, \ldots, x_m \allowbreak 
\in X$. Gromov in \cite{Gro} proposed a method to reduce 
this problem to the one of pseudoholomorphic curves
without such constraints. The idea is to blow up $X$ in the points $x_1, \ldots, 
x_m$ and consider the curves on the blown-up space $\ti X$. He has shown that a 
$C^\ell$-smooth almost complex structure $J$ on $X$ lifts to a $C^{\ell-1}$-smooth 
almost complex structure $\ti J$ on $\ti X$ such that the natural
projection $\pr: \ti X \to X$ is holomorphic and such that every $J$-holomorphic 
curve $C$ in $X$ passing through $x_1, \ldots, x_m$ lifts to a unique $\ti 
J$-holomorphic curve $\ti C$ in $\ti X$ with $C= \pr(\ti C)$.

\smallskip
Our aim here is to make an explicit construction for the moduli space of 
pseudoholomorphic curves passing through fixed points. Since the construction is 
simply a modification of the case $m=0$ where no points are marked, we shall 
mostly skip or merely indicate proof of claims.

We begin by introducing some notation. Denote by $\mbfx=(x_1, \ldots x_m)$ the
tuple of fixed points on $X$, which are supposed to be pairwise distinct. 
Also fix a tuple $\mbfz=(z_1, \ldots, z_m)$ of pairwise distinct points on the
surface $S$. Define 
$$
\scrs(\mbfz, \mbfx) \deff \{ u \in \scrs=L^{1,p}(S, X): u(z_i)=x_i \};
$$
$$
\scrp(S,\mbfz; X, \mbfx) \deff \{ (u, J_S, J) \in \scrp: u\in \scrs(\mbfz, 
\mbfx) \};
$$
$$
\scrp^*(S,\mbfz; X, \mbfx) \deff  \scrp(S,\mbfz; X, \mbfx) \cap \scrp^*(S,X).
$$
The linearization of the conditions $u(z_i)=x_i$ yields the equations $v(z_i)=0$
for $v \in T_u\scrs= L^{1,p}(S, E_u)$. Denote as above $E=E_u \deff u^*TX$, and
set $E_i= E_{u,i} \deff (E_u)_{z_i} = T_{u(z_i)}X$. Then we obtain the bundle
$E_\mbfz$ over $\scrs$ with a fiber $(E_\mbfz)_u\deff \oplus E_{u,i}$ equipped with
the natural {\sl evaluation homomorphism} $\ev_\mbfz: \scre \to E_\mbfz$
$$
\ev_{\mbfz}: v \in \scre_u=L^{1,p}(S, E) \mapsto (v(z_1), \ldots, v(z_m)) \in 
E_\mbfz.
$$
It is easy to see that $\ev_\mbfz: \scre \to E_\mbfz$ is surjective. This means
that the equations $u(z_i)=x_i$ are transversal and implies that $\scrs(\mbfz, 
\mbfx)$ is a Banach submanifold of $\scrs$ with the tangent space
\begin{equation}\eqqno(2.4.1)
T_u\scrs(\mbfz, \mbfx) = \{ v \in T_u\scrs= L^{1,p}(S, E_u): v(z_i)=0 \}.
\end{equation}
The same argument shows that $\scrp^*(S,\mbfz; X, \mbfx)$ is also a Banach 
submanifold of $\scrp^*(S,X)$ with the tangent space
$$
T_u\scrp^*(S,\mbfz; X, \mbfx) =  \{ (v, \dot J_S, \dot J) \in  
T_u\scrp^*(S,X): v(z_i)=0 \}.
$$

\smallskip
Let $\diff_+(S, \mbfz)$ be the subgroup of those $g\in \diff_+(S)$ which fix the 
marked points $z_1,\ldots, z_m$. Then $\diff_+(S, \mbfz)$ leaves the subsets 
$\scrs(\mbfz, \mbfx) \subset \scrs$ and $\scrp^*(S,\mbfz; X, \mbfx) \subset 
\scrp(S,X)$ invariant. So we can define the {\sl total moduli space of 
pseudoholomorphic curves through the given points $x_1, \ldots, x_m$} as the 
quotient $\scrm(\mbfx) \deff \scrp^*(S,\mbfz; X, \mbfx) /\diff_+(S, \mbfz)$. This
space is equipped with the natural projection $\pi_\scrj: \scrm(\mbfx) \to \scrj$ 
defined in an obvious way.

\smallskip
The smooth structure on $\scrm(\mbfx)$ is constructed in the same way as it was 
for $\scrm$. First, one constructs a global slice on the action of $\diff_+
 (S, \mbfz)$ on $\scrj_S$. To do this, we consider the action of the component 
of $\diff_0(S, \mbfz)$ the group $\diff_+(S, \mbfz)$ containing the identity. The 
quotient $\scrj_S/ \diff_0(S, \mbfz)$ is the {\sl Teichm\"uller space $\ttt
_{g,m}$ of complex structures on a Riemann surface of genus $g=g(S)$ with $m$ 
punctures}. The marked points $z_1,\ldots, z_m$ are the positions of punctures. 

As in the case $m=0$, one can imbed $\ttt_{g,m}$ in $\scrj_S$ in such a way that
the composition $\ttt_{g,m} \hook \scrj_S \twoheadrightarrow \scrj_S/ \diff_0(S, 
\mbfz) \cong \ttt_{g,m}$ is the identity map. This imbedding $\ttt_{g,m} \hook 
\scrj_S$ is the desired slice. The choice of such imbedding $\ttt_{g,m} \hook 
\scrj_S$ is equivalent to the choice of the complex structure $J_{\ttt, S, \mbfz}$
on the product $\ttt_{g,m} \times S$. One considers $\ttt_{g,m} \times S$ with
this complex structure and with the holomorphic projection $\pr: \ttt_{g,m} 
\times S \to \ttt_{g,m}$ as the universal family corresponding to $\ttt_{g,m}$.

After the choice of the slice $\ttt_{g,m} \hook \scrj_S$, $\scrm(\mbfx)$ is
obtained as the quotient of the space 
$$
\wh\scrm(\mbfx) \deff \{ (u, J_S, J) \in \scrp^*(S,\mbfz; X, \mbfx):
J_S \in \ttt_{g,m} \}
$$
by the group $\bfg$ of biholomorphisms of $\ttt_{g,m} \times S$ preserving the 
projection $\pr: \ttt_{g,m} \times S \to \ttt_{g,m}$. It is discrete except the
cases where $g=0$ and $m=1$ or $2$. In the case $m=1$, $S\bs \{z_1\}$ is the 
complex plane $\cc$ and $\bfg$ is its automorphism group $\cc^* \ltimes \cc$. 
Similarly, in  the case $m=2$, $S\bs \{z_1, z_2\}$ is the punctured complex plane
$\cc^*$ and $\bfg= \zz_2 \ltimes \cc^*$ is likewise its automorphism group. In 
either case one can construct a local slice for the action of $\bfg$ on $\wh\scrm
 (\mbfx)$ by repeating the arguments of \refsubsection{2.2}. So the quotient 
$\whcalm(\mbfx)/\bfg$ is a $C^\ell$-smooth Banach manifold.

\medskip
Now we define the notion of the normal sheaf of a pseudoholomorphic curve passing 
through fixed points on $X$. In this new situation, the linearization of 
$\dbar$-equations leads to the operator
\begin{equation}\eqqno(2.4.1a)
D=D_{u,J}: \{ v \in L^{1,p}(S,E_u): v(z_i)=0 \text{ for }i=1,\ldots,m\} 
\to L^p_{(0,1)}(S, E),
\end{equation}
which is the usual Gromov operator $D=D_{u,J}$, but now considered with 
a new domain of definition
\begin{equation*}
\scre_{u,\mbfx} \deff \{ v \in L^{1,p}(S,E_u): v(z_i)=0  
\text{ for }i=1,\ldots,m\}.
\end{equation*}
The space $\scre_{u,\mbfx}$ is the kernel of the evaluation homomorphism
$\ev_z: \scre_u \to E_\mbfz$ and is the tangent plane to $\scrs(\mbfz, \mbfx)$,
see \eqqref(2.4.1). 

We now describe the structure of the operator \eqqref(2.4.1a). Recall that we have
the 
decomposition $D_{u,J}= \dbar_{u,J} + R_{u,J}$, see \refsubsection{1.4}. Observe
that the sheaf $\scro(E_u)[-\mbfz]$ of holomorphic sections of $\scro(E_u)$
vanishing at the points $z_1, \ldots, z_m \in S$ is locally free and hence 
corresponds to a holomorphic bundle. Let us denote this bundle by $E_{u, -\mbfz}$.

\newlemma{lem2.4.1} \sli The (co)kernel of the operator 
\begin{equation}\eqqno(2.4.2a)
\dbar_{u,J}: \{ v \in L^{1,p}(S,E_u): v(z_i)=0 \text{ for }i=1,\ldots,m\} 
\to L^p_{(0,1)}(S, E),
\end{equation}
is canonically isomorphic to the cohomology groups $\sfh^0_\dbar(S, E_{u,-\mbfz})$
and $\sfh^1_\dbar(S, E_{u,-\mbfz})$. 

\slii The operator $D_{u,J}$ induces the operator
$$
D_{u,-\mbfz, J}: L^{1,p}(S, E_{u,-\mbfz}) \to L^p_{(0,1)}(S, E_{u,-\mbfz})
$$
which is of the form $D_{u,-\mbfz, J}= \dbar_{u,-\mbfz, J} + R_{u,-\mbfz, J}$, where
$\dbar_{u,-\mbfz, J}$ is the Cauchy-Riemann operator corresponding to the natural 
holomorphic structure in $E_{u,-\mbfz}$ and $R_{u,-\mbfz, J}$ is a $\cc$-antilinear
$L^\infty$-bounded bundle homomorphism, \ie 
$$
R_{u,-\mbfz, J} \in L^\infty\bigl(S, \barr\hom_\cc(E_{u,-\mbfz},
 E_{u,-\mbfz} \otimes \Lambda^{(0,1)}S)\bigr).
$$

\sliii The (co)kernel of the operator 
\begin{equation}\eqqno(2.4.3a)
D_{u,J}: \{ v \in L^{1,p}(S,E_u): v(z_i)=0 \text{ for }i=1,\ldots,m\} 
\to L^p_{(0,1)}(S, E)
\end{equation}
is canonically isomorphic to the cohomology groups $\sfh^0_D(S, E_{u,-\mbfz})$
and $\sfh^1_D(S, E_{u,-\mbfz})$ corresponding to the operator $D_{u,-\mbfz, J}$.
\end{lem}

\proof Fix local holomorphic coordinates $\zeta_i$ on $S$, each centered at
the corresponding marked point $z_i$. Consider the natural inclusions 
\begin{align}
\eqqno(2.4.4a)
j^0:L^{1,p}(S, E_{u,-\mbfz}) &\hook 
\{ v \in L^{1,p}(S,E_u): v(z_i)=0 \text{ for }i=1,\ldots,m\},
\\
\eqqno(2.4.5a)
j^1: L^p_{(0,1)}(S, E_{u,-\mbfz}) &\hook L^p_{(0,1)}(S, E_u).
\end{align}
Observe that $v\in L^{1,p}(S,E_u)$ with $v(z_i)=0$ belongs to $L^p_{(0,1)}(S, 
E_{u,-\mbfz})$ \iff locally near every $z_i$ it has the form $v(\zeta_i)= 
\zeta_i w(\zeta_i)$ for some (uniquely defined!) $L^{1,p}$-section $w(\zeta_i)$ 
of $E_{u,-\mbfz}$. This is equivalent to the condition $\zeta_i\inv \dbar_{u,J} 
v(\zeta_i) \in L^p$ as well as to the condition $\zeta_i\inv D_{u,J} v(\zeta_i) 
\in L^p$. Consequently, $D_{u,J}$ restricted to $L^{1,p}(S,E_{u,-\mbfz})$ takes 
values in $L^p_{(0,1)}(S, E_{u,-\mbfz})$. This yields the operator $D_{u,-\mbfz, 
J}$. Moreover,  $D_{u,-\mbfz, J}$ is of order 1 and has the Cauchy-Riemann symbol.
Consequently, it has the form $D_{u,-\mbfz, J}= \dbar_{u,-\mbfz, J} +R_{u,-\mbfz,
J}$, where $\dbar_{u,-\mbfz, J}$ is the Cauchy-Riemann operator corresponding
to the holomorphic structure in $E_{u,-\mbfz}$, and $R_{u,-\mbfz, J}$ is the
$\cc$-antilinear part of $D_{u,-\mbfz, J}$.

Let $v_1(\zeta_i), \ldots, v_n(\zeta_i)$ be a local holomorphic frame of $E_u$ in 
a neighborhood of $z_i$ and $R_{\alpha\beta}(\zeta_i)$ the matrix of $R_{u,J}$ in 
this frame. Then $\zeta_iv_1(\zeta_i), \ldots, \zeta_iv_n(\zeta_i)$ is a local 
frame of $E_{u,-\mbfz}$. From $\cc$-antilinearity of $R_{u,J}$ we obtain
$$
\textstyle
R_{u,J}(\zeta_i\,v_\alpha(\zeta_i))= 
\sum_\beta R_{\alpha\beta}(\zeta_i)\bar \zeta_i v_\beta(\zeta_i)=
\sum_\beta \msmall{\bar \zeta_i \over \zeta_i} R_{\alpha\beta}(\zeta_i)
\cdot \zeta_i v_\beta(\zeta_i).
$$
This shows that ${\bar \zeta_i \over \zeta_i} R_{\alpha\beta}(\zeta_i)$ is the 
matrix of $R_{u, -\mbfz, J}$ in the frame $\zeta_iv_1(\zeta_i), \ldots, \zeta_i
v_n(\zeta_i)$. Recall that $R_{u,J}$ is a continuous bundle homomorphism (see 
\lemma{lem1.4.1}, \slip). So we see that $R_{u, -\mbfz, J}$ is also continuous 
outside the marked points $z_i$ and has singularities of the form ${\bar 
\zeta_i \over \zeta_i} R_{u,J}$ at $z_i$. In particular, $R_{u, -\mbfz, J}$ is of
type $L^\infty$, but is not continuous in general.

\smallskip
The equality of the kernels of the operators in \sli and \slii with the
corresponding 0-cohomology groups follows directly from the definition of the 
operators $\dbar_{u,-\mbfz,J}$ and $D_{u,-\mbfz,J}$. The equality for 
1-cohomology groups will be shown only for the operator $D_{u,-\mbfz,J}$, the 
other one is carried out in the same manner. So let $\phi \in L^p_{(0,1)}(S, E_{u, 
-\mbfz})$. If $\phi = D_{u,-\mbfz,J}(v)$ for $v \in L^{1,p}(S, E_{u, -\mbfz})$,
then $v \in L^{1,p}(S, E_u)$ and $j^1\phi = D_{u,J}(v)$, or more precisely $j^1
\phi = D_{u,J}(j^0(v))$. This shows that the inclusion $j^1$ in \eqqref(2.4.5a) 
induces a well-defined homomorphism from $\sfh^1_D(S, E_{u, -\mbfz})$ to the 
cokernel of \eqqref(2.4.3a). Moreover, $\phi \in L^p_{(0,1)}(S, E_{u, -\mbfz})$ 
induces the zero class in the cokernel of \eqqref(2.4.3a) \iff $j^1(\phi) = 
D_{u,J}(v)$ for some $v \in L^{1,p}(S, E_u)$ with $v(z_i)=0$. But then locally
$$
D_{u,J}(\zeta_i\inv v(\zeta_i)) = 
\zeta_i\inv D_{u,-\mbfz,J}(v(\zeta_i)) =
\zeta_i\inv j^1 \phi(\zeta_i) \in L^p
$$
by the definition of the inclusion $E_{u, -\mbfz} \hook E_u$. This implies
that $v \in L^{1,p}(S, E_{u-\mbfz})$. Thus the homomorphism induces by $j^1$
is injective.

Further, for any $\phi \in L^p_{(0,1)}(S, E_u)$ there exists $v \in L^{1,p}(S, 
E_u)$ which vanishes at all $z_i$ and solves the equation $\phi = D_{u,J}(v)$
in a neighborhood of every $z_i$. Then $\phi - D_{u,J}(v)$ represents the same
class in the cokernel of \eqqref(2.4.3a) and is of the form $\phi - D_{u,J}(v)
= j^1(\psi)$ for some $\psi \in L^p_{(0,1)}(S, E_{u, -\mbfz})$. This finishes
the proof of the claim \sliiip. \qed

\medskip
Now we define the normal sheaf of a pseudoholomorphic curve passing through fixed 
points. The construction is completely analogous to that in the  case of no fixed
points. Here, instead of the tangent bundle $TS$ we use the bundle related 
the new situation. This is the bundle $TS_{-\mbfz}$ associated to the locally
free coherent sheaf $\scro(TS)[-\mbfz]$ of local holomorphic sections of $TS$ 
vanishing at the points $z_i$. One can prove the analog of \lemma{lem2.4.1}
for $TS_{-\mbfz}$. Observe however, that such a result follows immediately
from that lemma if we set $X=S$ and $u=\id_S$.

As in \refsubsection{1.5}, we obtain the sheaf homomorphism $du: \scro(TS
_{-\mbfz}) \to \scro(E_{u, -\mbfz})$, which is injective for non-constant $u:
S \to X$. Now, the {\sl normal sheaf to curve $C=u(S)$ passing through the 
points $\mbfx=(x_1,\ldots, x_m)$} is defined as the quotient $\scrn_{u,\mbfx} 
\deff \scro(E_{u, -\mbfz}) / du\bigl(\scro(TS_{-\mbfz}) \bigr)$ together with
the exact sequence
\begin{equation}\eqqno(2.4.7)
0 \lrar \scro(TS_{-\mbfz}) \buildrel du \over \lrar \scro(E_{u, -\mbfz}) 
\lrar \scrn_{u,\mbfx} \lrar 0.
\end{equation}
The sheaf $\scrn_{u,\mbfx}$ can be decomposed into its {\sl regular part} $\scrn
\reg_{u,\mbfx}$ and its {\sl singular part} $\scrn\sing_{u,\mbfx}$, where 
$\scrn \reg_{u,\mbfx}$ is locally free and $\scrn\sing_{u, \mbfx}$ is a torsion 
sheaf. Then $\scrn \reg_{u,\mbfx}$ is a sheaf of local holomorphic sections of 
the {\sl normal bundle $N_{u,\mbfx}$ to curve $C=u(S)$ passing through the points
$\mbfx=(x_1,\ldots, x_m)$}, so that $\scrn\reg_{u,\mbfx}= \scro(N_{u,\mbfx})$.

As in \refsubsection{1.5}, we also obtain the exact sequence
\begin{equation}\eqqno(2.4.8)
0 \lrar \scro(TS_{-\mbfz}) \otimes \scro([A])
\buildrel du \over \lrar \scro(E_{u, -\mbfz}) 
\scro(N_{u,\mbfx})
\lrar 0.
\end{equation}
where $[A]$ is the branching divisor of $du$ (see \refdefi{def1.5.1}). This
implies that the regular part $\scro(N_{u,\mbfx})$ is the quotient
$$
\scro(N_{u,\mbfx})= 
\scro(E_{u,-\mbfz}) / du\bigl(\scro(TS_{-\mbfz}) \otimes \scro([A])\bigr).
$$
From the definition of $E_{u,-\mbfz}$ and $TS_{-\mbfz}$ we obtain the isomorphism
$$
\scro(N_{u,\mbfx})\cong \scro(N_u) \otimes \scro([A]).
$$
On the other hand, the singular part remains the same as is the case without 
constraints:
$$
\scrn\sing_{u,\mbfx}\cong \scrn\sing_u \cong \scro /\scro(-[A]).
$$

Further, we observe that the operators $\dbar$ on $TS_{-\mbfz}$ and $D_{u,-\mbfz,
J}$ in $E_{u,-\mbfz}$ commute with the homomorphism $du: TS_{-\mbfz} \to E_{u,
-\mbfz}$. Consequently, $D_{u,-\mbfz, J}$ induces the operator 
$$
D^N_{u,-\mbfz, J}: L^{1,p}(S, N_{u, -\mbfx}) \to L^p_{0,1)}(S, N_{u, -\mbfx}) 
$$
with the properties similar to ones of \eqqref(1.5.3). Further, as in 
\lemma{lem1.5.3} and \refcorol{cor1.5.4} we obtain a long exact sequence of
$D$-cohomologies. 

\newprop{prop2.4.2} The~short exact sequence $\eqqref(2.4.7)$ induces
the~long exact sequence of $D$-cohomologies
\begin{equation*}
\def\normalbaselines{\baselineskip20pt\lineskip3pt \lineskiplimit3pt }
\def\mapright#1{\smash{\mathop{\longrightarrow}\limits^{#1}}}
\def\mapdown{\Big\downarrow}
\matrix
0& \mapright{}& \sfh^0(S, TS_{-\mbfz}) &\mapright{} 
 & \sfh^0_D(S, E_{u,-\mbfz})
 & \mapright{}& \sfh^0_D(S, N_{u,-\mbfz})\oplus \sfh^0(S, \scrn_u\sing)
 &\mapright{\delta} &\vphantom{\mapdown}\\
 & \mapright{}& \sfh^1(S, TS_{-\mbfz}) &\mapright{} 
 & \sfh^1_D(S, E_{u,-\mbfz})
 & \mapright{}& \sfh^1_D(S, N_{u,-\mbfz})         &\mapright{} &0.
\endmatrix
\end{equation*}
\end{prop}

\medskip
Finally, we note that the results of {\sl Paragraphs \ref{sec:2.2}\/} and
{\sl\ref{sec:2.3}\/} remain valid, after an appropriate modification,
also for curves passing through fixed points. We state without the proof 
the summary of results which will be used later.

\newthm{thm2.4.3} \sli The total moduli space $\scrm_\mbfx$ of pseudoholomorphic 
curves in a given homology class $[C] \in \sfh_2(X, \zz)$ passing through fixed 
pairwise distinct points $\mbfx=(x_1, \ldots, \allowbreak
x_m)$ on $X$ is a $C^\ell$-smooth 
Banach submanifold of $\scrm$ of real codimension $2m$. In particular, the 
projection $\pi_{\!\!\scrj}: \scrm_\mbfx \to \scrj$ is a $C^\ell$-smooth Fredholm 
map of index 
$$
2(c_1(X)[C] + (n-3)(1-g) -m).
$$

\slii For a generic $J \in \scrj$ and a generic $C^\ell$-smooth path $h:[0,1] \to
\scrj$ the fiber
$$
\scrm_{J,\mbfx} \deff \pi_\scrj\inv(J)
$$
and the {\sl relative moduli space}
$$
\scrm_{h,\mbfx} \deff [0,1] \times_\scrj \scrm_\mbfx 
$$
are $C^\ell$-smooth manifolds of expected dimension $2(c_1(X)[C] + (n-3)(1-g) -m)$
and $2(c_1(X)[C] + (n-3)(1-g) -m) +1$ respectively.
\end{thm}

\medskip
\newsection[4]{Cusp-curves in the moduli space.}
In this section we study the problem of deformation of pseudoholomorphic curves
with prescribed singularities and develop the techniques required for controlling 
their singularities under deformation. As the main result of this section we show 
that the locus of pseudoholomorphic curves with a prescribed type of singularity 
is a smooth 
Banach submanifold of expected codimension in the total moduli space of
pseudoholomorphic curves. This improvement of the result of 
Micallef and White (see \lemma{lem1.2.1}) plays a crucial role below in 
\refsection{3} in the proof of the {\sl saddle point property}.

Recall that our moduli space $\scrm$ consists of parameterized non-multiple
pseudoholomorphic curves, \ie pseudoholomorphic maps from a fixed real surface 
$S$ modulo reparameterizations.

\newdefi{def4.0.1} 
A point $z \in S$ on a $J$-holomorphic curve $u: S \to X$ is a {\sl cusp}, or a 
{\sl cuspidal point}, if $\ord_z du >0$.
The number $\ord_z du$ is called the {\sl order of the cusp} of $u$ at $z$.
A $J$-holomorphic curve $u: S \to X$ containing cuspidal points is called
a {\sl cusp curve}. 
\end{defi}

Note that in the literature on pseudoholomorphic curves the notion ``cusp curve''
has a different meaning. Our terminology agrees rather with the one used in
algebraic geometry where the notion ``cusp'' means a ``peak'', \ie an irreducible
singularity. This describes the situation at hand more accurately.

\newsubsection[4.1]{Deformation of pseudoholomorphic maps} 
Explicit construction of deformations is needed to obtain local charts for 
subspaces of curves with prescribed singularities.

\newlemma{lem4.1.1} Let $B \subset \rr^{2n} \cong \cc^n$ be the unit ball, 
$\scry$ a Banach manifold, $\{J_{\eta,t}\}_{\eta\in \scry, t\in  [0,1]}$ a 
family of homotopies of almost complex structures in $B$ with parameterized by $
\scry$ and depending $C^{\ell-1}$-smoothly on $(\eta,t)\in \scry\times [0,1]$.
Further, let $u_{\eta,0}: \Delta \to B$, $\eta\in \scry$, be a $C^{\ell-1}
$-smooth family of $J_{\eta,0}$-holomorphic map, such that $u_{\eta,0}(\Delta) 
\subset B(\half)$, and $\ord_0(du_{\eta,t})= \mu$. 

Then for any family $v_\eta\in \rr^{2n}$ depending $C^{\ell-1}$-smoothly on 
$\eta \in \scry$ and any $\nu \in \nn$ there exists $t^*=t^*(J_t, u_0,v, \mu, 
\nu)>0$, a neighborhood $U_\scry$ of a given $\eta^*\in \scry$, and a 
$C^{\ell-1}$-smooth family of homotopies $\{ w_{\eta,t} \}_{\eta\in U_\scry, 
t\in [0, t^*]}$ with $w_{\eta,t}\in L^{1,p}(\Delta, \rr^{2n})$ such that the 
maps $u_{\eta,t}: \Delta \to B$ given by
\begin{equation}
u_{\eta,t}(z) = u_{\eta,0}(z) + z^\nu (t\,v_\eta + w_{\eta,t}(z))
\eqqno(4.1.1)
\end{equation}

\noindent
\sli are $J_{\eta,t}$-holomorphic if $\nu \le 2\mu +1$, and

\noindent
\slii are $J_{\eta,0}$-holomorphic if $\nu >2\mu +1$.

Moreover, for $z\not =0$ the function $w_{\eta,t}(z)$ depends 
$C^{\ell-1}$-smoothly on $(\eta, t, z)$.
\end{lem}

\state Remarks.~1. 
In other words, there exists a pseudoholomorphic deformation $u_t$ of a given map
$u_0$ in a given direction ${d \over dt}u_t\ogran_{t=0}=z^\nu v+ O(|z|^{\nu+
\alpha})$; and moreover, for smaller $\nu$ it is possible to deform 
simultaneously the almost complex structure. Furthermore, if the initial data 
depend smoothly on the parameter $\eta$, then the corresponding constructions give 
a smooth dependence of the maps on $\eta$.

\state 2. The loss of smoothness from $C^\ell$ to $C^{\ell-1}$ is due to the fact 
that the Gromov operator $D_{u,J}$ depends only $C^{\ell-1}$-smoothly
on $u$. Indeed, $D_{u,J}$ is the derivative of the $\dbar$-operator $u \mapsto 
\dbar_J u$ in the $u$-direction, which is only $C^\ell$-smooth.

\proof We give only a sketch. First, we fix a family $\phi_{\eta,t}$ of affine 
transformations of $\rr^{2n}$ with depend $C^{\ell-1}$-smoothly on $(\eta,t)$ 
such that $\phi_{\eta,t} \scirc u_{\eta,0} (0) =0 \in B$ and $\phi_{\eta,t}\scirc 
J_{\eta,t} \scirc \phi_{\eta,t}\inv$ coincide with $J\st$ in $u_{\eta,0}(0)$. 
Setting $\ti u_{\eta,0} \deff \phi_{\eta,t} \scirc u_{\eta,0}$ and $\ti J_{\eta,t}
\deff \phi_{\eta,t} \scirc J_{\eta,t} \scirc \phi_{\eta,t}\inv$ we reduce the 
problem to the case where $\ti u_{\eta,0}(0) = 0$ and $\ti J_{\eta,t}(0) = J\st$. 

Now we assume that there is no dependence on the parameter and drop the index 
$\eta$. Using \eqqref(4.1.1) one writes the equation $\dbar_{J_t} u_t =0$ in the 
form
\begin{equation}
(x + y J\st)^{-\nu} \dbar_{J_t}
\bigl(u_0(z) + (x + y J\st)^\nu (t\,v + w_t(z)) \bigr)=0
\eqqno(4.1.2)
\end{equation}
with $x +\isl y =z$ the standard coordinates on $\Delta$, and considers 
\eqqref(4.1.2) as an equation for $w_t(z)$. Then one shows that under the 
hypotheses of the lemma the linearization of \eqqref(4.1.2) has the form
\begin{equation}
(\dbar^{(\nu)}_{u_t, J_t} + R^{(\nu)}_{u_t, J_t}) \dot w_t(z) =
\psi^{(\nu)}_{u_t, J_t}(\dot J_t)(z),
\eqqno(4.1.3)
\end{equation}
where $\dot w_t(z)= {d \over dt}w_t(z)$ and $\psi^{(\nu)}_t(\dot J_t)(z) \in 
L^\infty(\Delta, \cc^n)$. Thus it is sufficient to find a right inverse
$T^{(\nu)}_{u_t, J_t}$ of the Gromov type operator $D^{(\nu)}_{u_t, J_t} = 
\dbar^{(\nu)}_{u_t, J_t} + R^{(\nu)}_{u_t, J_t}$ with an additional condition 
$\dot w_t(0)=0$. We refer to \cite{Iv-Sh-1}, {\sl Lemma 3.3.1}, for the 
explicit construction of such a right inverse $T^{(\nu)}_{u, J}$. Moreover, 
the operator $T^{(\nu)}_{u, J}$ and the inhomogeneity term $\psi^{(\nu)}_{u,J}$
depend smoothly on $u$ and $J$. As a consequence, the solution $w_t$ of 
\eqqref(4.1.2) depends $C^{\ell-1}$-smoothly on the parameter $\eta \in \scry$.
\qed

\newdefi{def4.1.1} Let $B \subset \rr^{2n} \cong \cc^n$ be a ball, $J_0$ a 
$C^\ell$-smooth almost complex structure in $B$, $u_0: \Delta \to B$ a 
$J_0$-holomorphic map and $\nu\ge 1$ an integer exponent. Denote
by $\dfrm_\nu (u, J; v)$ the a map depending $C^{\ell-1}$-smoothly on

\begin{itemize}
\item
a $C^\ell$-smooth almost complex structure $J$ in $B$, sufficiently close to $J_0$;

\item a $J$-holomorphic map $u$, sufficiently close to $u_0$;

\item a vector $v \in \rr^{2n}$, sufficiently close to $0$;
\end{itemize}

\noindent 
such that $\ti u\deff \dfrm_\nu (u, J; v)$ is a $J$-holomorphic map of the form
$\ti u(z)= u(z) + z^\nu v + O(|z|^{\nu + \alpha})$. Note that the choice of such
a map $\dfrm_\nu$ is not unique.
\end{defi}

\smallskip
\lemma{lem4.1.1} allows us to construct local deformations of pseudoholomorphic 
maps with appropriate types of singularities. To obtain a global deformation,
we use 

\newlemma{lem4.1.2} Let $u_0: S \to X$ be a non-multiple $J_0$-holomorphic 
map, $z_1, \ldots, z_m$ fixed points on $S$, and $U_1, \ldots, U_m \subset S$
disjoint neighborhoods of these points. Further, let $\{J_t\}_{t\in [0,1]}$
be a given $C^{\ell-1}$-smooth homotopy of almost complex structures on $X$, and
$\{u_{i,t}\}_{t\in [0,1]}$ given $C^{\ell-1}$-smooth homotopies of
$J_t$-holomorphic maps $u_{i,t}: U_i \to X$.

Then there exist $t^*>0$, a $C^{\ell-1}$-smooth homotopy $\{\ti J_t\}_{t
\in [0,t^*]}$ of almost complex structures on $X$, a $C^{\ell-1}$-smooth homotopy
$\{\ti u_t\}_{t\in [0,t^*]}$ of $\ti J_t$-holomorphic maps $\ti u_t: S \to X$
such that $u_t$ coincides with each $u_{i,t}$ in some (possibly smaller)
neighborhood of $z_i$ and $\ti J_t$ coincides with $J_t$ in some
neighborhood of each $x_i\deff u_0(z_i)$.
\end{lem}

The proof of the lemma is left to the reader.

\smallskip
Refining the result of \lemma{lem4.1.1} we show that the condition $u_1(z)-
u_2(z)= o(|z|^k)$ of \lemma{lem1.2.4} defines a submanifold in the spaces of 
pairs of pseudoholomorphic maps.

\newdefi{def4.2b.1} Define {\sl the spaces of pairs of pseudoholomorphic maps 
coinciding up to order $k$ at $z=0$} as $\scrpp_k(\Delta, X) \deff $
\begin{equation}
\bigl\{ (u', u'', J)\in 
L^{1,p}(\Delta, X) \times L^{1,p}(\Delta, X) \times \scrj: 
\dbar_Ju'=0=\dbar_Ju'', u'(z) - u''(z) = o(|z|^k) \bigr\},
\eqqno(4.2b.1)
\end{equation}
where the condition $u'(z) - u''(z) = o(z^k)$ is related to any local 
coordinate system on $X$ in a neighborhood of the point $u'(0)= u''(0) \in X$.
\end{defi}

The structure of $\scrpp_m(\Delta, X)$ for the cases $k=0$ and $k=1$ is easily
obtained from transversality techniques. In general we have

\newthm{thm4.2b.1} Assume that $\scrj$ consists of $C^\ell$-smooth structures 
with $\ell\ge2$. Then the space $\scrpp_k(\Delta, X)$ is a $C^{\ell-1}
$-submanifold of the fiber product $\scrp(\Delta, X) \times_{\!\!\scrj} 
\scrp(\Delta, X)$ of codimension of $2n(k+1)$ with the the tangent space 
\begin{multline}
T_{(u',u'',J)}\scrpp_k(\Delta, X) = 
\\
\bigl \{ (v',v'',\dot J) \in 
T_{(u',u'',J)}\bigl(\scrp(\Delta, X) \times_{\!\!\scrj} \scrp(\Delta, X) 
\bigr) :j^k(v' -v'') = 0 \bigr\}.
\eqqno(4.2b.2)
\end{multline}

Moreover, for $k=0$ and $k=1$ the space $\scrpp_k(\Delta, X)$ is well-defined
and $C^\ell$-smooth also for $\ell \ge1$.
\end{thm}

\proof It follows from \lemma{lem2.2.4} that $\scrp(\Delta, X) 
\times_{\!\!\scrj} \scrp(\Delta, X)$ is a $C^\ell$-smooth Banach manifold with 
the tangent space 
\begin{multline}
T_{(u',u'',J)}\bigl(\scrp(\Delta, X) \times_{\!\!\scrj} \scrp(\Delta, X) 
\bigr)=
\\
\bigl\{ (v', v'', \dot J) : (v',\dot J) \in T_{(u',J)}\scrp(\Delta, X), 
(v'',\dot J) \in T_{(u'',J)}\scrp(\Delta, X) \bigr\},
\eqqno(4.2b.3)
\end{multline}
so that $D_{u',J}v' + \dot J \scirc du' \scirc J_\Delta =0$ and similarly for 
$v''$. 

Fix $(u'_0, u''_0, J_0) \in \scrpp_0(\Delta, X)$ and local coordinates $(w_i)$ in
a neighborhood $U \subset X$ of $x^* \deff u'_0(0)= u''_0(0) \in X$. Then there 
exists $r>0$ such that for any pair $(u', u'')$ of $L^{1,p}(\Delta, X)$-maps
sufficiently close to $(u'_0, u''_0)$ we  have $u'(\Delta(r)) \subset U$ and 
$u''(\Delta(r)) \subset U$. The coordinates in $U$ induce the linear structure.
Thus we can consider the difference $u'(z)- u''(z)$ having in mind that it 
is well-defined only for $z\in \Delta(r)$. 

The subspace $\scrpp_0(\Delta, X)$ is defined by the condition $u'(0)=u''(0)$ 
for $(u',u'',J) \in \scrp(\Delta, X) \times_{\!\!\scrj} \scrp(\Delta, X)$. Setting 
$F(u', u'', J) \deff u'(0) -u''(0)$ we obtain a $C^\ell$-smooth function, which
is well-defined in a neighborhood of $(u'_0, u''_0, J_0)$ and is a 
local defining function for $\scrpp_0(\Delta, X)$. The differential of $F$ in 
$(u',u'',J) \in \scrpp_0(\Delta, X)$,
$$
dF: T_{(u',u'',J)}\bigl(\scrp(\Delta, X) \times_{\!\!\scrj} \scrp(\Delta, X) 
\bigr) \to T_{u'(0)} X,
$$ 
is given by the formula $dF(v',v'', \dot J) = v'(0)-v''(0)$ and is
a surjective map. Thus $\scrpp_0(\Delta, X)$ is a $C^\ell$-smooth submanifold
of $\scrp(\Delta, X) \times_{\!\!\scrj} \scrp(\Delta, X)$ of codimension $2n =
\dimr X$. 

\smallskip
Considering the $C^\ell$-smooth map $\ev_0: \scrpp_0(\Delta, X) \to X$ with 
$$\ev_0(u',u'', J) \deff u'(0) = u''(0),$$ we obtain  a $C^\ell$-smooth bundle
$E^{(0)}$ over $\scrpp_0(\Delta, X)$ with fiber $E^{(0)}_{(u',u'', J)}= {u'}^*
T_{u'(0)}X$. The formulas $\sigma'(u',u'', J) \deff du'(0)$ and $\sigma''(u',u'',
 J) \deff du''(0)$ define $C^\ell$-smooth sections of the bundle $T^*_0\Delta 
\otimes E^{(0)}$ over $\scrpp_0(\Delta, X)$. Thus the condition $du'(0)= du''(0)$
is equivalent to the vanishing of $\sigma' -\sigma''$. Consequently, $\scrpp_1
 (\Delta, X)$ is a $C^\ell$-smooth submanifold of $\scrpp_0(\Delta, X)$ of 
codimension $2n$. 
 
\smallskip
We proceed further by induction using the case $k=0$ as the base. Our notation is
as follows. For a triple $(u',u'', J)\in \scrpp_0(\Delta,X)$ we 
consider the (integrable) complex structure $J\st$ in $U$ with coincides with 
$J$ at the point $u'(z)= u''(z)$ and is constant \wrt the coordinates in $U$. 
Note that $J\st$ depends $C^\ell$-smoothly on $(u',u'', J) \in \scrpp_k(\Delta,
X)$. Thus we can regard $U$ as an open subset in $\cc^n$.

For a pair $(u',u'')$ of $J$-holomorphic maps with values in $U \subset X$ we 
obtain
\begin{align}
0& = \dbar_J u' -\dbar_J u'' = \bigl( (\d_x u' - J(u')\cdot \d_y u') -
(\d_x u'' - J(u'')\cdot \d_y u'') 
\notag\\
&= \d_x(u'-u'') + J(u')\cdot \d_y(u' -u'') +
\bigl(J(u') -J(u'')\bigr)\cdot \d_y u'' 
\notag\\
&= \dbar_{J(u')}(u' -u'') +\bigl(J(u') -J(u'')\bigr)\cdot \d_y u'' .
\eqqno(4.2b.4)
\end{align}
Consequently,
\begin{align}
\llap{$\dbar_{J\st}$}(u'-u'')&
= \dbar_{J\st}(u'-u'') - (\dbar_J u' -\dbar_J u'')
\notag\\
&=\bigl(J\st - J(u')\bigr)\cdot \d_y (u'-u'')- 
\bigl(J(u') -J(u'')\bigr)\cdot \d_y u'' .  
\eqqno(4.2b.5)
\end{align}
Let us denote the last expression by $H_{u',u'',J}(z)$

\smallskip
Now suppose that $(u',u'', J)$ varies in $\scrpp_k(\Delta,X)$ with $k\ge1$.
We can assume by induction that $\scrpp_k(\Delta,X)$ is a $C^{\ell-1}$-smooth 
manifold. We claim that for any $p<\infty$
\begin{equation}
\eqqno(4.2b.7)
\llap{$f_k(z) \deff $}z^{-(k+1)}(u'(z) - u''(z)) 
\end{equation}
is a well-defined $L^{1,p}(\Delta(r), \cc^n)$-valued function depending $C^{\ell-
1}$-smoothly on $(u',u'',J) \in \scrpp_k(\Delta,X)$. The claim implies the 
theorem. Indeed, the function $F_k$ given by $F_k: (u',u'', J) \in \scrpp_k(
\Delta, X) \mapsto f_k(0) \in \cc^n$ is then a local defining function for
$\scrpp_{k-1}(\Delta,X)$ inside $\scrpp_k(\Delta,X)$, whereas non-degeneracy
of $dF_k$ can be easily obtained from \lemma{lem4.1.1}.

Again by induction, we can suppose that $f_{k-1}(z)= z^{-k}(u'(z) - u''(z))$ is a 
well-defined $L^{1,p}(\Delta(r), \cc^n)$-valued function depending $C^{\ell-1}
$-smoothly on $(u',u'',J) \in \scrpp_{k-1}(\Delta,X)$. Note that $f_{k-1}(0)$
vanishes identically on $\scrpp_k(\Delta,X)$. Further, for any exponents $p, p'$
with $2<p'<p <\infty$ the map $f(z)\in L^{1,p}(\Delta, \cc^n) \mapsto z\inv (f(z)
-f(0)) \in L^{p'}(\Delta, \cc^n)$ is linear and bounded. Consequently, for any
$p<\infty$ the function $f_k(z) = z\inv f_{k-1}(z)$ lies in $L^p(\Delta, \cc^n)$
and depends $C^{\ell-1}$-smoothly on $(u',u'',J) \in \scrpp_k(\Delta,X)$ \wrt
the $L^p$-topology.

\smallskip
Without loss of generality we may assume that $U$ is convex. The identity 
$$
J(w) = J(w^*) + \int_{t=0}^1 \d_t J(w^* +t(w-w^*)) dt
$$
for $(w,w^*) \in U \times U$ implies the relation $J(w)= J(w^*) +\sum_i (w_i-
w^*_i) S_i(w,w^*) = S(w,w^*; \allowbreak w-w^*)$ with the function $S(w,w^*; 
\ti w)$ depending $C^\ell$-smoothly on $J \in \scrj$, $C^{\ell-1}$-smoothly on 
$(w,w^*) \in U\times U$ and $\rr$-linearly on $\ti w\in \cc^n$. Substituting 
$u''(z) = u'(z) + z^{k+1} f_k(z)$ in $(J(u') -J(u''))\cdot \d_y u''$ we obtain
$$
(J(u'(z)) -J(u''(z)))\cdot \d_y u''(z) = S\bigl(u''(z), u'(z); z^{k+1}f_k(z) 
\bigr)\cdot \d_y u''(z) 
$$
By apriori regularity estimates, for $r<1$ we can consider $du''(z)$ as a map from 
$\scrpp_{\!k}(\Delta,X)$ to $C^0(\Delta(r),\cc^n)$ which depends $C^\ell$-smoothly 
on $(u',u'',J)$. Thus we have represented the term $(J(u') -J(u''))\cdot \d_y u''$
as a composition of the $C^{\ell-1}$-smooth map 
$$
(u',u'',J)\in \scrpp_k(\Delta,X) \;\mapsto\; S(u''(z),u'(z); f_k(z)) \cdot
\d_y u''(z) \in L^p(\Delta(r), \cc^n)
$$
and the linear bounded map
$$
S(u''(z),u'(z); f_k(z)) \d_y u''(z)\;\mapsto\; z^{-(k+1)} \cdot 
S(u''(z),u'(z); z^{k+1} \cdot f_k(z)) \d_y u''(z).
$$
Thus $(J(u') -J(u''))\cdot \d_y u''$ depends $C^\ell$-smoothly on $(u',u'',J)\in
\scrpp_k(\Delta,X)$ \wrt the norm topology in $L^p(\Delta(r), \cc^n)$. 
Consequently, the formula 
$$
(u',u'', J) \in \scrpp_k(\Delta,X) \;\mapsto \;
z^{-(k+1)} \cdot (J\st - J(u'(z))) \cdot \d_y (u'(z)-
u''(z))
$$
defines a $L^p(\Delta(r), \cc^n)$-valued map depending $C^{\ell-1}
$-smoothly on $(u',u'', J) \in \scrpp_k(\Delta,X)$. 

Similar estimates can be be carried out for the first term $(J\st - J(u'))\cdot 
\d_y(u'-u'')$ in \eqqref(4.2b.5). Together, this implies that $h_k(z) \deff z^{-k}
H_{u',u'',J}(z)$ lies in $L^p(\Delta(r),\cc^n)$ and depends $C^{\ell-1}$-smoothly
on $(u',u'', J) \in \scrpp_k(\Delta,X)$ \wrt $L^p$-topology. Now let $f_{\dbar, 
k}(z)$ be a solution of the equation $\dbar_{J\st}f_{\dbar, k}(z) = h_k(z)$
depending $C^{\ell-1}$-smoothly on $(u',u'', J) \in \scrpp_k(\Delta,X)$ \wrt 
the $L^{1,p}$-topology. Then $(u'(z) -u''(z)) - z^{k+1} f_{\dbar, k}(z)$ is 
a holomorphic $\cc^n$-valued function, depending $C^{\ell-1} $-smoothly on 
$(u',u'', J) \in \scrpp_k(\Delta,X)$ \wrt $L^{1,p}$-topology and vanishing in 
$z=0$ up to order $k+1$. Consequently, 
$$
(u'(z) -u''(z)) - z^{k+1} f_{\dbar, k}(z)= z^{k+1}f_{\scro, k}(z)
$$
and $f_k(z) = f_{\scro, k}(z) + f_{\dbar, k}(z)$ possesses the property claimed
above. \qed

\smallskip

\newsubsection[4.2a]{Curves with prescribed cusp order}
In this paragraph we show that $J$-curves with cusps of given order form a 
Banach submanifold of the moduli space and compute its codimension.

\newdefi{def4.2.1} For a given natural $m$ we denote by $\mbfk$ an $m$-tuple 
$(k_1, \ldots, \allowbreak
k_m)$ with $k_i \ge1$ and set $|\mbfk| \deff \sum_i
k_i$. The $m$-tuple $(1,\ldots,1)$ is denoted $\bfone_m$.
Define the {\sl moduli space $\scrm_\mbfk$ of pseudoholomorphic curves
with a given cusp order $\mbfk$} as the set of classes $[u,J, \mbfz]$ such
that $[u,J] \in \scrm$ and $u$ has $m$ (marked) cusp-points $\mbfz= \{z^*_1,
\ldots, z^*_m \}$ with $\ord_{z^*_i} \ge k_i$. Two triples $(u,J, \mbfz)$
and $(\ti u, \ti J, \ti\mbfz)$ define the same class $[u,J, \mbfz] = [\ti u,
 \ti J, \ti\mbfz] \in \scrm_\mbfk$ \iff there exists $g \in \bfg$ such that
$\ti u = u \scirc g$ and $\ti z^*_i = g(z^*_i)$.
\end{defi}

\smallskip
The main result of this paragraph is

\newthm{thm4.2.1} The set $\scrm_\mbfk$ is a $C^\ell$-smooth
manifold and the natural map $\scrm_\mbfk \lrar \scrm$ given by $[u,J, \mbfz]
\mapsto [u,J]$ of $\scrm$ is an immersion of codimension $2\, (n\, |\mbfk|
-m)$, where $n= \dimc X$ and $m$ is the number of marked cusp-points.
\end{thm}

\smallskip
We divide the proof in several steps. First we consider the corresponding
problem for $\whcalm$. The reason is that it is more convenient to work
with maps, \ie elements of
$\whcalm$, than with parameterized curves, \ie elements of $\scrm$. This
means that we are interested in the set
\begin{equation}
\whcalm_\mbfk \deff \left\{ (u,J_S, J; z^*_1, \ldots z^*_m) \in
\whcalm \times (S)^m :
\msmall{ \matrix 
 z^*_i\text{ are pairwise distinct, }\cr
\ord_{z^*_i} du \ge k_i
\endmatrix}
\right\},
\eqqno(4.2.1)
\end{equation}
where $(S)^m= S \times \cdots \times S$ is the $m$-fold product of $S$.
Obviously, the projection from $\whcalm \times (S)^m$ onto $\whcalm$ and then
onto $\scrm$ maps $\whcalm_\mbfk$ onto $\scrm_\mbfk$. In our proof of
\refthm{thm4.2.1} we shall show that this map $\whcalm_\mbfk \to
\scrm_\mbfk$ is a principle $\bfg$-bundle.

\newdefi{def4.2.1n} Set 
\begin{equation}
\whcalm^{(m)} \deff \{ (u,J_S, J; z^*_1, \ldots z^*_m) \in
\whcalm \times (S)^m : z^*_i \not = z^*_j \text{ for every } i\not=j\,\}
\eqqno(4.2.2n)
\end{equation}
denoting by $S_i$ the $i$-th factor in $(S)^m$. Equip $\whcalm^{(m)}$ with the 
maps $\ev_i: \whcalm^{(m)} \to X^m$ defined by $\ev_i(u,J_S, J; z^*_1, \ldots, 
z^*_m) \deff u(z^*_i)$. Denote by $E_i$ the pulled-back bundles $\ev_i^*TX$ and
$\ev^{(m)}{}^* T(X^m)$ over $\whcalm^{(m)}$. The fiber of $E_i$ over $(u, J_S, 
J; \mbfz)$ is $(E_i)_{(u, J_S, J; \mbfz)}= T_{u(z^*_i)}X$.
\end{defi}

Obviously, the space $\whcalm^{(m)}$ is a $C^\ell$-smooth Banach manifold, 
$\ev_i: \whcalm^{(m)} \to X^m$ are $C^\ell$-smooth maps, and $E_i$ are 
$C^\ell$-smooth bundles over $\whcalm^{(m)}$. Note that we also have line
bundles $TS_i$ and $T^*S_i$ over $\whcalm^{(m)}$ which are defined in an obvious 
way as the (co)tangent bundles to each $S_i$.

\newlemma{lem4.2.2n} The formula $\yps(u,J_S, J; z^*_1, \ldots z^*_m) \deff
 (du(z^*_1), \ldots, du(z^*_m)) \in \bigoplus_i T^*S_i \otimes E_i$ defines
a $C^\ell$-smooth section of\/ $\bigoplus_i T^*S_i \otimes E_i$ over 
$\whcalm^{(m)}$, transversal to the zero section. The zero-set of $\yps$
coincides with the space $\whcalm_{\bfone_m}$ of maps having cups in each
marked $z^*_i$. Thus $\whcalm_{\bfone_m}$ is a $C^\ell$-smooth Banach 
submanifold of $\whcalm^{(m)}$ of codimension $2nm$.
\end{lem}

Before starting the proof we introduce some new notation.

\newdefi{def4.2.2n} Let $\scry$ be a $C^\ell$-smooth Banach manifold and $f: 
\scry \to \ttt_g \times S$ a $C^\ell$-smooth map of the form $f(y)= (J_S(y), 
z^*(y))$. Set $F(y) \deff (y, z^*(y))$ so that $F: \scry \to \scry \times S$ is 
an imbedding. A {\sl local $J_S(y)$-holomorphic coordinate (or simply 
a {\sl$J_S$-holomorphic coordinate}) on $\scry \times S$ centered at $z^*$} 
is a $C^\ell$-smooth $\cc$-valued function $z$ defined in some neighborhood $U
\subset \scry \times S$ of $F(\scry)$ which vanishes along $F(\scry)$ and is
$J_S(y)$-holomorphic along each $\{y\} \times S$. One can use \lemma{lem4.1.1}
for a proof of the existence of such a local holomorphic coordinate. 
\end{defi}

\statep Proof of. \lemma{lem4.2.2n}. It is obvious that $\yps$ is 
well-defined. To show the $C^\ell$-smoothness of $\yps$, for any $i=1,\ldots,m$,
we fix some local coordinate $z_i$ on $\whcalm^{(m)}$ which is 
$J_S$-holomorphic along $S_i$ and centered at $z^*_i \in S_i$. Now we can
find a local frame $\bfxi=(\xi_1, \ldots, \xi_n)$ of $T^*S_i \otimes E_i$ which
depends $C^\ell$-smoothly on $(u,J_S, J) \in \whcalm$ and
holomorphically on the coordinate $z_i$. The existence of such a frame follows
from \refdefi{def1.4.1} and a parametric version of \lemma{lem1.4.1A}.
The coefficients of $du \in T^*S_i \otimes E_i$ \wrt such a frame 
$\bfxi$ depend $C^\ell$-smoothly on $(u,J_S, J) \in \whcalm$ and holomorphically
on $z_i$. Consequently, the $du(z^*_i)$ depend $C^\ell$-smoothly on $(u,J_S, J; 
\mbfz) \in \whcalm^{(m)}$. Thus $\yps$ is $C^\ell$-smooth.

The transversality of $\yps$ to the zero-section of $\bigoplus_i T^*S_i \otimes
E_i$ follows immediately from results of \refsubsection{4.1}. In particular,
$\whcalm^{(m)}_{\bfone_m}$ is the $C^\ell$-smooth Banach submanifold of
$\whcalm^{(m)}$. The corresponding codimension is $\rank_\rr\left(
\bigoplus_i T^*S_i \otimes E_i\right) =2nm$. \qed

\newdefi{def4.2.3n} For (finite-dimensional) complex vector spaces $V$, $W$, 
and $k\in \nn$ denote by $j^k(V,W)$ the vector space of polynomial maps 
$f:V \to W$ of degree $\deg f\le k$ with $f(0)=0$, considered as the space of 
$k$-jets of holomorphic maps $F: V \to W$. For $l\ge k$ the natural projection
$\pr: j^l(V,W) \to j^k(V,W)$ is well-defined. Let $j^{k,l}(V,W)$ denote 
its kernel. Similar notation for complex bundles is used. Note that 
$j^1(V,W) = \hom(V,W) = V^* \otimes W$.
\end{defi}

\newlemma{lem4.2.3n} \sli For any $(u, J_S, J; \mbfz) \in \whcalm_\mbfk$ the
jet $j^{2k_i+1}u(z^*_i)$ is a well-defined element of $j^{2k_i+1} (T_{z^*_i}, 
T_{u(z^*_i)}X) = j^{2k_i+1} (TS_i, E_i)_{(u, J_S, J; \mbfz)}$. 

\slii Moreover, $j^{2k_i+1}u(z^*_i) \in j^{k_i+1, 2k_i+1} (T_{z^*_i}, 
T_{u(z^*_i)}X)= j^{k_i+1, 2k_i+1} (TS_i, E_i)_{(u, J_S, J; \mbfz)}$.

\sliii Set $\yps_\mbfk(u, J_S, J; \mbfz) \deff \bigl(j^{k_1+1, 2k_1+1} u(z^*_1), 
\ldots, j^{k_m+1, 2k_m+1}u(z^*_m) \bigr)$. Then $\yps_\mbfk: \whcalm_\mbfk \to 
\bigoplus_i j^{k_i+1, 2k_i+1} (TS_i, E_i)$ is a section which is $C^\ell$-smooth 
and transversal to the zero-section.
\end{lem}

\proof Assertions \sli and \slii follow essentially from \lemma{lem1.2.4}. 
The nontrivial points here are the following. First, the jet 
$j^{2k_i+1}u(z^*_i)$ is defined even if the structure $J$ is $C^\ell$-smooth 
with $\ell<2k_i$ and the map $u$ is $C^{\ell+1}$-smooth, since in general there
are no higher smoothness for $u$. Second, the jet $j^{2k_i+1}
u(z^*_i)$ is a {\sl complex} polynomial. Finally, the jet $j^{2k_i+1}u(z^*_i)$
is independent of the choice of the integrable structure $J\st$ and 
$J\st$-holomorphic coordinates in a neighborhood of $u(z^*_i)$ used in 
\lemma{lem1.2.4} for definition of the jet. Let us give a proof of the latter
property.

Let $J'$ and $J''$ be integrable complex structures in a neighborhood of 
$u(z^*_i)$ such that $J'(u(z^*_i)) = J''(u(z^*_i)) = J(u(z^*_i))$. Find local
complex coordinate systems $\mbfw'= (w'_1, \ldots,  \allowbreak
w'_n)$ and $\mbfw''= (w''_1, \ldots,w''_n)$ which are centered in $u(z^*_i)$ and 
holomorphic \wrt $J'$ and
$J''$ respectively. Without loss of generality we may assume that the frames 
$({\d \over \d w'_1}, \ldots,  {\d \over \d w'_n})$ and $({\d \over \d w''_1}, 
\ldots,  {\d \over \d w''_n})$ coincide in $u(z^*_i)\in X$. Consequently, 
we can express one system by another using the formula $\mbfw''= \mbfw' + 
F(\mbfw')$ with 
\begin{equation}
F(\mbfw') = O(|\mbfw'|^2) \qquad \text{and} \qquad 
dF(\mbfw') = O(|\mbfw'|).
\eqqno(4.2.3n)
\end{equation}
Let $u'(z)$ and $u'(z)$ be the local expressions of $u: S \to X$ in the local 
coordinate systems $\mbfw'$ and $\mbfw''$ respectively. Then $u''(z)= u'(z) + 
F(u'(z))$. So from \eqqref(4.2.3n) and $u'(z)= O(|z|^{k_i+2})$ we see that
coefficients of polynomials $j^{2k_i+1}u'(z)$ and $j^{2k_i+1}u''(z)$ coincide. 

\smallskip
To show the smoothness of the section $\yps_\mbfk$ we fix an element 
$(u_0,J_{S,0}, J_0; \mbfz_0)\in \whcalm_\mbfk$, $\mbfz_0=(z^*_{1,0}, \ldots,
z^*_{m,0})$, and a sufficiently small neighborhood $\scry \subset \whcalm
_\mbfk$ of $y_0 \deff (u_0, J_{S,0}, J_0; \mbfz_0)$.
In what follows, for any $i=1, \ldots,m$, we fix families of certain structures 
on various spaces. We assume that the members of the families are parameterized
by and depend $C^\ell$-smoothly on $y=(u,J_S, J; \mbfz)\in \scry$. The 
families are:
\begin{enumerate}
\item integrable complex structures $J'_i$ in a neighborhood 
of each $u(z^*_{i,0})$ such that each $J'_i$ coincides with $J$ in $u(z^*_i)$;
\item local complex coordinate systems $\mbfw'_i= (w'_{i,1}, \ldots, 
w'_{i,n})$ on $X$ centered in $u(z^*_i)$ and holomorphic \wrt $J'_i$;
\item local frames $\bfxi_i= (\xi_{i,1},\ldots ,\xi_{i,n})$ of the bundles 
$E_i$ which are defined in a neighborhood of $z^*_i$ and holomorphic along 
$S_i$; 
\item  local $J_S$-holomorphic coordinates $z_i$ on $S_i$ centered in 
$z^*_i$.
\end{enumerate}
Further, we assume that every coordinate $z_i$ has image the whole disc $\Delta$.
Note that pulling back the frames $({\d \over \d w'_{i,1}}, \ldots, {\d \over 
\d w'_{i,n}} )$ we obtain local frames $\left(u^*({\d \over \d w'_{i,1}}), \ldots, 
u^*({\d \over \d w'_{i,n}}) \right)$ of $E_i$ which depend $C^\ell$-smoothly 
on $y=(u,J_S, J; \mbfz)\in \scry$. Now, the expression of $u(z)$ 
in the local coordinate system $\mbfw'_i$ yields an element $u'_i(z_i) \in
L^{1,p}(\Delta, \cc^n)$ which depends $C^\ell$-smoothly on $y\in \scry$ 
\wrt the standard smooth structure in $L^{1,p}(\Delta, \cc^n)$. Deriving, we 
obtain an element $du'_i(z_i)\in L^p(\Delta, \cc^n\otimes_\rr T^*\Delta)$ which
depends $C^\ell$-smoothly on $y\in \scry$ \wrt the standard smooth structure 
in $L^p(\Delta, \cc^n\otimes_\rr T^*\Delta)$. 

Consider now $du'_i$ as a 
section of $E_i\otimes T^*S_i$, and its coefficients of $du'_i$ in the frame 
$\left(u^*({\d \over \d w'_{i,1}}) \otimes dz_i, \ldots, u^*({\d \over 
\d w'_{i,n}}) \otimes dz_i \right)$ as $L^p(\Delta, \cc)$-functions. Thus we
can conclude that the coefficients of $du'_i$ depend $C^\ell$-smoothly on 
$y\in \scry$ \wrt the standard smooth structure in $L^p(\Delta, \cc)$.
Consequently, the same is true for the coefficients of $du'_i$ in the frame 
$(\xi_{i,1}\otimes dz_i,\ldots, \xi_{i,n}\otimes dz_i)$. Since the latter frame
is holomorphic, the coefficients of the jet $j^{2k_i}du(z^*_i)$ depend 
$C^\ell$-smoothly on $y\in \scry$. This provides the desired smoothness 
property of $\yps_\mbfk$.

\smallskip
Finally, note that the transversality of $\yps_\mbfk$ to the zero-section 
follows from results of \refsubsection{4.1}. \qed

\newcorol{cor4.2.4n} For any $\mbfk=(k_1,\ldots, k_m)$ with $k_i \ge 1$ 
the space $\whcalm_\mbfk$ is a $C^\ell$-submanifold of $\whcalm^{(m)}$ of 
codimension $2|\mbfk|n$.
\end{corol}

\proof Assume that for a given $\mbfk=(k_1,\ldots, k_m)$ with 
$k_i \ge 1$ the claim holds. Fix some $\mbfk^+=(k^+_1,\ldots, k^+_m)$ with 
$k_i \le k^+_i \le 2 k_i$ and consider truncated section $\yps_{\mbfk,\mbfk^+}
: \whcalm_\mbfk \to \bigoplus_i j^{k_i+1, k^+_i+1} (TS_i, E_i)$ given
by 
$$
\yps_{\mbfk,\mbfk^+}(u, J_S, J; \mbfz) \deff \bigl(j^{k_1+1, k^+_1+1} 
u(z^*_1), \ldots, j^{k_m+1, k^+_m+1}u(z^*_m) \bigr).
$$ 
Then $\whcalm_{\mbfk^+}$ is identified with the zero set of $\yps_{\mbfk,
\mbfk^+}$. By \lemma{lem4.2.3n}, $\yps_{\mbfk,\mbfk^+}$ is transversal to 
the zero-section. Thus $\whcalm_{\mbfk^+}$ is a $C^\ell$-smooth submanifold
of $\whcalm_\mbfk$ of codimension equal to $\rank_\rr \bigoplus_i j^{k_i+1, 
k^+_i+1} (TS_i, E_i) = 2n(|\mbfk^+| -|\mbfk|)$. So we can apply the induction.
\qed

\smallskip
\newlemma{lem4.2.5n} The natural projection $\wh\pr_\mbfk : \whcalm_\mbfk 
\to \whcalm$ given by the formula $\wh\pr_\mbfk(u, J_S, J; \mbfz) 
\allowbreak
\deff  (u, J_S, J)$ is an immersion of codimension $2(|\mbfk|n -m)$. 
\end{lem}

\proof The differential of the projection $\wh\pr_\mbfk$ is given by 
$$
d\wh\pr_\mbfk: (v, \dot J_S, \dot J; \dot \mbfz)\in 
T_{(u, J_S, J; \mbfz)}\whcalm_\mbfk \mapsto 
(v, \dot J_S, \dot J) \in T_{(u, J_S, J)}\whcalm.
$$
Thus the kernel $\ker d\wh\pr_\mbfk$ consists of vectors of the form $(0,0,0; 
\dot \mbfz)$ with $\dot \mbfz_i = (\dot z^*_1, \ldots, \dot z^*_m) \in 
\bigoplus T_{z^*_i}S_i$ and we must show that $\ker d\wh\pr_\mbfk$ is trivial.
Intuitively this is obvious, since elements of the kernel correspond to 
deformations leaving $(u,J_S,J)$ unchanged but moving cusp-points $z^*_i$ on 
$S$ and this is impossible.

For a rigorous proof we use conclusions of the proof of \lemma{lem4.2.3n}. 
Consider $du(z_i)$ as a holomorphic section of $T^*S_i 
\otimes E_i$. Then $du(z_i)$ vanishes in $z^*_{i,t}$ up to the order $\ge k_i$ 
and there are no other zeros of $du(z_i)$ in a neighborhood of $z^*_i$. Thus 
we can locally express $z^*_i$ as the zero set of $du(z_i)$. This implies
that locally there exists $C^\ell$-smooth functions $F_i$ of $(u,J_S,J)\in 
\whcalm$ such that $F_i(u,J_S,J) =z^*_i$ for $(u,J_S,J)\in \whcalm_\mbfk$.
Thus $\wh\pr_\mbfk : \whcalm_\mbfk \to \whcalm$ is an immersion.

To compute the codimension of $\wh\pr_\mbfk: \whcalm_\mbfk \hook \whcalm$
one represents $\wh\pr_\mbfk$ as the composition $ \whcalm_\mbfk \hook 
\whcalm^{(m)} \buildrel \pr \over \lrar \whcalm$.
\qed

\medskip
Now we can finish

\nobreak
\statep Proof of. \refthm{thm4.2.1}. Consider the action of $\bfg$ on 
$\whcalm$ and the diagonal action of $\bfg$ on $\whcalm \times (S)^m$. The 
both actions are $C^\ell$-smooth, free, and commute with the projection $\pr: 
\whcalm \times (S)^m \to \whcalm$. Moreover, for every $\mbfk=(k_1,\ldots, 
k_m)$ with $k_i\ge 1$ the submanifold $\whcalm_\mbfk \hook \whcalm \times 
 (S)^m$ is $\bfg$-invariant \wrt the diagonal action of $\bfg$. 
For the quotient $\scrm_\mbfk= \whcalm_\mbfk /
\bfg$ one can construct a $C^\ell$-smooth atlas in the same way as it was
done for $\scrm= \whcalm /\bfg$. The construction shows that the map 
$\scrm_\mbfk \to \scrm$ is a $C^\ell$-smooth immersion of codimension
equal to the codimension of $\wh\pr_\mbfk: \whcalm_\mbfk \hook \whcalm$.
\qed

\medskip
Summarizing the results and notation of this paragraph, we obtain

\newcorol{cor4.2.6n} The maps $\ev_\mbfk : \scrm_\mbfk \to X^m$ and $\ev_i:
\scrm_\mbfk \to X$ given by $\ev_\mbfk([u,J, \mbfz]) \deff (u(z^*_1), \ldots , 
u(z^*_m))$ and $\ev_i([u,J, \mbfz]) \deff u(z^*_i)$ are well-defined and 
$C^\ell$-smooth. This yields $C^\ell$-smooth bundles $E_i \deff \ev_i^* TX$ 
with a fiber $(E_i)_{[u,J, \mbfz]} = T_{u(z^*_i)}X$. The bundles $T_{z^*_i}
S_i$ over $\whcalm_\mbfk$ induce $C^\ell$-smooth bundles $L_i$ over $\scrm
_\mbfk$ with the fiber $(L_i) _{[u,J, \mbfz]}= T_{z^*_i} S_i$.
The section $\yps_\mbfk: \whcalm_\mbfk \to \bigoplus_i j^{2k_i+1}(TS_i, E_i)$
induces the section $\yps_\mbfk: \scrm_\mbfk \to \bigoplus_i j^{2k_i+1}
 (L_i, E_i)$ with $\yps_\mbfk([u,J, \mbfz]) \deff \yps_\mbfk(u,J_S,J; \mbfz)$. 
\end{corol}

\proof The claim follows from the fact that all the constructions are 
compatible with $\bfg$-action.
\qed

\smallskip
\newdefi{def4.2.4n} 
For a given $\mbfk=(k_1, \ldots, k_m)$ we set
\begin{align}
\whcalm_{=\mbfk} &\deff \bigl\{ (u, J_S, J; \mbfz) \in  \whcalm_\mbfk
\;:\; \ord_{z^*_i} du =k_i \,\bigl\};
\\
\scrm_{=\mbfk} &\deff \bigl\{ [u, J; \mbfz] \in  \scrm_\mbfk
\;:\; \ord_{z^*_i} du =k_i \,\bigl\}.
\end{align}
\end{defi}

\newlemma{lem4.2.7n} 
\sli The set $\whcalm_{=\mbfk}$ is an open $C^{\ell-1}$-smooth submanifold 
of $\whcalm_\mbfk$ invariant \wrt the natural action of\/ $\bfg$ on 
$\whcalm_\mbfk$.

\slii The image of the projection of $\whcalm_{=\mbfk}$ to $\whcalm$ is an
{\sl imbedded} submanifold of $\whcalm$, and the projection is a non-ramified
covering over the image.

\sliii There exists a $C^{\ell-1}$-smooth bundle $N$ over $\whcalm_{=\mbfk}
\times S$ whose restriction onto  $\{(u, J_S, 
\allowbreak
J; \mbfz)\} \times S$ coincides with $N_u$. The diagonal action 
of\/ $\bfg$ on $\whcalm_{=\mbfk} \times S$ lifts canonically to the action 
on the bundle $N$.

\sliv The bundle $N$ induces Banach bundles $L^{1,p}(S, N)$ and $L^p_{(0,1)}
 (S, N)$ over $\whcalm_{=\mbfk}$ with fibers $L^{1,p}(S, N_u)$ and $L^p_{(0,1)}
 (S, N_u)$ over $(u, J_S, J, \mbfz) \in \whcalm_{=\mbfk}$ respectively. The
operators $D^N_{u,J}: L^{1,p}(S, N_u) \to L^p_{(0,1)}(S, N_u)$ induce a 
$C^{\ell-1}$-smooth bundle homomorphism $D^N: L^{1,p}(S, N) \to L^p_{(0,1)}
 (S, N)$.
\end{lem}

\proof \sli The complement $\whcalm_\mbfk \bs \whcalm_{=\mbfk}$ is of the 
union of (the projections of) the spaces $\whcalm_{\mbfk'}$ such that either 
$\mbfk' =(k_1, \ldots, k_m, 1)$, or $\mbfk' =(k'_1, \ldots, k'_m)$ with $k'_i 
\ge k_i$ and $k'_{i_0} > k_{i_0}$ for some $i_0$. In other words, we have 
either at least one additional cusp-point or a higher order cusp in at least 
one point. Obviously, these conditions define closed subsets in 
$\whcalm_\mbfk$. The $\bfg$-invariance of $\whcalm_{=\mbfk}$ follows from 
the definition.

\slii The set $\whcalm_{=\mbfk}$ admits a finite transformation group $\aut(
\mbfk)$ generated by transpositions of marked cusp-points $z^*_i$ and
$z^*_j$ with $k_i = k_j$. The rest of part \slii follows.

\sliii Let $z_i$ be a local $J_S$-holomorphic coordinate on $\whcalm_{=\mbfk}
\times S$ centered at $z^*_i$ as in \refdefi{def4.2.2n}. It follows from the
proof of \lemma{lem4.2.3n} that $z_i^{-k_i} du(z_i)$ is a well-defined 
{\sl non-vanishing} local section of $\hom(TS, E_u)$, 
which depends $C^{\ell-1}$-smoothly on $(u,J_S, J)$ and holomorphically on
$z_i$. This provides the existence on $N$ with the stated property,
at least locally in a neighborhood of $(u,J_S, J; z^*_i)$. The globalization
of $N$ is trivial. Since the constructions involved are natural, the
$\bfg$-action admits the desired lift.

\sliv One uses the fact that the constructions of the bundles $L^{\!1,p}
\!(S, N)$, $L^p_{(0,1)}\!(S, N)$, and the operator $D^N$ are natural. This
implies $C^{\ell-1}$-smoothness of the obtained objects. \qed

\smallskip
\state Remark. One could explain the meaning of \lemma{lem4.2.7n} as follows.
First, we note that for the globalization of normal bundles $N_u$ to $\scrm$
we should use not the Cartesian product $\scrm \times S$, but the
$\bfg$-twisted product $\scrm \ltimes S$, \ie
$\whcalm \times_\bfg S \deff \bigl( \whcalm \times S\bigr)/\bfg$. Second, we
must choose a stratification of $\scrm$ by strata where $N_u$ does not
``jump''. By the definition of $N_u$ such strata are exactly $\scrm_{=\mbfk} 
= \whcalm_{=\mbfk} / \bfg$ where there is no ``jump'' of
the cusp-order.

\medskip
Another application of the techniques used in the proof of \refthm{thm4.2.1}
is a local version of the theorem. Below $\scrp(\Delta, X)$ denotes the Banach 
space of pseudoholomorphic maps between the unit disc $\Delta$ with the standard 
structure $J\st$ and $X$, \ie $\scrp(\Delta, X) = \{ (u,J) \in L^{1,p}(\Delta,
 X) \times \scrj: \dbar_{J\st,J} u=0 \}$. 

\newlemma{lem4.2.8n} \sli For any given integer $k\ge1$ the set 
\begin{equation}
\eqqno(4.2.4n)
\llap{$\scrp_k($}\Delta,0; X) \deff
\{ (u,J) \in \scrp(\Delta, X) : \ord_{z=0}(du) \ge k \}
\end{equation}
is a $C^\ell$-smooth submanifold of $\scrp(\Delta, X)$ of 
real codimension $2kn$, $n\deff \dimc X= \half \dimr X$, with tangent
space
\begin{equation}
\eqqno(4.2.5n)
T_{(u,J)}\scrp_k(\Delta,0; X) = 
\{ (v,\dot J) \in T_uL^{1,p}(\Delta, X)\times T_J\scrj : 
D_{u,J}v=0, j^k(v(z)-v(0)) =0 \}.
\end{equation}
\end{lem}

\newsubsection[4.2s]{Curves with prescribed secondary cusp index}
Recall that by \lemma{lem1.2.3} for a pseudoholomorphic map $u:(S,J_S) \to (X,J)$ 
with cusp order $k$ at $z^*\in S$ the jet $j^{2k+1}u(z^*)$ is well-defined. As we 
shall see, the part of the jet $j^{2k+1}u(z^*)$ invariant under reparameterization 
plays an important role for determining the type of critical points on
moduli spaces (see \refsubsection{3.3}).

\newdefi{def4.2s.1} \sli Let $u:(S,J_S) \to (X,J)$ be pseudoholomorphic map with
a cusp of order $k\deff \ord_{z^*} du$ at $z^*\in S$, $\pr_N: E_u \to N_u$
the projection to the normal bundle, and $z$ a local holomorphic coordinate
on $S$ centered at $z^*$. Define the {\sl secondary cusp index 
$l$ of $u$ at $z^*\in S$} by setting $l\deff k$ if $\pr_N \scirc j^{2k+1}
u(z^*)$ is zero polynomial and $l\deff \ord_{z=0} \pr_N \scirc j^{2k+1} 
u(z^*) -k-1$ otherwise.

\slii For a given $m$-tuple $\mbfk=(k_1, \ldots, k_m)$ of prescribed
orders of cusps we consider $m$-tuples $\mbfl= (l_1, \ldots, l_m)$ 
with $0 \le l_i \le k_i$ and set $|\mbfl| \deff \sum_i l_i$. 
Define the {\sl moduli space $\scrm_{\mbfk, \mbfl}$ of pseudoholomorphic maps
with cusps of given order and secondary index $(\mbfk, \mbfl)$} as the set of 
$[u,J, \mbfz] \in \scrm_\mbfk$ such that $\ord_{z^*_i} du =k_i$
and the secondary cusp index of $u$ at $z^*_i$ is at least $l_i$. Set
\begin{equation}
\eqqno(4.2s.1)
\whcalm_{\mbfk, \mbfl} \deff \{ (u, J_S, J, \mbfz) \in \whcalm_{=\mbfk}
\;:\; [u, J, \mbfz] \in \scrm_{\mbfk, \mbfl} \}.
\end{equation}
\end{defi}

\smallskip
\newthm{thm4.2s.1} The space $\scrm_{\mbfk, \mbfl}$ is a closed $C^{\ell-1}
$-smooth submanifold of $\scrm_{=\mbfk}$ of codimension $2(n-1)|\mbfl|$.
\end{thm}

\state Remark.
The meaning of the notion of secondary cusp index can be explained as follows.
One expects that for a $J$-holomorphic map $u:S \to X$ with a cusp of order $k=
\ord_{z^*}du$ at $z^* \in S$ the polynomial $\pr_N \scirc j^{2k+1} u(z^*)$ has 
vanishing order $k+1$. Thus the secondary cusp index $l$ is the order of deviation 
from this condition. The content of \refthm{thm4.2s.1} is that for generic 
$[u,J;\mbfz] \in \scrm_{=\mbfk}$ there is no deviation and that the space of 
curves with cusps of prescribed degeneration order is of expected codimension.

We note also that the range $0\le l_i \le k_i=\ord_{z^*_i}du$ is the maximal 
one where the secondary cusp index is well-defined: The higher order terms
of $\pr_N \scirc du$, as well as the coefficients of $du$ (considered as a 
holomorphic section of $T^*S \otimes E_u$), depend on the choice of the local
holomorphic coordinate $z_i$ centered at $z^*_i\in S$.

\proof We maintain the notation of \lemma{lem4.2.3n}. Now, for any $(u,J_S,J; 
\mbfz)\in \whcalm_{=\mbfk}$,  $\mbfz = (z^*_1,\ldots, z^*_m)$, the jets 
$j^{2k_i+1} u(z^*_i) \in j^{2k_i+1} (TS_i, E_i)_{(u,J_S, J; \mbfz)}$ are 
well-defined and depend $C^{\ell-1}$-smoothly on $(u,J; \mbfz)$. By 
\lemma{lem4.2.7n}, for any $i=1,\ldots,m$ the formula $(N_i)_{(u,J_S, J; 
\mbfz)} \deff (N_u)_{z^*_i}$ defines a $C^{\ell-1}$-smooth bundle $N_i$ over 
$\whcalm_{=\mbfk}$ with the projection $\pr_N : E_i \to N_i$. 
This yields the compositions $\pr_N \scirc j^{2k_i+1}u(z^*_i) \in j^{2k_i+1} 
 (TS_i, N_i)_{(u,J_S, J; \mbfz)}$ which depend $C^{\ell-1}$-smoothly on 
$(u,J_S,J; \mbfz)$. Thus we obtain a $C^{\ell-1}$-smooth bundle 
$$
\bigoplus\nolimits_{i=1}^m j^{k_i+1, k_i+l_i+1}(TS_i, N_i)
_{(u,J_S, J; \mbfz)}
$$
over $\whcalm_{=\mbfk}$ of rank $2(n-1)|\mbfl|$ and a $C^{\ell-1}$-smooth 
section 
$$
\yps^N_{\mbfk,\mbfl} \deff
(\pr_N \scirc j^{k_i+1, k_i+l_i+1} u(z^*_i))_{i=1}^m.
$$
Observe that $\whcalm_{\mbfk,\mbfl}$ is defined in $\whcalm_{=\mbfk}$ as the 
zero set of $\yps^N_{\mbfk,\mbfl}$. It follows from \lemma{lem4.2.3n} that
$\yps^N_{\mbfk,\mbfl}$ is transversal to the zero section. Consequently,
$\whcalm_{\mbfk,\mbfl}$ is a submanifold of $\whcalm_{=\mbfk}$ of codimension 
$2(n-1)|\mbfl|$. The claim of the theorem follows now by taking the 
$\bfg$-quotient.\qed

\newsubsection[4.3a]{Curves with cusps of prescribed type}
In this paragraph we give a construction of $J$-curves of any given cusp type,
completing the result of Micallef and White. In particular, we obtain a
more direct and constructive proof of \lemma{lem1.2.1} without referring
to local structure of minimal surfaces, as  is done in \cite{Mi-Wh}. Then 
we show that the set of cusp-curves with prescribed cusp type is a Banach 
submanifold of the total moduli space and compute its codimension.

\smallskip
Let $J$ be an almost complex structure on the ball $B \subset \cc^n$ such that 
$J(0) =J\st(0)$. We assume that $J$ is $C^\ell$-smooth with $\ell \ge 2$.
First we consider the local structure of multiple maps.

\newlemma{lem4.3a.1} Let $u: \Delta \to B$ be a non-constant $J$-holomorphic 
map with $u(0)=0 \in B$. Then there exist a radius $r>0$, a uniquely defined 
$\nu\in \nn$, and a non-multiple $J$-holomorphic map $u': \Delta(r^\nu) \to B$
such that $u(z) = u'(z^\nu)$ for $z\in \Delta(r)$. 
\end{lem}

\proof By \lemma{lem1.2.4}, $u(z) = v \cdot z^\mu + O(|z|^{\mu + \alpha})$
with some $v\in T_0B=\cc^n$, positive $\mu \in\nn$, and $\alpha>0$. If $\mu 
=1$, then $u$ is already non-multiple in some $\Delta(r)$ and there are
nothing to prove. Thus we may assume that $\mu \ge 2$. 

Take a sufficiently small $\rho_0>0$ and consider $U\deff u\inv(B(\rho_0))$.
By the first part of \lemma{lem1.2.4}, $U$ is a disc and $u$ is an immersion 
in $U\bs \{0\}$. Using the second part of \lemma{lem1.2.4} it is not difficult
to show that $u(U \bs \{0\})$ is an immersed $J$-holomorphic punctured disc 
in $B$. Therefore the restriction $u\ogran_U$ is a composition 
of a non-multiple $J$-holomorphic map and a covering branched only in 
$0\in U$.\qed

\medskip
It is known that any non-multiple holomorphic map $u: \Delta \to \cc^n$, 
in appropriate holomorphic coordinates on $\Delta$ and $\cc^n$,
has locally the form
\begin{equation}
\notag
u(z) = \sum_{i=0}^l v_i z^{p_i},
\end{equation}
with the following properties. $p_0 = \ord_0(du) +1$, $v_0 \not=0$, the vectors
$v_i \in \cc^n$ are linearly independent of $v_0$ for $i>0$, and $\gcd(p_0,\ldots, 
p_l)=1$. We want to establish a similar result for pseudoholomorphic curves, 
replacing the operation $u_{i-1}(z) \mapsto u_{i-1}(z) + v_i z^{p_i}$ by 
$u_{i-1}(z) \mapsto \dfrm_{p_i}(u_{i-1}, v_i)$. 

\newlemma{lem4.3a.2} Let $B \subset \cc^n$ be the unit ball, $J$ a $C^\ell$-smooth
almost complex structure on $B$ with $J(0)= J\st$, and $u: \Delta \to B$ a 
non-multiple $J$-holomorphic map such that $u(z)= v_0 z^{p_0}  +o(z^{p_0})$ for 
some $p_0>1$ and $v_0 \not=0 \in \cc^n$. Take a divisor $d>1$ of $p_0$ and 
denote by $\eta$ a primitive $d$-th root of unity. Then there exist an integer 
$q>0$, a vector $v\in \cc^n$, and a complex polynomial $\psi(z)$ such that

\sli $q$ is {\sl not} a multiple of $d$;

\slii $v$ is $\cc$-linearly independent of $v_0;$ in particular, $v \not=0;$

\sliii $\psi(z)= z + o(z)$ and $\deg\psi(z)\le q;$ 

\sliv $u(\eta z) = u( \psi(z))  + z^q \cdot v + o(z^q).$ 
\lineeqqno(4.3a.2a)

\end{lem}

\proof Denote by $v_0^\perp\subset \cc^n$ a complex orthogonal complement to 
$v_0$, and by $B^\perp(\rho)$ the ball of radius $\rho$ in $v_0^\perp$. 
Note that we can canonically identify the space $v_0^\perp$ with the fiber 
$(N_u)_{z=0}$ of the normal bundle $N_u$ of $u$ (see \refdefi{def1.5.1}). 
Fixing a holomorphic frame $w_1(z), \ldots w_{n-1}(z)$ of $N_u$ we can
identify $\Delta \times v_0^\perp$ with the total space of $N_u$ over $\Delta$
and use $(z, w_1, \ldots w_{n-1})$ as coordinates in $\Delta \times v_0^\perp$.
Denote by $J\st$ the standard integrable complex structure in $\Delta  \times 
v_0^\perp$. It coincides with the canonical holomorphic structure in $N_u$. Set
$$
U_{r,\rho} \deff \Delta(r) \times B^\perp(\rho)
$$
Fix a holomorphic splitting $F_0: N_u \to E_u$ of the projection $\pr_N : E_u 
\to N_u$. We shall identify $N_u$ as a subbundle of $E_u$ by means of $F_0$.
Define the map $F: U_{r,\rho}  \to \cc^n$ as the composition
$$
(z,w) \mapsto F_0(z)(w) \in (E_u)_z = T_{u(z)}B=\cc^n  
\mapsto F(z,w) \deff u(z) + F_0(z)(w).
$$
It is not difficult to see that for sufficiently small $r$ and $\rho$ the map 
$F=F(z,w)$ takes values in $B$ and has the following properties:

\begin{itemize}
\item $F(z,w)$ is $C^1$-smooth;
\item $F(z,0) = u(z)$ and $\nabla_{\dot w}F(z,0)=\dot w$; or more precisely,
$\nabla_{\dot w}F(z,0)=F_0(z)(\dot w)$;
\item the pulled-back structure $\ti J \deff F^*J$ coincides with $J\st$ along
the set $\check\Delta \times \{0\}$, \ie 
\begin{equation}
\ti J(z,0) = J\st(z,0).
\eqqno(4.3a.3a)
\end{equation}
\end{itemize}

From \eqqref(4.3a.3a) we obtain a uniform estimate 
\begin{equation}
\eqqno(4.3a.4a)
\bigl|\ti J(z,w) - J\st(z,w)\bigr| \le C \cdot |w|.
\end{equation}
Further, $\eta^{p_0}=1$ obviously gives $u(\eta z) -u(z)= o(z^{p_0}) = o(z^d)$. 
This implies that for sufficiently 
small $r'$ we can represent $u(z)$ in the form $u(z) = F(\zeta(z), \ti w(z))$ 
with uniquely defined $C^1$-smooth $\zeta(z): \Delta(r') \to \Delta(r)$ and 
$\ti w(z): \Delta(r') \to B^\perp(\rho)$ fulfilling the condition $\zeta(z)= z 
+o(z)$. Further, $\ti w(z)= o(z^d)$.

Set $\ti u(z) \deff \bigl( \zeta(z), \ti w(z) \bigr)$. We obtain a $C^1$-smooth 
map $\ti u(z): \Delta(r') \to U_{r,\rho}$, for which
$$
|J\st(\ti u(z)) - \ti J(\ti u(z)) | = 
|J\st (\zeta(z), \ti w(z))- \ti J(\zeta(z), \ti w(z)) | 
\le C'\cdot |\ti w(z)|.
$$
Consequently
$$
\bigl|\dbar_{J\st} \ti u(z) \bigr|= 
\bigl|\dbar_{J\st} \ti u(z) - \dbar_{\ti J} \ti u(z) \bigr|= 
\bigl|\bigl(J\st(\ti u(z)) - \ti J(\ti u(z)) \bigr) \d_y\ti u(z)\bigr|
\le C''\cdot |\ti w(z)|,
$$
or explicitly for components $\zeta(z)$ and $\ti w(z)$
\begin{align}
\bigl|\dbar_{J\st} \ti w(z) \bigr| & \le C''\cdot |\ti w(z)|;
\eqqno(4.3a.5a)\\
\bigl|\dbar_{J\st} \zeta(z) \bigr| & \le C''\cdot |\ti w(z)|.
\eqqno(4.3a.6a)
\end{align}
Observe that $\ti w(z)$ is not identically zero. Otherwise we would obtain
that $u(\eta z) = u(\zeta(z))$, which would contradict the condition of 
non-multiplicity of $u(z)$.

Hence, by \lemma{lem1.2.1}, $\ti w(z) = z^q v + o(z^q)$ and $\zeta(z)= \psi(z)
+ o(z^q)$ for some $q>0$, non-zero $v\in v_0^\perp$, and a complex polynomial 
$\psi(z)$ of degree $\le q$. Substituting these relations in $u(z)= 
F(\ti u(z))$ we obtain \eqqref(4.3a.2a).

Finally, the identity $\sum_{j=1}^d \bigl(u(\eta^j z) - u(\eta^{j-1} z) \bigr)
\equiv 0$ together with \eqqref(4.3a.2a) implies that $\sum_{j=1}^d (\eta^{j-1}
z)^q \cdot v =0$. Thus $\sum_{j=1}^d \eta^{jq} =0$ which is possible \iff $q$
is not a multiple of $d$.  \qed

\medskip
Iterating the construction of \lemma{lem4.3a.2}, we obtain

\newcorol{cor4.3a.3} Let $B \subset \cc^n$ be the unit ball, $J$ a $C^\ell
$-smooth almost complex structure in $B$ with $J(0)= J\st$, and $u: \Delta \to B$ 
a non-multiple $J$-holomorphic map with $u(0)=0$.

Then there exist uniquely defined sequences $(p_0, p_1, \ldots, p_l)$ and 
$(d_0, d_1, \ldots, d_l)$ of positive integers with the following properties:

\sli $p_0 = \ord_0 du +1$, so that $u(z) = z^{p_0} v_0 + o(z^{p_0})$ with 
non-zero $v_0 \in \cc^n$;

\slii $d_i =  \gcd(p_0, \ldots, p_i);$

\sliii $p_i< p_{i+1}$, $d_i > d_{i+1}$, and $d_l=1;$ 
in particular, $p_{i+1}$ is not a multiple of $d_i$;
\lineeqqno(4.3a.7a)

\sliv if $\eta_i$ is the primitive $d_i$-th root of unity, then 
$$
u(\eta_i z) = u(\psi_i(z)) \cdot v_0 +  z^{p_{i+1}} \cdot v_{i+1}
+ o(z^{p_{i+1}})
$$
for appropriate complex polynomials $\psi_i(z)$ with $\psi_i(z)=z+o(z)$, and 
vector $v_{i+1}\in \cc^n$, $\cc$-linearly independent of $v_0$.
\end{corol}


\newdefi{def4.3a.1} \sli To any increasing sequence of positive integers
$1\le p_0 < p_1 < \cdots < p_l$ we associate the {\sl sequence of divisors}
$d_i \ge d_1 \ge \cdots \ge d_l$ defined by $d_i =  \gcd(p_0, \ldots, p_i)$.
In particular, $d_0 =p_0$. 

\slii A sequence $\vec p= (p_0,p_1,\ldots, p_l)$ of positive integer exponents
is called a {\sl cusp type} if $p_i$ and the associate divisors $d_i = \gcd
 (p_0, \ldots, p_i)$ satisfy the condition \eqqref(4.3a.7a). In the situation of 
\refcorol{cor4.3a.3}, the sequence $\vec p= (p_0,p_1, \ldots, p_l)$ is called
the {\sl cusp type of $u$ at $z=0$}, $p_i$ the {\sl critical exponents of $u$
at $z=0$}, and $\vec d =(d_i)$ the {\sl sequence of divisors of $u$ at $z=0$}. 

\sliii For a given cusp type $\vec p= (p_0,p_1, \ldots, p_l)$, an integer $p'$ 
is called an {\sl admissible exponent} if $p'$ equals $p_l$ or is of the 
form $p'= p_i+ j\cdot d_i$ for some $i=0,\ldots,l-1$ and $j=0,\ldots, l_i$, 
$l_i \deff \left[ {p_{i+1} - p_i \over d_i } \right]$. Thus all critical 
exponents are admissible and there are exactly $l_i$ non-critical admissible 
exponents between $p_i$ and $p_{i+1}$. Denote by $\vpp = (p'_0, \ldots, 
p'_{l'})$ the sequence of the admissible exponents ordered by growth. Its 
length  is $l' = l + \sum_{i=0}^{l-1} l_i = l + \sum_{i=0}^{l-1}\left[ 
 {p_{i+1} - p_i \over d_i } \right]$.

Note that the corresponding sequence of divisors $d'_j \deff \gcd(p'_0, 
\ldots, p'_j)$ consists of divisors $d_i$ of critical exponents, repeated $l_i+
1$ times. Vice versa, an admissible exponent $p'_j>p'_0=p_0$ is critical \iff 
$d'_j < d'_{j-1}$.
\end{defi}

\smallskip
\newthm{thm4.3a.4} Let $B \subset \cc^n$ be the unit ball, $J$ an almost complex 
on $B$ with $J(0)= J\st$, and $u: \Delta \to B$ a non-multiple $J$-holomorphic map
such that $u(0)= 0$. Further, let $\vec p= (p_0, \ldots, p_l)$ and $\vpp= (p_0, 
\ldots, p_{l'})$ be the sequences of critical and resp.\ admissible exponents 
of\/ $u$ at\/ $z=0$, and $\vdp =(d'_0, \ldots, d'_{l'})$ the corresponding 
sequence of divisors.

Then there exist a sequence $(v_0, \ldots, v_{l'})$ of vectors in $\cc^n$ (one 
$v_j$ for each $p'_j$), a complex polynomial $\phi(z)$, and a radius $r>0$, such 
that the following holds.

\sli $u(z)= z^{p_0}\cdot v_0 + o(z^{p_0})$; $v_0 \not=0$, $v_1, \ldots, v_{l'}$ 
are complex orthogonal to $v_0;$
\lineeqqno(4.3a.2)

\slii $\phi(z)=z + o(z)$ and $\deg \phi(z) \le p_l-p_0 +1;$
\lineeqqno(4.3a.3)

\sliii for appropriately chosen maps $\dfrm$, the recursive formula
\begin{equation}
\eqqno(4.3a.4)
u_j(z) \deff \dfrm_{p'_j/d'_j}\bigl(u_{j-1}(z^{d'_{j-1}/d'_j}), J; v_j\bigr) 
\quad j=0,1,\ldots,l'
\end{equation}
beginning from $u_{-1}(z) \equiv 0$ yields a sequence of well-defined 
$J$-holomorphic maps $u_j: \Delta(r^{d'_j}) \to B$ with the property
\begin{equation}
\eqqno(4.3a.5)
u(\phi(z)) - u_j(z^{d'_j}) = v_{j+1} z^{p'_{j+1}} + o(z^{p'_{j+1}}).
\end{equation}

\noindent
Moreover, such $v_j$ and $\phi(z)$ are uniquely defined. Further, $v_j$ is 
non-zero if $p'_j$ is  critical.
\end{thm}

\proof The choice of the maps $\dfrm_d $ ensuring that at each recursive step 
the right hand side of \eqqref(4.3a.5) is well-defined will be made below. Now
we assume that for given $j<l'$ we have constructed a $J$-holomorphic map $u_j:
\Delta(r^{d'_j}) \to X$ and a polynomial $\phi(z)$ such that $\phi(z)= z+o(z)$
and $u(\phi(z)) = u_j(z^{d'_j}) +o(z^{p'_j})$. Then by \lemma{lem1.2.4},
\begin{equation}
\eqqno(4.3a.6)
u(\phi(z)) = u_j(z^{d'_j}) + z^q w +  o(z^q)
\end{equation}
for some non-zero $w \in \cc^n$ and $q>p'_j$. Represent $w\in \cc^n$ in the form 
$w= a \cdot v_0 + w'$ and replace $\phi$ by $\phi'(z) \deff \phi(z)- {a\over p_0}
\cdot z^{q-p_0 +1}$. The relations $u(z) = z^{p_0} v_0 + o(z^{p_0})$ and 
\eqqref(4.3a.6) yield
\begin{equation}
\eqqno(4.3a.6')
u(\phi'(z)) = u_j(z^{d'_j}) + z^q w' +  o(z^q).
\end{equation}
If $w' =0$, we can consider \eqqref(4.3a.6) with some $q'>q$. Thus we may assume
that $w' \not=0$. Moreover, we see that $\phi(z)$ is defined uniquely by 
\eqqref(4.3a.5) up to degree $p'_j -p_0 +1$.

Denote by $\eta_j$ the primitive $d'_j$-th root of unity. Then 
by \lemma{lem4.3a.2},
\begin{equation}
\eqqno(4.3a.7)
u(\eta_j z) = u(\psi_j(z)) + v z^p + o(z^p)
\end{equation}
for an appropriate polynomial $\psi_j(z) =z +o(z)$ and $v$ linearly independent 
of $v_0$. Moreover, $p$ is the first critical exponent after $p'_j$ in the 
sequence $\vpp$ of the admissible exponents of $u(z)$ at $z=0$. In particular, $p$ 
is not a multiple of $d'_j$. Set $\hat\phi_{j+1}(z) \deff \eta_j\inv\phi_{j+1}
 (\eta_j z)$. Then we obtain $\hat\phi_{j+1}(z) = z+ o(z)$ and
\begin{align}
u(\phi_{j+1}(\eta_j z)) &= u(\eta_j\hat\phi_{j+1}(z)) = 
u(\psi_j(\hat\phi_{j+1}(z))) + z^p v + o(z^p)
\notag\\
&= u(\phi_{j+1}(\hat\psi_j(z))) + z^p v + o(z^p),
\eqqno(4.3a.8)
\end{align}
where $\hat\psi_j(z)$ is a polynomial with $\hat\psi_j(z)= z +o(z)$ and 
$\hat\psi_j(\hat\phi_{j+1}(z)) = \phi_{j+1}(\hat\psi_j(z)) +o(z^p)$.
Substitution of \eqqref(4.3a.6') in \eqqref(4.3a.8) together with the identity
$\eta_j^{d'_j}=1$ yields
\begin{equation}
\eqqno(4.3a.9)
u_j( z^{d'_j})  + \eta_j^q z^q w' = 
u_j\bigl( \hat\psi_j^{d'_j}(z) \bigr) 
+ z^q w' + v z^p + o(z^{\min(p,q)}).
\end{equation}
Further, since $\hat\psi_j(z)= z +o(z)$, we can find a polynomial $\ti\psi_j(z)$ 
with the properties $\ti\psi_j(z)= z +o(z)$ and $\hat\psi_j^{d'_j}(z)= \ti\psi_j
 (z^{d'_j})$. For such $\ti\psi_j(z)$, the relation \eqqref(4.3a.9) transforms to
\begin{equation}
\eqqno(4.3a.10)
u_j( z^{d'_j})  =
u_j\bigl( \ti\psi_j(z^{d'_j}) \bigr) 
+ (1 - \eta_j^q) z^q w' + v z^p + o(z^{\min(p,q)}).
\end{equation}
Assume that $q<p$. Then $q$ is a multiple of $d'_j$. In particular, $q\ge 
p'_{j+1}$. In the case $q> p'_{j+1}$ we simply set $v_{j+1} \deff 0$ and obtain
the relation \eqqref(4.3a.5). In the case $q= p'_{j+1}$ we set $v_{j+1} \deff w'$
and come to the relation \eqqref(4.3a.5) again. The case $p<q$ is impossible since
$p$ is not a multiple of $d'_j$. 

In the remaining case $q=p$ we have two subcases, $p'_{j+1} <p$ and $p'_{j+1} 
=p$. Then we set $v_{j+1} \deff 0$ or respectively $v_{j+1} \deff w'$ and obtain 
\eqqref(4.3a.5) from \eqqref(4.3a.6').

\medskip
Now we construct the maps $\dfrm_{p'_j/d_j}$ with the desired properties. The idea
is to rescale the maps $u_j$ making the norms $\norm{du_j}_{L^2}$ sufficiently 
small and obtaining a recursive apriori estimate on $v_j$. For this fix some 
$r \in\; ]0,1[$ and maps $\wt \dfrm_p$ with the properties listed in 
\refdefi{def4.1.1}. Then the substitutions $\ti u(z) \deff u(r z)$, 
$\ti u_j(z) \deff u_j(r^{d'_j} z)$, 
$\ti v_j \deff r^{p'_j}v_j$, and $ \ti \phi_j(z) \deff r\inv \phi_j(r z)$
transform \eqqref(4.3a.4) and \eqqref(4.3a.5) into recursive relations
\begin{equation}
\eqqno(4.3a.a6)
\llap{$\ti u_j(z)= $}\;\wt\dfrm_{p'_j/d'_j}
\bigl(\ti u_{j-1}(z^{d'_{j-1}/d'_j}), J; 
\ti v_j\bigr),
\end{equation}
\begin{equation}
\eqqno(4.3a.a7)
\ti u(\ti \phi(z)) - \ti u_j(z^{d'_j})= 
\ti v_{j+1} z^{p'_{j+1}} + o(z^{p'_{j+1}}).
\end{equation}
for $J$-holomorphic maps $\ti u_j: \Delta \to B$. Note that $\norm{d\ti u}
_{L^2(\Delta)} = \norm{du}_{L^2(\Delta(r))}$ will be arbitrarily small for $r$ 
small enough. Choosing an appropriate $r\ll1$ and using induction, on can obtain 
sufficiently small upper bounds on $\ti v_j$, ensuring that \eqqref(4.3a.a6) is
well-defined for $j=0,\ldots, l$. For such $r$, we define $\dfrm_{p'_j/d_j}$
by the reverse substitutions in \eqqref(4.3a.a6).
\qed

\smallskip
\state Remark. For almost complex surface, \ie in the case $n=2$, the critical 
exponents determine a topological type of a cusp. In particular, under 
hypotheses of \refthm{thm4.3a.4}, the intersection of the image $u(\Delta)$ 
with the sphere $S^3_r$ of a sufficiently small radius $r>0$ is an iterated 
toric knot $\gamma$ transversal to the 2-plane distribution $\xi$ on $S^3_r$ 
given by $\xi_x \deff T_xS^3_r \cap J(x) T_xS^3_r$. Thus the Bennequin index 
$\beta(\gamma, \xi)$ is well-defined. We refer to \cite{Iv-Sh-1} for the proof
of the formula $\beta= 2\delta -1$ relating the Bennequin index $\beta$ of 
$\gamma$ and the nodal number $\delta$ of $u(\Delta)$ in $0\in B$. On the 
other hand, $\delta$ can be computed by the formula 
\begin{equation}
\eqqno(4.3a.11a)
\delta = \sum_{i=1}^m (d_{i-1}- d_i)(p_i -1),
\end{equation}
see \cite{Rf} or \cite{Mil}. In the higher dimensional setting, 
\ie for $n\ge3$, the topological type of the 
cusp $u(\Delta)$ is not determined by the critical exponents and depends on 
additional information encoding further linear relations between $v_j$. For 
example, the condition {\sl $v_2$ and $v_1$ are linearly dependent} defines a 
proper subset in $\scrp_{\vec p}(\Delta, 0;B)$. Moreover, using the techniques
of this paragraph one can show that this subset is a $C^{\ell-1}$-smooth 
submanifold in $\scrp_{\vec p}(\Delta, 0;B)$. Details can be recovered by an 
interested reader.

\smallskip
\newthm{thm4.3a.5} Let $B \subset \cc^n$ be the unit ball, $\vec p= (p_0,\ldots,
p_l)$ a cusp type, $\vpp= (p_0, \ldots, p_{l'})$ the corresponding sequences of 
admissible exponents, and $\vdp =(d'_0, \ldots, d'_{l'})$ the sequence of 
divisors associated with $\vpp$. Then the set
\begin{equation}
\eqqno(4.3a.11)
\llap{$\scrp_{\vec p} (\Delta, 0$};B) \deff 
\bigl\{ (u,J) \in \scrp_{p_0-1} (\Delta, 0; B):
u \text{ \rm has a cusp type $\vec p$ in $z=0$}\; \bigr\}
\end{equation}
is a $C^{\ell-1}$-smooth submanifold of $\scrp_{p_0-1} (\Delta, 0;B)$ of real
codimension $2(n-1)(p_l -p_0-l')$.
\end{thm}

Note that $p_0=p'_0$ and $p_l= p'_{l'}$.

\proof Let $u:\Delta(r) \to B$ be a $J$-holomorphic map, $r>0$, and let
\begin{equation}
\scrp_k(\Delta, u; B,J) \deff 
\{ u'\in \scrp(\Delta; B, J) : u(z) -u'(z) = o(z^k) \} 
\notag
\end{equation}
By \refthm{thm4.2b.1}, $\scrp_k(\Delta, u; B,J)$ is a $C^{\ell-1} $-smooth 
submanifold of $\scrp(\Delta; B, J)$ of codimension $2n(k+1)$. For $l>k$ it 
follows that $\scrp_l(\Delta, u; B,J)$ has codimension $2n(l-k)$ in $\scrp_k(
\Delta, u; B,J)$. Moreover, if $J_y$ is a $C^\ell$-smooth family of almost complex 
structures in $B$ parameterized by a (Banach) manifold $\scry$ and $u_y\in 
L^{1,p}(\Delta(r), B)$ a $C^{\ell-1}$-smooth family of $J_y$-holomorphic maps, 
then $\cup_{y\in \scry} \scrp_k(\Delta, u_y; B,J_y)$ is a $C^{\ell-1} $-smooth 
manifold.

\smallskip
For a given $(u^*, J^*) \in \scrp_{\vec p} (\Delta, 0;B)$, let $v^*_0, \ldots, 
v^*_{l'}$ and $\phi^*(z) = z + \phi^*_2 z^2 + \phi^*_3 z^3+ \cdots$  be the
parameters of $u^*$ constructed in \refthm{thm4.3a.4}. Define $\scry$ to be
the space of small deformations of $v^*_j$ and $\phi^*_i$. This means that $y
\in \scry$ is a tuple $(v_0, \ldots, v_{l'}; \phi_2, \ldots, \phi_{p_l-p_0+1})$
with $v_j \in \cc^n$ and $\phi_i \in \cc$ satisfying $|v_j - v^*_j| <\epsi$ 
and $|\phi_i - \phi^*_i| < \epsi$ with $\epsi$ sufficiently small. Further, 
let $U$ denote a sufficiently small neighborhood of $J^*$ in the space of 
$C^\ell$-smooth almost complex structures in $B$. For $y=(v_0, \ldots;\phi_2, 
\ldots) \in \scry$ and $J\in U$ we construct the maps $u_{y,J; j}$, $j=0, 
\ldots, l'$, using the recursive relation \eqqref(4.3a.4) and set $\phi_y(z) 
\deff z + \phi_2 z^2 + \cdots +\phi_{p_l -p_0 +1} z^{p_l -p_0 +1}$. Then for 
$|z| < r' \ll1$ the inverse map $\phi_y\inv(z)$ is well-defined and 
holomorphic. Define $u_{y,J}(z) \deff u_{y,J; l'}(\phi\inv(z))$. We obtain a 
$C^{\ell-1}$-smooth family of pseudoholomorphic maps $u_{y,J}: \Delta(r') \to B$ 
parameterized by $(y,J) \in \scry\times U$. 

Note that by \refthm{thm4.3a.4} every $(u,J) \in \scrp_{\vec p} (\Delta, 0;B)$
sufficiently close to $(u^*, J^*)$ lies in $\scrp_{p_l}(\Delta, u_{y,J}; B,J)$ 
for an appropriate $y\in \scry$, and such $y\in \scry$ is uniquely defined. 
Thus the union $\cup_{(y,J)\in \scry\times U} \scrp_{p_l}(\Delta, u_{y,J}; B,
J)$ is a local $C^{\ell-1}$-smooth chart for $\scrp_{\vec p} (\Delta, 0;B)$. 

Finally, note that the union $\cup_{(y,J)\in \scry\times U} \scrp_{p_0-1}
 (\Delta, u_{y,J}; B,J)$ is naturally isomorphic to $\scrp_{p_0}(\Delta, 0; B)
\times \scry$. Computing the number of parameters we obtain the codimension of
the imbedding $\scrp_{\vec p} (\Delta, 0;B) \hook \scrp_{p_0} (\Delta, 0;B)$. 
\qed

Globalizing \refthm{thm4.3a.5} we obtain

\newcorol{cor4.3a.6} Let $\vec \mbfp =(\vec p_1, \ldots, \vec p_m)$ be a
sequence of cusp types, $\vec p_i= (p_{i,0}, \ldots p_{i, l_i})$. Set 
$k_i \deff p_{i,0} -1$ and $\mbfk\deff (k_1, \ldots, k_m)$. Then the space
$$
\scrm_{\vec \mbfp} \deff \{ [u,J; z^*_1, \ldots z^*_m] \in \scrm_\mbfk:
u \text{ \sl has cusp type $\vec p_i$ in }z^*_i\; \}
$$
is a $C^{\ell-1}$-smooth submanifold of $\scrm_\mbfk$ of codimension
$2(n-1)\sum_{i=0}^m (p_{i,l_i} -p_{i,0}-l'_i)$.
\end{corol}

\medskip
\newsection[3]{Saddle points in the moduli space}

\newsubsection[3.1]{Critical and saddle points in the moduli space}
In application of the continuity method for constructing $J$-holomorphic
curves two main difficulties occur. The first one appears in the
proof of the ``closedness part'', when one tries to extend a deformation $[u_t, 
J_t] \in \scrm, t\in [0,t'[\;$ into the endpoint $t'$. This difficulty is 
connected with the fact that the projection $\pi_{\!\!\scrj}: \scrm \to \scrj$ 
is, in general, not proper. In particular, for a path $J_t\in \scrj, t\in [0,t']$ 
there may not exist a lift $[u_t, J_t]$ to $\scrm$, and the fibers
$\scrm_{\!\!\scrj} = \pi_J\inv (J)$ can be non-compact. Gromov's compactness
theorem (\cite{Gro}, see also  \cite{Iv-Sh-3}) gives a fiberwise topological
compactification of $\scrm$ by adding certain degenerate curves $C$.
However, for the moment we neglect this difficulty assuming we can avoid it in 
our case.

The second main difficulty appears in the proof of the ``openness part'' when one
tries to extend a lift $[u_t, J_t] \in \scrm$, $t \in[0, t']$, of a
path of $J_t\in \scrj, t\in [0, 1]$ to a bigger interval $t \in[0, t''[$ with
some $t''>t'$. Obviously, this difficulty can appear only if $[u_{t'}, J_{t'}]$ 
is a {\sl critical point} of $\pi _{\!\!\scrj}$, \ie when the differential
of the projection $d\pi_{\!\!\scrj}$ is not surjective in $[u_{t'}, J_{t'}]$.
Thus it is desirable to find conditions on the critical points of
$\pi_{\!\!\scrj}$ which ensure the existence of such a lift.

\smallskip
Assume additionally that the given path $h: [0,1] \to \scrj$, $h(t)\deff
J_t$, is $C^2$-smooth and transversal to $\pi_h: \scrm \to \scrj$, \ie
$\scrm_h$ is a manifold. Let $\pi_h: \scrm_h \to I$, $I\deff [0,1]$, be the
projection. Then we have a well-defined $C^1$-smooth bundle homomorphism
$d\pi_h: \scrm_h \to \pi_h^*(TI) \cong \rr$. Further, $d\pi_h$ vanishes
exactly at critical points of $\pi_h$ and, by \lemma{lem1.3.1}, at each such
point $p\deff [u,h(t)]$ we have a well-defined quadratic form $\nabla
d\pi_h(p): T_p\scrm_h \to \rr\cong T_tI$. For our purpose it is
sufficient to show that each critical point $p$ is a {\sl saddle}, \ie
the quadratic form $\nabla d\pi_h(p)$ has at least one positive and one
negative eigenvalue.

\smallskip
It turns out that this condition depends only on the geometry of the projection
$\pi_{\!\!\scrj}:\scrm \to \scrj$ at $p=[u,h(t)]$, and not on the particular choice
of a transversal map $h:I \to \scrj$. In more detail, the situation is as
follows.

\smallskip
First, since $\scrm$ is $C^\ell$-smooth with $\ell\ge2$, the map $\pi_{\!\!\scrj}:
\scrm \to \scrj$ defines a $C^1$-smooth homomorphism of Banach bundles $d\pi_{\!\!
\scrj}: T\scrm \to \pi_{\!\!\scrj}^*(T\scrj)$. \refcorol{cor2.2.2} relates
the (co)kernel of $d\pi_{\!\!\scrj}$ for a given $[u,J] \in \scrm$ with $\sfh^i(S,
\scrn_u)$, and \lemma{lem1.3.1} provides a well-defined bilinear map
\begin{equation}
\Phi=\Phi_{[u,J]} \deff \nabla d \pi_{\!\!\scrj}:
T_{[u,J]}\scrm \times \sfh^0(S, \scrn_u)
\to \sfh^1(S, \scrn_u).
\eqqno(3.1.1)
\end{equation}

The situation remains essentially the same if we consider a relative moduli
space $\scrm_h= Y\times_h\scrm$ with a $C^\ell$-smooth map $h: Y \to \scrj$
transversal to $\pi_{\!\!\scrj}$. Indeed, one can easily see that for the natural
projection $\pi_h: \scrm_h \to Y$ and a point $[u,y] \in \scrm_h$ with $h(y)
\ddef J$ one has the natural isomorphisms
\begin{equation}
\matrix
\ker(d\pi_h :T_{[u,y]}\scrm_h \to T_yY)
&\cong& \ker(d\pi_{\!\!\scrj} :T_{[u,J]}\scrm \to T_J\scrj)
&\cong \sfh^0_D(S, \scrn_u),
\cr
\coker(d\pi_h :T_{[u,y]}\scrm_h \to T_yY)
&\cong& \coker(d\pi_{\!\!\scrj} :T_{[u,J]}\scrm \to T_J\scrj)
&\cong \sfh^1(S, \scrn_u).
\endmatrix
\eqqno(3.1.2)
\end{equation}
Further, the relation between $\Phi=\nabla d \pi_{\!\!\scrj}$ and $\nabla d \pi_h$
is given by the following

\newlemma{lem3.1.1}
\sli The isomorphism $\coker(d\pi_h)\cong \sfh^1(S,
\scrn_u)$ is induced by composition $T_yY \buildrel dh \over \lrar T_J\scrj
\buildrel \barr \Psi_{u,J} \over {\relbar\!\!\relbar\!\!\lrar}
\sfh^1(S, \scrn_u)$.

\smallskip
\baselineskip=14pt
\slii The bilinear map $\nabla d \pi_h: T_{[u,y]}\scrm_h \times \sfh^0_D(S,
\scrn_u) \to \sfh^1(S, \scrn_u)$ is induced by the composition $T_{[u,y]}
\scrm_h \hook T_{[u,J]}\scrm \oplus T_yY \twoheadrightarrow T_{[u,J]}\scrm$
and the bilinear map $\Phi: T_{[u,J]}\scrm \times \sfh^0_D(S,  \scrn_u) \to
\sfh^1(S, \scrn_u)$.
\end{lem}

\medskip
Summing up, we obtain the following situation in the most important case $Y=I$.

\newlemma{lem3.1.2}
For a map $h: I \to \scrj$ transversal to $\pi_J$,
the singular points of the projection $\pi_h: \scrm_h \to I$ are exactly those
$[u,t] \in \scrm_h$ for which $\sfh^1(S, \scrn_u)=\rr$.

For such $[u,t] \in \scrm_h$ with $J\deff h(t)$ one has the equality
$T_{[u,t]} \scrm_h = \sfh^0_D(S, \scrn_u)$ and the isomorphism $\barr\Psi
_{[u,J]}: T_tI \buildrel \cong \over \lrar \sfh^1(S, \scrn_u)$. Moreover,
the quadratic form $\Phi_{[u,J]}: \sfh^0_D(S, \scrn_u)\to \sfh^1(S, \scrn_u)$
equals to the composition of the quadratic form $\nabla d\pi_h: T_{[u,t]}
\scrm_h  \to T_tI$ with $\barr\Psi_{[u,J]}: T_tI \to \sfh^1(S, \scrn_u)$.
\end{lem}

\newcorol{cor3.1.3} The nullity, rank and signature of $\Phi_{[u,J]}$ and
$\nabla d\pi_h$ coincide.
\end{corol}

\newdefi{def3.1.1} Let $Q$ be a quadratic form defined on a
(finite-dimen\-sional) vector space $V$ and taking values in a vector space $W$
with $\dimr W=1$. Define the {\sl saddle index of $Q$} by $\sind Q \deff
\min \{ \ind_+ Q, \ind_- Q \}$, where $\ind_\pm Q$ are respectively the
positive and negative indices of $Q$ \wrt some (in fact, any) orientation of
$W$. For a critical point $[u,J]\in \scrm$ with $\sfh^1(S, \scrn_u) \cong \rr$,
call $\sind \Phi_{[u,J]}$ the {\sl saddle index of $[u,J]$}. A point
$[u,J]\in \scrm$ is a {\sl saddle point} of the moduli space $\scrm$ \iff
$\sind \Phi_{[u,J]}$ is strictly positive.
\end{defi}

\newsubsection[3.2]{Second variation of the $\dbar$-equation}
To find saddle points of $\scrm$ we need to find an explicit formula for the
form $\Phi$ in \eqqref(3.1.1). Note that, since the space $\scrp$ appears as 
the zero-set of the $\dbar$-equation \eqqref(1.1.1), the description of the 
tangent space $T\scrp$ is given by the variation of the $\dbar$-equation. 
Similarly, we show that the form $\Phi$ is essentially the part of the second 
variation of the $\dbar$-equation invariant \wrt the choice of a connection 
being used.

\smallskip
Let $[u,J]\in\scrm$ be represented by $(u,J_S,J)\in\whcalm$. Recall the
description of $T_{(u,J_S,J)}\whcalm$ given in \eqqref(2.2.4). Moreover, since
$du$ is non-vanishing at a generic point, $\dot J_S$ is determined by $v$ and
$\dot J$. Note that the tangent space to an orbit $\bfg\cdot(u,J) \subset\whcalm$
can be identified with  $du(\sfh^0(S,TS)) \subset \scre_{(u,J_S,J)}$. This
defines a subbundle of $\scre$ which we also denote by $du(\sfh^0(S,TS))$.
Thus we obtain the isomorphism
\begin{equation}
T_{[u,J]}\scrm\cong T_{(u,J_S,J)}\whcalm/du(\sfh^0(S,TS)).
\eqqno(3.2.1)
\end{equation}
Explicitly, the tangent space $T_{[u,J]}\scrm$ consists of triples
$([v], \dot J_S, \dot J)$ for which $(v, \dot J_S,\dot J) \in T_{(u,J_S,J)}
\whcalm$ and $[v]$ is the equivalence class $v + du\bigl(\sfh^0(S,TS) \bigr)$.

\newdefi{def3.2.1} Set
\begin{equation}
 \wh\scre_{(u,J_S,J)} \deff
\bigl(\scre_{(u,J_S,J)} /du(\sfh^0(S,TS) \bigr) \oplus \sfh^1(S,TS).
\eqqno(3.2.2)
\end{equation}
Recall the canonical isomorphism $T_{J_S}\ttt_g \cong \sfh^1(S,TS)$ given
by \eqqref(2.2.2). For $(u,J_S,J)\allowbreak\in \whcalm$ define the operator
\begin{equation}
\wh D=\wh D_{u,J}: \wh\scre_{(u,J_S,J)} \to  \scre'_{(u,J_S,J)}
\qquad
\wh D ([v], [I_S]) \deff Dv + J\scirc du \scirc I_S.
\eqqno(3.2.3)
\end{equation}
\end{defi}

\newlemma{lem3.2.1}
Formula \eqqref(3.2.2) defines a $C^\ell$-smooth Banach bundle $\wh\scre$ 
over $\scrm$ with the fiber $\wh\scre_{(u,J_S,J)}$
over $[u,J]\in \scrm$. The tangent bundle $T\scrm$ can be included in the
following exact sequence of bundles over $\scrm$
\begin{equation}
0\to T\scrm \buildrel \alpha \over\lrar
\wh\scre \oplus \pi_{\!\!\scrj}^*T\scrj
\buildrel \beta \over\lrar \scre' \to0,
\eqqno(3.2.4)
\end{equation}
where the homomorphisms $\alpha=(\alpha_1,\alpha_2)$ and $\beta=(\beta_1,
\beta_2)$ are given by
\begin{equation}
\matrix
\alpha_1([v],\dot J_S, \dot J) &\deff&
([v], [\dot J_S]) \in \wh\scre=
\scre/du\bigl(\sfh^0(S,TS)\bigr)\bigoplus \sfh^1(S,TS)
\cr
\alpha_2 &\deff& d\pi_{\!\!\scrj}: T_{[u,J]}\scrm \to T_J\scrj
\cr
\beta_1 &\deff& \wh D_{u,J}: \wh\scre_{(u,J_S,J)}\to\scre'_{(u,J_S,J)}
\cr
\beta_2 &\deff& \Psi_{u,J}: T_J\scrj\to\scre'_{(u,J_S,J)}
\endmatrix
\end{equation}
\end{lem}

\proof It is easy to show that $\sfh^0(S,TS)$ and $\sfh^1(S,TS)$ can
be considered as smooth bundles over $\whcalm$ equipped with the natural
$\bfg$-action. Then $du$ defines a $\bfg$-equivariant homomorphism between the
bundles $\sfh^0(S, TS)$ and $\scre$. Hence, using formula \eqqref(3.2.2), we 
can construct a bundle $\wh\scre_{\whcalm}$ over $\whcalm$ with the induced
$\bfg$-action. By \lemma{lem2.2.1}\.\slip, this is equivalent to the first
assertion of the lemma.

The exactness of \eqqref(3.2.4) follows from relations \eqqref(2.2.4) and 
(\ref{eq3.2.1}--\ref{eq3.2.3}).
\qed

\smallskip
\newlemma{lem3.2.2}
The homomorphisms $\alpha_1$ and $\beta_2$ yield isomorphisms
\begin{equation}
\sfh^0_D(S, \scrn_u) \cong \ker \alpha_2 \buildrel
\alpha_1 \over \cong \ker \beta_1
\quad\text{and}\quad
\sfh^1_D(S, \scrn_u) \cong \coker \alpha_2 \buildrel
\beta_2 \over \cong \coker \beta_1,
\eqqno(3.2.5)
\end{equation}
inducing the identity
\begin{equation}
\Phi_{u,J}= -\nabla \wh D: T_{[u,J]}\scrm \times \sfh^0_D(S, \scrn_u)
\to \sfh^1_D(S, \scrn_u).
\eqqno(3.2.6)
\end{equation}
\end{lem}

\proof The isomorphisms \eqqref(3.2.5) follow from definitions and
\refcorol{cor2.2.2}. Moreover, we can identify $\sfh^0_D(S, \scrn_u)$ with
$\ker\bigl( \wh D_{u,J}: \wh\scre_{(u, J_S, J)} \to \scre'_{(u, J_S, J)}
\bigr)$.

Let $i: \sfh^0_D(S, \scrn_u) \to T_{[u,J]}\scrm$ and $p:\scre'_{(u,J_S,J)} \to
\sfh^1_D(S, \scrn_u)$ denote the corresponding inclusion and projection. Fix
some connections on $T\scrm$, $\pi^*_{\!\!\scrj} T\scrj$, $\wh\scre$, and $\scre'$,
and denote all of them simply by $\nabla$. Covariant differentiation of the
relation $\beta_1 \scirc \alpha_1 + \beta_2 \scirc \alpha_2=0$ gives
\begin{equation}
\nabla\beta_1 \scirc \alpha_1 + \nabla\beta_2 \scirc \alpha_2 +
\beta_1 \scirc \nabla\alpha_1 + \beta_2 \scirc \nabla\alpha_2 = 0,
\end{equation}
which together with $\alpha_2\scirc i=0$ and $p\scirc \beta_1=0$ yields
\begin{equation}
  p\scirc \nabla\beta_1 \scirc \alpha_1\scirc i=
- p\scirc \beta_2 \scirc \nabla \alpha_2\scirc i.
\end{equation}
\qed

\newdefi{def3.2.2} Using the isomorphisms $\ker \bigl( \wh D_{u,J}: \wh
\scre_{(u,J_S,J)} \to \scre'_{(u,J_S,J)} \bigr) \cong \sfh^0_D(S, \scrn_u)$
from \eqqref(3.2.5) and $\sfh^1(S, TS) \cong T_{J_S} \ttt_g$ from \eqqref(2.2.2),
redefine
\begin{equation}
\sfh^0_D(S, \scrn_u) \deff \bigl\{\, ([v], I_S)\in \wh\scre_{(u,J_S,J)}
\oplus  T_{J_S} \ttt_g\;:\; Dv + J \scirc du \scirc I_S =0 \;\bigr\}.
\eqqno(3.2.7)
\end{equation}
Then the projection $\sfh^0_D(S, \scrn_u) \to \sfh^0_D(S, N_u)$ is given by
the formula $([v], I_S) \mapsto \pr_N(v)$ with $\pr_N: E_u \to N_u$
defined by \eqqref(1.5.2).
\end{defi}

\smallskip
Now assume that some symmetric connections on $TX$ and $TS$ are fixed.
They induce connections on $\scre$ and $\scre'$, on the tangent bundles
$T\whcalm$ and $T\scrm$, and so on. We shall use the same notation $\nabla$
for all these connections. Further, denote by $R^X(\cdot,\cdot;\cdot)$ the
curvature operator of the connection $\nabla$ on $X$.

\newlemma{lem3.2.3}
For $([v], \dot J_S, \dot J) \in T_{[u,J]}\scrm$
and $([w], I_S) \in \sfh^0_D(S, \scrn_u) \subset \wh\scre_{(u, J_S, J)}\!\!$

\smallskip
$\qquad
\bigl(\nabla_{([v],\dot J_S, \dot J)}\wh D\bigr) ([w],I_S)=
\term1{R^X(v,du ;w)} + \term2{J \scirc R^X(v, du\scirc J_S;w)} +
$
\lineeqqno(3.2.8)

\smallskip
$\qquad
+ \term3{\nabla_vJ\scirc \nabla w \scirc J_S} +
\term4{\nabla^2_{v,w}J \scirc du \scirc J_S} +
\term5{\nabla_wJ\scirc \nabla v \scirc J_S} +
\term6{\dot J\scirc \nabla w \scirc J_S} +
$

\smallskip
$\qquad
+ \term7{\nabla_w\dot J\scirc du \scirc J_S} +
\term8{J\scirc \nabla w \scirc \dot J_S} +
\term9{\nabla_wJ\scirc du \scirc \dot J_S} +
\term{10}{\nabla_vJ \scirc du \scirc I_S} +
$

\smallskip
$\qquad
+ \term{11}{J \scirc \nabla v \scirc I_S} +
\term{12}{\dot J \scirc du \scirc I_S}.
$
\end{lem}

\state Remark. The numerical subscripts on the various terms are for future 
reference.  

\proof Consider the bundle $\wt\scre\deff \scre \oplus \sfh^1(S,TS)$
over $\whcalm$ with the bundle homomorphism $\wt D: \wt\scre \to \scre'$,
$\wt D(w, [I_S]) \deff Dw + J \scirc du \scirc I_S$. We claim that for the
covariant derivative $\bigl(\nabla_{([v],\dot J_S, \dot J)}\wt D\bigr)(w,
I_S)$ we obtain the same expression as in the statement of the lemma. Obviously,
this would imply the lemma.

The only nontrivial point here is to compute the derivative of the operator of
covariant differentiation $\nabla^\op\deff \nabla: L^{1,p}(S, u^*TX) \to L^p(S,
u^*TX \otimes T^*S)$ in the direction given by some $v \in T_u L^{1,p}(S, X)=
L^{1,p}(S, u^*TX)$. To do this, we fix a smooth vector field $\xi$ on $S$ and
a local section $\mbfw$ of $\scre$. Then $\bigl( \nabla^\op\mbfw \bigr) (\xi)=
\nabla_\xi \mbfw$ is a local section of a Banach bundle with the fiber $L^p(S,
u^*TX)$ over $u\in L^{1,p}(S, X)$.

Differentiation in the direction $v$ yields
\begin{equation}
\nabla _v \bigl( \nabla^\op \mbfw  \bigr) (\xi)
= \nabla^2 _{v,\xi} \mbfw=
R^X(v, du(\xi); \mbfw) + \nabla^2_{\xi,v}\mbfw =
\end{equation}
\begin{equation}
R^X(v, du; \mbfw)(\xi) +  \bigl(\nabla^\op_\xi(\nabla  \mbfw) \bigr)(v).
\end{equation}
Thus we obtain the formula $\nabla_v (\nabla^\op)= R^X(v, du; \cdot)$. Besides,
we have the relation $\nabla_v du= \nabla v$, which was already used for
deriving \eqqref(1.3.3) from \eqqref(1.1.1). Now, the proof of the lemma 
can be completed by explicit calculations. \qed

Using \eqqref(3.2.8) we can describe in more detail the structure of $\pi
_{\!\!\scrj} :\scrm \to \scrj$ at critical points with $\sfh^1_D(S, N_u) \cong \rr$.
Note that the term $[4]$ in \eqqref(3.2.8) is the only one that depends on 
second order derivatives of $J$. Further, the operator $D=D_{u,J}$ is also 
independent of second order derivatives. Thus, deforming $J$ and preserving 
the order one jet $j^1\!J|_{u(S)}$, the map $u:S \to X$ remains $J$-holomorphic
with same the $D$-cohomology groups $\sfh^i(S, \scrn_u)$. The result of such
changes of $J$ is given by

\newlemma{lem3.2.4} Let $[u,J] \in \scrm$ with $\sfh^1_D(S, N_u)\cong
\rr$ and a quadratic form $\ti \Phi: \sfh^0_D(S, N_u)\to \sfh^1_D(S, N_u)$ be
given. Then there exists a $C^1$-small perturbation $\ti J \in \scrj$ of $J$
such that $j^1\! J|_{u(S)} = j^1\! \ti J|_{u(S)}$ and the restriction of
$\Phi_{u, \ti J}$ to $\sfh^0_D(S, N_u)$ equals the given $\ti \Phi$.
Moreover, such a perturbation $\ti J$ of $J$ can be realized in an arbitrarily
small neighborhood of a given point $x \in u(S)$.
\end{lem}

\proof Let $U\subset X$ be a neighborhood of the given $x$. Find
$U' \subset U$ such that $u\inv(U') \not= \emptyset$ and $u$ is an imbedding
on $u\inv(U')$. Obviously, it is sufficient to find an appropriate jet $j^2\!
\ti J|_{u(S)}$ which differs from $j^2\! J|_{u(S)}$ only in $U'\cap u(S)$.
Then $j^2\! \ti J|_{u(S)}$ can be extended to $\ti J$ with the desired
properties.

Covariant differentiation of the identity $J^2 = -\id$ gives the relations
$\nabla_v J \scirc J + J \scirc \nabla_v J=0$ and
\begin{equation}
\nabla^2 _{v_1,v_2}J \scirc J + \nabla_{v_1} J \scirc \nabla_{v_2} J +
\nabla_{v_2} J \scirc \nabla_{v_1}J  + J \scirc \nabla^2 _{v_1,v_2}J =0,
\qquad v_1,v_2 \in T_xX.
\end{equation}
Consequently, we have the following description of the possible choice for 
$j^2\! \ti J|_{u(S)}$ with $j^1\! \ti J|_{u(S)}= j^1\!  J|_{u(S)}$. The tensor
field $u(S) \ni x \mapsto \Theta_x$ defined by
\begin{equation}
v_1, v_2, w \in T_xX \mapsto \Theta_x(v_1, v_2; w)
\deff \nabla^2 _{v_1,v_2}(\ti J -J) (w)\in  T_xX
\end{equation}
must be supported in $U'$, symmetric%
\footnote{\.Obviously, $\nabla^2 _{v_1,v_2}(\ti J -J)- \nabla^2 _{v_2,
v_1} (\ti J -J)$ can be expressed via $\ti J -J$ and the curvature tensor
$R^X(\cdot,\cdot; \cdot)$ of $\nabla$. But $\ti J -J$ vanishes on $u(S)$.}
in $v_1$ and $v_2$, $J$-antilinear in $w$, and zero for $v_1,v_2\in T_x(u(S))
\subset T_xX$. Vice versa, any tensor field $\Theta$ with these properties has
the form $\Theta(v_1, v_2; w) = \nabla^2 _{v_1,v_2}(\ti J -J) (w)$ for an
appropriate $\ti J$ with $j^1\! \ti J|_{u(S)}= j^1\! J|_{u(S)}$.

The condition that $\Theta_x(v_1, v_2; w)$ vanishes for $v_1, v_2\in T_x(u(S)
)$  means that for $x\in U'\cap u(S)$ we can consider $\Theta_x(v_1, v_2;w)$
as a tensor with arguments $v_1, v_2$ varying in the normal bundle $N_u$.
Now, for $\ti J$ as above, $v\in \sfh^0_D(S, N_u)$, and $\psi \in \sfh^0_D(S,
N_u \otimes K_S) \allowbreak \cong \sfh^1(S, N_u)^*$ we obtain the relation
\begin{equation}
\la\psi, \Phi_{u,\ti J}(v,v) \ra =  \la\psi, \Phi_{u,J}(v,v) \ra +
\re \int_S \psi \scirc \Theta(v,v; du).
\end{equation}
Finally, observe that any quadratic form on a {\sl finite dimensional}
space $\sfh^0_D(S, N_u)$ can be realized as $\re \int_S \psi \scirc \Theta(v,v;
du)$ with $\Theta$ satisfying the conditions stated above. \qed

\newsubsection[3.3]{Second variation at cusp-curves}
Our aim in this paragraph is to find conditions ensuring that a critical 
point $[u,J] \in \scrm$ with $\sfh^1(S, N_u)\cong \rr$ is a saddle point. 
\lemma{lem3.2.4} shows that such critical points with $\scrn\sing_u \cong 0$ 
are ``hopeless'' from this point of view. Hence, we need to understand in more
detail the structure of the bilinear operator $\Phi$ on the component $\sfh^0
 (S, \scrn_u \sing) \subset \sfh^0_D(S, \scrn_u)$.

Recall that by the definition of the normal sheaf the stalk $(\scrn_u\sing)_z$
at $z \in S$ is non-trivial exactly when $z$ is a cusp-point of $u: S \to X$
and in this case $\dimc (\scrn_u\sing)_z = \ord_z du$. Thus we want to 
understand the structure of the moduli space at critical points corresponding
to cusp-curves. The following two lemmas contain technical results needed for
this purpose. Recall that the holomorphic line bundle $\scro([A])$ was
introduced in \refdefi{def1.5.1}. We maintain the notation $\nabla$ and
$R^X( \cdot{,} \cdot{;} \cdot)$ from \lemma{lem3.2.3}. In particular, we
have $\nabla_\xi J= \nabla_{du(\xi )}J$, $\nabla^2_{\xi, \eta}v - \nabla^2
_{\eta,\xi} v = R^X(du(\xi), du(\eta); v)$, and other similar relations.
Further, we assume that $\nabla J_S=0$.

\newlemma{lem3.3.1} An element $([w], I_S) \in \sfh^0_D(S, \scrn_u)$
lies in $\sfh^0(S, \scrn\sing_u)$ \iff $w=du(\ti w)$ for some $\ti w \in
L^{1,p}_\loc (S\bs \supp(\scrn_u\sing), TS)$ that extends to
$\ti w \in L^{1,p}(S, TS \otimes \scro([A])\mkern .5mu)$.

\smallskip
In this case, outside the zero-set of $du$ one has
\begin{equation}
\dbar \ti w \equiv  \nabla\ti w +J_S \scirc \nabla\ti w \scirc J_S
= - J_S \scirc I_S
\eqqno(3.3.1)
\end{equation}
and
\begin{align}
0\equiv &\nabla_{\ti w}(D_{u,J}v + \dot J \scirc du \scirc J_S
+ J \scirc du \scirc \dot J_S)  =\term{D}{D_{u,J}(\nabla_{\ti w}v) }+
\eqqno(3.3.2) \\
&\term{1'}{R^X(w,du ;v)} + \term{2'}{J \scirc R^X(w, du\scirc J_S;v)} +
\term5{\nabla_wJ\scirc \nabla v \scirc J_S} +
\term{4'}{\nabla^2_{w,v}J \scirc du \scirc J_S} +
\notag\\
&\term{3'}{\nabla_vJ\scirc \nabla_{\ti w}du \scirc J_S} +
\term{6'}{\dot J\scirc \nabla_{\ti w}du \scirc J_S} +
\term7{\nabla_w\dot J\scirc du \scirc J_S} +
\term{8'}{J\scirc \nabla_{\ti w}du \scirc \dot J_S} +
\notag\\
&\term9{\nabla_wJ\scirc du \scirc \dot J_S} +
\term{13}{J\scirc du \scirc \nabla_{\ti w} \dot J_S} -
\term{14}{\nabla v \scirc \nabla \ti w} -
\term{15}{J \scirc \nabla v \scirc \nabla \ti w  \scirc J_S}.
\notag
\end{align}
\end{lem}

\proof Relation \eqqref(3.3.1) follows from the equality
\begin{align*}
0&= D_{u,J}(du(\ti w))+  J \scirc du \scirc I_S
= du(\dbar\ti w)+  du \scirc J_S \scirc I_S =
\cr
 &=du \Bigl( ( \nabla \ti w + J_S \scirc \nabla \ti w\scirc J_S)
+ J_S \scirc I_S \Bigr).
\end{align*}
Using $J_S^2 = -\id$ we can write the relation in the form
$I_S = J_S \scirc \nabla \ti w - \nabla \ti w \scirc J_S$.

\medskip
To show \eqqref(3.3.2) we start with the computation of 
$D_{u,J}(\nabla _{\ti w} v)$:
\begin{gather}
\eqqno(3.3.2a)
D_{u,J}(\nabla_{\ti w}v) = \nabla(\nabla_{\ti w}v) +
J \scirc \nabla(\nabla_w v) \scirc J_S +
\nabla J({\nabla_w v}, du \scirc J_S)=
\\
\notag
\term{16}{\nabla^2_{\cdot, \ti w}v} \!+
\term{14}{\nabla v \scirc \nabla \ti w \mathstrut} \!+
\term{17}{J \scirc \nabla^2_{\cdot,\ti w}v \scirc J_S} \!+
\term{15}{J \scirc \nabla v \scirc \nabla \ti w \scirc J_S \mathstrut} 
\!+
\term{18}{\nabla J({\nabla_{\ti w} v}, du \scirc J_S)} .
\end{gather}
\pagebreak[3]
Similarly,
\begin{align}
\eqqno(3.3.2b)
&\nabla_{\ti w}(D_{u,J}v + \dot J \scirc du \scirc J_S
+ J \scirc du \scirc \dot J_S)=
\\
&\nabla_{\ti w}( \nabla v + J\scirc \nabla v \scirc J_S +
\nabla_v J \scirc du \scirc J_S +
\dot J \scirc du \scirc J_S + J \scirc du \scirc \dot J_S)=
\notag\\
&\term{16'}{\nabla^2_{\ti w,\cdot}v} +
\term{17'}{J \scirc \nabla^2_{\ti w,\cdot}v \scirc J_S} +
\term{5}{\nabla_wJ \scirc \nabla v \scirc J_S} +
\term{4'}{\nabla^2_{w,v}J \scirc du \scirc J_S} +
\notag\\
&\term{18}{\nabla J(\nabla_{\ti w}v;  du \scirc J_S)} +
\term{3'}{\nabla_v J \scirc \nabla_{\ti w}du \scirc J_S} +
\term{7}{\nabla_w\dot J \scirc du \scirc J_S} +
\notag\\
&\term{6'}{\dot J \scirc \nabla_{\ti w}du \scirc J_S} +
\term{9}{\nabla_w J \scirc du \scirc \dot J_S} +
\term{8'}{J \scirc \nabla_{\ti w}du \scirc \dot J_S} +
\term{13}{J \scirc du \scirc \nabla_{\ti w}\dot J_S}.
\notag
\end{align}
Comparing the terms [16] and $[16']$, it follows that $\nabla^2_{\cdot, \ti w}v -
\nabla^2_{\ti w,\cdot}v = R^X(du, w;v)$. A similar relation holds for the
terms [17] and $[17']$. The equality \eqqref(3.3.2) of the lemma is obtained
by subtracting \eqqref(3.3.2a) from \eqqref(3.3.2b).
\qed

\newlemma{lem3.3.2}
\sli In the situation of \lemma{lem3.3.1}, let $z^*
\in S$ be a cusp-point. Consider $\ti w$ as a section of $TS$ with poles. Set
$k \deff \ord_{z^*}du = \dimc (\scrn_u\sing)_{z^*}$ and choose a local complex
coordinate $z$ on $S$ centered in $z^*$. Fix additionally $([v], \dot J_S, \dot
 J) \in T_{[u, J]} \scrm$ and $\psi \in \sfh^0_D(S, N_u^* \otimes K_S)$. Then 
locally in a neighborhood of $z^*$
\begin{equation}
\matrix
z^k \cdot \ti w(z) &=&
w_0 & + z\cdot w_1 &+ \cdots + z^k\cdot w_k  &+ z^k\cdot w^*(z),
\cr
v(z) &=&
v_0 & + z\cdot v_1 &+ \cdots +z^k \cdot v_k &+ z^k \cdot v^*(z),
\cr
\psi(z) &=&
\psi_0 & + z\cdot \psi_1 &+ \cdots +z^k \cdot \psi_k &
+ z^k \cdot \psi^*(z),
\endmatrix
\eqqno(3.3.2a1)
\end{equation}
where $w^*(z)$, $v^*(z)$, and $\psi^*(z)$ are $L^{1,p}$-smooth local sections
of the corresponding bundles vanishing at $z=0$.

\slii The polynomials in \eqqref(3.3.2a1) can be considered as the order 
$k$ jets of
the following local holomorphic objects: a section of $TS$ for $w_0 + \cdots +
z^k\cdot w_k$, a $(E_u)_{z^*}$-valued function for $v_0+ \cdots +z^k \cdot
v_k$, and resp.\ $(N_u^*)_{z^*}$-valued $(0,1)$-form for $\psi_0 +\cdots +z^k
\cdot \psi_k$. In particular, the coefficients can be considered as
well-defined elements
\begin{equation}
\matrix
w_i &=& \left({\d \over \d z}\right)^i(z^k \ti w(z))\ogran_{z=0}
&\in& (T^*_{z^*} S)^{\otimes i-k} \otimes T_{z^*} S,
\cr
v_i &=& \nabla^i(v(z))\ogran_{z=0}  &\in&
(T^*_{z^*} S)^{\otimes i} \otimes (E_u)_{z^*},
\cr
\psi_i &=& \nabla^i(\psi(z))\ogran_{z=0}  &\in&
(T^*_{z^*} S)^{\otimes i} \otimes (N^*_u \otimes K_S)_{z^*}.
\cr
\endmatrix
\eqqno(3.3.2b1)
\end{equation}
\end{lem}

\proof \sli It follows from \refdefi{def1.5.1} that $du$, considered as a 
{\sl holomorphic} section of the bundle $\hom_\cc(TS, E_u)$, locally has the 
form $du(z)= z^k s(z)$ for some local holomorphic non-vanishing section $s$. 
Consequently, $\scro([A]) = \scro( k\cdot [z^*])$ in a neighborhood of $z^*$. 
Thus by \lemma{lem3.3.1} $\ti w$ can be locally
represented in the form $\ti w(z) = z^{-k} \cdot \hat w(z)$ for some local
$L^{1,p}_\loc$-smooth section of $TS$. Equation \eqqref(3.3.1) is equivalent
to $\dbar \hat w = -\isl\cdot z^k \cdot I_S$ and implies the estimate $|\dbar
\hat w(z)| \le |z^k|\cdot |I_S(z)|$. Now we use \lemma{lem1.2.3}.

The same argument applies to $v(z)$ and $\psi(z)$. Indeed, equation 
\eqqref(2.2.4) on $v$ and relation \eqqref(1.4.1) imply the inequality 
$|\dbar v(z)| \le c\cdot |z^k|$ with some constant $c$. A similar inequality for 
$\psi(z)$ follows from \eqqref(1.5.3).

\smallskip
\slii This part of the lemma can be reformulated in terms of the transformation
of coefficients $w_i$, $v_i$, and $\psi_i$ under the change of a local
holomorphic coordinate $z$ on $S$ and local coordinates on $X$. The claim
concerning the change of $z$ is obvious.

Considering changes of coordinates on $X$ we make the following observation.
If $\mbfx'=(x'_1, \ldots, x'_{2n})$ and $\mbfx''=(x''_1, \ldots,
x''_{2n})$ are two systems of coordinates on $X$ centered in $u(z^*)$, then
$\mbfx''=L(\mbfx') + Q(\mbfx') + \cdots$, where $L: \rr^{2n} \to
\rr^{2n}$ (resp. $Q: \rr^{2n} \to \rr^{2n}$) is an appropriate linear (resp.\
quadratic) map, and so on. In particular, $\mbfx'' - L(\mbfx')= O(
|\mbfx'|^2)$. Consequently, for local frames $\d_{\mbfx'}= (\d_{x'_1},
\ldots,\d_{x'_{2n}})$ and $\d_{\mbfx''}= (\d_{x''_1}, \ldots,\d_{x''_{2n}}
)$ of $TX$ we obtain the relation $\d_{\mbfx''}(\mbfx') - L^{\sf t}(
\d_{\mbfx''}) (\mbfx') = O(|\mbfx'|)$. Thus, for the pulled-back
frames $u^*\d_{\mbfx'}$ and $u^*\d_{\mbfx''}$ of $E_u$ we have
$u^*\d_{\mbfx'}(z) - L^{\sf t}(u^*\d_{\mbfx''})(z)= O(|z|^{k+1})$.
This implies that the change of local coordinates on $X$ induces only a linear
transformation of the $k$-jet of $v$, \ie the $k$-jet of $v$ behaves like a
$k$-jet of a $T_{u(z^*)}$-valued function. The same argumentations
can be applied to the $k$-jet of $\psi$.
\qed

\newlemma{lem3.3.3}
For $[u,J] \in \scrm\!\!$, $([w], I) \in \sfh^0_D(S,
\scrn_n\sing)$, $([v], \dot J_S, \dot J) \in T_{[u,J]}\scrm\!\!$, and
$\psi \in \sfh^0_D(S, N_u \otimes K_S) \cong \sfh^1_D(S, N_u)^*$
it follows that 
\begin{equation}
\Bigl\la \psi ,
\Phi_{u,J}\bigl( ([v], \dot J_S, \dot J),\; ([w],I_S) \bigr) \Bigr\ra=
\re \res_S(\psi\scirc \nabla_{\ti w}v),
\eqqno(3.3.3)
\end{equation}
where $\res_S(\psi\scirc \nabla_{\ti w}v)$ denotes the residual type sum
\begin{equation}
\res_S(\psi\scirc \nabla_{\ti w}v) \deff
\mathop{\textstyle\sum}\limits_{du(z^*_i)=0}\;
\lim\limits_{\epsi \lrar 0}\;
\int_{|z - z^*_i|=\epsi} \psi\scirc \nabla_{\ti w}v
\eqqno(3.3.4)
\end{equation}
over all cusp-points $z^*_i\in S$ of $u$.

Moreover, if $\dot J=0$ and $([v], \dot J_S) \in \sfh^0(S, \scrn_u\sing)$,
then $v=du(\ti v)$ with $\ti v \in L^{1, p}(S, TS \otimes \scro([A])\,)$
and
\begin{equation}
\Bigl\la \psi ,
\Phi_{u,J}\bigl( ([v], \dot J_S, 0),\; ([w],I_S) \bigr) \Bigr\ra=
\re \res_S(\psi\scirc \nabla du(\ti w, \ti v)).
\eqqno(3.3.5)
\end{equation}
\end{lem}

\proof First, we note that by \lemma{lem3.3.2} the formulas
(\ref{eq3.3.3}--\ref{eq3.3.5}) are well-defined.

Now, compute the subtraction of \eqqref(3.3.2) from \eqqref(3.2.8). The terms
$[5]$, $[7]$, and $[9]$ cancel. To simplify further terms we use the Bianchi
identity and antisymmetry of $R^X(\cdot {,} \cdot{;}\cdot)$ in the first two
arguments. The difference of terms $[1] +[2] +[4] -[1'] -[2'] -[4']$ is zero:
\begin{equation}
\term1{R^X( v, du; w)} - \term{1'}{R^X( w, du; v)} +
\term2{J \scirc R^X( v, du \scirc J_S; w)} -
\end{equation}
\begin{equation}
\term{2'}{J \scirc R^X( w, du\scirc J_S; v)} +
\term4{\nabla^2_{v,w}J \scirc du \scirc J_S} -
\term{4'}{\nabla^2_{w,v}J \scirc du \scirc J_S} =
\end{equation}
\begin{equation}
= R^X( v, du; w) + R^X( du, w; v) +
J \scirc R^X( v, du \scirc J_S; w) +
\end{equation}
\begin{equation}
J \scirc R^X( du\scirc J_S ,w ; v) +
R^X(v,w; J \scirc du \scirc J_S) -
J \scirc R^X(v,w; du \scirc J_S)  =
\end{equation}
\begin{equation}
=R^X( v, du; w) + R^X( du, w; v) + R^X(w,v; du) +
\end{equation}
\begin{equation}
J \scirc R^X( v, du \scirc J_S; w) +
J \scirc R^X( du\scirc J_S, w; v) +
J \scirc R^X( w, v; du \scirc J_S)  = 0.
\end{equation}

\smallskip
In the differences $[3] -[3']$, $[6] -[6']$, and $[8] -[8']$ respectively,
we use the relation
\begin{equation}
\nabla(w)=\nabla(du(\ti w))= \nabla_{\ti w}du + du\scirc
\nabla \ti w.
\eqqno(3.3.6)
\end{equation}
This yields
\begin{equation}
\term3{\nabla_vJ \scirc \nabla w \scirc J_S} -
\term{3'}{\nabla_vJ \scirc \nabla_{\ti w}du \scirc J_S} =
\term{3''}{\nabla_vJ \scirc du \scirc \nabla \ti w \scirc J_S},
\end{equation}
and similar equalities for $[6] -[6']$ and $[8] -[8']$. Thus we obtain
\begin{equation}
\nabla_{([v],\dot J_S, \dot J)}(\wh D)([w],I_S) =
\term{3''}{\nabla_vJ\scirc du\scirc \nabla\ti w \scirc J_S} +
\term{6''}{\dot J\scirc du\scirc \nabla\ti w \scirc J_S} +
\end{equation}
\begin{equation}
\term{8''}{J\scirc du\scirc \nabla\ti w \scirc \dot J_S}
+ \term{10}{\nabla_vJ \scirc du \scirc I_S}
+ \term{11}{J \scirc \nabla v \scirc I_S}
+ \term{12}{\dot J \scirc du \scirc I_S}
\end{equation}
\begin{equation}
- \term{13}{J\scirc du \scirc \nabla_{\ti w} \dot J_S} +
\term{14}{\nabla v \scirc \nabla \ti w} +
\term{15}{J \scirc \nabla v \scirc \nabla \ti w  \scirc J_S}
- \term{D}{D_{u,J}(\nabla_{\ti w}v) }.
\end{equation}

\smallskip
Further simplification uses the relation $I_S = J_S \scirc \nabla \ti w -
\nabla \ti w \scirc J_S$. This gives
\begin{equation}
\term{3''}{\nabla_vJ \scirc du\scirc \nabla\ti w \scirc J_S} +
\term{10}{\nabla_vJ \scirc du \scirc I_S} =
\nabla_vJ\scirc du\scirc \nabla\ti w \scirc J_S +
\end{equation}
\begin{equation}
\nabla_vJ \scirc du \scirc
(J_S \scirc \nabla \ti w - \nabla \ti w \scirc J_S)
= \term{3'''}{\nabla_vJ \scirc du\scirc J_S \scirc \nabla \ti w},
\end{equation}
and similarly,
\begin{equation}
\term{6''}{\dot J\scirc du\scirc \nabla\ti w \scirc J_S} +
\term{12}{\dot J \scirc du \scirc I_S} =
\term{6'''}{\dot J\scirc du\scirc J_S \scirc \nabla \ti w},
\end{equation}
\begin{equation}
\term{11}{J \scirc \nabla v \scirc I_S} +
\term{15}{J \scirc \nabla v \scirc \nabla \ti w  \scirc J_S}
= \term{15'}{J \scirc \nabla v \scirc J_S \scirc \nabla \ti w}.
\end{equation}

\smallskip
Now we put together the terms $[3''']$, $[6''']$, $[14]$, and $[15']$.
Because of the relation
\begin{equation}
\nabla v + J \scirc \nabla v \scirc J_S + \nabla_vJ \scirc du\scirc J_S +
\dot J \scirc du\scirc J_S + J\scirc du\scirc \dot J_S =0
\end{equation}
this yields
\begin{equation}
\term{3'''}{\nabla_vJ \scirc du\scirc J_S \scirc \nabla \ti w}+
\term{6'''}{\dot J\scirc du\scirc J_S \scirc \nabla \ti w} +
\term{14}{\nabla v \scirc \nabla \ti w} +
\term{15'}{J \scirc \nabla v \scirc J_S \scirc \nabla \ti w}=
\end{equation}
\begin{equation}
\bigl( \nabla_vJ \scirc du\scirc J_S +
\dot J\scirc du\scirc J_S +
\nabla v +
J \scirc \nabla v \scirc J_S \bigr)\scirc \nabla \ti w =
- \term{8'''}{J\scirc du\scirc \dot J_S \scirc \nabla \ti w}.
\end{equation}

\smallskip
Finally, we conclude that outside the zero-set of $du$ one has
\begin{equation}
\nabla_{([v],\dot J_S, \dot J)}(\wh D)([w],I_S) =
\end{equation}
\begin{equation}
J\scirc du\scirc \bigl( \term{8''}{\nabla\ti w \scirc \dot J_S}
- \term{8'''}{\dot J_S \scirc \nabla \ti w}
- \term{13}{\nabla_{\ti w} \dot J_S} \bigr)
- \term{D}{D_{u,J}(\nabla_{\ti w}v) }.
\end{equation}
Note that $\psi \scirc J \scirc du = \psi \scirc du \scirc J_S =0$, since
$\psi$ vanishes on $du(TS)$. Consequently,
\begin{equation}
\Bigl\la \psi , \Phi_{u,J}\bigl( ([v], \dot J_S, \dot J),
\; ([w],I_S) \bigr)  \Bigr\ra=
\re \int_S \psi \scirc
\Bigl( - \nabla _{([v], \dot J_S, \dot J)} \wh D\Bigr)([w], I_S)=
\end{equation}
\begin{equation}
\re \lim\limits_{\epsi \to 0} \int_{S \bs \cup \Delta(z_i, \epsi)}
\psi \scirc \Bigl( -\nabla _{([v], \dot J_S, \dot J)} \wh D\Bigr)([w], I_S)
=
\re \lim\limits_{\epsi \to 0} \int_{S \bs \cup \Delta(z_i, \epsi)}
\psi \scirc D_{u,J}(\nabla_{\ti w}v).
\end{equation}
Integrating by parts and using $D \psi=0$ we obtain the desired formula 
\eqqref(3.3.3).

To obtain \eqqref(3.3.5) we use \eqqref(3.3.6) and relation $\psi \scirc du =0$.
\qed

\smallskip
Now we can describe the structure of $\Phi_{u,J}$ for cusp-curves. Here we 
restrict ourselves to the case when $(X,J)$ is an almost complex surface. The 
point is that, unlike to the higher dimensional situation, in this dimension 
there are {\sl topological} reasons for the existence of cusp-curves, see
\lemma{lem2.3.4}. 

Other than the restriction on dimension, our setting is as follows. $[u,J] \in 
\scrm$ is a $J$-holomorphic curve with $\sfh^1(S, N_u) \cong \rr$, $z^*\in S$ a 
cuspidal point, $k\deff \ord _{z^*} du$, $z$ a local complex coordinate centered 
at $z^*$, $J^*$ a local (integrable) complex structure in a neighborhood of
$u(z^*)$, and $(w^1, w^2)$ a local system of $J^*$-holomorphic coordinates on
$X$ centered at $u(z^*)$. Finally, we fix some {\sl non-zero} $\psi \in
\sfh^0(S, N_u^* \otimes K_S) \cong \sfh^1(S, N_u)^*$ and denote by $\ord_{z^*}
\psi$ the order of vanishing of $\psi$ at $z^*$.

\newlemma{lem3.3.4}
\sli After a polynomial transformation of the coordinates $(w^1, 
\allowbreak
w^2)$, chosen above, the map $u$ will have the form
\begin{equation}
u(z)=\bigr( z^{k+1} P_1(z), z^{k+l +2} P_2(z) \bigl) + z^{2k+1} g(z),
\eqqno(3.3.7)
\end{equation}
such that $0\le l \le k$, $P_1$ is a polynomial of degree $\le k$ with $P_1(0)
\not =0$, $P_2$ is a polynomial of degree $\le k-l-1$, trivial if $l=k$ or with
$P_2(0) \not =0$ otherwise, and $g(z)$ is an $L^{1,p}$-smooth $\cc^2$-valued
function.

\slii The integers $\ord_{z^*} \psi$ and $l$ do not depend on the particular
choice of coordinates $(w^1, w^2)$ and $\psi\in \sfh^0(S, N_u^* \otimes K_S)$.
For the restriction of $\Phi_{u,J}$ on the stalk $(\scrn_u\sing)_{z^*}
\subset \sfh^0 (S, \scrn\sing_u)$, it follows that
\begin{equation}
\matrix 
\ind_+ \bigl( \Phi_{u,J}\ogran_{(\scrn_u\sing)_{z^*} } \bigr)&=&
\ind_- \bigl( \Phi_{u,J}\ogran_{(\scrn_u\sing)_{z^*} } \bigr)=\cr
\sind \bigl( \Phi_{u,J}\ogran_{(\scrn_u\sing)_{z^*} } \bigr)&=&
\max\,( 0, k-l -\ord_{z^*} \psi).
\endmatrix
\eqqno(3.3.8)
\end{equation}

\sliii If $z^*_1$ and $z^*_2$ are distinct cusp-points of $u:S \to X$,
then the stalks $(\scrn_u\sing)_{z^*_1}$ and $(\scrn_u\sing)_{z^*_2}$ are
$\Phi$-orthogonal, \ie
\begin{equation}
\Phi_{u,J}\bigl( (\scrn_u\sing)_{z^*_1}, (\scrn_u\sing)_{z^*_2} \bigr)=0.
\eqqno(3.3.9)
\end{equation}
\end{lem}

\proof {\sl Part} \sli follows immediately from \lemma{lem1.2.4}. It
simply says that if $\ord_{z^*}du=k$, then the jet $j^{2k+1}u$ is well-defined
and {\sl holomorphic}, \ie can be represented by a complex polynomial.
Note that the theorem of \cite{Mi-Wh} (see \lemma{lem1.2.1}) says that {\sl
topologically} one can also define higher terms which determine the whole
behavior of $u$ at $z^*$.

{\sl Part} \sliii can be easily obtained from \eqqref(3.3.4) and 
\eqqref(3.3.5). It remains to consider

{\sl Part \sliip}. First, we observe that the integer $l$ is the secondary cusp 
index of $u$ at $z^*$ (see \refdefi{def4.2s.1}). It follows then from the
results of \refsubsection{4.2s} that this integer is well defined and 
independent of the choice of $(w^1, w^2)$. The independence of $\ord_{z^*} \psi$ of
the choice of $(w^1, w^2)$ and $\psi$ is obvious. 

Let $J^*$ and $(w^1, w^2)$
be a complex structure and $J^*$-holomorphic coordinates in a neighborhood
of $u(z^*)$, such that $J^*(u(z^*)) = J(u(z^*))$ and $u$ has the local form
\eqqref(3.3.7). Differentiating \eqqref(3.3.7) we see that in the coordinates 
$(w^1, w^2)$
\begin{equation}
du(z)=\bigr( z^k P'_1(z), z^{k+l +1} P'_2(z) \bigl) + z^{2k} g'(z),
\eqqno(3.3.10)
\end{equation}
with polynomials
\begin{equation}P'_1(z)= (k+1)P_1(z) + z {\textstyle{d \over dz}} P_1(z)
\quad\text{and}\quad
P'_2(z)= (k+l+2)P_2(z) + z {\textstyle{d \over dz}} P_2(z),
\eqqno(3.3.11)
\end{equation}
of degree $\le k$ and $\le k-l-1$ respectively
and with $g'(z) = (2k+1) g(z) dz + z dg(z)$ being $L^p$-bounded.

From the definition of the Nijenhuis torsion tensor $N_J$ of $J$ we obtain a
pointwise estimate $|\dbar_J w_\alpha| \le |N_J|$. Further,
\begin{equation}
\dbar(w_\alpha \scirc u) = (dw_\alpha \scirc du)^{(0,1)} =
\dbar_J w_\alpha  \scirc du,
\end{equation}
since $u$ is $J$-holomorphic. Consequently, we obtain a pointwise  estimate
\begin{equation}
|\dbar(w_\alpha \scirc u)(z)| \le c\cdot |z^k|
\eqqno(3.3.12)
\end{equation}
with some constant $c$. Let $\{ e^*_\alpha\}_{\alpha=1,2}$ be the local
$J^*$-complex frame of $T^*X$ dual to the frame $\{ dw_\alpha\}_{\alpha=1,2}$.
Then there exists a local $J$-complex frame $\{ e_\alpha\}_{\alpha=1,2}$ of
$T^*X$ with pointwise estimates
\begin{equation}
|e^*_\alpha(w) -e_\alpha(w)| \le c\cdot |w|
\quad\text{and}\quad
|\nabla e^*_\alpha(w) - \nabla e_\alpha(w)| \le c,
\eqqno(3.3.13)
\end{equation}
where $|w|^2 = |w_1|^2 + |w_2|^2$ and $c$ is some constant.
Using (\ref{eq3.3.10}--\ref{eq3.3.13}) and the estimates 
$|u(z)| \le c\cdot |z^{k+1}|$ and $|du(z)|
\le c\cdot |z^k|$ we conclude that

{\sl a)\.} $\mbfe_\alpha \deff u^* e_\alpha$ is a local complex frame of
$E_u = u^*TX$ with a pointwise estimate
\begin{equation}
|\dbar_{u,J} \mbfe_\alpha(z)| \le c\cdot |z^k|;
\end{equation}

{\sl b)\.} $du$, considered as a section of $E_u\otimes T^*S$ with the frame
$\mbfe_\alpha\otimes dz$, has local form \eqqref(3.3.10), possibly with another
$g'(z) \in L^p$. Moreover, since $du$ is a holomorphic section and
$\mbfe_\alpha$ are sufficiently regular, this new $g'(z)$ is $C^1$-smooth.
Further, since $z dg(z)= g'(z) - (2k+1)zg(z) dz$ is continuous and $dg(z) \in
L^p$ with $p>2$, we conclude that $z dg(z)$ vanishes at $z=0$. This gives
an additional relation $g'(0)=0$.

Differentiating \eqqref(3.3.10) we obtain that in the frame 
$\mbfe_\alpha\otimes dz^2$
\begin{equation}
\nabla du(z)=\bigr( z^{k-1} P''_1(z), z^{k+l } P''_2(z) \bigl) +
z^{2k-1} g''(z),
\eqqno(3.3.14)
\end{equation}
with polynomials
\begin{equation}
P''_1(z)= k\,P'_1(z) + z {\textstyle{d \over dz}} P'_1(z)
\quad\text{and}\quad
P''_2(z)= (k+l+1)P_2(z) + z {\textstyle{d \over dz}} P'_2(z),
\eqqno(3.3.15)
\end{equation}
of degree $\le k$ and $\le k-l-1$ respectively and with
\begin{equation}
g''(z) = (2k+1) g'(z) \otimes dz + z \nabla g'(z)
\end{equation}
continuous and vanishing at $z=0$. By our construction, $P''_1(0)= (k+1)k
P_1(0) \not =0$ and $P''_2(0)= (k+l+2)(k+l+1) P_1(0)$ vanishes \iff $l=k$.

Since the projection $\pr: E_u \to N_u$ is obtained as the quotient \wrt
the image of $z^{-k} du \sim (P'_1(z), z^{l+1} P'_2(z))$, we have the following
form for the composition:
\begin{equation}
\pr \scirc \nabla du(z)= P'''(z) + g'''(z),
\end{equation}
where $P'''(z)$ is a polynomial $P'''(z)$ of degree $\le k-l-1$ given by 
the relation
\begin{equation}
z^{k+l}P'''(z) = z^{k+l} P''_2(z) -
{z^{k+l+1} P'_2(z) \cdot z^{k-1} P''_1(z) \over
z^k P'_1(z)} + o(z^{2k-1}).
\end{equation}
In particular, $P'''(0)= (k+l+1)(l+1)P_2(0)$ vanishes \iff $l=k$.

Denoting $\nu \deff \ord_{z^*} \psi$ we obtain that
\begin{equation}
\psi \scirc \nabla du(z) = a z^{k+l+\nu} + o(z^{k+l+\nu})
\end{equation}
with $a$ vanishing \iff $l=k$. The proof of part \slii of the lemma can be now
finished using the following algebraic result. \qed

\newlemma{lem3.3.5} For a given polynomial $P(z)= a_0 + a_1 \,z
+ \cdots + a_{ k-l-1} z^{ k-l-1}$ with $a_0 \not=0$ and $0\le l <k$
the quadratic form
\begin{equation}
(w_0,\ldots,w_k) \in \cc^{k+1}
\mapsto \re\res_{z=0} \left( { z^{k+l}\, P(z)\,
\bigl(\sum_{i=0}^k w_i\, z^i\bigr)^2  \over z^{2k} } dz \right) \in \rr
\end{equation}
is equivalent to the quadratic form
\begin{equation}
(w_0,\ldots,w_k) \in \cc^{k+1}
\mapsto \re\res_{z=0} \left( { z^{k+l}\, a_0\,
\bigl(\sum_{i=0}^k w_i\, z^i\bigr)^2  \over z^{2k} } dz \right)\in \rr
\end{equation}
and satisfies the index relations
\begin{equation}
\ind_+ Q = \ind_- Q =  \sind Q = k-l.
\eqqno(3.3.16)
\end{equation}
\end{lem}


\newsubsection[4.4]{Critical points and cusp-curves in the moduli space}
Recall that in \lemma{lem3.3.4} we found two obstructions for existence of
saddle points. They are encoded in the secondary cusp-indices $l_i$ of
cusp-points $z^*_i$ of $u$ (see \refdefi{def4.2s.1}) and the vanishing
order at $z^*_i$ of a generic $\psi \in \sfh^0_D(S, N_u \otimes K_S) \cong
\bigl(\sfh^1_D(S, N_u) \bigr) ^*$. The behavior $l_i$ under deformation was
studied in \refsubsection{4.2s}. 
In this paragraph we describe the behavior of $\sfh^0_D(S, N_u \otimes K_S)$.
Our main interest is, of course, $[u,J] \in \scrm$ with $\dimr\sfh^1_D (S, N_u)
=1$, because these are candidates for saddle points. We start with

\newlemma{lem4.4.1} Let $\mbfk=(k_1, \ldots, k_m)$ and $h^1\in \nn$ 
be given. Then the set
\begin{equation}
\whcalm_{=\mbfk, h^1} \deff \bigl\{ (u, J_S, J; \mbfz) \in \whcalm
_{=\mbfk} \;:\; \dimr \sfh^1_D(S, N_u) = h^1\,\bigr\}
\subset \whcalm_{=\mbfk}
\end{equation}
is a $C^{\ell-1}$-smooth submanifold of codimension $h^0{\cdot} h^1$ where
$h^0 \deff \dimr \sfh^0_D(S, N_u)$.

The set $\whcalm_{=\mbfk, h^1}$ is $\bfg$-invariant and the projection $\pr:
\whcalm_{=\mbfk, h^1} \lrar \scrm_{=\mbfk, h^1} \deff \whcalm_{=\mbfk, h^1}
/\bfg$ is a $C^{\ell-1}$-smooth principle $\bfg$-bundle.
\end{lem}

\state Remark. The definition \eqqref(1.5.1) of $\scrn_u$ and the index formula
\eqqref(2.2.3) imply that $h^0= h^1 +2(\mu +(g-1)(3-n) -|\mbfk|)$, where $n=\half
\dimr X$ and $\mu \deff \la c_1(X), [u(S)] \ra$. So $h^0=\dimr \sfh^0_D(S,
N_u)$ is constant along $\whcalm_{=\mbfk, h^1}$ and
\begin{equation}
\whcalm_{=\mbfk} = \mathop{\textstyle\bigsqcup\,}\nolimits
_{h^1 =0} ^\infty \;\whcalm_{=\mbfk, h^1}
\end{equation}
is a stratification of $\whcalm_{=\mbfk}$ indexed by $h^1 =\dimr \sfh^1_D(S,
N_u)$. Taking the $\bfg$-quotients, we obtain a similar stratification of $\scrm
_{=\mbfk}$. Another stratification, more interesting for our purpose, is
\begin{equation}
\scrm_{h^1} = \mathop{\textstyle\bigsqcup\,}\nolimits
_{\mbfk}\scrm_{=\mbfk, h^1}
\end{equation}
with $h^1=1$.
Note that if for given $\mbfk$ and $h^1$ the expected value of $h^0$ is
negative, then $\whcalm_{=\mbfk, h^1}$ is empty.

\proof Consider the Banach bundles $L^{1,p}(S, N)$, $L^p_{(0,1)}(S, N)$
over $\whcalm_{=\mbfk}$, and the bundle homomorphism $D^N: L^{1,p}(S, N) \to
L^p_{(0,1)}(S, N)$ constructed in \lemma{lem4.2.7n}. Then $\sfh^i_D(S, N_u)$
is the (co)kernel of $D^N$. From \lemma{lem1.3.1} we obtain the map
\begin{equation}
\nabla_{(v, \dot J_S, \dot J)}D^N: \sfh^0_D(S, N_u) \to \sfh^1_D(S, N_u)
\eqqno(4.4.1),
\end{equation}
which is bilinear in $(v, \dot J_S, \dot J) \in T_{(u, J_S, J)}\whcalm_{=\mbfk}$ 
and $w \in \sfh^0_D(S, N_u)$. It is not difficult to see that the map 
\eqqref(4.4.1) can be computed using \eqqref(3.2.8) and that it coincides with 
the restriction of $\Phi$ from \eqqref(3.1.1) to the corresponding spaces.

The key point of the proof is to show the surjectivity of the induced map
\begin{equation}
\Phi: T_{(u, J_S, J)}\whcalm_{=\mbfk} \;\lrar\; \homr \bigl(
\sfh^0_D(S, N_u),\;  \sfh^1_D(S, N_u) \bigr)
\end{equation}
Then the claim of the lemma will follow from  the implicit function theorem.

Fix bases $(w_1, \ldots, w_{h^0})$ of $\sfh^0_D(S, N_u)$ and $(\psi_1, \ldots,
\psi_{h^1})$ of $\sfh^0_D(S, N_u \otimes K_S) \cong \bigl(\sfh^1(S,N_u)
\bigr)^*$. The last isomorphism is the Serre duality from \lemma{lem1.5.1}.
We must find tangent vectors $(v_{ij}, \dot J_{S,ij}, \dot J_{ij}) \in T
_{(u, J_S, J)}\whcalm_{=\mbfk}$, $i=1,\ldots, h^0$, $j=1,\ldots, h^1$
obeying the relation
\begin{equation}
\bigl\la \psi_{j'},\,
\Phi\bigl( (v_{ij}, \dot J_{S,ij}, \dot J_{ij}), w_{i'} \bigr) \bigr\ra
= \delta_{ii'} \delta_{jj'}
\eqqno(4.4.2)
\end{equation}
with $\la\cdot\,, \cdot \ra$ denoting the pairing from \eqqref(1.5.5).

The main idea is to find solutions of \eqqref(4.4.2) in the special form such that
$v_{ij}$ and $\dot J_{S,ij}$ are identically zero, and $\dot J_{ij}$ vanish
along $u(S)$ and in a neighborhood of all cusp-points on $u(S)$. This
assumption implies that all the terms in \eqqref(3.2.8) except $[7]$ vanish.
Thus \eqqref(4.4.2) reduces to
\begin{equation}
\re \int_S \psi_{j'} \scirc \nabla_{w_{i'}} \dot J_{ij}\scirc du
\scirc J_S = \delta_{ii'} \delta_{jj'}.
\end{equation}
From this point we can use the arguments either of \lemma{lem2.1.2} or {\sl
Lemma 3.2.4}. Note that we can arrange $\dot J_{ij}$ to have support in
any given open subset $U \subset X$ with $U\cap u(S) \not = \emptyset$. \qed

\smallskip
The big freedom in the choice of $\dot J_{ij}$ implies the following

\newcorol{cor4.4.2} Let $\dimr X =4$, \ie $X$ is an almost complex
surface. Then the intersection of $\whcalm_{=\mbfk, \mbfl}$ and
$\whcalm_{=\mbfk, h^1}$ is transversal, so that the set
\begin{equation}
\whcalm_{=\mbfk, \mbfl, h^1}\deff \whcalm_{=\mbfk, \mbfl} \cap
\whcalm_{=\mbfk, h^1}
\end{equation}
is a $C^{\ell-1}$-smooth submanifold of $\whcalm_{=\mbfk}$ of codimension
$2|\mbfl| + h^0{\cdot}h^1$. A similar result also holds for
$\scrm_{=\mbfk, \mbfl, h^1} \deff \whcalm_{=\mbfk, \mbfl, h^1} /\bfg =
\scrm_{=\mbfk, \mbfl} \cap \scrm_{=\mbfk, h^1}$.
\end{corol}

\medskip
Now we will study the behavior of zeros of a non-trivial $\psi \in \sfh^0_D(S,
N_u \otimes K_S)$ for $[u, J] \in \scrm_{=\mbfk, h^1=1}$. Note that, modifying
the construction from \lemma{lem4.4.3}, we obtain a bundle $N^* \otimes K_S$
over $\whcalm_{=\mbfk} \times S$, $C^{\ell-1}$-smooth Banach
bundles $L^{1,p}(S, N^* \otimes K_S)$ and $L^p_{(0,1)}(S, N^* \otimes K_S)$
over $\whcalm_{=\mbfk}$, and a $C^{\ell-1}$-smooth bundle homomorphism
\begin{equation}
(D^N)^*: L^{1,p}(S, N^* \otimes K_S) \to L^p_{(0,1)}(S, N^* \otimes K_S).
\end{equation}
Since the kernel of $(D^N)^*$ is of constant dimension on each $\whcalm_{=\mbfk,
h^1}$, we obtain a $C^{\ell-1}$-smooth bundle $\sfh^0_D(S, N^* \otimes K_S)$ of
$\rank_\rr =h^1$ on $\whcalm_{=\mbfk, h^1}$. This means that there exists a
(local) frame $\psi_1, \ldots, \psi_{h^1}$ of $\sfh^0_D(S, N^* \otimes K_S)$
which depends $C^{\ell-1}$-smoothly on $(u, J_S, J) \in \whcalm_{=\mbfk, h^1}$.

In the particular case $h^1=1$ we obtain a (local) $C^{\ell-1}$-smooth family of
non-trivial $\psi \in \sfh^0_D(S,N^* \otimes K_S)$ such that for every $(u,
J_S, J)$ the corresponding $\psi$ is defined uniquely up to a constant factor.
\lemma{lem1.2.3} ensures that the zero-divisor of such $\psi$ is well-defined
and has degree
$c_1(N^* \otimes K_S)$. By \lemma{lem2.3.4}, the possible range for $c_1(N^*
\otimes K_S)$ is the interval between 0 and $g-1$. We are interested in the
distribution of the zeros of $\psi$, especially at cusp-points of $u$. For
a given $\mbfk=(k_1, \ldots, k_m)$ we consider  $m$-tuples $\bfnu=(\nu_1,
\ldots, \nu_m)$ with $0 \le \nu_i \le k_i$, $i=1, \ldots, m$. Denote
$|\bfnu| \deff \sum_{i=1}^m \nu_i$.

\newlemma{lem4.4.3}
\sli The set $\whcalm_{=\mbfk, \bfnu} \deff \bigl\{ (u, J_S, J) \in
\whcalm_{=\mbfk, h^1=1} \;:\; \ord_{z^*_i} \psi \ge \nu_i \,\bigr\} \subset
\whcalm_{=\mbfk, h^1=1}$ is a $C^{\ell-1}$-smooth submanifold of codimension
$2(n-1)|\bfnu|$, $n=\half\,\dimr X$.

\slii  Let $n=2$, \ie $X$ is an almost complex surface. Then the intersection of
 $\whcalm_{=\mbfk, \mbfl}$ and $\whcalm_{=\mbfk, \bfnu}$ is transversal, so
that the set
\begin{equation}
\whcalm_{=\mbfk, \mbfl, \bfnu} \deff \whcalm_{=\mbfk, \mbfl, h^1=1} \cap
\whcalm_{=\mbfk, \bfnu} \;\subset\; \whcalm_{=\mbfk, h^1=1}
\end{equation}
is a $C^{\ell-1}$-smooth submanifold of codimension $2|\mbfl| + |\bfnu|$.

\sliii Similar results hold for $\scrm_{=\mbfk, \mbfl, \bfnu} \deff
\whcalm_{=\mbfk, \mbfl, \bfnu}/\bfg= \scrm_{=\mbfk, \mbfl,h^1=1} \cap
\scrm_{=\mbfk, \bfnu}$.
\end{lem}

\proof \slip.
Fix $(u_0, J_{S,0}, J_0) \in \whcalm_{=\mbfk, \bfnu}$. Let $z_i$
be a local $J_S$-holomorphic coordinate on $\whcalm_{=\mbfk, h^1=1}$ in the
sense of \refdefi{def4.2.2n}, centered at the cusp-point $z^*_i$ of $(u,
J_S, J)$, $i=1, \ldots, m$. Further, let $\psi$ be a local $C^{\ell-1}$-smooth
family of non-trivial elements of $\sfh^0_D(S,N^* \otimes K_S)$. Then by
\lemma{lem3.3.2}, for each $i=1, \ldots, m$ and each $(u, J_S, J) \in
\whcalm_{=\mbfk, h^1=1}$ we can construct the jets $j^{k_i\!}\psi =
\sum_{j=0} ^{k_i} \psi_{i,j} \cdot z_i^{\;j}$ of $\psi$ at $z^*_i$.

Repeating the arguments used in the proof of \lemma{lem4.2.3n} 
we can show
that the coefficients $\psi_{i,j} \in (T^*_{z^*_i})^{\otimes j} \otimes (N^*
\otimes K_S)_{z^*_i}$ depend $C^{\ell-1}$-smoothly on $(u, J_S, J) \in \whcalm
_{=\mbfk, h^1=1}$. This means that $\whcalm_{=\mbfk, \bfnu}$ is the zero set
of the (locally defined) function $\yps^\psi_\bfnu$ on $\whcalm_{=\mbfk,
h^1=1}$ given by the first $\nu_i$ coefficients of each $j^{k_i\!}\psi$,
$i=1, \ldots,  m$, \ie
\begin{equation}
\yps^\psi_\bfnu(u, J_S, J) =
(\psi_{1,0}, \ldots, \psi_{1, \nu_1 -1}, \ldots, \psi_{m,0}, \ldots,
\psi_{m, \nu_m -1}).
\end{equation}

Consequently, it is sufficient to show the surjectivity of the differential 
$d\yps^\psi_\bfnu$ at the fixed $(u_0, J_{S,0}, J_0)$. But first we must compute 
$d\yps^\psi_\bfnu$ for a given $(v, \dot J_S, \dot J) \in T_{(u_0, J_{S,0}, J_0)} 
\whcalm_{=\mbfk, h^1=1}$. Let $\gamma(t)= (u_t, J_{S,t}, J_t)$ be a curve in 
$\whcalm_{=\mbfk, h^1=1}$ which starts at $(u_0, J_{S,0}, J_0)$ and has the 
tangent vector $(v, \dot J_S, \dot J)$ at $t=0$. Then we obtain a family $\psi_t$ 
of non-trivial $\psi_t \in \sfh^0_D(S, N_{u_t}^* \otimes K_S)$. In particular, 
for each $t$ we obtain the relation $D^*_t \psi_t=0$, where  $D^*_t$ denotes 
the operator $D^{N^* \otimes K_S}$ corresponding to $(u_t, J_{S,t}, J_t)$.

Fix some symmetric connections on $X$ and $S$. As in \refsubsection{3.2}, we 
obtain induced connections for all (usual and Banach) bundles involved in our
computations. We use the same notation $\nabla$ for all these connections,
in particular, for the connection in the bundle  $L^{1,p}(S, N_u)$ with
the fiber $L^{1,p}(S, N_{u_t})$ over $(u_t, J_{S,t}, J_t)$. Hence for any
$w_0 \in L^{1,p}(S, N_{u_0})$ we can construct a family  $w_t \in L^{1,p}(S,
 N_{u_t})$ which is covariantly constant. This yields a covariantly constant
trivialization of the Banach bundle $L^{1,p}(S, N_u)$ along $\gamma$.

For every such family $w_t$ we have the relation
\begin{equation}
\la w_t, D^*_t \psi_t \ra=0.
\eqqno(4.4.3)
\end{equation}
Vice versa, a family $\psi_t \in L^{1,p}(S, N^*_{u_t} \otimes K_S)$ lies
in $\sfh^1_D(S, N^*_{u_t} \otimes K_S)$ if \eqqref(4.4.3) holds. Rewrite 
\eqqref(4.4.3) in the form
\begin{equation}
\la \psi_t, D_t w_t \ra=0
\eqqno(4.4.4)
\end{equation}
with $D_t$ denoting the operator $D^N_{u_t, J_t}: L^{1,p}(S, N_{u_t}) \to
L^p_{(0,1)}(S, N_{u_t})$. After covariant differentiation in $t$ we obtain
$\la \dot \psi_t, D_t w_t \ra + \la \psi_t, (\nabla_{(v_t, \dot J_{S,t},
\dot J_t)} D_t) w_t \ra =0$. The latter is equivalent to
\begin{equation}
\la  D^*_t \dot\psi_t,  w_t \ra +
\la \psi_t, (\nabla_{(v_t, \dot J_{S,t}, \dot J_t)} D_t) w_t \ra =0.
\eqqno(4.4.5)
\end{equation}

Now we can give the description of $d\yps^\psi_\bfnu$ at $(u_0, J_{S,0}, J_0) \in 
\whcalm_{=\mbfk, h^1=1}$. For a given tangent vector $(v, \dot J_S, \dot J)$ we 
find $\dot \psi \in L^{1,p}(S, N_{u_0} \otimes K_S)$ such that \eqqref(4.4.5)
holds for every $w \in L^{1,p}(S, N_{u_0} \otimes K_S)$. The existence of such
$\dot \psi$ is equivalent to the condition that $(v, \dot J_S, \dot J)$ is
tangent to $\whcalm_{=\mbfk, h^1=1}$. Such $\dot \psi$ is unique up to
addition of $\psi \in \sfh^0_D(S, N_{u_0} \otimes K_S)$. The
jets $j^{k_i\!} \dot \psi= \sum_{j=0} ^{k_i} \dot \psi_{i,j} \cdot z_i^{\;j}$
of such $\dot \psi$ at $z^*_i$ are well-defined and
\begin{equation}
d\yps^\psi_\bfnu(v, \dot J_S, \dot J) =
(\dot \psi_{1,0}, \ldots, \dot \psi_{1, \nu_1 -1}, \ldots,
\dot \psi_{m,0}, \ldots,
\dot \psi_{m, \nu_m -1}).
\end{equation}
$d\yps^\psi_\bfnu$ is independent of the choice of $\dot \psi$ provided
$(u_0, J_{S,0}, J_0) \in \whcalm_{=\mbfk, \bfnu}$.

\smallskip
To show the surjectivity of $d\yps^\psi_\bfnu$ we must invert the construction
above. Let $j^{\nu_i-1\!} \dot \psi$ be given jets. Extend them to jets
$j^{k_i\!} \dot \psi$. Note that by definition the operator $D^*_0 = D^{N^*
\otimes K_S} _{u_0, J_0}$ has the form $\dbar^{N^* \otimes K_S}_{u_0, J_0} +
R^{N^* \otimes K_S}_{u_0, J_0}$ where
\begin{equation}
R^{N^* \otimes K_S}_{u_0, J_0}: N^* \otimes K_S \to
N^* \otimes K_S \otimes \Lambda^{(0,1)}
\end{equation}
is a continuous bundle homomorphism. Consider the equations
\begin{equation}
z_i^{-k_i}\bigl(\dbar^{N^* \otimes K_S}_{u_0, J_0} +
R^{N^* \otimes K_S}_{u_0, J_0} \bigr)
\bigl(j^{k_i\!} \dot \psi + z_i^{k_i}  \phi_i(z_i)\bigr)=0
\eqqno(4.4.6)
\end{equation}
for unknown $\phi_i(z_i)$ defined in a neighborhood of $z^*_i$. Using
\lemma{lem1.4.1} we obtain pointwise estimates $|R^{N^* \otimes K_S}_{u_0,
J_0}(z_i)| \le C\cdot |z_i|^{k_i}$. Thus equation \eqqref(4.4.6) is equivalent to
\begin{equation}
\left(\dbar^{N^* \otimes K_S}_{u_0, J_0} +
\left(\msmall{\bar z_i \over z_i}\right)^{k_i}
R^{N^* \otimes K_S}_{u_0, J_0} \right) \phi_i(z_i)
+ z_i^{- k_i} R^{N^* \otimes K_S}_{u_0, J_0}j^{k_i\!}=0.
\eqqno(4.4.7)
\end{equation}
The existence of solutions of \eqqref(4.4.7) can be deduced from the
surjectivity of the operator $\dbar + R: L^{1,p}(\Delta,\cc^n) \to L^p(\Delta,
\cc^n)$ with $R \in L^p$, $p>2$.
We refer to \cite{Iv-Sh-1} for the construction of a right inverse for such
$\dbar + R$. This implies the local existence of solutions $\phi_i(z_i)$
of \eqqref(4.4.6).

The regularity property of $R^{N^* \otimes K_S}_{u_0, J_0}$ implies that the $z_i
^{k_i} \phi_i(z_i)$ are $C^{\ell-1}$-smooth. Thus we can construct a $\dot \psi 
\in C^{\ell-1}(S, N^*_{u_0} \otimes K_S)$ which locally near $z^*_i$ has the form 
$\dot \psi(z_i) = j^{k_i\!} \dot \psi + z_i^{k_i} \phi_i(z_i)$ and satisfies 
\eqqref(4.4.6). Now, the surjectivity of $\yps^\psi_\bfnu$
will follow from the existence of $(v, \dot J_S, \dot J)\in T_{(u_0, J_{S,0},
J_0)}  \whcalm_{=\mbfk, h^1=1}$ such that for the constructed $\dot \psi$ and
a fixed non-zero $\psi_0 \in \sfh^0_D(S, N^*_{u_0} \otimes K_S)$
the relation \eqqref(4.4.5) holds for any $w \in L^{1,p}(S, N_{u_0})$.

Now observe that we can use \eqqref(3.2.8) to compute $\nabla_{ (v, \dot J_v, 
\dot J)} D^N_{u_0, J_0}$. This implies that we can use the trick from the proof 
of \lemma{lem4.4.1}.  Namely, we look for the desired $(v, \dot J_S, \dot J)$
in the special form, such that $v$ and $\dot J_S$ vanish identically, and
$\dot J$ vanishes  along $u_0(S)$ and in some neighborhoods of cusp-points of
$u(S)$. Now all terms in \eqqref(3.2.8) except $[7]$ vanish, and 
\eqqref(4.4.5) is equivalent to
\begin{equation}
D^{N^* \otimes K_S}_{u_0, J_0} \dot \psi +
\psi_0 \scirc \nabla \dot J \scirc du_0 \scirc J_S =0.
\end{equation}
To finish the construction of $\dot J$ we use the fact that $D^{N^* \otimes
K_S}_{u_0, J_0} \dot \psi$ vanishes in a neighborhood of each cusp-point
$z^*_i$. This yields the surjectivity of $\yps^\psi_\bfnu$ and
the first assertion of the lemma.

\smallskip
The second and third assertions follow from previous considerations. \qed

\newsubsection[4.5]{(Non)existence of saddle points in the moduli space}
The results obtained above in this section allow us to prove the main technical
result of the paper. Let $X$ be a manifold of dimension $2n$,
$\scrj$ an open connected set in the space of $C^\ell$-smooth almost complex
structures on $X$ with $\ell >2$ non-integer, $S$ a closed surface of genus
$g\ge1$, and $[C] \in \sfh_2(X,\zz)$ a homology class.

\newdefi{def4.5.1} A pseudoholomorphic  map $u:S \to X$ has an {\sl
ordinary cusp} at $z^*\in S$ if for appropriate coordinates $z$ on $S$ and
$(w_1, w_2)$ on $X$
\begin{equation}
u(z)=\bigl(z^2 +O(|z|^3), z^3 + O(|z|^{3+\alpha}) \bigr).
\end{equation}
This property is equivalent to the condition that $u$ has a cusp of order 1
and the secondary cusp-index 0 at $z^*$.
\end{defi}

\newthm{thm4.5.1} Let $h(t)=J_t$, $t\in I=[0,1]$, be a {\sl generic}
path in $\scrj$ and $\scrm_h$ the corresponding relative moduli space of
parameterized pseudoholomorphic curves of genus $g\ge1$ in the homology class
$[C]$.

\sli If $n\ge 3$, then every critical point of the projection $\pi_h: \scrm_h
\to I$ is represented by an {\sl imbedded} curve $C=u(S)$, $u: S \to X$;

\slii If $n=2$, then every critical point of the projection $\pi_h: \scrm_h
\to I$ is represented by a curve $C=u(S)$ such that:
\begin{itemize}
\item the only singularities on $C$ are nodes or ordinary cusps;
\item the possible number of cuspidal points $\vkappa$ on $C$ is 
\begin{equation}\eqqno(4.5.2)
\mu \le \vkappa \le \mu +g-1,
\end{equation}
where $\mu \deff \la c_1(X), [C] \ra$ and $g$ is the (geometric) genus of $C$,
$g=g(S)$;
\item the saddle index of $d^2\pi_h$ at $C$ is at least $\vkappa$, \ie
$$
\sind_C d^2\pi_h \ge \vkappa \ge \mu .
$$
\end{itemize}

\sliii In the case when the inequality \eqqref(4.5.2) is a contradiction, 
the claim \slii has the following meaning:
\begin{itemize} 
\item If $g=0$, then $\pi_h$ has no critical points;
\item If $\mu +g -1<0$, then the space $\scrm_h$ is empty for generic $h$.
\end{itemize}
\end{thm}

Before giving the proof we must specify the meaning of the notion {\sl generic
path}. One of the most reasonable conditions is that any two {\sl regular}
almost complex structures $J_0,\; J_1 \in \scrj$ (see \S\.{\sl2.3}) can be
connected by a path $\{J_t\}_{t\in I=[0,1]}$ with the property stated in
the theorem. To ensure this we need the following easy

\newprop{prop4.5.2} Let $F: \scrx \to \scry$ be a $C^1$-smooth
Fredholm map between separable Banach manifolds. Assume that $\scry$ is
connected and that the index of $F$ is at most $-2$. Then the set $\scry \bs
F(\scrx)$ is path-connected.
\end{prop}

\state Remark. The proposition generalizes the obvious fact that submanifolds 
of codimension at least 2 do not divide the ambient manifold. Note that one can 
have at most countably many connected components of $\scrx$ and that on these 
components the index of $F$ can vary from component to component.

\statep Proof. \refthm{thm4.5.1}. We already know from \refsection{2} that for 
a generic path $h(t) =J_t$ in $\scrj$ the set $\scrm_h$ is a manifold. In previous
paragraphs of this section we have showed that critical points $[u,J]$ of the
projection $\pi_h: \scrm_h \to I$ have an intrinsic description independent of the
particular choice of the path $J_t$. Moreover, the quadratic form $d^2\pi_h$ at 
these points also admits a similar intrinsic description. Furthermore, we have 
found a stratification of the set of ``suspicious''
points $[u,J] \in \scrm$ by submanifolds and estimated their codimension. It
remains to find the strata with the Fredholm index $\le -2$ over $\scrj$ and
apply \propo{prop4.5.2}.

The ``suspicious'' points $[u,J]$ on $\scrm$ are those with $\sfh^1_D(S,
\scrn_u) \cong \rr$. They can be separated into classes according to the
structure of the normal sheaf $\scrn_u$. Since the singular part $\scrn_u\sing$
of $\scrn_u$ reflects the cusp-curves we are led to the spaces $\scrm_{=\mbfk}$
of curves with prescribed order of cusps.

Denote by $\ind$ the index of the projection $\pr _{\!\!\scrj}: \scrm \to \scrj$,
so that $\ind = 2(\la c_1(X), [C] \ra + (g-1)(3-n))$. If $\ind < 0$, then
for a generic path $h(t)=J_t$ the set $\scrm_h$ is empty and the claim of the
theorem holds. Thus we may assume that $\ind \ge0$. By \refthm{thm4.2.1}, we
must ``pay'' at least $2(n-1)|\mbfk|$ dimensions to go to $\scrm_{=\mbfk}$. By
\lemma{lem4.4.1}, we must ``pay'' further $\ind-2|\mbfk| +1$ dimensions to 
obtain the condition $\sfh^1_D(S, \scrn_u) \cong \rr$. Note that $2(n-1)|\mbfk|
\ge 2|\mbfk| +2$ if $n\ge 3$ and $\mbfk$ is non-trivial. Thus in the case $n\ge3$ 
we ``overdraw'' our ``credit'' $\ind$ at least by $3$. This means that for 
non-trivial $\mbfk$ the index of the projection from $\scrm_{=\mbfk,h^1=1}$ to 
$\scrj$ is at most $-3$ and we can apply \propo{prop4.5.2}. Thus for $n\ge3$ any 
critical point of $\scrm_h$ is represented by an immersion $u:S \to X$. 

\medskip
In the case $n=2$ we can ``strike the balance'' in a similar way. Indeed, we
come to the ``overdraw'' of at least $3$ dimensions in each of the following
cases:

{\sl a)} $\sfh^1_D(S, \scrn_u) \cong \rr$ and there exists at least one
cusp-point of cusp-order$\ge2$;

{\sl b)} $\sfh^1_D(S, \scrn_u) \cong \rr$ and there exists at least one
cusp-point the secondary cusp-index $\ge1$;

{\sl c)} $\sfh^1_D(S, \scrn_u) \cong \rr$ and a non-trivial $\psi \in
\sfh^1_D(S, N_u^* \otimes K_S )$ vanishes in at least one cusp-point.

\noindent
Thus for generic $h(t)=J_t$ we can exclude all these possibilities. The
remaining case admits only cusps of order 1 with the secondary cusp-index 0.
This means that $u$ has only ordinary cusps. Since possibility
{\sl c)} is excluded, each such cusp gives input 1 into the saddle index by
\lemma{lem3.3.4}. Finally, we estimate the number of such cusps using
\lemma{lem2.3.4}. 

\medskip
Now we show that for $n\ge3$ and generic $h$ any critical point of $\scrm_h$ is 
represented by an {\sl imbedding} $u:S \to X$. Since this result will be not used
in the sequel, we give only a sketchy proof. 

Denote by $\scrm_{\sf{imm}}$ 
the total moduli space of {\sl immersed pseudoholomorphic curves} with the same 
topological data $g=g(S)$ and $[C] \in \sfh_2(X, \zz)$ as usual. In other words,
$\scrm_{\sf{imm}}= \scrm_{=\mbfk}$ with trivial $\mbfk$. It follows easily from 
\refsection{4} that this is an open set in the whole space $\scrm$. The space 
$\scrm_{\sf{imm}}$ admits a natural stratification in which every stratum contains
curves with the same number and type of multiple point on the image $C=u(S)$. 
Obviously, the biggest stratum is the subspace of {\sl imbedded} curves, and this 
is an open subset in $\scrm_{\sf{imm}}$. The next biggest stratum consists of 
curves with exactly one transversal double point on $C=u(S)$. Let us denote it by 
$\scrm_{\sf{imm}}^\times$ with the character $\times$ symbolizing a transversal 
self-intersection of exactly 2 branches of $C=u(S)$. 

Locally, $\scrm_{\sf{imm}}^\times$ is defined by the condition $u(z_1) =u(z_2)$ 
for some $z_1 \not= z_2 \in S$. Linearization of this condition is the equation 
$$
\pr_{N^\times}(v(z_1) -v(z_2)) =0
$$
on $[v, \dot J_S, \dot J] \in T_{[u,J]}\scrm_{\sf{imm}}$, where $N^\times$ denotes
the plane in $T_{x^\times}X$ normal to both branches of $C=u(S)$ at the 
point $x^\times=u(z_1)=u(z_2)$, \ie $N^\times \deff T_{x^\times}
X /\bigl(du(T_{z_1}S) \oplus du(T_{z_2} S) \bigr)$. It is easy to see that this 
condition is transversal. Thus $\scrm_{\sf{imm}}^\times$ is a $C^\ell$-smooth 
submanifold of real codimension $2(n-2)$. Moreover, it follows from the proof of 
\lemma{lem4.4.1} that biggest stratum is transversal to the subspace $\scrm
_{\sf{imm}, h^1=1}$ of immersed curves with $\sf{h}^1(S, \scrn_u)=1$. Consequently,
the space $\scrm^\times_{\sf{imm}, h^1=1} \deff \scrm_{\sf{imm}}^\times \cap \scrm
_{\sf{imm}, h^1=1}$ of immersed curves with $\sf{h}^1(S, \scrn_u)=1$ and with 
exactly one transversal self-intersection point is a $C^\ell$-smooth submanifold 
of real codimension $2(n-2)$ in $\scrm_{\sf{imm}, h^1=1}$, and of real codimension
$\ind + 1+ 2(n-2)$ in $\scrm$. Since $2(n-2) \ge2$ for $n\ge3$, we can apply 
\propo{prop4.5.2}. 

The complementary strata $\scrm _{\sf{imm}}^{\mib{a}}$ of $\scrm_{\sf{imm}}$ 
consist of curves having either several double points, or one double point with 
tangency of higher degree, or even more complicated multiple points, with the index
$\mib{a}$ encoding the number and the type of multiple points. A similar argument 
shows that these strata $\scrm_{\sf{imm}}^{\mib{a}}$ are transversal to $\scrm
_{\sf{imm},h^1 =1}$ and that the intersections $\scrm _{\sf{imm}, h^1=1}^{\mib{a}} 
\deff \scrm _{\sf{imm}}^{\mib{a}} \cap \scrm _{\sf{imm}, h^1=1}$ are transversal. 
The computation of the number of conditions shows that these strata have even 
higher codimension in $\scrm_{\sf{imm}}$. So \propo{prop4.5.2} still applies. 
This finishes the proof of the theorem. \qed

\medskip
In applications, one needs a version of \refthm{thm4.5.1} for the case of curves
passing through given fixed points $\mbfx=(x_1,\ldots,x_m)$ on $X$. Recall that
for a $C^\ell$-smooth map $h: I\deff[0,1] \to \scrj$ we denote by $\scrm_{h,\mbfx}$
the relative moduli space of $J_t=h(t)$-holomorphic curves passing through $\mbfx=
 (x_1,\ldots,x_m)$ (see \refsubsection{2.4}). Let $\pi_{h,\mbfx}: \scrm_{h,\mbfx}
\to I$ be the corresponding projection. We also assume that $\dimr X=4$.

\newthm{thm4.5.3} For a generic $h$ every critical point of the projection 
$\pi_{h, \mbfx}: \scrm_{h,\mbfx} \to I$ is represented by a curve $C$ such that:
\begin{itemize}
\item the only singularities on $C$ are nodes or ordinary cusps;
\item the marked points $(x_1,\ldots, x_m)$ are smooth points of $C=u(S)$;
\item the possible number of cuspidal points $\vkappa$ on $C$ is 
\begin{equation}\eqqno(4.5.3)
\mu -m \le \vkappa \le \mu -m +g-1
\end{equation}
where $g$ is the (geometric) genus of $C$;
\item the saddle index of $d^2\pi_h$ at $C$ is at least $\vkappa$, \ie
$$
\sind_C d^2\pi_h \ge \vkappa \ge \mu -m.
$$
\end{itemize}

In the case when the inequality \eqqref(4.5.3) is a contradiction
the claim has the following meaning:
\begin{itemize} 
\item If $g=0$, then $\pi_{h, \mbfx}$ has no critical points;
\item If $\mu -m +g -1 <0$, then the space $\scrm_h$ is empty for generic $h$.
\end{itemize}

\end{thm}

\proof The main observation in the proof is that after an appropriate modification
all the results of this section remain valid also for curves passing through 
fixed points. In particular, the most important formulas \eqqref(3.3.3) and 
\eqqref(3.3.5) from \lemma{lem3.3.3} holds after replacing $\scrn\sing_u$ by 
$\scrn\sing_{u,\mbfx}$. To show this we note first that {\sl Lemmas 
\ref{lem3.2.3}, \ref{lem3.3.1}}, and {\sl \ref{lem3.3.2}} can be applied without
any modification. After this, the proof of \lemma{lem3.3.3} applies with
the only difference that the usual Gromov operator $D_{u,J}$ acting in $E$
should be replaced by the operator $D_{u,-\mbfz,J}$ acting in $E_{-\mbfz}$.
The validation of such a replacement is justified in \refsubsection{2.4}. Indeed,
by the very definition, $D_{u,-\mbfz,J}$ is the restriction of $D_{u,J}$ to
the subspace of sections of the subbundle $E_{u,-\mbfz} \subset E_u$.
In a similar way one modifies the argumentation of \refsubsection{4.4}.

Finally, we note that the condition of coincidence of some cusp point of $C=u(S)$
with some of marked points $x_1,\ldots,x_m$ defines a subset in $\scrm_\mbfx$
which has a natural stratification into submanifolds of codimension $\ge2$.
Every such stratum is defined by the cusp order $\mbfk$ of $C=u(S)$ and indication
of the those cuspidal points which pass through the marked points $x_1,\ldots,
x_m$. This means that every such stratum is a submanifold of the space $\scrm
_{=\mbfk}$. Moreover, the codimension of every such stratum in $\scrm_{=\mbfk}$ is 
$4a$, where $a$ is the number of cusps lying in the marked points. As in the case 
$m=0$ above, one can show that the intersection of such a stratum with the space
$\scrm_{=\mbfk, h^1=1}$ is transversal and has the expected codimension. 
Hence we may conclude that for generic $h$ such a coincidence can not occur in 
the critical points of $\pi_{h,\mbfx}$. The same argument is applied to show that 
for generic $h$ there are no coincidence of the marked points $x_1,\ldots,x_m$ 
with nodal points of the curve $C=u(S)$ representing a critical point of $\pi_h$. 
\qed

\newsection[5]{Deformation of nodal curves}

\newsubsection[5.2]{Nodal curves and Gromov compactness theorem}
The total moduli space $\scrm$ constructed in {\sl Section 2} is not complete.
More precisely, the projection $\pi_{\!\!\scrj}: \scrm \to \scrj$ is, in 
general, not proper. This means that there exists a sequence $[u_i, J_i] \in 
\scrm$ such that $J_i$ converges to $J_\infty \in \scrj$ but no subsequence of
$\{ u_i \}$ converges in $L^{1,p}(S, X)$-topology, even after 
reparameterization. Gromov compactness theorem ensures that there still exists 
subsequence of $\{ u_i \}$ which converges \wrt the Gromov topology, which is 
weaker that the Sobolev $L^{1,p}$-topology.

\smallskip
In the literature one can find several non-equivalent definitions for Gromov
topo\-logy. In this paper  we shall use that one which is equivalent to the
original definition of Gromov (\cite{Gro}). However, our version is more
detailed in the sense that it is based on the notion of {\sl stable maps}.
This notion for curves in a complex algebraic manifold $X$ was introduced by
Kontsevich in \cite{K}, see also \cite{K-M}. Our definition of stable maps over 
$(X,J)$ is simply a translation of this notion to almost complex manifolds.

\newdefi{def5.2.1} The {\sl standard node} is the complex analytic
set
\begin{equation}
\scra_0 \deff \{ (z_1,z_2)\in \Delta^2 : z_1\cdot z_2 =0\}.
\eqqno(5.2.1)
\end{equation}
A point on a complex curve is called a {\sl nodal point}, if has a
neighborhood biholomorphic to the standard node. A {\sl nodal curve} $C$ is a
complex analytic space of pure dimension 1 with only nodal points as
singularities.
\end{defi}

\newdefi{def5.2.1a} An annulus $A$ with a complex structure $J$ has
{\sl conformal radius} $R>1$ if $A$ is biholomorphic to $A(1,R) \deff \{ z\in
\cc \,:\, 1 <|z| < R \}$. Define a {\sl cylinder} $Z(a,b) \deff S^1
\times [a,b] = \{ (\theta, t) : 0 \le \theta \le 2\pi,\; a\le t \le b \}$, $a<
b$, with the complex structure $J_Z ({\d\over \d\theta}) \deff {\d\over \d t}$.
Obviously, $Z(a,b)$ is also an annulus $A$ of conformal radius $R= e^{b-a}$.
Also denote $Z_k\deff Z(k, k+1)$.
\end{defi}

\smallskip
In other terminology, nodal curves are called {\sl prestable}. We shall
always suppose that $C$ is connected and has a ``finite topology", \ie $C$ has
finitely many irreducible components, finitely many nodal points, and that
$C$ has a smooth boundary $\d C$ consisting of finitely many smooth circles
$\gamma_i$, such that $\barr C \deff C \cup \d C$ is compact.

\newdefi{def5.2.2} A real oriented surface with boundary $(\Sigma,
\d\Sigma)$ {\sl parameterizes} a complex nodal curve $C$ if there
is a continuous map $\sigma :\barr\Sigma \to \barr C$ such that:

\sli if $a\in C$ is a nodal point, then $\gamma_a = \sigma\inv(a)$ is a
smooth imbedded circle in $\Sigma \bs \d \Sigma $, and if $a\not= b$ then
$\gamma_a \cap \gamma_b= \emptyset$;

\slii $\sigma :\barr\Sigma \bs \bigcup_{i=1}^N\gamma_{a_i}\to \barr C \bs \{
a_1,\ldots ,a_N\} $ is a diffeomorphism, where $a_1,\ldots ,a_N$ are the
nodes of $C$.
\end{defi}

\bigskip
\vbox{\xsize=.54\hsize\nolineskip\rm
\putm[.17][-.01]{\gamma_1}%
\putm[.31][.27]{\gamma_2}%
\putm[.53][0.275]{\gamma_3}
\putm[.71][.265]{\gamma_4}%
\putm[.87][.22]{\gamma_5}%
\putt[1.1][0]{\advance\hsize-1.1\xsize%
\centerline{Fig.~1}\smallskip
Circles $\gamma_1,..., \gamma_5$ are contracted by the parameterization
map $\sigma$ to nodal points  $a_1, \ldots a_5$.
}%
\putm[.56][.37]{\bigg\downarrow\sigma}%
\epsfxsize=\xsize\epsfbox{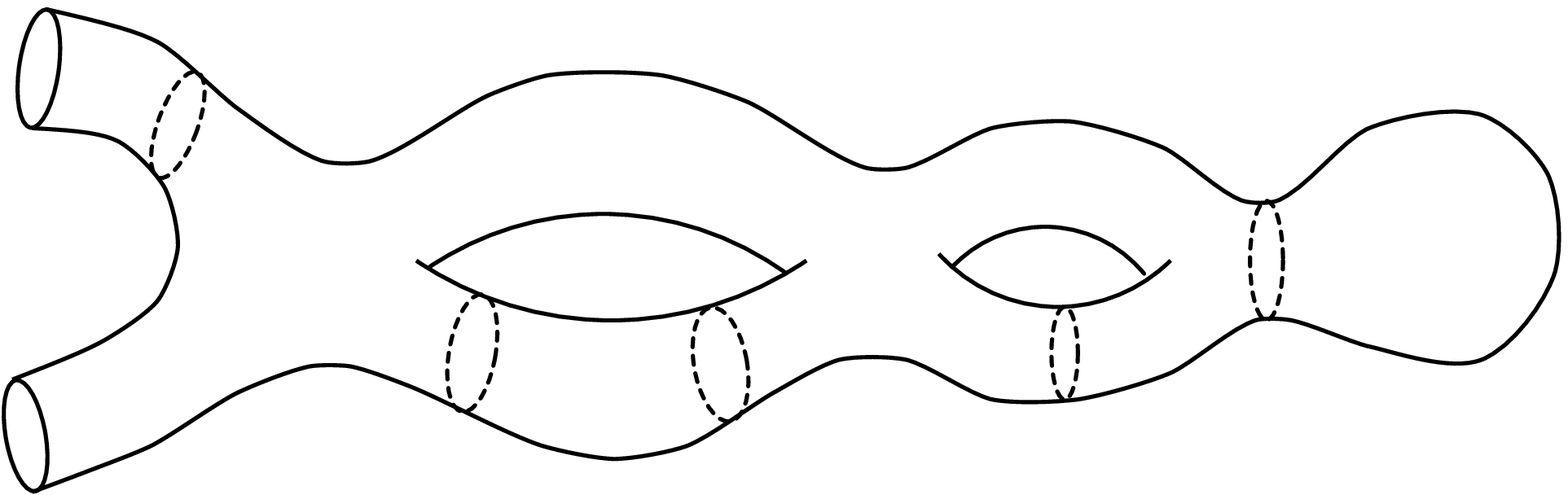}
\vskip.15\xsize
\putm[.20][.05]{a_1}%
\putm[.36][.255]{a_2}%
\putm[.49][0.27]{a_3}%
\putm[.725][.23]{a_4}%
\putm[.86][.19]{a_5}%
\epsfxsize=\xsize\epsfbox{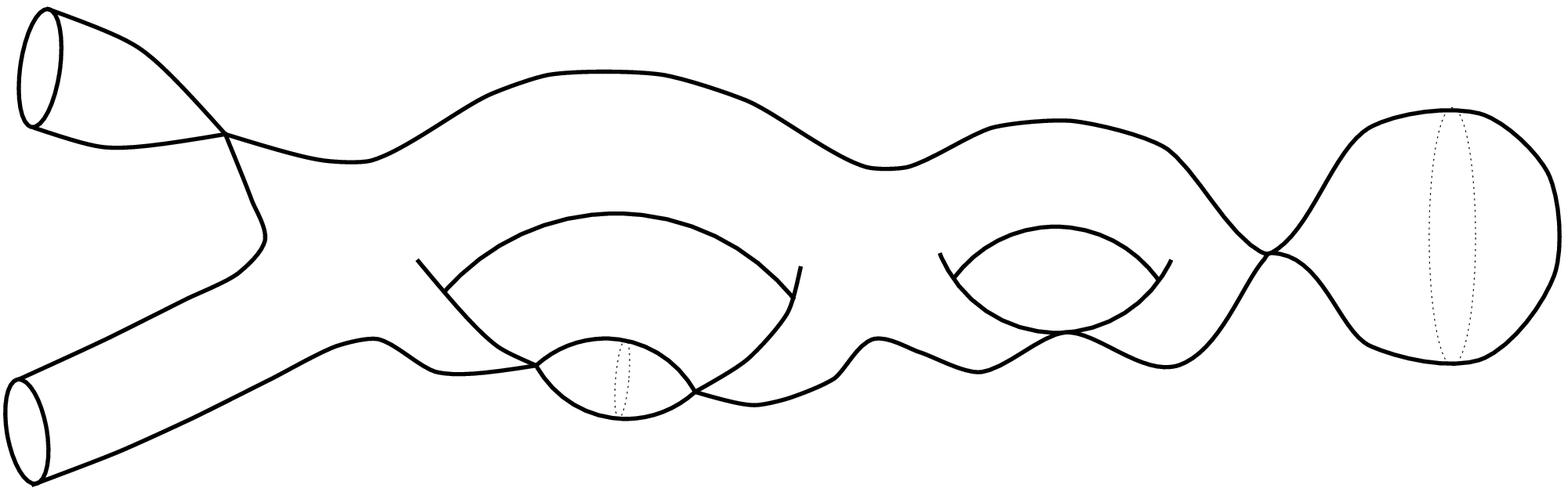}
}

\smallskip
Note that such a parameterization is not unique: if $g:\barr\Sigma \to
\barr\Sigma$ is any orientation preserving diffeomorphism then $\sigma \scirc
g: \barr\Sigma \to \barr C$ is again a parameterization.

A parameterization of a nodal curve $C$ by a real surface can be considered as
a method of ``smoothing'' of $C$. We shall also use an alternative method of
``smoothing'', the normalization.
Consider the normalization $\hat C$ of $C$. Mark on each component of this
normalization the pre-images (under the normalization map $\pi_C: \hat C \to
C$) of nodal points of $C$. Let $\hat C_i$ be a component of $\hat C$. We can
also obtain $\hat C_i$ by taking an appropriate irreducible component $C_i$,
replacing nodes contained in $C_i$ by pairs of discs with marked points, and
marking the remaining nodal points. Since it is convenient to consider components
in this form, we make the following

\newdefi{def5.2.3} A {\sl component $C'$} of a nodal curve $C$ is the 
normalization of an irreducible component of $C$ with marked points selected
as above.
\end{defi}

\smallskip
This definition allows us to introduce Sobolev and H\"older spaces of
functions and (continuous) maps of nodal curves.

\newdefi{def5.2.3a} A continuous map $u: C \to X$ is Sobolev $L^{1,p}
$-smooth, $u\in  L^{1,p}(C, X)$ if the induced maps $u_i \deff u\ogran_{C_i}
: C_i \to X$ of all of its components $C_i$ are $L^{1,p}$-smooth.
The notion of $J$-holomorphic maps $u: C \to X$ is similarly defined. 
For $u \in L^{1,p}(C, X)$ define $E_u \deff u^*TX$. Thus $E_u$ is determined 
by restricting to each component and by identifying fibers over pairs $(z', z'')$ 
of marked points corresponding to nodal points. An $L^{1,p}$-smooth section $v$ of 
$E_u$ over $C$ is given by a collection of sections $v_{C_i} \in L^{1,p}(C_i, 
E_u)$, one for every component $C_i$ of $C$, such that $v(z')=v(z'')$ for each pair
$(z', z'')$ of marked points chosen as above. Denote by $L^{1,p}(C, E_u)$ the space
of $L^{1,p}$-sections of $E_u$.
\end{defi}

\smallskip
\newdefi{def5.2.4} The {\sl energy} or the {\sl area} of a
continuous $L_\loc ^{1,2}$-smooth map \wrt a metric $h$ on $X$ is defined as
\begin{equation}
\area_h(u) \deff \norm{du}^2_{L^2(C)} = \int_C |du|^2_h
\eqqno(5.2.2)
\end{equation}
\end{defi}

This definition depends only on the complex structure on $C$ but not on the
choice of a metric on $C$ in the given conformal class. If an $\omega$-tame
almost complex structure $J$ is given, there is a prefered choice of a metric
$h$ on $X$ defined by $h(v,v) \deff \omega(v, Jv)$ for $v\in TX$.

\state Remark. Our definition of the area uses the following fact. Let $g$ be
a Riemannian
metric on $C$ compatible with $j_C$, $h$ a Riemannian metric on $X$, and $u:C
\to X$ a $J$-holomorphic immersion. Then $\norm{du}^2_{L^2(C)}$ is independent
of the choice of $g$ and coincides with the area of the image $u(C)$ \wrt
the metric $h_J(\cdot, \cdot) \deff \half(h(\cdot, \cdot) + h(J\cdot, J\cdot)
)$. The metric $h_J$ here can be seen as a ``Hermitization'' of $h$ \wrt
$J$. It is well-known that $\norm{du}^2_{L^2(C)}$ is independent of the choice
of a metric $g$ on $C$ in the same conformal class, see \eg \cite{S-U}.
Thus we can use the flat metric $dx^2 + dy^2$ to compare area and energy.
For a $J$-holomorphic map we obtain
\begin{equation}
\norm{du}^2_{L^2(C)}= \int_C |\d_x u|_h^2 + |\d_y u|_h^2 =
\int_C |\d_x u|_h^2 + |J \d_x u|_h^2 = \int_C |du|_{h_J}^2 = \area_{h_J}(u(C)),
\end{equation}
where the last equality is another well-known result, see \eg \cite{Gro}. Since we
consider varying almost complex structures on $X$, it is useful to know that
we can use any Riemannian metric on $X$ having a reasonable notion of area.

\smallskip
\newdefi{def5.2.5} A {\sl stable curve over $(X,J)$} is a pair
$(C,u)$, where $C$ is a nodal curve and $u: C\to X$ is a $J$-holomorphic map
satisfying the following condition: If $C'$ is a closed component of $C$ such
that $u$ is constant on $C'$, then there exist only finitely many biholomorphisms
of $C'$ which preserve the marked points of $C$. In this case $u$ is
called a {\sl stable map}.
\end{defi}

\state Remark. One can see that stability condition is nontrivial only in the
following cases:

\begin{itemize}
\item[{\sl 1)}] some component $C'$ is biholomorphic to $\cp^1$ with 1 or 2
marked points; in this case $u$ should be non-constant on any such component $C'$;
\item[{\sl 2)}] some irreducible component $C'$ of $C$ is $\cp^1$ or a torus
without nodal points.
\end{itemize}

\noindent
Since we consider only connected nodal curves, case {\sl 2)} can occur only
if $C$ irreducible, \ie $C'=C$. In this case $u$ must be non-constant on
$C$.

\newdefi{def5.2.6} A component $C'$ of a nodal curve $C$ is called
{\sl non-stable} in the following cases:

{ 1)} $C'$ is $\cp^1$ and has one or two marked points;

{ 2)} $C'$ is $\cp^1$ or a torus and has no marked points.

\noindent
Let $u: C \to X$ be a pseudoholomorphic  map. An irreducible component $C'$ of $C$
is a {\sl ghost} component (\wrt $u$) if $u$ is constant on $C'$. The {\sl ghost 
part $C^{gh}$ of $C$ (\wrt $u$)} is the union of all ghost components. In this 
paper we shall deal only with the case when all ghost components are closed.

A map $u$ is {\sl non-multiple} if, except finitely many points $z\in C$, 
one has $u\inv(u(z)) =\{z \}$. Note that this condition excludes also ghost 
components. 
\end{defi}

\smallskip
Now we are going to describe the Gromov topology on the space of stable curves 
over $X$ introduced in \cite{Gro}. Let $\{J_n\}$ be a sequence of continuous 
almost complex structures on $X$ which converges to $J_\infty$ in the $C^0
$-topology. Furthermore, let $(C_n, u_n)$ be a sequence of stable curves over 
$(X, J_n)$, such that all $C_n$ are parameterized by the same real surface $S$.

\newdefi{def5.2.7} We say that $(C_n,u_n)$ {\sl converges in the Gromov 
topology to a stable $J_\infty$-holomorphic curve $(C_\infty,u_\infty)$ over 
$X$} if the parameterizations $\sigma_n: \barr S \to \barr C_n$ and 
$\sigma_\infty: \barr S \to \barr C_\infty$ can be chosen in such a way 
that the following holds:

\sli $u_n\scirc \sigma_n$ converges to $u_\infty\scirc \sigma_\infty$ in the
$C^0(S, X)$-topology;

\slii if $\{ a_k \}$ is the set of nodes of $C_\infty$ and $\{\gamma_k\}$ are
the corresponding circles in $S$, then on any compact subset $K\comp
S \bs \cup_k \gamma_k$ the convergence $u_n\scirc \sigma_n \to u_\infty
\scirc \sigma_\infty$ is $L^{1,p}(K, X)$ for all $p< \infty$;

\sliii for any compact subset $K\comp \barr S \bs \cup_k\gamma_k$ there
exists $n_0=n_0(K)$ such that $ \sigma_n^{-1}(\{ a_k \}) \cap K= \emptyset$
for all $n\ge n_0$ and the complex structures $\sigma_n^*j_{C_n}$ converge
smoothly to $\sigma_\infty^*j_{C_\infty}$ on $K$;

\sliv the structures $\sigma_n^*j_{C_n}$ are independent of $n$ near the
boundary $\d S$.
\end{defi}

\medskip
Condition \sliv is trivial if $S $ is closed, but it is useful when
one considers the ``free boundary case'', \ie when $S$ (and thus all
$C_n$) are not closed and no boundary condition is imposed.

The reason for introducing the notion of a curve stable over $X$ is similar
to the one for the Gromov topology. We are looking for a completion of the
space of smooth imbedded pseudoholomorphic  curves which has ``nice''
properties, namely: \.1) such a completion should contain the limit of a
subsequence of every sequence of smooth curves which is bounded in an appropriate
sense; \.2) such a limit should also exist for every sequence in the
completed space; \.3) such a limit should be unique. The Gromov's compactness
theorem ensures us that the space of curves stable over $X$ has these nice
properties.

\newdefi{def5.2.8} Let $C_n$ be a sequence of nodal curves,
parameterized by the same real surface $S$. We say that the complex
structures on $C_n$ {\sl do not degenerate near boundary}, if there exist
$R>1$, such that for any $n$ and any boundary circle $\gamma_{n, i}$ of $C_n$
there exist an annulus $A_{n,i} \subset C_n$ adjacent to $\gamma_{n, i}$,
such that all $A_{n,i}$ are mutually disjoint, do not contain nodal points of
$C_n$, and have the same conformal radius $R$.
\end{defi}

Since the conformal radii of all $A_{n, i}$ are all the same, we can identify them
with $A(1,R)$. This means that all changes of complex structures of $C_n$
take place away from boundary. The condition is trivial if $C_n$ and $S$
are closed, $\d S = \d C_n = \emptyset$.

\state Remark. Changing our parameterizations $\sigma_n: S \to C_n$, we
may suppose that for any $i$ the pre-image $\sigma_n\inv (A_{n,i} )$ is the
same annulus $A_i$ independent of $n$.

\medskip
Now we state Gromov's compactness theorem for stable curves. Assume that
$X$ is a compact manifold and fix some Riemannian metric $h$ on $X$.

\newthm{thm5.2.1} Let $(C_n,u_n)$ be a sequence of stable
$J_n$-holo\-mor\-phic curves over $X$ with parameterizations $\delta_n: 
S \to C_n$. Suppose that:

\begin{itemize}
\item[{\sl a)}] $\{J_n\}$ is a sequence of continuous almost complex structures 
on $X$, which converges to $J_\infty$ in the $C^0$-topology;

\item[{\sl b)}] there is a  constant $M$ such that $\area_h [u_n (C_n)]\le M$
for all $n$;

\item[{\sl c)}] complex structures on the $C_n$ do not degenerate near the
boundary.

\end{itemize}

Then there is a subsequence $(C_{n_k},u_{n_k})$ and parameterizations $\sigma
_{n_k}: S \to C_{n_k}$, such that $(C_{n_k}, u_{n_k}, \sigma_{n_k})$
converges to a $J_\infty$-ho\-lo\-mor\-phic curve $(C_\infty,
u_\infty, \sigma_\infty)$ stable over $X$. 

Moreover, the limit curve $(C_{n_k}, u_{n_k})$ is unique up to the choice
of the parameterization $\sigma_\infty$.

Furthermore, if the structures $\delta_n^*j_{C_n}$ are constant on the fixed
annuli $A_i$, each adjacent to a boundary circle $\gamma_i$ of $S$, then
the new parameterizations $\sigma_{n_k}$ can be taken equal to $\delta_{n_k}$
on some subannuli $A'_i \subset A_i$, also adjacent to $\gamma_i$.
\end{thm}

A detailed proof of the theorem in the stated form can be found in \cite{Iv-Sh-3}.
We also refer to the original proof of Gromov in \cite{Gro}.

\smallskip
Gromov's compactness theorem induce a natural completion of the moduli space 
$\scrm$. Let a closed real surface $S$ of genus $g$ and a homology class $A \in 
\sfh_2(X,\zz)$ be given.

\newdefi{def5.2.9} Nodal $J$-holomorphic curves $u': C' \to X$ and $u'': C'' \to 
X$ are {\sl equivalent} if there exists a biholomorphism $\phi: C' \to C''$ with 
$u' = u'' \scirc \phi$. The {\sl total moduli space $\barm^{st}$ of stable nodal 
curves over $X$} is the set of equivalence classes $[C,u, J]$ with $J \in \scrj$ 
and $u: C \to X$ a stable $J$-holomorphic curve representing a given class $A\in 
\sfh_2(X,\zz)$. The space $\barm^{st}$ is equipped with the {\sl Gromov topology}
in which a sequence $[C_n,u_n, J_n]$ converges to $[C_\infty,u_\infty, J_\infty]$ 
if $J_n$ converges to $J_\infty$ in the $C^\ell$-topology and $(C_n, u_n)$ to $(C
_\infty, u_\infty,)$ in the sense of \refdefi{def5.2.7}. Denote by $\pr_{\!\!
\scrj}^{st}$ the natural projection $\pr_{\!\!\scrj}^{st}:  [C, u,J] \in \barm
^{Gr} \mapsto J \in \scrj$.

Define the {\sl Gromov compactification} $\barm^{Gr}$ of the total moduli space 
$\scrm$ of pseudoholomorphic curves $X$ as the closure of $\scrm$ in $\barm^{st}$.
\end{defi}

Note that every fiber $\barm^{Gr}_J \deff \barm^{Gr} \cap \bigl(\pr_{\!\!\scrj}
^{st} \bigr)\inv (J)$ is compact. Note also that in general $\barm^{st} \not=
\barm^{Gr}$. This means that there are stable curves $[C,u,J] \in \barm^{st}$
which can not be reached from $\scrm$.

\newsubsection[5.2c] {The cycle topology for pseudoholomorphic curves} The 
Gromov compactness theorem gives a precise description of the behavior of {\sl 
parameterized} pseudoholomorphic curves at ``infinity'' of the total moduli
space. However, what we are really interested in is not a pseudoholomorphic map
$u:C \to X$ itself but rather the image $u(C)\subset X$, \ie a {\sl 
non-pa\-ra\-me\-te\-rized} pseudoholomorphic curve. The natural space
where non-pa\-ra\-me\-te\-rized curves ``live'' is the space $\scrz_2(X)$ of 
2-currents on the ambient manifold $X$. Recall that $\scrz_2(X)$ is the dual space
to the space $C^\infty(X, \Lambda^2X)$ of smooth 2-forms on $X$ (see \eg
\cite{Gr-Ha}, {\sl Chapter 3}).

\newdefi{def5.2c.1} Let $X$ be a manifold, $C$ an abstract nodal curve with the 
smooth boundary $\d C$ such that $\barr C \deff C \cup \d C$ is compact, and $u: 
C \to X$ a map which is $L^{1,p}$-smooth up to boundary. Define the {\sl cycle 
$u[C]$ associated with the map $u: C \to X$} as the current whose pairing with 
a smooth $2$-form $\phi$ on $X$ equals $\< u(C), \phi \> \deff \int_C u^*\phi$. In
this case we also say that $u[C]$ is {\sl represented by the map $u: C \to X$}.

If additionally $u: C \to X$ is $J$-holomorphic \wrt some an almost complex 
structure on $X$, we call $C' \deff u[C]$ an {\sl $J$-holo\-mor\-phic curve $C$ 
\underline{\it in} $X$}. In this case we say that the $J$-holomorphic curve $(C, 
u)$ \underline{\it over} $X$ and the map $u:C\to X$ {\sl represent the curve 
$C'= u[C]$}. The set $u(C)$ is called the {\sl support of $C'$} and denoted 
by $\supp(C')$.

A curve $C'$ in $X$ is {\sl non-multiple} if it can be represented by a 
non-multiple pseudoholomorphic map $u:C\to X$ (see {\sl Definitions 
\ref{def1.2.2}} and {\sl\ref{def5.2.6}}). In this case we identify the set $u(C)$ 
and the current $u[C]$ and use the same notation $u(C)$.

A sequence of cycles $u_n[C_n]$  {\sl converges to a cycle $u_\infty[C_\infty]$}
if $\< u_n(C_n), \phi \>$ converges to $\< u_\infty(C_\infty), \phi \>$ for any 
smooth smooth $2$-form $\phi$ on $X$. In other words, the cycle topology is the 
topology induced from the space of currents $\scrz_2(X)$.
\end{defi}

\smallskip
\newlemma{lem5.2c.1} \sli Let $(X,J)$ be an almost complex manifold
and $(C, u)$, $(C', u')$ closed $J$-holomorphic curves over $X$. Assume that
\begin{itemize}
\item $C$ and $C'$ are parameterized by the same closed surface $S$;
\item $(C, u)$ contains no multiple and ghost components;
\item the associated cycles $u[C]$ and $u'[C']$ coincide. 
\end{itemize}
Then $(C, u)$ are $(C', u')$ equivalent.

\slii Let $J_n$ be a sequence of continuous almost complex 
structures on $X$ which converges to an almost complex structure $J_\infty$ in the
$C^0$-topology, and $(C_n, u_n)$ a sequence of stable $J_n$-holomorphic closed 
curves over $X$ which converges to $(C_\infty, u_\infty)$ in the Gromov topology. 
Then $u_n[C_n]$ converges to $u_\infty[C_\infty]$ in the cycle topology.

\sliii Let $J_n$ be a sequence of continuous almost complex 
structures on $X$ which converges to an almost complex structure $J_\infty$ in the
$C^0$-topology, $(C_n, u_n)$ a sequence of stable $J_n$-holomorphic closed 
curves over $X$, and $(C_\infty, u_\infty)$ a parameterized $J_\infty$-holomorphic 
curve. Assume that 
\begin{itemize}
\item $u_n[C_n]$ converges to $u_\infty[C_\infty]$ in the cycle topology;
\item $C_n$ and $C_\infty$ are parameterized by the same closed surface $S$;
\item $(C_\infty, u_\infty)$ contains no multiple and ghost components;
\item $J_\infty$ is $C^1$-smooth.
\end{itemize}
Then $(C_n, u_n)$ converges to $(C_\infty, u_\infty)$ in the Gromov topology. 
\end{lem}

\proof {\sl Part \slip}. The hypotheses on $(C, u)$ and $(C', u')$ imply
that $(C', u')$ also contains no multiple and ghost components. The claim then 
follows from the unique continuation property of pseudoholomorphic curves
(see \lemma{lem1.2.4}).

{\sl Part \sliip}. This follows from the definition of the Gromov 
topology and the description of the convergence at nodes given in 
{\sl Step 0)} of the proof of \lemma{lem5.2.2} below.

{\sl Part \sliiip}. 
Fix a $J_\infty$-Hermitian metric $h$ on $X$. Let $\omega$ be the associated 
2-form, $\omega(v,w) \deff h(J_\infty v,w)$. Then the structures $J_n$ are 
$\omega$-tame for $n\gg1$.
Note that even if $\omega$ is apriori only continuous and not closed, the notion 
of $\omega$-tameness is still meaningful. Moreover, the $\omega$-tameness provides
a uniform bound of $h$-area of $u_n[C_n]$. Consequently, some subsequence
$(C_{n'}, u_{n'})$ of $(C_n, u_n)$ converges in the Gromov topology to 
a stable $J_\infty$-holomorphic curve $(C'_\infty, u'_\infty)$. The hypotheses
of the corollary imply that $(C'_\infty, u'_\infty)$ is equivalent to 
$(C_\infty, u_\infty)$ and the result follows.\qed

\state Remark. The meaning of \lemma{lem5.2c.1} is that, in the absence of 
multiple and ghost components, the notions of pseudoholomorphic curves {\sl over} 
$X$ and {\sl in} $X$ essentially coincide. The same also holds for the Gromov 
and the cycle topologies. Note also that several authors (\cite{Ye}, 
\cite{Pa-Wo}, \cite{Hum}) considered a weaker version of the Gromov compactness
theorem where the cycle topology is used instead of the Gromov one.

\newdefi{def5.2c.2} Define the {\sl cycle compactification $\barm$} of the total 
moduli space $\scrm$ as the set of pairs $(C,J)$, where $J \in \scrj$ and $C$ is a
$J$-holomorphic curve {\it in} $X$ which considered as the cycle $u[C']$ 
represented 
by some $J$-holomorphic map $u: C' \to X$. Equip the space $\barm$ with the 
{\sl cycle topology} in which a sequence $(C_n, J_n)$ converges to $(C_\infty, 
J_\infty)$ if $J_n$ converges to $J_\infty$ in the $C^\ell$-topology and $C_n$ 
converges to $C_\infty$ in the sense of \refdefi{def5.2.7}. Denote by $\pr_{\!\!
\scrj}$ the natural projection $\pr_{\!\!\scrj}:  (C,J) \in \barm \mapsto J \in 
\scrj$. Define the natural projection $\pr^{Gr}:\barm^{Gr} \to \barm$ by
$\pr^{Gr}: (C, u,J) \in \barm^{Gr} \mapsto u[C]\in \barm$.
\end{defi}

\newdefi{def5.2c.3} A {\sl normal parameterization} of a $J$-holomorphic curve 
$C$ in $X$ is given by a Riemann surface $S$ (possibly not connected) and a map $u:
S \to X$ such that 
\begin{itemize}
\item[1)] $u[S]=C$;
\item[2)] $u$ is $J$-holomorphic and $L^{1,p}$-smooth up to boundary;
\item[3)] the restriction of $u$ to every connected component of $S$ is 
non-multiple; in particular, there are no ghost components,  \ie $u$ is 
non-constant on every connected component of $S$;
\item[4)] the number of boundary circles of $S$ is as small as possible.
\end{itemize}
\end{defi}

\state Remark. Without condition (3) one could add new ghost spheres to $S$ 
and make the Euler characteristic $\Chi(S)$ arbitrarily large. Condition (4) 
excludes the possibility of dividing components of $S$ into pieces which also 
allows to increase $\Chi(S)$.

\newlemma{lem5.2c.2} Let $J$ be a $C^1$-smooth almost complex structure on $X$,
$C$ an abstract nodal curve, and $u: C\to X$ a $J$-holomorphic map which is an 
imbedding near the boundary $\d C$. Then, up to a diffeomorphism, there exists 
a unique normal parameterization $\ti u: S \to X$ of $u[C]$.
\end{lem}

\proof Let $C= \cup_i C_i$ be the decomposition of $C$ into irreducible 
components. Denote by $m_i$ the degree $u$ on $C_i$. This means that
\begin{itemize}
\item $m_i=0$ if $C_i$ is a ghost component, \ie $u$ is constant on $C_i$;
\item $m_i$ is the number of points in the preimage $u\inv(x)$ for a
generic $x\in u(C_i)$ otherwise.
\end{itemize}
For every non-zero $m_i$, denote by $S_i$ the normalization of the image $u(C_i)$.
Denote by $\ti u_i: S_i \to u(C_i)$ the corresponding normalization maps. In 
particular, $m_i=1$ for every non-closed component $C_i$ and in this case $S_i$ 
is the normalization of $C_i$. Define $S$ as the disjoint union of the surfaces
$S_i$, each taken $m_i$ times. Let $\ti u: S \to X$ be the map which coincides 
with the composition $\ti u_i: S_i \to u(C_i) \hook X$ on every copy of $S_i$.
One can see that $\ti u: S \to X$ is a normal parameterization of $u[C]$.

The uniqueness of such a normal parameterization follows from \lemma{lem1.2.4}
\qed

\state Remark. Let us give an example showing that the condition on the behavior 
of $u$ at the boundary imposed in \lemma{lem5.2c.2} is necessary. Define curves
$C'$ and $C''$ as the disjoint unions $C'\deff \{ z\in \cc: |z|<2\} \sqcup \{ z\in
\cc: 1<|z|<3\}$ and $C''\deff \{ z\in \cc: |z|<3\} \sqcup \{z\in \cc:
1<|z|<2\}$. Let $u': C' \to \cc$ and $u'': C'' \to \cc$ be the maps which are the 
standard imbeddings on every component of $C'$ and $C''$. Then
obviously $(C',u')$ and $(C'',u'')$ are not equivalent in the sense of 
\refdefi{def5.2.9} and define non-equivalent normal parameterizations
of $u'[C'] =u''[C'']$.

\newcorol{cor5.2c.3} Under the hypotheses of \lemma{lem5.2c.2}, the curve $u[C]$,
considered as a current $u[C]\in \scrz_2(X)$, admits a unique representation in 
the form 
$$
\textstyle
u[C]= \sum_i m_i u_i[C_i], 
$$
where the $u_i: C_i \to X$ are $J$-holomorphic maps and the $u_i(C_i)$ are 
the irreducible components of $\supp(u[C])$. 
\end{corol}

The corollary ensures that the notions of an irreducible component and the 
multiplicity of a closed pseudoholomorphic curve {\it in} $X$ are well-defined.

\medskip
The importance of the notion of a normal parameterization lies in the fact that it 
allows us to define a natural stratification of the cycle compactification $\barm$
of the total moduli space. Let $S$ be a given connected real surface $S$ of genus 
$g$ and $[C] \in \sfh_2(X, \zz)$ a homology class, and $\barm=\barm(S, X, [C])$ 
the cycle compactification of the total space $\scrm=\scrm(S, X, [C])$ of 
irreducible pseudoholomorphic curves of genus $g$ in the homology class $[C]$. 
Take $(C,J) \in \scrm(S, X, [C])$ and consider a normal parameterization $u': S' 
\to X$ of $C$. Let $C= \sum_i m_i C_i$ be the decomposition into irreducible 
components in the sense of \refcorol{cor5.2c.3}. Restricting $u'$ to 
appropriate connected components we obtain normal parameterizations $u'_i: S'_i
\to X$ of the corresponding $C_i$. 

\newdefi{def5.2c.4} The {\sl topological type of a component $C_i$} is the triple
$(S'_i, m_i, [C_i])$, where $[C_i]$ denotes the homology class of $C_i$. The {\sl
topological type $\bftau$ of a curve $(C,J) \in \scrm(S, X, [C])$} is the sequence
of all topological types of components $(S'_i, m_i, [C_i])$ defined up to 
permutation. 
\end{defi}

\newlemma{lem5.2c.4} \sli The space $\scrm_\bftau$ of pseudoholomorphic curves 
$(C,J) \in \barm(S, X, [C])$ of a given topological type $\tau$ is a 
$C^\ell$-smooth Banach manifold. The natural projection $\pr_{\!\!\scrj}: \scrm
_\bftau \to \scrj$ is a $C^\ell$-smooth Fredholm map.

\slii The space $\barm=\barm(S, X, [C])$ is the union of subspaces
$\scrm_\bftau$. 
\end{lem}

\proof 
The decomposition $C= \sum_i m_i C_i$ of every $(C,J) \in \scrm_\bftau$
shows that $\scrm_\bftau$ is the fiber product of the spaces $\scrm(S'_i, X, 
 [C_i])$ over all triples $(S'_i, m_i, [C_i]) \in \bftau$ taken over the space
$\scrj \!\!$, 
$$
\textstyle
\scrm_\bftau=\prod_{(S'_i, m_i, [C_i]) \in \bftau } 
\scrm(S'_i, X, [C_i]) \bigm/ \!\!\!\scrj\!\!.
$$
Checking the transversality condition, one obtains the desired  differentiable 
structure on $\scrm_\bftau$.

The second assertion of the lemma is obvious. \qed 

\smallskip

\newsubsection[5.2a]{Fine apriori estimates for convergence at a node}
For the purpose of this paper we need a refined version of the {\sl Second
apriori estimate} given in \cite{Iv-Sh-3}, {\sl Lemma 3.4}. This gives a 
precise description with estimates of the Gromov convergence in neighborhoods 
of the contracted circles.

\newlemma{lem5.2.2} Let $X$ be a compact manifold $X$, $J^*$ a
$C^{0,s}$-smooth almost complex structure on $X$ with $s>0$, and $h$ a metric on 
$X$. Then there exist constants $\epsi=\epsi(X,h, J^*, s)>0$ and $C <\infty$ 
such that for any $C^{0,s}$-smooth almost complex structure $J$ with 
\begin{equation}
\norm{J -J^*}_{C^{0,s}(X)} \le \epsi
\eqqno(5.2a.1)
\end{equation}
and any $J$-holomorphic map $u: Z(0, l) \to X$ the condition
\begin{equation}
\norm{du}_{L^2(Z_k)} \le \epsi \quad \text{for any }k\in[0, l-1]
\eqqno(5.2.3)
\end{equation}
implies the uniform estimate
\begin{equation}
\norm{du}^2_{L^2(Z_k)} \le
C \cdot e^{-2k} \cdot \norm{du}^2_{L^2(Z(0, 2))} +
C \cdot e^{-2 (l-k)} \cdot\norm{du}^2_{L^2(Z(l-2, l))}
\eqqno(5.2.4)
\end{equation}
for any $k\in[1, l-2]$.
\end{lem}

\proof \step0.\. {\sl Lemma 3.3} in \cite{Iv-Sh-3} states that
under hypotheses of the lemma one has a ``local'' estimate
\begin{equation}
\norm{du}^2_{L^2(Z_k)} \le \msmall{\gamma \over 2}
\left( \norm{du}^2_{L^2(Z_{k-1})} +
\norm{du}^2_{L^2(Z_{k+1})}  \right ) \qquad \text{for any $k\in[1, l-2]$}
\eqqno(5.2.5)
\end{equation}
with a universal constant $\gamma <1$. Then in {\sl Corollary 3.4} in
\cite{Iv-Sh-3} it is shown that \eqqref(5.2.5) implies the estimate
\begin{equation}
\norm{du}^2_{L^2(Z_k)} \le
e^{-2\alpha (k-1)} \cdot \norm{du}^2_{L^2(Z(0, 2))} +
e^{-2\alpha (l-2-k)} \cdot\norm{du}^2_{L^2(Z(l-2, l))}
\eqqno(5.2.6)
\end{equation}
for any $k\in[1, l-2]$ with a constant $\alpha >0$ related to $\gamma$
by $\gamma = {1\over \cosh(2\alpha)}$.

\state Remark. Note that in the proof of the estimates \eqqref(5.2.5) and 
\eqqref(5.2.6) are proven in \cite{Iv-Sh-3} under the following assumption:
It is supposed that $J^*$ and $J$ in question are only continuous and that
$\norm{J- J^*}_{C^0(X)} \le \epsi'$ for some $\epsi'=\epsi'(X,J^*,h)>0$
independent of $J$. 

\medskip
From the relation $\gamma = {1\over \cosh(2\alpha)}$ we see that the smaller the $
\gamma$ we have the bigger the $\alpha$ in 
\eqqref(5.2.6) we obtain. For our purpose it would be sufficient to prove 
estimate \eqqref(5.2.5) with the parameter $\gamma^* \deff {1\over \cosh2}$. 
Note however that in the ``ideal'' case when $(X,J,h)$ is $\cc^n$ with the 
standard complex and Hermitian structures, $\gamma^*$ is exactly the best 
possible constant, see the proof of {\sl Lemma 3.3} in \cite{Iv-Sh-3}. Thus one 
can not expect that estimate \eqqref(5.2.5) holds with {\sl uniform} $\gamma \le 
\gamma^*$. The idea is to consider \eqqref(5.2.5) with parameters $\gamma_k$ 
depending on $k$ and to estimate the difference $\gamma_k -\gamma^*$.

\medskip\noindent
\step1.\. {\sl Under hypotheses of the lemma, for any $k\in[1, l-2]$,
one has the estimate
\begin{equation}
\norm{du}^2_{L^2(Z_k)} \le 
\msmall{\gamma_k \over 2} 
\cdot
\left( \norm{du}^2_{L^2(Z_{k-1})} +
\norm{du}^2_{L^2(Z_{k+1})}  \right )
\eqqno(5.2.7)
\end{equation}
for 
\begin{equation}
\gamma_k \deff \gamma^* + 
C_1\cdot\left(e^{-\alpha s k} + e^{-\alpha s (l-k)} \right)
\eqqno(5.2a.2)
\end{equation}
with the parameter $\alpha >0$ as in \step0 and some constant $C_1$
depending only on $X$, $h$, $J^*$, and $s$.}

\smallskip
While proving this estimate we shall denote by $C$ a constant whose particular
value is not important and which may not be the same in different formulas.
The main condition is that these constants are {\sl uniform}, \ie independent
of $J$, $u$, and $l$, and depend only on $X$, $h$, $J^*$, and $s$.

Estimates  \eqqref(5.2.3) and  \eqqref(5.2.6) together with apriori 
estimates show that
\begin{equation}
\diam(u(Z(k-1, k+2)) \le C \cdot 
\left(e^{-\alpha k} + e^{-\alpha (l-k)} \right).
\end{equation}
Consequently, due to a uniform H\"older $C^{0,s}$-estimate on $J$,
for the oscillation of $J$ on the image $u(Z(k-1, l-k))$ we obtain
\begin{equation}\eqqno(5.2a.2a)
\osc(J, u(Z(k-1, k+2))) \le C \cdot 
\left(e^{-\alpha s k} + e^{-\alpha s (l-k)} \right).
\end{equation}
This implies that in a neighborhood of each $u(Z(k-1, k+2))$ there exist an 
integrable structure $J\st$ and a flat (\ie Euclidean) metric $h\st$ such that
\begin{equation}
\norm{J -J\st}_{L^\infty(u(Z_k))} + \norm{h -h\st}_{L^\infty(u(Z_k))}
\le C\cdot \left(e^{-\alpha s k} + e^{-\alpha s (l-k)} \right).
\end{equation}
Using this we obtain estimates
\begin{align}
\norm{\dbar\st u}_{L^2(Z_k)}
= \norm{\dbar\st u - \dbar_J u}_{L^2(Z_k)}
&\le \norm{J -J\st}_{L^\infty(u(Z_k))} \cdot \norm{du}_{L^2(Z_k)} \le
\notag
\\
&\le C\cdot \left(e^{-\alpha s k} + e^{-\alpha s (l-k)} \right);
\end{align}
\begin{equation}
\Bigl| \norm{du}_{L^2(Z_k), h} 
-\norm{du}_{L^2(Z_k), h\st} \Bigr| \le
C\cdot \left(e^{-\alpha s k} + e^{-\alpha s (l-k)} \right)
\cdot \norm{du}_{L^2(Z_k), h}.
\end{equation}
In particular, we can use $h\st$ instead of $h$ in our estimates.

Now consider $U$ as a subset on $\cc^n$ with the standard $J\st$ and $h\st$.
Then we can find $u_\dbar\in L^{1,2}(Z(k-1, k+2), \cc^n)$ such that $\dbar\st
u_\dbar= \dbar\st u$ and
\begin{equation}
\norm{du_\dbar}_{L^2(Z(k-1, k+2))} \le C \norm{\dbar\st u}
_{L^2(Z(k-1, k+2))}.
\eqqno(5.2a.3)
\end{equation}
Set $u_\scro \deff u -u_\dbar$, so that $u_\scro$ is $J\st$-holomorphic.
It follows that 
\begin{equation}
\norm{du_\scro}^2_{L^2(Z_k)} \le \msmall{\gamma^* \over2}
\left( \norm{du_\scro}^2_{L^2(Z_{k-1})} +
\norm{du_\scro}^2_{L^2(Z_{k+1})}  \right).
\eqqno(5.2.8)
\end{equation}
Together with the estimates on $u_\dbar$, \eqqref(5.2.8) implies 
\eqqref(5.2.7).
          
\medskip\noindent
\step2.\. {\sl There exist a uniform $k_0=k_0(X, h, J^*,s)$ and $A^\pm_k$, 
$k=k_0, \ldots,l-k_0 $ with the properties

\sli $A^\pm_k$ are ``supersolutions'' of \eqqref(5.2.7), \ie
\begin{equation}
\eqqno(5.2a.4a)
A^\pm_k \ge \msmall{\gamma_k \over2} (A^\pm_{k-1} + A^\pm_{k+1})
\end{equation}

\slii  $A^\pm_k$ have the desired exponential decay}
\begin{equation}
\eqqno(5.2a.5a)
A^+_k \le C\cdot e^{-2k},
\qquad\qquad
A^-_k \le C\cdot e^{-2(l-k)}.
\end{equation}

\sliii $\gamma_k <1$ for $k \in [k_0, l-k_0]$.

\smallskip
Fix $k^*\in \zz$ such that $l-1 \le k^* <l+1$, so that $k^* \approx {l\over2}$.
Set 
\begin{align}
A^+_k & \deff \cases 
      e^{-2k -{1\over k}} & 0\le k\le k^* \\
      e^{-2k -{1\over k^*} + {1\over l-k}- {1\over l-k^*}} 
   \hskip7pt& k^*\le k \le l
\endcases  \\
A^- _k& \deff \cases 
      e^{-2(l-k) -{1\over l-k^*} + {1\over k}-{1\over k^*}}
    & 0\le k\le k^* \\
      e^{-2(l-k) -{1\over l-k}} & k^*\le k \le l
\endcases  
\end{align}
Making the Taylor expansion in $k\inv$ we obtain 
$$
\msmall{ 2 A^\pm_k \over A^\pm_{k-1} + A^\pm_{k+1}}= 
\cases
   \msmall { 1\over \cosh(2)}+  
   \msmall { \sinh(2) \over \cosh^2(2)} \cdot k^{-2} + O(k^{-3})
 & \text{for }0< k \le k^*;
\\
   \msmall { 1\over \cosh(2)}+  
   \msmall { \sinh(2) \over \cosh^2(2)} \cdot (l-k)^{-2} + O((l-k)^{-3})
 & \text{for }k^*\le k <l;
\endcases
$$
So the existence of the desired $k_0(s)$ follows from the asymptotic behavior
$C_1 e^{-\alpha s k} = o(k^{-2})$ for $k\lrar \infty$.

\medskip\noindent
\step3. {\sl There exists a constant $C_2= C_2(X,h, J^*,s)$ such that 

$\qquad
\norm{du}^2_{L^2(Z_k)} \le C_2\cdot \left(
A^+_k \cdot \norm{du}^2_{L^2(Z(0,2))} +
A^-_k \cdot\norm{du}^2_{L^2(Z(l-2,l))} 
\right)
$
\lineeqqno(5.2.9)

\smallskip\noindent
for any $k\in[k_0, l-k_0]$ with the uniform constant $k_0=k_0(s)$ chosen as 
above.}

\smallskip
Obviously, \eqqref(5.2.9) implies the claim of the lemma. Set 
$$
A^*:= C_2 \cdot \left(A^+_k \cdot \norm{du}^2_{L^2(Z(0,2))} +
A^-_k \cdot\norm{du}^2_{L^2(Z(l-2,l))} \right)
$$
and choose a constant $C_2$ so that 
$$
A^*_{k_0} \ge \norm{du}^2_{L^2(Z_{k_0})}
\qquad \text{and} \qquad
A^*_{l-k_0} \ge
\norm{du}^2_{L^2(Z_{l-k_0})}. 
$$
Then by \eqqref(5.2.7) and \eqqref(5.2a.4a) 
\begin{equation}\notag
\norm{du}^2_{L^2(Z_k)} - A^*_k \le 
\msmall{\gamma_k \over2} \cdot 
\left(\norm{du}^2_{L^2(Z_{k-1})} - A^*_{k-1}
+\norm{du}^2_{L^2(Z_{k+1})} - A^*_{k+1} \right).
\end{equation}
Find $k_\max\in [k_0, l-k_0]$ realizing the maximum of
$\norm{du}^2_{L^2(Z_k)} - A^*_k$. Then
\begin{equation}
\matrix
\norm{du}^2_{L^2(Z_{k_\max})} - A^*_{k_\max}
&\le& \msmall{\gamma_{k_\max} \over2}
\left(\norm{du}^2_{L^2(Z_{k_\max-1})} - A^*_{k-1}
+\norm{du}^2_{L^2(Z_{k+1})} - A^*_{k_\max+1} \right)
\cr\noalign{\vskip3pt}
&\le&
\gamma_{k_\max} \cdot( \norm{du}^2_{L^2(Z_{k_\max})} - A^*_{k_\max}).
\endmatrix
\end{equation}
Since $\gamma_{k_\max}<1$, the last inequality holds only if
$\norm{du}^2_{L^2(Z_{k_\max})} \le A^*_{k_\max}$. Thus
\begin{equation}
\norm{du}^2_{L^2(Z_k)} \le A^*_k
\qquad\text{for any }k\in [k_0, l-k_0].
\end{equation}
This finishes the proof. \qed

\smallskip
\newthm{thm5.2.3} Let $J^*$ be a $C^{0,s}$-smooth almost complex structure on the 
ball $B\subset \rr^{2n}$ with $0<s<1$. Then there exists $\epsi^*= \epsi^* (J^*,
s)$ with the following property. For any almost complex  structure $J$ on $B$ with 
$\norm{J-J^*}_{C^{0,s}(B)} \le \epsi^*$ and any $J$-holomorphic map 
$u: Z(0,l) \to B(\half)$ with $l\ge3$ satisfying the condition 
$$
\norm{du}_{L^2(Z_k)} \le \epsi^*\qquad\text{for any }k\in[1,l]
$$
there exist a {\sl linear} complex structure $J\st$ in $\rr^{2n}$ and vectors 
$v^+, v^0, v^- \in \rr^{2n}$ such that 
\begin{multline}\eqqno(5.2.10)
\left\Vert
u - \left( e^{-t+J\st\theta}v^+ + v^0 + e^{t-J\st\theta}v^- \right)
\right\Vert^2_{L^{1,2}(Z_k)} \le
\\
\le C^*  \cdot\bigl(k \;e^{-2(1+s)k}\norm{du}^2_{L^2(Z(0,2))} + 
(l-k)\;e^{-2(1+s)(l-k)}\norm{du}^2_{L^2(Z(l-2,l))} 
\bigr)\ 
\end{multline}
for any $k=1,\ldots, l-1$ with a constant $C^*=C^*(J^*,s) <\infty$ independent of 
$J$, $l$, and $u$.
\end{thm}

\proof In fact, we prove that for the constant $\epsi^*$ one can take the 
$\epsi$ from \lemma{lem5.2.2}. The proof also exploits the same ideas which were 
used in the proof of that lemma. 

\smallskip\noindent {\sl Step 1}. 
Let $J$ and $u: Z(0,l) \to B(\half)$ be as in the hypotheses of the theorem. For 
$k= 1,\ldots,l$, let $x_k \in B(\half)$ be the average value of $u$ on $Z_k$ \wrt
the cylinder metric, \ie
$$
x_k \deff \bint_{Z_k} u \deff
\msmall{1\over 2\pi}\int_{(t,\theta) \in Z_k} u(t,\theta) 
dt\,d\theta.
$$
Define the complex structures $J_k$ by $J_k\deff J(u(k,0))$. We consider every 
$J_k$ as a {\sl linear} complex structure in $\rr^{2n}$, \ie constant in $x\in 
\rr^{2n}$. Further, any $k=1,\ldots,l$ we define the metric $g_k$ setting $g_k(v, 
w) \deff \half (g\st(v,w) + g\st(J_kv, J_kw))$, where $g\st$ denotes the standard 
Euclidean metric in $\rr^{2n}$. Then $g_k$ are linear in the same sense as $J_k$. 
In  computing various norms related to $Z_k$ or $Z(k-2, k+1)$, we shall use the 
metric $g_k$ without indicating this in the notation. Observe that all $g_k$ are 
equivalent since the $J_k$ are uniformly bounded. Further, convergence of 
$J_{k_\nu}$ implies convergence of $g_{k_\nu}$.

For any $k=1,\ldots,l$ there exist uniquely defined vectors 
$v^+_k, v^0_k, v^-_k \in \rr^{2n}$ such that for the function 
$$
v_k(t, \theta) \deff  e^{-t+J_k\theta}v^+_k  
+v^0_k+  e^{t-J_k\theta}v^-_k 
$$
the norm $\norm{u -v_k}_{L^{1,2}(Z_k)}$ (computed with $g_k$) attains the minimum.

We claim that under the hypotheses of the theorem there exist a constant $C_1=
C_1(J^*, s)$ and and an integer $k_0=k_0(J^*, s)$ such that for any integer 
$k=k_0,\ldots, l-k_0$
\begin{multline}\eqqno(5.2.11)
\norm{u -v_k}^2_{L^{1,2}(Z_k)} +
\norm{v_{k-1} -v_k}_{L^{1,2}(Z_k)}^2 +
\norm{v_{k+1} -v_k}_{L^{1,2}(Z_k)}^2 \le
\\
\le \gamma_s \cdot\left(
   \norm{u -v_{k-1}}^2_{L^{1,2}(Z_{k-1})} +
   \norm{u -v_{k+1}}^2_{L^{1,2}(Z_{k+1})} \right)  +
\\
+C_1  \cdot \bigl(e^{-2(1+s)k} \norm{du}^2_{L^2(Z(0,2))} + 
e^{-2(1+s)(l-k)}\norm{du}^2_{L^2(Z(l-2,l))}\bigr) \ 
\end{multline}
with the parameter $\gamma_s \deff {1\over \cosh(2+2s)}\cdot\;$ 
Assuming the contrary, there must exist sequences of 
\begin{itemize}
\item integers  $l_\nu\lrar \infty$;
\item integers $k_\nu \lrar \infty$ with $l_\nu -k_\nu \lrar \infty$;
\item structures $J_\nu$ in $B$ with $\norm{J_\nu -J^*}_{C^{0,s}(B)} \le \epsi^*$;
\item $J_\nu$-holomorphic maps $u_\nu: Z(0,l_\nu) \to B(\half)$ with 
$\sup_{k=0,\ldots,l_\nu} \norm{du_\nu}_{L^2(Z_k)} \lrar 0$
\end{itemize}
with the following property. For the points $x_{\nu,k} \deff \bint_{Z_k} u_\nu$, 
the linear complex structure $J_{\nu,k} \deff J_\nu(x_{\nu,k})$, the corresponding 
metrics $g_{\nu,k}$, and vectors $v^+_{\nu,k}, v^0_{\nu,k}, v^-_{\nu,k} \in \rr
^{2n}$ constructed as above for every $u_\nu \ogran_{Z_k}$ with $k=1,\ldots,
l_\nu$, at the position $k=k_\nu$ we obtain the inequality in the opposite 
direction: 
\begin{multline}\eqqno(5.2.12)
\norm{u_\nu -v_{\nu,k_\nu}}^2_{L^{1,2}(Z_{k_\nu})} +
\norm{v_{k_\nu-1} -v_{k_\nu}}^2_{L^{1,2}(Z_{k_\nu})} +
\norm{v_{k_\nu+1} -v_{k_\nu}}^2_{L^{1,2}(Z_{k_\nu})} \ge
\\
\ge \gamma_s \cdot\left(
   \norm{u_\nu -v_{\nu,k_\nu-1}}^2_{L^{1,2}(Z_{k_\nu-1})} +
   \norm{u_\nu -v_{\nu,k_\nu+1}}^2_{L^{1,2}(Z_{k_\nu+1})} \right)  +
\\
+\nu  \cdot \bigl(e^{-2(1+s)k_\nu} \norm{du_\nu}^2_{L^2(Z(0,2))} + 
e^{-2(1+s)(l-k_\nu)}\norm{du_\nu}^2_{L^2(Z(l_\nu-2,l_\nu))}\bigr). \ 
\end{multline}

\smallskip
Let us estimate the behavior of $u_\nu -v_{\nu,k}$ in $Z_k$ for $k\approx k_\nu$.
Set
$$
A_{\nu, k} \deff e^{-k} \norm{du_\nu}_{L^2(Z(0,2))} + 
e^{-(l-k)}\norm{du_\nu}_{L^2(Z(l_\nu-2,l_\nu))}.
$$
Then by \lemma{lem5.2.2} we have $\norm{du_\nu}_{L^2(Z_k)} \le C\cdot A_{\nu,k}$. 
This yields a similar estimate on the diameter: $\diam(u_\nu(Z_k)) \le C \cdot 
A_{\nu,k}$, possibly with a new constant $C$. Further, for a linear complex 
structure $J'$ with the corresponding operator $\dbar' \deff \dbar_{J'}$ we obtain
the pointwise estimate
\begin{multline*}
\left|\dbar' u_\nu \right| = 
\left|\dbar' u_\nu - \dbar_{J_\nu} u_\nu \right| = 
\left| (\d_x u_\nu - J' \cdot \d_y u_\nu) - 
 (\d_x u_\nu - J_\nu(u_\nu) \cdot \d_yu_\nu)\right| 
\\
\le \left|J' - J_\nu \scirc u_\nu\right| \cdot \left|du_\nu \right|.
\qquad
\end{multline*}
For $J_{\nu,k}$ this yields the estimate
\begin{equation}\eqqno(5.2.14)
\norm{\dbar_{J_{\nu,k}} u_\nu}_{L^2(Z(k-2,k+1))} \le 
C \bigl(\diam(u_\nu(Z(k-2,k+1))) \bigr)^s \cdot 
\norm{du_\nu}_{L^2(Z(k-2,k+1))} 
\le C'\cdot A^{1+s}_{\nu,k}.
\end{equation}

By construction, $J_{\nu,k}$ are uniformly bounded. This implies that we can 
represent $u_\nu\ogran_{Z_k}$ in the form $u_\nu\ogran_{Z_k} = w_{\nu,k} + 
f_{\nu,k}$, where $w_{\nu,k}$ is $J_{\nu,k}$-holomorphic and $f_{\nu,k}$ 
is estimated as 
$$
\norm{f_{\nu,k}} _{L^{1,2}(Z_k)} \le C\cdot A^{1+s}_{\nu,k}.
$$

\smallskip
Define the positive $\eta_\nu$ by the relation 
$$
\eta_\nu^2= \norm{u_\nu -v_{\nu,k_\nu}}^2_{L^{1,2}(Z_{k_\nu})} +
\norm{v_{k_\nu-1} -v_{k_\nu}}^2_{L^{1,2}(Z_{k_\nu})}  + 
\norm{v_{k_\nu+1} -v_{k_\nu}}^2_{L^{1,2}(Z_{k_\nu})} 
$$
and set
$$
\ti u_\nu(t, \theta) \deff 
\msmall{1\over \eta_\nu} u_\nu(t+k_\nu, \theta), 
\qquad\qquad 
\ti w_{\nu,k}(t, \theta) \deff 
\msmall{1\over \eta_\nu} w_{\nu,k+k_\nu}(t+k_\nu, \theta),
$$
$$
\ti f_{\nu,k}(t, \theta) 
\deff \msmall{1\over \eta_\nu} f_{\nu,k+k_\nu}(t+k_\nu, \theta),
\qquad
\ti J_{\nu, k} \deff J_{\nu, k+k_\nu},
\qquad 
\ti v^\iota_{\nu,k} \deff 
\msmall{1\over \eta_\nu} v^\iota_{\nu,k+k_\nu},\; \iota=+,0,-,
$$
In other words, we shift all the picture from $Z_{k_\nu}$ to $Z_0$ and
rescale the maps $u_\nu$, the vectors $v^\iota_{\nu,k}$, $\iota=+,0,-$, and so on
in a way as to make the left hand side of \eqqref(5.2.12) equal to $1$. 

It follows from \eqqref(5.2.12) that $A^{1+s}_{\nu, k+k_\nu} \le C \nu^{-\half} 
\eta_\nu =o(\eta_\nu)$ for any fixed $k$. Consequently, $\norm{\ti f_{\nu, k}}
_{L^{1,2} (Z_k)} \lrar 0$ for any fixed $k$ and $\nu\lrar \infty$. This implies 
that the norms $\norm{\ti w_{\nu, k} - \ti v_{\nu, k}}_{L^{1,2}(Z_k)}$ remain 
uniformly bounded in $\nu$ for any fixed $k$. 

Represent every $\ti w_{\nu, k}$ as the Laurent series 
\begin{equation}\eqqno(5.2.14a)
\textstyle
\ti w_{\nu, k}(t, \theta)= \sum_{m=-\infty}^{+\infty}  
e^{m(-t + J_{\nu, k}\theta)} w^m_{\nu, k}
\end{equation}
and denote by $\ti w'_{\nu, k}$ the sum of terms with $m=0, \pm1$, \ie 
$$
\ti w'_{\nu, k}(t, \theta)\deff 
e^{t - J_{\nu, k}\theta} w^{-1}_{\nu, k} +  w^0_{\nu, k} +
e^{-t + J_{\nu, k}\theta} w^1_{\nu, k}.
$$
It follows from the construction of $\ti v_{\nu,k}$ that 
\begin{equation}\eqqno(5.2.15)
\norm{\ti w'_{\nu, k} - \ti v_{\nu,k} }_{L^{1,2}(Z_k)}= 
O(\norm{\ti f_{\nu, k}}_{L^{1,2}(Z_k)}) \lrar 0
\end{equation}
for any fixed $k$. Indeed, $\ti v_{\nu,k}$, considered as a function in $\theta$, 
is a linear combination of a constant and the trigonometric functions $\cos\theta$
and $\sin\theta$. So it is orthogonal to the remaining terms $e^{m(-t + J_{\nu, k}
\theta)} w^m_{\nu, k}$, $|m|\ge2$. Thus $\ti w'_{\nu, k}$ is the best 
approximation of $\ti w_{\nu, k}$ by such linear combinations, whereas the 
difference $\ti w'_{\nu, k} - \ti v_{\nu,k}$ appears as the best approximation of 
$\ti f_{\nu, k}$. 

Since $\ti J_{\nu, 0}$ is bounded uniformly in $\nu$, there exists a subsequence,
still indexed by $\nu$, which converges to a linear complex structure $\ti J$.
It follows from the definition of $\ti J_{\nu, k}$ and the estimate on the
diameter of $u_\nu(Z_k)$ for $k\approx k_\nu$ that for any fixed $k$
the structures $\ti J_{\nu, k}$ also converge to $\ti J$. 

Now we show that, after going to a subsequence, $\ti u_{\nu, 0} - \ti v_{\nu, 0}$ 
converges weakly in the $L^{1,2}(Z(-2,1))$-topology to a $\ti J$-holomorphic 
function, and that this convergence is strong in the $L^{1,2}(Z_0)$-topology. The 
inequality \eqqref(5.2.12) together with the choice of $\eta_\nu$ and the
construction of $\ti u_{\nu, k}$ gives boundedness of the norms $\norm{\ti u_{\nu,
0} - \ti v_{\nu, 0} }_{L^{1,2}(Z(-2,1))}$ uniform in $\nu$. So the weak 
convergence follows. From \eqqref(5.2.14) and $A^{1+s}_{\nu,k} =o(\eta_\nu)$ we 
obtain the vanishing $\norm{\dbar_{\ti J _{\nu, 0}} \ti u_{\nu, 0} }_{L^{1,2} 
 (Z(-2,1))} \lrar0$. The estimate \eqqref(5.2.15) and $\norm{\ti f_{\nu, k}}
_{L^{1,2} (Z_k)} \lrar 0$ yield $\norm{\dbar_{\ti J _{\nu, 0}} \ti v_{\nu, 0} }
_{L^{1,2} (Z(-2,1))} \lrar0$. Now the desired strong $L^{1,2}(Z_0)$-convergence 
follows from elliptic regularity of $\ti J _{\nu, 0}$. In the same way for $k= 
\pm1$ we obtain the weak $L^{1,2}$-convergence of $\ti u_{\nu, k} -\ti v_{\nu, k}$
in $Z_k$.  

For $k=0,\pm1$, let $\ti u_k \deff \lim \ti u_{\nu, k} -\ti v_{\nu, k}$ be the 
limit functions obtained above. Since $\norm{\ti f_{\nu, k}} _{L^{1,2} (Z_k)} 
\lrar 0$, these are $\ti J$-holomorphic functions in $Z_k$, $\ti J = \lim 
\ti J_{\nu,k}$.

Observe that the functions $\ti v_{\nu,\pm1} -\ti v_{\nu,0}$ are linear 
combinations of constants and the function $e^{\pm t}\cos\theta$, $e^{\pm t}\sin
\theta$ which are uniformly bounded in the $L^{1,2}(Z_0)$-norm. Consequently,
after taking a subsequence, we also obtain the strong $L^{1,2}$-convergence
in $Z_0$. This implies the strong $L^{1,2}$-convergence in $Z(-2,1)$. Finally, from
\eqqref(5.2.15) we conclude that the Laurent series for each $\ti u_k$ does not
contain terms of degree $m=0,\pm1$, \ie  a constant term and a multiple of
$e^{\pm(-t + \ti J \theta)}$. This, in turn, implies that, first, the $\ti u_k(t, 
\theta)$ are restrictions to $Z_k$ of the same $\ti J$-holomorphic function 
$\ti u$, and second, $\lim \ti v_{\nu,-1} -\ti v_{\nu,0} = \lim \ti v_{\nu,1} -
\ti v_{\nu,0} =0$.

Substituting into \eqqref(5.2.12), we see that $\ti u$ satisfies the inequality
\begin{equation}\eqqno(5.2.16)
\norm{\ti u}^2_{L^{1,2}(Z_0)} 
\ge \gamma_s \cdot\left(
   \norm{\ti u}^2_{L^{1,2}(Z_{-1})} +
   \norm{\ti u}^2_{L^{1,2}(Z_1)} \right).
\end{equation}
On the other hand, the absence of the terms of degree $m=0,\pm1$ in the Laurent 
decomposition of type \eqqref(5.2.14a) for $\ti u$ implies the inequality
\begin{equation}\eqqno(5.2.17)
\norm{\ti u}^2_{L^{1,2}(Z_0)} 
\le \gamma_2 \cdot\left(
   \norm{\ti u}^2_{L^{1,2}(Z_{-1})} +
   \norm{\ti u}^2_{L^{1,2}(Z_1)} \right).
\end{equation}
with $\gamma_2 \deff {1\over \cosh(4)}$. This inequality is easily obtained for 
mutually orthogonal terms $\ti u^m e^{m(-t+\ti J \theta)}$. However, since $s<1$,
$\gamma_2 < \gamma_s= {1\over \cosh(2 +2s)}$, which is a contradiction.

This implies the validity of \eqqref(5.2.11) for all $k= k_0,\ldots, l-k_0$
with $k_0$ independent of $J$, $l$, and $u$.

\smallskip\noindent {\sl Step 2}. 
We now turn back to the proof of the theorem. To show that \eqqref(5.2.11) 
implies \eqqref(5.2.10), we set for $k=0, \ldots, l$
\begin{multline}
A'_k \deff k\cdot \msmall{\cosh(2+2s) \over \sinh(2+2s)} \cdot C_1 \cdot
e^{-2(1+s)k} \norm{du}^2_{L^2(Z(0,2))}+ 
\\ 
+ (l-k)\cdot \msmall{\cosh(2+2s) \over \sinh(2+2s)} \cdot C_1 \cdot
e^{-2(1+s)(l-k)} \norm{du}^2_{L^2(Z(l-2,l))}.\ 
\end{multline}
Then $A'_k$ satisfies the equality
$$
A'_k = \msmall{\gamma_s \over 2} \cdot (A'_{k-1} + A'_{k+1}) 
+ C_1  \cdot \bigl(e^{-2(1+s)k} \norm{du}^2_{L^2(Z(0,2))} + 
e^{-2(1+s)(l-k)}\norm{du}^2_{L^2(Z(l-2,l))}\bigr). 
$$
Consequently,
\begin{multline}\eqqno(5.2.18)
\norm{u -v_k}^2_{L^{1,2}(Z_k)} - A'_k +
\norm{v_{k-1} -v_k}_{L^{1,2}(Z_k)}^2 +
\norm{v_{k+1} -v_k}_{L^{1,2}(Z_k)}^2 \le
\\
\le \gamma_s \cdot\left(
   \norm{u -v_{k-1}}^2_{L^{1,2}(Z_{k-1})} - A'_{k-1} +
   \norm{u -v_{k+1}}^2_{L^{1,2}(Z_{k+1})} - A'_{k+1})\right).\quad
\end{multline}
As in the proof of \lemma{lem5.2.2}, \eqqref(5.2.18) implies the estimate
\begin{multline}\eqqno(5.2.19)
\norm{u -v_k}^2_{L^{1,2}(Z_k)} - A'_k \le 
\\
\le C_2  \cdot \bigl(e^{-2(1+s)k} \norm{du}^2_{L^2(Z(0,2))} + 
e^{-2(1+s)(l-k)}\norm{du}^2_{L^2(Z(l-2,l))}\bigr)
\end{multline}
for all $k= k_0,\ldots, l-k_0$ with $k_0$ independent of $J$, $l$, and $u$.
Substitution the definition of $A'_k$ yields
\begin{multline}\eqqno(5.2.20)
\norm{u -v_k}^2_{L^{1,2}(Z_k)} +
\norm{v_{k-1} -v_k}_{L^{1,2}(Z_k)}^2 +
\norm{v_{k+1} -v_k}_{L^{1,2}(Z_k)}^2 \le 
\\
\le C_3  \cdot \bigl(k\cdot e^{-2(1+s)k} \norm{du}^2_{L^2(Z(0,2))} + 
(l-k) \cdot e^{-2(1+s)(l-k)}\norm{du}^2_{L^2(Z(l-2,l))}\bigr)
\end{multline}

\smallskip\noindent {\sl Step 3}. For concrete $J$, $l$, and $u$ as in the 
hypotheses of the theorem, find $k^*$ for which the right hand side of 
\eqqref(5.2.20) takes its minimum. Set $J\st \deff J_{k^*} = J(x_{k^*})$, 
and $v^\iota \deff v^\iota_{k^*}$, $\iota=-,0,+$, and $v(t,\theta) \deff v_{k^*}(t,
\theta)$. In view of \eqqref(5.2.20), for the proof of the theorem it is 
sufficient to estimate $\norm{v_k - v}_{L^{1,2}(Z_k)}$.

We do this by descending recursion starting from $k=k^*$. Assume that we have 
shown that 
\begin{equation}\eqqno(5.2.21)
\norm{v_k - v}_{L^{1,2}(Z_k)} \le 
C_4 \, k^{\half}\, e^{-(1+s)k}\,\norm{du}_{L^2(Z(0,2))}
\end{equation}
for all $k=k_1+1,\ldots, k^*$ with the constant $C_4$ to be chosen below. 
By our choice of $k^*$, for $k=1,\ldots, k_1$ we obtain from 
\eqqref(5.2.20)
$$
\norm{v_k - v_{k+1}}_{L^{1,2}(Z_k)} \le 
2\,C_3\, k^{\half}\, e^{-(1+s)k}\,\norm{du}_{L^2(Z(0,2))}.
$$
Observe that for any function $w(t,\theta)$ of the form 
$$
w(t,\theta)= w^0 
+ (e^t w^+_c + e^{-t} w^-_c) \cos(\theta) 
+ ( e^t w^+_s + e^{-t} w^-_s) \sin(\theta) 
$$
with constant vectors $w^0, w^\pm_c, w^\pm_s \in \rr^{2n}$---so are all our 
differences $v_k - v_{k'}$---we have the estimate
$$
\norm{w}_{L^{1,2}(Z_k)} \le e \cdot \norm{w}_{L^{1,2}(Z_{k+1})}.
$$
Applying this, we obtain
\begin{align*}
\norm{v_{k_1} - v}_{L^{1,2}(Z_{k_1})} 
&\le \norm{v_{k_1} - v_{k_1+1}}_{L^{1,2}(Z_{k_1})} + 
\norm{v_{k_1+1} -v }_{L^{1,2}(Z_{k_1})} 
\\
\le& 2\,C_3\, k_1^{\half}\, e^{-(1+s)k_1} \,\norm{du}_{L^2(Z(0,2))}+ 
e \cdot \norm{v_{k_1+1} -v }_{L^{1,2}(Z_{k_1+1})}
\\
\le& 2\,C_3\, k_1^{\half}\, e^{-(1+s)k_1} \,\norm{du}_{L^2(Z(0,2))}+ 
C_4\, (k_1+1)^{\half}\, e^{1-(1+s)(k_1+1)}\,\norm{du}_{L^2(Z(0,2))}
\\
=&2\,C_3\, k_1^{\half}\, e^{-(1+s)k_1} \,\norm{du}_{L^2(Z(0,2))}+ 
e^{-s} \,C_4\, (k_1+1)^{\half}\, e^{-(1+s)k_1} \,\norm{du}_{L^2(Z(0,2))}.
\end{align*}
Assume additionally that $e^{-s/2}(k_1+1)^{\half} \le k_1^{\half}$. Then 
setting $C_4 \deff {2C_3 \over 1- e^{-s/2}}$ we can conclude that \eqqref(5.2.21) 
also holds for $k=k_1$. Since the condition $e^{-s/2}(k_1+1)^{\half} 
\le \cdot k_1^{\half}$ is equivalent to $k_1 \ge { 1 \over e^s-1 }$, our recursive 
construction implies \eqqref(5.2.21) for all $k\in \left[{ 1 \over e^s-1 }, k^*
\right]$. For the remaining $k\in \left[1, { 1 \over e^s-1 }\right]$ the 
estimate \eqqref(5.2.21) follows from \lemma{lem5.2.2}.

Making a similar recursive construction for $k=k^*, \ldots, l$
we obtain the estimate 
\begin{equation}\eqqno(5.2.22)
\norm{v_k - v}_{L^{1,2}(Z_k)} \le 
C_4 \, (l-k)^{\half}\, e^{-(1+s)(l-k)}\,\norm{du}_{L^2(Z(l-2,l))}
\end{equation}
for all $k=k^*,\ldots, l-1$ with the the same constant $C_4$. Now \eqqref(5.2.21),
\eqqref(5.2.22), and \eqqref(5.2.20) imply the desired estimate 
\eqqref(5.2.10). \qed

\newsubsection[5.3]{Deformation of a node and gluing}
The cycle topology on $\barm$, introduced in \refsubsection{5.2c},
has the nice property that $\pr_{\!\!\scrj}: \barm \to \scrj$ is continuous and 
{\sl proper}. The last property is follows from the Gromov compactness for 
closed curves. However, it is desirable to have a better understanding of the 
topological structure of $\barm$. Recall that in \refsubsection{5.2c} we obtained
a natural stratification of $\barm$ in which the strata are distinguished by a 
topological type of curves. Moreover, every stratum $\scrm_\bftau$ has a natural 
structure of a $C^\ell$-smooth Banach manifold such that the restricted
projection $\pr_{\!\!\scrj}: \scrm_\bftau \to \scrj$ is Fredholm. 
So to understand of the topology of $\barm$
means to describe how different strata are attached to each other.
The most important problem is to describe deformations of the standard node. 
Let us formulate the question as follows:

\state Gluing problem. {\sl Let $J_0\in \scrj$ be an almost complex structure and 
$u_0: \scra_0 \to X$ a $J_0$-holomorphic map. Describe possible $J$-holomorphic 
maps $u: Z(0,l) \to X$ with $J\in \scrj$ sufficiently close to $J_0$ and $l\gg0$ 
which are sufficiently close to $u_0$ \wrt the Gromov topology}. In other words, 
we try to reverse the bubbling and construct a single map $u$ of a long cylinder 
$Z(0,l)$ by gluing together the components $u'_0, u''_0: \Delta \to X$ of the map 
$u_0$. 

\smallskip
Moreover, one would like to have a smooth structure on the set of such deformation,
so that the transversality techniques could be applied. This means that one seeks
a family of deformations of a given $u_0:  \scra_0 \to (X, J_0)$ depending
smoothly on the parameter.

\smallskip
As the main result of this paragraph we give a satisfactory solution to
the {\sl Gluing problem}. Let us start with introducing some notation.

\newdefi{def5.3.1} For a fixed sufficiently small $\epsi>0$, let
\begin{equation}
\scra \deff \{\, (z^+, z^-) \in \Delta^2\;:\;
|z^+| \cdot |z^-| < \epsi \,\}
\end{equation}
with the projection
\begin{equation}
\pr_\scra : \scra \to \Delta(\epsi),
\quad \pr_\scra(z^+, z^+) = \lambda(z^+, z^+) \deff z^+ \cdot z^- .
\end{equation}
Further, for $\lambda\in \Delta(\epsi)$ define the analytic sets
\begin{equation}
\scra_\lambda \deff \{\, (z^+, z^-) \in \Delta^2\;:\;
z^+ \cdot z^- = \lambda\,\} = \pr_\scra \inv(\lambda).
\end{equation}
For $\lambda=0$ this is the standard node and for $\lambda \not=0$ a cylinder of 
conformal radius $R= \log{1\over|\lambda|}$. Define 
$$
\scra^\pm_\lambda \deff
\{(z^+,z^-) \in \scra_\lambda: |\lambda| \le z^\pm <1 \}.
$$
Then $\scra^\pm_\lambda$ are subannuli for $\lambda \not=0$ and $\scra^\pm_0$ are
discs $\Delta^\pm$, the irreducible components of $\scra_0$. In any case, $\scra
_\lambda = \scra^+_\lambda \cup \scra^-_\lambda$. 

To describe a Hermitian metric on a complex manifold $X$, it is sufficient to 
indicate only the corresponding K\"ahler form $\omega$. In this case we shall say
that $\omega$ {\sl induces} a metric on $X$ or even that $\omega$ {\sl is} a 
metric on $X$. The author begs the reader's pardon for such informality in the
terminology. The advantage of such notation is that the restriction of a metric on 
a complex 
submanifold is given by the restriction of the corresponding K\"ahler form.
In this notation, the standard metric on the disc $\Delta$ with the coordinate
$z$ is given by the form ${\isl \over 2}dz \wedge d\bar z$. 

We equip $\scra_\lambda$ with the Riemannian metric induced from $\Delta^2$. This
gives the standard metric ${\isl \over 2}dz^\pm \wedge d\bar z^\pm$ on each 
component $\Delta^\pm$ of $\scra_0$ and the {\sl hyperbola metric} 
\begin{equation}\eqqno(5.3.0.1)
{\isl \over 2}\Bigl(1 + {|\lambda|^2 
\over |z^+|^4}\Bigr) dz^+ \wedge d\bar z^+ = 
{\isl \over 2}\Bigl(1 + {|\lambda|^2 
\over |z^-|^4}\Bigr) dz^- \wedge d\bar z^-
\end{equation}
on $\scra_\lambda$ with $\lambda \not=0$. 

Set
\begin{equation*}
\check\scra \deff  \{\, (z^+, z^-) \in \scra\;:\;
z^+ \cdot z^- \not=0 \,\} = 
\sqcup_{\lambda \in \check\Delta(\epsi)} \scra_\lambda 
\end{equation*}
and
\begin{equation*}
V^\pm \deff \{ 1-\epsi < |z^\pm | < 1 \}, \qquad V \deff V^+ \sqcup V^-.
\end{equation*}
For a given $\lambda$ we have the canonical imbedding $V \to \scra_\lambda$,
defined by the coordinate functions $z^\pm$ on $\scra_\lambda$ and on $V^\pm$.
This imbedding defines the restriction map
\begin{equation*}
u \in L^{1,p}(\scra_\lambda, X) \mapsto u\ogran_V \in L^{1,p}( V, X)
\qquad
u \mapsto (u(z^+)\ogran_{V^+}, \;  u(z^-)\ogran_{V^-}).
\end{equation*}
\end{defi}

\newdefi{def5.3.2} For a nodal curve $C$ with smooth boundary $\d C=
\sqcup_i \gamma_i$, $\gamma_i \cong S^1$, let $\scrp(C)$ be the set of
{\sl stable pseudoholomorphic maps between $C$ and $X$},
\begin{equation}
\scrp(C) \deff \{\, (u, J) \in L^{1,p}(C, X) \times \scrj
\;:\; \dbar_J u=0,\; \text{ $u$ is stable}\;\}.
\end{equation}
Equip $\scrp(C)$ with the topology induced from $L^{1,p}(C, X) \times \scrj$.
In particular, $\scrp(V)$ consists of triples $(u^+, u^-, J)$, where
$u^\pm: V^\pm \to X$ is $L^{1,p}$-smooth $J$-holomorphic map.
Denote by $\scrp^*(C)$ the subset of $(u,J) \in \scrp(C)$ for which
$u$ is {\sl non-multiple on the union of compact components of $C$}.
Further, we define
\begin{equation}
\scrp(\scra) \deff \sqcup_{\lambda \in \Delta(\epsi)}
\scrp(\scra_\lambda) ,
\qquad
\scrp(\check\scra) \deff \sqcup_{\lambda \in \check\Delta(\epsi)}
\scrp(\scra_\lambda)
\end{equation}
and equip this spaces with the topology induced by the Gromov convergence in the 
interior of $\scra_\lambda$ and $L^{1,p}$-convergence near boundary. This means 
that $(u_n, J_n, \lambda_n)$ converges to $(u_\infty, J_\infty, \lambda_\infty)$ 
if $(J_n, \lambda_n)$ converges to $(J_\infty,  \lambda_\infty)$ in $\scrj \times
\Delta (\epsi)$, the restrictions $u_n\ogran_V$ converges to $u_\infty \ogran_V$
\wrt $L^{1,p}$-norm, and $u_n$ converges to $u_\infty$ in the sense of {\sl
Definition 5.2.7}. Elements of $\scrp(\scra)$ will be denoted by $(u, J,
\lambda)$. As usual, $\pr_{\!\!\scrj}$ stands for the natural projections from
$\scrp(C)$ or $\scrp(\scra)$ to $\scrj$.
\end{defi}

\newthm{thm5.3.1} The natural map $\pr_V: \scrp(\scra) \to \scrp(V)
\times \Delta(\epsi)$, defined by
\begin{equation}
\pr_V(u, J, \lambda) \deff (u(z^+)\ogran_{V^+}, u(z^-)\ogran_{V^-}, J;
\lambda)
\end{equation}
is an imbedding of a {\sl topological} Banach submanifold. 

Moreover, for every $(u_0,J_0) \in \scrp(\scra_0)$ there exists a neighborhood 
$\scru \subset \scrp(\scra_0)$ of $(u_0,J_0)$, an $\epsi'>0$, and a map
$\Phi: \scru \times \Delta(\epsi') \to \scrp(\scra)$ such that
\begin{itemize}
\item $\Phi$ is a homeomorphism onto the image;
\item for every $\lambda \in \Delta(\epsi')$ the restricted map $\Phi_\lambda 
\deff \Phi\ogran_{\scru \times \{\lambda\}}: \scru \to \scrp(V) \times \Delta(
\epsi)$takes values in $\scrp(\scra_\lambda) \subset \scrp(\scra)$ and is a 
$C^1$-diffeomorphism;
\item the family of maps $\Psi_\lambda \deff \pr_V \scirc \Phi_\lambda: \scru \to 
\scrp(V)$ depends continuously on $\lambda \in \Delta(\epsi')$ \wrt the 
$C^1$-topology. 
\end{itemize} 
\end{thm}

In other words, the theorem statess that $\scrp(\scra) =\sqcup \scrp(\scra_\lambda)$
is a continuous family of $C^1$-submanifolds. Before proving this result, we state 
and prove a corollary which provides a technique which allows one to smooth nodal 
points on pseudoholomorphic curves.

\newthm{thm5.3.1a} Let $C^*$ be a closed connected nodal curve parameterized by
a real surface $S$, $J^*\in \scrj$, $u^*: C^* \to X$ a $J^*$-holomorphic map, and 
$(C^*,u^*,J^*) \in \barm^{Gr}=\barm^{Gr}(S, X, [C^*])$ the corresponding
element of the Gromov compactification of the total moduli space. Assume
that the map $u^*:C^* \to X$ is non-multiple.

Then there exist $(C', u', J') \in \scrm(S, X, [C^*])$ arbitrarily close to
$(C^*,u^*,J^*)$ \wrt the Gromov topology such that $C'$ is a smooth curve.
\end{thm}

The notation used in this theorem  was introduced in {\sl Definitions 
\ref{def5.2.6}} and {\sl\ref{def5.2c.2}}. Note that the condition of 
non-multiplicity of $u^*:C^* \to X$ is equivalent to the absence of ghost and
multiple components. 

\proof Let $\{ z^*_1, \ldots z^*_k\}$ be the set of nodal points of $C^*$.
For every $z^*_i$ fix a neighborhood $V_i$ isomorphic to the standard node.
We may also assume that the sets $V_i$ are pairwise disjoint. 
Let $u^*_i$ denote the restriction of $u^*$ to $V_i$. Applying \refthm{thm5.3.1} 
we can perturb $V_i\cong \scra_0$ to an annulus $V'_i = \scra_{\lambda_i}$ and 
$u^*_i: V_i \to X$ to a $J^*$-holomorphic map $u'_i: V'_i \to X$. If these
perturbations $(V'_i, u'_i)$ are made small enough, then we can adjust the
structure $J^*$ and the map $u^*$ on the remaining part of the curve $C^*$ in a way
yielding the desired $(C', u', J') \in \scrm(S, X, [C^*])$. \qed

Modifying the proof of \refthm{thm5.3.1a} one can also obtain 

\newprop{prop5.3.1b} Let $C^*$ be a closed connected nodal curve parameterized by
a real surface $S$, $J^*\in \scrj$, $u^*: C^* \to X$ a $J^*$-holomorphic map, and 
$(C^*,u^*,J^*) \in \barm^{Gr}=\barm^{Gr}(S, X, [C^*])$ the corresponding
element of the Gromov compactification of the total moduli space. Assume
that the map $u^*:C^* \to X$ is non-multiple. 

Then in a neighborhood of $(C^*,u^*,J^*)$ the space $\barm^{Gr}$ is a topological
Banach manifold and  the natural projection $\pr^{Gr}:\barm^{Gr} \to \barm$ 
a homeomorphism. 
\end{prop}

Since the result of \propo{prop5.3.1b} is not needed for the 
purposes of this paper, we leave it without a proof.

\medskip
The proof of \refthm{thm5.3.1} is divided in the subsequent lemmas. The first two 
are simple but useful technical results.

\newlemma{lem5.3.2} Let $C$ be a nodal curve without closed compact components,
and $E$ a holomorphic vector bundle over $C$. Then any operator $D: L^{1,p}(C,
E) \to L^p_{(0,1)}(C, E)$ of the form $D= \dbar_E + R$ with $R \in L^p$ is
surjective and its kernel $\sfh^0_D(C, E)$ admits a closed complement.
\end{lem}

\proof Imbed $C$ into a compact nodal curve $\wt C$ and extend $E$ to
a holomorphic vector bundle $\wt E$ over $\wt C$. Without loss of generality
we may assume that the Chern numbers $\la c_1(\wt E), \wt C_i\ra$ are
sufficiently large for each component $\wt C_i$ of $\wt C$. Now extend
$R \in L^p(C, \homr(E, E\otimes \Lambda^{(0,1)}))$ to $\wt R \in L^p(\wt C,
\homr(\wt E, \wt E\otimes \Lambda^{(0,1)}))$ and set $\wt D \deff \dbar
_{\wt E} + \wt R$. Adjusting $\wt R$ on the complement $\wt C \bs C$ we may
assume that $\wt D: L^{1,p}(\wt C, \wt E) \to L^p_{(0,1)}(\wt C, \wt E)$
is surjective. The sufficient condition for existence of such an adjustment
is provided by the condition on the Chern numbers of $\wt E$. Since any $\eta
\in L^p_{(0,1)}(C, E)$ extends to $\wt \eta \in L^p_{(0,1)}(\wt C, \wt E)$,
the surjectivity of $\wt D$ implies the surjectivity of $D$.

The existence of a closed complement to the kernel of $D$ is shown in 
\cite{Iv-Sh-2} in the case when $D= \dbar$. This proof applies also in our case 
with only minor adjustments.
\qed

\state Remark. The existence of a closed complement to the kernel of 
$D$ allows us 
to apply the implicit function theorem.

\smallskip 
\newlemma{lem5.3.3}
\sli The space $L^{1,p}(C, X)$ is a smooth Banach
manifold with tangent space
\begin{equation}
T_uL^{1,p}(C, X)=L^{1,p}(C, E_u).
\end{equation}

\slii The space $\scrp^*(C)$ is a $C^\ell$-smooth
submanifold of $L^{1,p}(C, X) \times \scrj$ with tangent space
\begin{equation}
T_{(u,J)}\scrp^*(C)= \{ (v, \dot J)\in L^{1,p}(C, E_u) \times T_J\scrj\,:\,
D_{u,J}v + \dot J \scirc du \scirc J_C =0 \}.
\end{equation}

\sliii $\pr_V: \scrp(\scra_0) \to \scrp(V)$ and $\pr_V: \scrp(\check\scra)
\to \scrp(V) \times \check\Delta(\epsi)$ are $C^\ell$-smooth imbeddings on
Banach submanifolds.
\end{lem}

The definitions of the spaces which are involved here are given in 
\refdefi{def5.2.3a} and \refdefi{def5.3.2}.

\proof Let $C= \cup C_i$ be the decomposition of $C$ into components
and $\{(z'_a, z''_a)\}$ the set of pairs of points on the normalization
$\hat C$ corresponding to the nodal points.

\sli The space $L^{1,p}(C, X)$ is a subset of a smooth Banach
manifold $\prod_i L^{1,p}(C_i, X)$ defined by equations $u(z'_a) =u(z''_a)$.
One checks the transversality condition and computes the tangent space.

\slii One can use the same arguments as in part \slip.

\sliii First we note that \lemma{lem1.2.4} \slii implies the following unique 
continuation property: Any
$J$-holomorphic map $u$, defined on an open set $U$ of a nodal curve $C$ admits
at most one $J$-holomorphic extension to an irreducible component $C'$ of $C$ 
provided $C'\cap U \not =\emptyset$. This shows that the restriction maps $F_0: 
\scrp(\scra_0) \to \scrp(V)$ and $\check F: \scrp(\check\scra) \to \scrp(V)\times
\check\Delta(\epsi)$ are set-theoretically injective. Note that we have
introduced new notation, $F_0$ and $\check F$, for the restrictions of
the map $\pr_V$ to the corresponding definition domains.

Further,  $F_0: \scrp(\scra_0) \to \scrp(V)$ is obviously $C^\ell$-smooth
and the differential $dF_0: T_u\scrp(\scra_0) \to T_u\scrp(V)$ is simply
the restriction map $dF_0: (v, \dot J) \mapsto (v\ogran_V, \dot J)$.
\lemma{lem5.3.2} shows that the restriction map
\begin{equation}
\{\, v \in L^{1,p}(\scra_0, E_u)\;:\; D_{u,J}v=0\, \}\mapsto
\{\, v \in L^{1,p}(V, E_u)\;:\; D_{u,J}v=0\, \}
\eqqno(5.3.2)
\end{equation}
is a closed imbedding and splits, \ie admits a closed complement.

\smallskip
The claim about $\check F: \scrp(\check\scra) \to \scrp(V) \times
\Delta(\epsi)$ is proven in a similar way. Details are left to the reader.
\qed

\smallskip
A crucial point in the proof of \refthm{thm5.3.1} is to find an apriori estimate
for the operator $D_{u,J,\lambda}$ which is uniform as $\lambda\lrar0$. Because of
local nature of the estimate it is sufficient to work with the ball $B\subset 
\rr^{2n} \cong \cc^n$ equipped with the standard complex structure $J\st$. We 
start with introducing a chart for the space $L^{1,p}(\scra, \cc^n) \deff \sqcup
_{\lambda \in \Delta(\epsi)} L^{1,p}(\scra_\lambda, \cc^n)$.

\newdefi{def5.3.3} For a nodal complex curve $C$ and a complex manifold $X$ let
$\scrh(C, X)$ denote the space of holomorphic maps $f: C \to X$ which are
$L^{1,p}$-smooth up to the boundary $\d C$. In the case $X=\cc$ we abbreviate
the notation to $\scrh(C)$.
\end{defi}

\newlemma{lem5.3.3a}%
\footnote{\.The results presented in the lemma have been obtained jointly with
S.~Ivashkovich.} There exist families of homomorphisms $T_\lambda: L^p_{(0,1)} 
 (\scra_\lambda{,} \cc)
\to L^{1,p}(\scra_\lambda, \cc)$ and isomorphisms $\sfl_\lambda: \scrh(\scra_0) 
\to \scrh(\scra_\lambda)$ and $Q_\lambda: L^p_{(0,1)} (\scra_0, \cc) \to 
L^p_{(0,1)} (\scra_\lambda, \allowbreak
\cc)$ with the following properties:

\sli the homomorphisms $T_\lambda$ are right inverses of $\dbar: L^{1,p}(\scra
_\lambda, \cc) \to L^p_{(0,1)} (\scra_\lambda, \cc)$;

\slii the norms of $T_\lambda$, $\sfl_\lambda$, and $Q_\lambda$, as well as  
$\sfl_\lambda\inv$ and $Q_\lambda\inv$, are uniformly bounded;

\sliii the homomorphisms $T_\lambda$, $\sfl_\lambda$, and $Q_\lambda$ depend
smoothly on $\lambda\not=0$. 
\end{lem}

\proof 
By the definition, the hyperbola metric ${\isl \over 2}\left(1 + {|\lambda|^2 \over
|z^+|^4} \right) dz^+ \wedge d\bar z^+$ on $\scra_\lambda$ is the sum of the 
standard flat metrics ${\isl \over 2}dz^+ \wedge d\bar z^+$ and ${\isl \over 2}dz^-
\wedge d\bar z^-$. Note also that the function $\left(1 + {|\lambda|^2 \over |z^+|
^4}\right)$, restricted to the subannulus $\scra^+_\lambda= \{ |\lambda|^\half <
|z^+| <1\} \subset \scra_\lambda$, takes values in the interval $[1, 2]$. This 
implies that in every subannulus $\scra^\pm_\lambda$ the metric ${\isl \over 2}
\left(1 + {|\lambda|^2 \over |z|^4}\right) dz^+ \wedge d\bar z^+$  is equivalent to
the disc metric ${\isl \over 2}dz^\pm \wedge d\bar z^\pm$. In particular, the norm 
$\norm{v}_{L^{1,p}(\scra_\lambda)}$ is equivalent to the norm
\begin{equation}
\left(\int_{|\lambda|^ < |z^+|^2 <1} (|v| + |dv|)^p\;
\msmall{\isl \over 2}dz^+ \wedge d\bar z^+
\,+\,
\int_{|\lambda|^ < |z^-|^2 <1} (|v| + |dv|)^p\;
\msmall{\isl \over 2}dz^- \wedge d\bar z^-
\right)^{1\over p}.
\end{equation}

We start with construction of $T_\lambda$. For the discs $\Delta^\pm$ with the 
coordinates $z^\pm$ respectively we define $\ti T^\pm: L^p_{(0,1)} (\Delta^\pm,
\cc) \to L^{1,p}(\Delta^\pm, \cc)$ to be the Cauchy-Green operators, \ie 
$$
\ti T^+(\phi^+)(z^+) \deff 
\msmall{1\over 2\pi \isl} \int_{\zeta \in \Delta} 
\msmall{d\zeta \wedge \phi^+(\zeta) \over \zeta -z^+}
$$
and similarly for $\ti T^-$. Then we set
$$
T^\pm(\phi^\pm)(z^\pm) \deff \ti T^\pm(\phi^\pm)(z^\pm) - 
\ti T^\pm(\phi^\pm)(0)    
$$
So $T^\pm$ are normalizations of $\ti T^\pm$ respectively to the condition 
$T^\pm(\phi^\pm)(0) =0$. For a form $\phi \in L^p_{(0,1)} (\scra_\lambda, \cc)$ we 
denote by $\phi^\pm (z^\pm)$ its restriction to $\scra^\pm \subset \Delta^\pm$ 
extended by $0$ to the whole discs $\Delta^\pm$, and set  
$$
T_\lambda(\phi) \deff T^+(\phi^+)(z^+) + T^-(\phi^-)(z^-) 
$$
For the special case $\lambda=0$ this construction should be modified follows:
Every form $\phi \in L^p_{(0,1)} (\scra_0, \cc)$ has two components $\phi^\pm
 (z^\pm)$ corresponding to the decomposition 
$$
L^p_{(0,1)} (\scra_0, \cc)= L^p_{(0,1)}  (\Delta^+, \cc) \oplus 
L^p_{(0,1)} (\Delta^-, \cc),
$$
and the operator $T_0$ transforms $\phi^\pm$ into functions $f^\pm(z^\pm)\deff
T^\pm_0(\phi^\pm)(z^\pm)$ which satisfy $f^+(0)= f^-(0) =0$. Then $f\deff (f^+,f^-)
\in L^{1,p}(\scra_0, \cc)$. The desired properties of $T_\lambda$ can be seen 
in a straightforward way.

\medskip
Defining the operators $\sfl_\lambda$, we recall the identification
$$
\scrh(\scra_0) =  \{ (g^+(z^+), g^-(z^-)) \in 
\scrh(\Delta^+) \oplus \scrh(\Delta^-) : g^+(0) =g^-(0) \}.
$$
We set $\sfl_0 = \id: \scrh(\scra_0) \to \scrh(\scra_0)$ and 
$$
\sfl_\lambda: g=(g^+(z^+), g^-(z^-)) \in \scrh(\scra_0) \mapsto 
g^+(z^+) + g^-(z^-) - g(0)
$$
for $\lambda \not=0$. The desired properties of $\sfl_\lambda$ are obvious as well.

Note also that for $\lambda \not=0$ the inverse operator $\sfl_\lambda\inv$
essentially gives the Laurent decomposition of functions $g \in \scrh(\scra
_\lambda)$.

\medskip
The definition of $Q_\lambda$ is more subtle. For $\lambda \not=0 \in 
\Delta(\epsi)$ we set 
\begin{equation}\eqqno(5.3.a1)
\rho_\lambda(r) \deff 
\msmall{ r^2 - {|\lambda|^2 \over r^2} \over 1-|\lambda|^2}\cdot
\end{equation}
Then every $\rho_\lambda$ induces a diffeomorphism of the interval $[|\lambda|,1]$ 
onto  $[-1,1]$, such that $[|\lambda|, |\lambda|^\half]$ and $[|\lambda|^\half,1]$
are mapped onto the intervals $[-1,0]$ and $[0,1]$ respectively. The inverse map 
is given by 
\begin{equation}\eqqno(5.3.a2)
r= R_\lambda(\rho) = \sqrt{
\msmall{ \rho(1-|\lambda|^2) + \sqrt{ \rho^2(1-|\lambda|^2)^2 +4\,|\lambda|^2}
          \over2}} \cdot
\end{equation}
For $\lambda \not=0\in \Delta(\epsi)$ we define the maps $\sigma^\pm_\lambda: 
Z(-1,1) = [-1,1] \times S^1 \to \scra_\lambda$ which are given in the coordinates
$z^\pm= r^\pm \, e^{\isl\theta^\pm}$ by relations $r^\pm = R_\lambda(\rho)$ and 
$\theta^\pm =\theta$ respectively, so that
$$
\sigma^\pm_\lambda: (\rho, \theta)\in Z(-1,1)  \mapsto 
z^\pm= R_\lambda(\rho) e^{\isl\theta} \in \scra_\lambda.
$$
Now we can explain the reason for the choice of \eqqref(5.3.a1), which at first 
glance probably seems rather wild. 
Our choice is made to insure that the pull-back 
by $\sigma_\lambda$ of the natural volume form of $\scra_\lambda$ is a constant 
multiple of the standard volume form on $Z(-1,1)$, \ie
\begin{equation}\eqqno(5.3.a3)
\sigma_\lambda^*\left( \left(1+ {|\lambda|^2 \over (r^+)^4}  \right)
\msmall{\isl\over2} dz^+ \wedge d\bar z^+ \right)
= (1-|\lambda|^2) d\rho \wedge d\theta .
\end{equation}
Moreover, $1-|\lambda|^2$ remains uniformly bounded as $\lambda$ varies
in $\Delta(\epsi)$, so that the volume forms in the right hand side of 
\eqqref(5.3.a3) are equivalent. 

The behavior of $\sigma^\pm_\lambda$ by $\lambda$ close to $0$, which is rather 
delicate, can be described as follows. Denote $Z^+ \deff Z(0,1)$ and $Z^- \deff 
Z(-1,0)$. Then for $\lambda \lrar 0$ we have convergence of $\sigma^+_\lambda$ on 
$Z^+$ to a map $\sigma^+_0: Z^+ \to \scra^+_0 =\Delta^+$, and resp.\ convergence 
of $\sigma^-_\lambda$ on $Z^-$ to $\sigma^-_0: Z^- \to \scra^-_0$, which are given
by 
\begin{align*}
\sigma^+_0: &(\rho, \theta) \mapsto z^+=\sqrt\rho \, e^{\isl\theta},
\\
\sigma^-_0: &(\rho, \theta) \mapsto z^-=\sqrt\rho \, e^{\isl\theta}.
\end{align*}
Observe also that we obtain the map $R_0(\rho) = \sqrt\rho$ in the limit of
\eqqref(5.3.a2) as $\lambda\lrar0$. The convergence $\sigma^\pm_\lambda \lrar
\sigma^\pm_0$ is in the $C^\infty$-sense in the interiors of $Z^\pm$, and in
the $C^0$-topology up to boundary of $Z^\pm$.

On the other hand, there is no convergence of $\sigma^\pm_\lambda$ on $Z^\mp$. 
The topological reason for the absence of the convergence of $\sigma^\pm_\lambda
\ogran_{Z^\mp}$ is that, making with $\lambda$ a small bypass around $0$, we 
perform a {\sl Dehn twist} with $\scra_\lambda$.\footnote{\.The author is
indebted to Bernd Siebert for this remark.}

Now, we define $Q_\lambda$ representing the components $\phi^\pm$ of every $\phi\in
L^p_{(0,1)} (\scra_\lambda, \cc) = L^p_{(0,1)} (\Delta^+, \cc) \oplus L^p_{(0,1)} 
 (\Delta^-, \cc)$ in the form $\phi^\pm(z^\pm) = f^\pm(r^\pm, \theta^\pm) d\bar
z^\pm$ with $L^p$-functions $f^\pm(r^\pm, \allowbreak
\theta^\pm)$ and setting 
\par\noindent
$
Q_\lambda (\phi) \deff \begin{cases}
f^+\Bigl(\! R_\lambda\bigl( (r^+)^2\bigr), \theta^+\Bigr)d\bar z^+ 
& \!\!\!\!\text{at the point $z^+= r^+ e^{\isl\theta^+}$ with 
$r^+ \in [|\lambda|^\half, 1]$};
\\
f^-\Bigl(\! R_\lambda\big( (r^-)^2 \bigr), \theta^-\Bigr)d\bar z^-
& \!\!\!\!\text{at the point $z^-= r^- e^{\isl\theta^-}$ with 
$r^- \in [|\lambda|^\half, 1]$}. 
\end{cases}$ \lineeqqno(5.3.a4)

\smallskip
Let us explain the meaning of the construction of $Q_\lambda$ given in
\eqqref(5.3.a4). The first point is that we essentially transform every form
$\phi \in L^p_{(0,1)}(\scra_\lambda, \cc)$ into a function $f\in L^p(\scra
_\lambda, \cc)$. This is done by representing $\phi$ in the form $\phi(z^+) = f^+
 (z^+) d\bar z^+$ on $\scra^+_\lambda$ and in the form $\phi(z^-) = f^-(z^-) 
d\bar z^-$ on $\scra^-_\lambda$. Observe that we use different coordinates $z^\pm$
on the different parts $\scra^\pm_\lambda$ of $\scra_\lambda$. Independently of
this, we obtain a well-defined $L^p$-function $f$ on $\scra_\lambda$, and hence a
family of well-defined maps $F_\lambda: L^p_{(0,1)}(\scra_\lambda, \cc) \to L^p(
\scra_\lambda, \cc)$. It is not difficult to see that the $F_\lambda$ are complex 
linear isomorphisms and that the operator norms $\norm{F_\lambda}_{\sf{op}}$
and $\norm{F_\lambda\inv}_{\sf{op}}$ are bounded uniformly in $\lambda$.

The second point is the observation that for the definition of a space $L^p(Y)$ 
only the involved measure $\mu$ on $Y$ is essential. In particular, if a measurable
map $g: (Y_1, \mu_1) \to (Y_2, \mu_2)$ induces an equivalence of measures, \ie 
$g^*\mu_2 = e^\rho(y) \mu_1$ for some bounded $\rho \in L^\infty(Y_1, \mu_1)$,
then the induced map $g^*: L^p(Y_2, \mu_2) \to  L^p(Y_1, \mu_1)$ is an isomorphism
of Banach spaces. Thus our construction exploits the fact that the measures
$$
(\sigma^\pm_\lambda)^* \left( \left(1+ {|\lambda|^2 \over (r^\pm)^4}  \right)
\msmall{\isl\over2} dz^\pm \wedge d\bar z^\pm \right)
$$
on $Z^\pm$ are equivalent.

\newcorol{cor5.3.3b} The Banach spaces $L^{1,p}(\scra_0, \cc^n)$ and $L^{1,p}
 (\scra_\lambda, \cc^n)$ with $\lambda \not=0$ are isomorphic.
\end{corol}

\proof One uses $T_\lambda$ to split the exact sequences
$$
0 \lrar \scrh(\scra_\lambda, \cc^n) \lrar L^{1,p}(\scra_\lambda, \cc^n)
\buildrel \dbar \over \lrar L^p_{(0,1)}(\scra_\lambda, \cc^n)
\lrar 0
$$
for $\lambda=0$ and $\lambda \not=0$. Then one applies $\sfl_\lambda$ and 
$Q_\lambda$ to identify $\scrh(\scra_0, \cc^n)$ with $\scrh(\scra_\lambda, \cc^n)$
and respectively $L^p_{(0,1)}(\scra_0, \cc^n)$ with $L^p_{(0,1)}(\scra_\lambda, 
\cc^n)$. \qed

\medskip
Let us explain now the main difficulty in the proof of \refthm{thm5.3.1}. Several 
authors (see \eg \cite{Sie}, \cite{Li-T},  \cite{Ru}) have approached the 
{\sl Gluing problem} 
by making an appropriate local imbedding $\scrp(\scra) \hook L^{1,p} (\scra, 
\cc^n)= \sqcup_{\lambda \in \Delta(\epsi)} L^{1,p}(\scra_\lambda, \cc^n)$
and showing that, roughly speaking, the $\dbar$-equation induces an operator which 
defines $\scrp(\scra)$ and has a certain (rather weak) smoothness property. 
This property is sufficient for showing that for any fixed $J \in \scrj$ the 
compactification $\barr\scrm_J$ has a well-defined ``virtual fundamental class'',
the basic object to define in the theory of Gromov-Witten invariants. 

The difficulty with this approach is that the smooth structure in $L^{1,p}(\scra,
\cc^n)$ given by the isomorphism $L^{1,p}(\scra, \cc^n) \cong L^{1,p}(\scra_0, 
\cc^n) \times \Delta(\epsi)$ depends heavily on both the linear and the complex 
structures in $\cc^n$. Moreover, one can show that for a generic smooth 
diffeomorphism $g: \cc^2 \to \cc^2$ the induced map $u \in L^{1,p}(\scra, \cc^n) 
\mapsto u\scirc g \in L^{1,p}(\scra, \cc^n)$ is not even Lipschitzian at generic 
$u \in L^{1,p}(\scra_0, \cc^n)$. 

To the contrary, our approach of imbedding $\scrp(\scra)$ into $\scrp(V)\times
\Delta(\epsi)$ by means of tracing the restriction $u\ogran_V$ has the
advantage that the smooth structure in $\scrp(V)$ is natural and canonical. 
Here we shall use \lemma{lem5.3.3a} rather in a different way: It 
serves for us as an approximative description of the behavior of the Gromov 
operator $D_{u,J,\lambda}$ as $\lambda \lrar 0$. 

\newdefi{def5.3.4} Let $X$ be a manifold with a fixed symmetric connection
$\nabla$. Then for $\lambda \in \Delta(\epsi)$, $u\in L^{1,p}(\scra_\lambda, X)$, 
and a $C^1$-smooth almost complex structure on $X$ we denote by $J_\lambda$ the
complex structure on $\scra_\lambda$ and by 
$$
D_{u,J,\lambda}: L^{1,p}(\scra_\lambda, E_u) \to 
L^p_{(0,1)}(\scra_\lambda, E_u) 
$$
the operator given by 
\begin{equation}\eqqno(5.3.a5)
D_{u,J,\lambda}(v) \deff \nabla v + J \scirc \nabla v \scirc J_\lambda
+ \msmall{1\over2} \bigl( \nabla_vJ \scirc du \scirc J_\lambda -
J \scirc \nabla_vJ \scirc du \bigr).
\end{equation}
Observe that in the definition of the space $L^p_{(0,1)}(\scra_\lambda, E_u)$ we 
equip $E_u \deff u^*TX$ with the structure $u^*J$. The construction of 
$D_{u,J,\lambda}$ is an extension of the definition of the Gromov operator, because
for $(u, J)\in \scrp(\scra_\lambda)$ the definition \eqqref(5.3.a5) coincides
with the original one in \eqqref(1.3.4).
\end{defi}

\newlemma{lem5.3.4} Let $B\subset \rr^{2n}$ be the ball and $J^*$ a  $C^1$-smooth 
almost complex structure in $B$. Then there exist constants  $\epsi=\epsi(J)>0$ 
and $C =C(J) < \infty$ such that for any 
almost complex structure $J$ with
$$
\norm{J -J^*}_{C^1(B)} \le \epsi,
$$
any $\lambda \in\Delta(\epsi)$, and any $J$-holomorphic map $u: \scra_\lambda \to
B(\half)$ with
$$
\norm{du}_{L^p(\scra_\lambda)} \le \epsi,
$$

\smallskip
\sli 
one has a uniform estimate
\begin{equation}
\norm{v}_{L^{1,p}(\scra_\lambda)} \le C \cdot 
\bigl(\norm{v}_{L^{1,p}(V)} + 
\norm{D_{u,J,\lambda}v}_{L^p(\scra_\lambda)} \bigr)
\eqqno(5.3.4)
\end{equation}
for any $v\in L^{1,p}(\scra_\lambda, E_u)$;

\slii there exists an operator $T_{u,J,\lambda}: L^p_{(0,1)}(\scra_\lambda, E_u)
\to L^{1,p}(\scra_\lambda, E_u)$ with $D_{u,J,\lambda} \scirc T_{u,J,\lambda}= id$;

\sliii moreover, the family of operators $T_{u,J,\lambda}$ depends continuously 
on $(u,J,\lambda)$. 
\end{lem}

\proof Recall that on each half-annulus $\scra^\pm_\lambda$ the hyperbola metric
\eqqref(5.3.0.1) is equivalent to the (flat) disc metric ${\isl\over2} dz^\pm
\wedge \bar z^\pm$. Thus we can apply the Morrey estimate
$$
\diam(u(\scra_\lambda^\pm)) \le 
C \cdot \norm{du}_{L^p(\scra^\pm_\lambda)}
$$
which is uniform in $\lambda$.  This gives
\begin{equation}\eqqno(5.3.a6)
\diam(u(\scra_\lambda)) \le C \cdot \norm{du}_{L^p(\scra_\lambda)} 
\end{equation}
again uniformly in $\lambda$. 

\smallskip
Consequently, $\diam(u(\scra_\lambda))$ is sufficiently small. Let $J\st$ be a 
linear complex structure in $\rr^{2n}$ with coincides with $J$ at some point $x_0=
u(z_0)\in u(\scra_\lambda)$. Then $\norm{J\scirc u -J\st}_{C^0(\scra_\lambda)}$ is 
also small enough. 

The canonical trivialization of the tangent bundle $TB$ yields the canonical
trivialization of $E_u = u^*TB$ and the identification $L^{1,p}(\scra_\lambda, E_u)
=L^{1,p}(\scra_\lambda, \cc^n)$. By \refcorol{cor5.3.3b} we obtain a structure of 
a continuous Banach bundle on the union $\sqcup_{(u,\lambda)\in L^{1,p}(\scra,B)}
L^{1,p}(\scra_\lambda, E_u)$.

A similar identification for $L^p_{(0,1)}(\scra_\lambda, E_u)$ requires some 
modification, since in the definition of this space one
involves the complex structure $J_u \deff u^*J = J\scirc u$. Therefore we must fix 
a {\sl complex} isomorphism $\phi$ of $(TB, J)$ with the trivial complex bundle 
$(TB, J\st)$ over $B$. This means that $\phi$ is an $\rr$-linear endomorphism of 
$TB$ satisfying $\phi \scirc J = J\st \scirc \phi$. Since $J$ coincides with 
$J\st$ at $x_0=u(z_0)$ with $z_0 \in \scra_\lambda$, we may also 
assume that $\norm{\phi - \id_{TB}}_{C^0(u(\scra_\lambda))}$ is small enough. Using
$\phi$, we obtain an isomorphism $\phi_*: L^p_{(0,1)}(\scra_\lambda, E_u)\buildrel 
\cong \over \lrar L^p_{(0,1)} (\scra_\lambda, \cc^n)$. Moreover, the norms 
$\norm{\phi_*}_{\sf{op}}$ and $\norm{\phi_*\inv}_{\sf{op}}$ are bounded uniformly 
in $\lambda \in \Delta(\epsi)$. Now the composition $\phi_* \scirc D_{u,J,\lambda}$
is a homomorphism between $L^{1,p}(\scra_\lambda, \cc^n)$ and $L^p_{(0,1)} 
 (\scra_\lambda, \cc^n)$. 

Using the trivialization $(u^*TB, u^*J\st)$ we define the operator 
$$
\dbar: L^{1,p}(\scra_\lambda, \cc^n) \to L^p_{(0,1)} (\scra_\lambda, \cc^n).
$$
By construction, $\norm{ \dbar - \phi_* \scirc D_{u,J,\lambda}}_{\sf{op}}$
is small enough. Consequently, the restriction of $D_{u,J,\lambda}$ to the image
of the operator $T_\lambda$ constructed in \lemma{lem5.3.3a} is an isomorphism.
Now the existence of $T_{u,J,\lambda}$ follows from the ``closed graph theorem'' 
of Banach. 

Note that the choice of the point $x_0$ and the isomorphism $\phi$ used in the 
construction of $T_{u,J,\lambda}$ can be made continuous on $(u,J,\lambda)$.
This means that there exist families $x_0(u,J,\lambda)$ and $\phi(u,J,\lambda)$
which depend continuously on  $(u,J,\lambda)$ and have the properties needed above.
For example, one can set $x_0(u,J,\lambda) \deff u(1,\lambda)$, the image of
the point $(1,\lambda) \in \scra_\lambda$. For such a choice of $x_0(u,J,\lambda)$
and $\phi(u,J,\lambda)$ the family $T_{u,J,\lambda}$ also depends 
continuously on $(u,J,\lambda)$.

\smallskip
As a consequence of \sliip, it is sufficient to prove \eqqref(5.3.4) 
under the additional condition $D_{u,J,\lambda}v =0$. But then 
\begin{multline*}
\norm{v}_{L^{1,p}(\scra_\lambda)} \le C_2 \cdot 
\bigl(\norm{v}_{L^{1,p}(V)} + 
\norm{\dbar v}_{L^p(\scra_\lambda)} \bigr)=
C_2 \cdot \bigl(\norm{v}_{L^{1,p}(V)} + 
\norm{\dbar v-  \phi_* \scirc D_{u,J,\lambda}v}_{L^p(\scra_\lambda)} \bigr)
\\
\le 
C_2 \cdot \bigl(\norm{v}_{L^{1,p}(V)} 
  + \norm{\dbar -  \phi_* \scirc D_{u,J,\lambda}}_{\sf{op}} 
  \cdot \norm{v}_{L^{1,p}(\scra_\lambda)}\bigr). 
\end{multline*}
This yields \eqqref(5.3.4) provided $C_2 \cdot \norm{\dbar -  \phi_* \scirc 
D_{u,J,\lambda}}_{\sf{op}} \le \half$.

\smallskip
Note once more that all the estimates in the proof are uniform in $\lambda \in
\Delta(\epsi)$. \qed

\smallskip
\newlemma{lem5.3.7} The union $\sqcup_{\lambda \in \Delta(\epsi)} T\scrp(
\scra_\lambda)$ is a continuous locally trivial Banach bundle over 
$\scrp(\scra)$.
\end{lem}

\proof We use the map $\pr_V$ to identify $\sqcup_{\lambda \in \Delta(\epsi)} 
\scrp(\scra_\lambda)$ with its image in $\scrp(V) \times \Delta(\epsi)$. By 
\lemma{lem5.3.4} \sliiip, we can we consider $\sqcup_{\lambda \in \Delta(\epsi)} 
T\scrp(\scra_\lambda)$ as subbundle of the bundle $T\scrp(V) \times \Delta(\epsi)$
over $\scrp(V) \times \Delta(\epsi)$. This defines the canonical topology on 
$\sqcup_{\lambda \in \Delta(\epsi)} T\scrp(\scra_\lambda)$. Moreover, 
\lemma{lem5.3.4} \sliii implies the claim at all $(u,J,\lambda)$ with
$\lambda\not=0$. Hence it remains to show the local triviality of the bundle
in question in a neighborhood of a given $(u_0, J_0, 0) \in \scrp(\scra_0)$.

Consider first the special case when $\norm{du_0}_{L^p(\scra_0)}$ is small 
enough. Under this assumption, \lemma{lem5.3.4} provides existence of 
a continuous family of splittings $L^{1,p}(\scra_\lambda, E_u) = \ker( D_{u,J,
\lambda}) \oplus \im(T_{u,J,\lambda})$ defined in a neighborhood of $(u_0,J_0,0)$ 
in $\sqcup_{\lambda \in \Delta(\epsi)} \scrp(\scra_\lambda)$. Moreover, it follows
from the construction of $T_{u,J,\lambda}$ that the map $\dbar:\im(T_{u,J,\lambda})
\to L^p_{(0,1)}(\scra_\lambda, \cc^n)$ is an isomorphism. By \refcorol{cor5.3.3b},
this implies local triviality of $\im(T_{u,J,\lambda})$. In a similar way one
shows that for a continuous family of isomorphisms $\phi_{u,J,\lambda}: (TB, J) 
\to  (TB, J\st)$ the operators
$$
v \in \ker( D_{u,J, \lambda}) \subset L^{1,p}(\scra_\lambda, E_u) 
\mapsto \bigl(\id - T_\lambda \scirc \dbar\bigr)\phi_{u,J,\lambda} v
\in \scrh(\scra_\lambda, \cc^n)\subset L^{1,p}(\scra_\lambda,  \cc^n)
$$
are isomorphisms. This gives the local triviality of $\ker( D_{u,J, \lambda})$.

Now observe that the bundle
\begin{equation}\eqqno(5.3.f1)
\{ (v, \dot J) \in T_{(u,J)}\scrp(\scra_\lambda):
v = -T_{u,J,\lambda}(\dot J \scirc du \scirc J_\lambda) \}
\end{equation}
is a complement to $\ker( D_{u,J, \lambda})$ in $T_{(u,J)}\scrp(\scra_\lambda)$.
Since the bundle \eqqref(5.3.f1) is isomorphic to the lift of $T\!\!\scrj(B)$ onto
$\scrp(\scra_\lambda)$, it is locally trivial over the whole union $\sqcup_\lambda
\scrp(\scra_\lambda)$. 

\smallskip
Turn back to the general case of the lemma. For a given $(u_0, J_0, 0) \in 
\scrp(\scra_0)$, find a radius $r>0$ such that for the subnode 
$$
\scra'_0 \deff \{ (z^+,z^-) \in \scra_0 : |z^+|<r,\; |z^-|<r\,\}
= \Delta^+(r) \cup \Delta^-(r)
$$
the norm $\norm{du_0}_{L^p(\scra'_0)}$ is small enough. Define 
$$
\textstyle
A^+\deff \{ z^+ \in \Delta^+: {r \over2} <|z^+| <1\,\},
\qquad
A^-\deff \{ z^- \in \Delta^-: {r \over2} <|z^-| <1\,\},
$$ 
$$ 
\scra'_\lambda \deff \{ (z^+,z^-) \in \scra_\lambda: |z^+|<r,\; |z^-|<r\,\}.
$$
and 
$$
V' \deff A^+ \cup A^-,
\qquad
V'' \deff V' \cap \scra'_0.
$$
For $|\lambda|< {r \over2}$ we identify $A^\pm$ with the corresponding subsets in 
$\scra_\lambda$ by means of the coordinates $z^\pm$. Then for $|\lambda|< {r \over
4}$ $A^+$ is disjoint from $A^-$. This gives a family of coverings $\scra_\lambda=
\scra'_\lambda \cup V'$ parameterized by $\lambda$ with $|\lambda|<  {r \over4}$. 
Moreover, we can consider $V'$ and $V''$ as constant (\ie independent of $\lambda$)
complex curves. Now, for any $(u,J,\lambda) \in \scrp(\scra)$ sufficiently
close to $(u_0,J_0,0)$ we obtain the sequence
\begin{equation}\eqqno(5.3.f2)
0 \to T_{(u,J,\lambda)} \scrp(\scra_\lambda) 
\buildrel \alpha_{u,J,\lambda} \over \lrar
T_{(u,J,\lambda)} \scrp(V') 
\oplus T_{(u,J,\lambda)} \scrp(\scra'_\lambda) 
\buildrel \beta_{u,J,\lambda} \over \lrar
T_{(u,J,\lambda)} \scrp(V'') 
\to 0,
\end{equation}
where we set
$$
\alpha_{u,J,\lambda}(v) \deff 
\bigl(v\ogran_{V'}, v\ogran_{\scra'_\lambda} )
\qquad
\beta_{u,J,\lambda}(v, w) 
\deff v\ogran_{V''} - w\ogran_{V''}.
$$
We claim that the sequence \eqqref(5.3.f2) is exact and splits.
Since $\alpha_{u,J,\lambda}$ is obviously injective, it is sufficient
to construct a right inverse to each $\beta_{u,J,\lambda}$, \ie a homomorphism
$$
\gamma_{u,J,\lambda}: T_{(u,J,\lambda)} \scrp(V'') \to
T_{(u,J,\lambda)} \scrp(V') 
\oplus T_{(u,J,\lambda)} \scrp(\scra'_\lambda)
$$
such that $\beta_{u,J,\lambda} \scirc \gamma_{u,J,\lambda}= \id$. Furthermore,
since $\beta_{u,J,\lambda}$ depends continuously on $(u,J,\lambda) \in \scrp(\scra
)$ close to the $(u_0,J_0,0)$, it is sufficient to construct a right inverse only 
for the $\beta_{u_0,J_0,0}$. 

\smallskip
Consider the restriction map between $L^{1,p}(\scra_0)$ and $L^{1,p}(V'')$ induced
by the imbedding $V'' \hook \scra_0$. It is well-known that this map admits a left 
inverse, \ie a map $Q:L^{1,p}(V'') \to L^{1,p}(\scra_0)$ such that $Q(v'')\ogran
_{V''} =v''$ for any $v''\in L^{1,p}(V'')$. Let us use the same notation for the 
restriction of $Q$ onto $T_{(u_0,J_0,0)}\scrp(V'') \subset L^{1,p}(V'')$. Recall
that in \lemma{lem5.3.4} we have constructed the operator $T_{u_0,J_0,0}: L^p_{
 (0,1)}(\scra_0) \to L^{1,p}(\scra_0)$ which is a right inverse to the operator 
$D_{u_0,J_0,0}: L^{1,p}(\scra_0) \to L^p_{(0,1)}(\scra_0)$. Denote by $\Chi
_{\scra'_0}: L^p_{(0,1)}(\scra_0, E_{u_0}) \to L^p_{(0,1)}(\scra_0, E_{u_0})$ the 
operator given the multiplication on the characteristic function of the subset 
$\scra'_0 \subset \scra_0$. Now, we obtain a left inverse $\gamma_{u_0,J_0,0}$ to 
the operator $\beta_{u_0,J_0,0}$ setting 
$$
\gamma_{u_0,J_0,0} (v'') \deff
\biggl(\!\!
\bigl( T_{u_0,J_0,0} \scirc \Chi_{\scra'_0} \scirc D_{u_0,J_0,0}
\scirc  Q(v'') \bigr)\ogran_{V'},
\bigl((T_{u_0,J_0,0} \scirc \Chi_{\scra'_0} \scirc D_{u_0,J_0,0} - \id)
\scirc  Q(v'') \bigr)\ogran_{\scra'_0}
\biggr)
$$
for $v'' \in T_{(u_0,J_0,0)} \scrp(V'')$. Let us check that $\gamma_{u_0,J_0,0}$ 
has the desired properties. For a given $v'' \in T_{(u_0,J_0,0)} \scrp(V'')$,
set $\ti v'' \deff Q(v'')$. Then $D_{u_0,J_0,0}\ti v''$ vanishes on $V''$.
Consequently, $D_{u_0,J_0,0}\ti v''$ is the sum of the forms $\phi, \psi
\in L^p_{(0,1)}(\scra_0, E_{u_0})$ with the supports in $V'\bs V''$ and $\scra'_0
\bs V''$ respectively. Moreover, $\psi = \Chi_{\scra'_0} \scirc D_{u_0,J_0,0} 
\ti v''$. From the relation $D_{u_0,J_0,0}\scirc T_{u_0,J_0,0} =\id$ we obtain
$$
D_{u_0,J_0,0} \bigl( T_{u_0,J_0,0} \scirc \Chi_{\scra'_0} \scirc D_{u_0,J_0,0}
\scirc  Q(v'') \bigr)\ogran_{V'} =
\bigl( \Chi_{\scra'_0} \scirc D_{u_0,J_0,0}
\scirc  Q(v'') \bigr)\ogran_{V'} = 0.
$$
This means that the first component $\bigl( T_{u_0,J_0,0} \scirc \Chi
_{\scra'_0} \scirc D_{u_0,J_0,0} \scirc  Q(v'') \bigr)\ogran_{V'}$ of $\gamma
_{u_0,J_0,0} (v'')$ satisfies $D_{u_0,J_0,0}v=0$. Thus the first component
of $\gamma_{u_0,J_0,0}$ takes values in $T_{(u_0,J_0,0)}\scrp(V')$. In the same
way one can show that the second component of $\gamma_{u_0,J_0,0}$ takes values 
in $T_{(u_0,J_0,0)}\scrp(\scra'_0)$. Finally, it is obvious that the difference
of the components of $\gamma_{u_0,J_0,0} (v'')$ is $v''$. 
The lemma follows.\qed

\statep Proof. of \refthm{thm5.3.1}. Take some $(u_0, J_0) \in \scrp(\scra_0)$.
Using \lemma{lem5.3.3} \sliiip, we identify $\scrp(\scra_0)$ with its image 
$\pr_V\bigl( \scrp(\scra_0) \bigr) \subset \scrp(V)$ and $T_{(u_0, J_0)}\scrp
 (\scra_0)$ with its image  $\pr_V\bigl( T_{(u_0, J_0)}\scrp (\scra_0) \bigr) 
\subset T_{(u_0, J_0)}\scrp(V)$.

We claim that there exists a closed complement to $\pr_V\bigl( T_{(u_0, J_0)}
\scrp (\scra_0)\bigr)$ 
in $T_{(u_0, J_0)}\scrp(V)$. To show this, let us fix a closed nodal complex curve
$C$ and an imbedding $\scra_0 \hook C$. Then there exists the unique open set $\wt
V \subset C$ such that $\wt V \cap \scra_0 =V$ and $\wt V \cup \scra_0 =C$. Extend
the bundle $E_{u_0} = u_0^*\cc^n$ to an $L^{1,p}$-smooth complex bundle $\wt E$ 
over $C$. Also extend the operator $D_{u_0,J_0,0}: L^{1,p}(\scra_0, E_{u_0}) \to 
L^p_{(0,1)} (\scra_0, E_{u_0})$ to an operator $\wt D: L^{1,p}(C, \wt E) \to L^p
_{(0,1)} (C, \wt E)$ of the form $\wt D= \dbar_{\wt E} + \wt R$ where $\dbar_{\wt 
E}$ is the Cauchy-Riemann operator corresponding to some holomorphic structure in 
$\wt E$ and $\wt R$ is a $\cc$-antilinear $L^p$-integrable homomorphism between 
$\wt E$ and $\wt E \otimes \Lambda^{(0,1)}C$, \ie $\wt R \in L^p\bigl(C, \barr
\hom_\cc (\wt E, \wt E \otimes \Lambda^{(0,1)}C) \bigr)$. It is not difficult to 
show that the extensions $\wt E$ and $\wt D$ can be made in such a way that the 
operator $\wt D$ is an isomorphism. In particular, there exists the inverse 
operator $\wt T: L^p_{(0,1)}(C, \wt E) \to L^{1,p}(C, \wt E)$. For an open set 
$U \subset C$ define
$$
\scrh(U) \deff \{ v \in L^{1,p}(U, \wt E) : \wt Dv=0\}.
$$
In particular, we have $\scrh(\scra_0)=\pr_V\bigl( T_{(u_0,J_0)}\scrp(\scra_0)
\bigr)$ and
similar identifications for $V$ and $V'$. Define the operator $\beta: \scrh(\wt V)
\oplus \scrh(\scra_0) \to \scrh(V)$ setting $\beta(v,w) \deff v\ogran_V -
w\ogran_V$. Then there exists an operator $\gamma: \scrh(V)\to \scrh(\wt V)
\oplus \scrh(\scra_0)$ which is left inverse to $\beta$. Indeed, 
the construction of the operator $\gamma_{u_0,J_0,0}$ made in the proof of
\lemma{lem5.3.7} can be applied with appropriate modifications. In particular,
one uses $\wt T$ instead of $T_{u_0,J_0,0}$. 

Observe that the kernel $\ker(\beta)$ can be identified in a natural way with 
the kernel of $\wt D: L^{1,p}(C, \wt E) \to L^p_{(0,1)} (C, \wt E)$. This implies
that $\beta$ is an isomorphism and $\gamma$ its inverse. Consequently, every
$v\in \scrh(V)$ can be uniquely represented in the form $v=v_1 + v_2$ such that 
$v_1$ extends to $\ti v_1 \in \scrh(\wt V)$ and $v_1$ extends to $\ti v_2 \in 
\scrh(\scra_0)$. The set of all $v_1\in \scrh(V)= T_{u_0,J_0}\scrp(V)$ obtained in
this way forms the desired closed compliment to $\pr_V\bigl( T_{(u_0, J_0)}\scrp 
 (\scra_0) \bigr)$ in $T_{(u_0, J_0)}\scrp(V)$.

\medskip
The existence of such a complement implies that there exists a small ball $\scru 
\subset \scrp(\scra_0)$ centered at $(u_0, J_0)$ and an open imbedding $F: \scru 
\times \scrb \hook \scrp(V)$ with the following properties:
\begin{itemize}
\item $\scrb$ is a small ball in a closed complement of $\pr_V \bigl(T_{(u_0, 
J_0)}\scrp(\scra_0) \bigr)$ in $T_{\pr_V(u_0, J_0)} \scrp(V)$;
\item the map $F$ is a $C^1$-diffeomorphism onto image;
\item the restricted map $F\ogran_{\scru \times \{0\}}: \scru \to \pr_V(\scru)
\subset \pr_V(\scrp(\scra_0)) \subset \scrp(V)$ coincides with $\pr_V$.
\end{itemize} 
In other words, $\scru \times \scrb$ appears as a chart for $\scrp(V)$ in
which $\pr_V(\scru)= F(\scru \times \{0\})$ is a linear subspace.

Now consider the natural projection $\pi_\scru: \scru \times \scrb \to \scru$
restricted to $F\inv\bigl(\pr_V(\scrp(\scra_\lambda))\bigr)$. By 
\lemma{lem5.3.7} this is a $C^1$-diffeomorphism for all sufficiently small 
$\lambda$, provided the ball $\scru$ is small enough. Denote by 
$$
\Phi_\lambda: \scru \to \pr_V\inv\bigl(\pr_V(\scrp(\scra_\lambda)) \cap 
 F(\scru \times \scrb)  \bigr) \subset \scrp(\scra_\lambda)
$$ 
the inverse map. We claim that the family $\Phi_\lambda$, parameterized by $\lambda
\in\Delta(\epsi')$ with a sufficiently small $\epsi'>0$, has the properties
stated in \refthm{thm5.3.1}. In fact, it remains to check only the fact that
the whole map 
$$
\Phi: \scru \times \Delta(\epsi') \to \sqcup_{\lambda
\in\Delta(\epsi')} \scrp(\scra_\lambda)
$$
is continuous. However, this follows from the construction of $\Phi$. \qed

\newsection[6]{Symplectic isotopy problem in $\cp^2$}

\newsubsection[6.1]{Symplectic isotopy problem}
Let $\Sigma, \; \Sigma'$ be two (connected) symplectically imbedded surfaces in a 
symplectic 4-fold $(X, \omega)$. Assume that they have the same homology class. 
Then they have the same genus, see \lemma{lem1.1.2}. Thus one can ask whether or 
not there exists an isotopy $\{ \Sigma_t \}_{t \in [0, 1]}$ from $\Sigma$ to 
$\Sigma'$ such that all $\Sigma_t$ are also symplectically imbedded.
This is refered to as the symplectic isotopy problem.

The example of Fintushel and Stern \cite{Fi-St} shows that there is no hope 
to obtain a results of this type in the case when $\la c_1(X), [\Sigma] \ra \le 0$.
Namely, they proved that under certain conditions on a symplectic 4-fold $(X, 
\omega)$ there exists an infinite 
collection of symplectic imbeddings $\Sigma_i \hook X$, such that $\Sigma_i$ 
represent the same homology class $[C] \in \sfh_2(X, \zz)$ but are pairwise 
non-isotopic, even smoothly. Moreover, the class of symplectic 4-folds with these 
conditions is sufficiently wide, so that one has enough examples of this type.

\smallskip
On the other hand, \refthm{thm4.5.1} hints that a satisfactory solution for the 
symplectic isotopy problem in the case $\la c_1(X), [C] \ra \ge 1$ is possible. 
We state the problem in a more precise form.

\state Conjecture 1. {\sl (Symplectic isotopy problem). \it Let $(X,\omega)$ be 
a compact symplectic $4$-dimensional manifold and $[C]\in \sfh_2(X,\zz)$ a 
homology class with $\la c_1(X), [C] \ra \ge1$. Then every two symplectically {\sl
 immersed} surfaces $\Sigma$ and $\Sigma'$ in the class $[C]$ are symplectically 
isotopic provided they have the same genus $g$ and the only singularities are 
{\sl positive} nodal points.}

Recall, there exists a complete classification of compact symplectic 4-folds
$X$ which come in question. Namely, {\bf Corollary 1.5} in \cite{McD-Sa-3}, 
claims

\newprop{prop6.1.0a}
Let $X$ be a symplectic manifold and $\Sigma
\subset X$ a symplectically imbedded surface which is not an exceptional
sphere. Then $X$ is the blow-up of a rational or ruled manifold.
\end{prop}

The complete description of possible symplectic structures on such $X$ was done in 
\cite{McD-4}, \cite{La-McD}, and \cite{McD-Sa-3}, see also \cite{Li-Liu}, 
\cite{Liu}.

As the main result of this paper we give a positive solution of the symplectic
isotopy problem for imbeddings of low degree in $\cp^2$.

\newthm{thm6.1.1} Any two symplectically imbedded surfaces $\Sigma,\; \Sigma'
\subset \cp^2$ of the same degree $d\le6$ are symplectically isotopic.
\end{thm}

The case $d=1$ and $2$ of the theorem has been proven by Gromov in \cite{Gro},
the case $d=3$ by Sikorav \cite{Sk-3}.

\smallskip
In this connection a result of S.~Finashin about (non-symplectic) isotopy problem 
in $\cp^2$ should be mentioned. He proves in \cite{Fin} that for any
even degree $d=2k\ge6$ there exist infinitely
many isotopy classes of imbedded real surfaces in $\cp^2$ having the degree $d$
and the genus $g$ given by the genus formula, \ie $g={(d-1)(d-2) \over2}\cdot\;$
Note that \refthm{thm6.1.1} claims that for $d=6$ only one of these isotopy 
classes is realizable by a {\sl symplectic} imbedding.

\medskip
Let us explain main ideas of the proof of \refthm{thm6.1.1}. First we observe
that existence of a symplectic isotopy $\{ \Sigma_t \}_{t \in [0,1]}$ between
surfaces $\Sigma,\; \Sigma'$ in a symplectic manifold $(X,\omega)$ implies  
existence of an ``accompanying'' homotopy $\{ J_t \}_{t \in [0, 1]}$ of tame 
almost complex structures, such that the imbeddings $\Sigma_t \hook X$ are 
$J_t$-holomorphic. Conversely a homotopy of $\omega$-tame $J_t$-holomorphic 
imbeddings is necessarily a symplectic isotopy. So given $\Sigma _0$ and 
$\Sigma_1$, the natural thing to do is to outfit them with compatible structures 
$J_0$ and $J_1$, take a generic curve $J_t$ and attempt to find appropriate 
liftings $\Sigma _t$. We do this using the following theorem of Harris \cite{Ha}
for an intermediate construction.

\newprop{prop6.1.2} Any two irreducible nodal algebraic curves $C_0$ and $C_1$ in 
$\cp^2$ of the same degree $d$ and the same geometric genus $g$ are {\sl 
holomorphically isotopic}, \ie can be connected by an isotopy $\{C_t\}_{t\in [0,
1]}$ consisting of nodal algebraic curves.
\end{prop}

By this result, in order to construct the symplectic isotopy, it is enough
to construct a lifting as above for the case where $J_1$ is the standard
integrable structure on $\cp^2$ and $\Sigma_1$ is some smooth algebraic curve.

\medskip
Obviously, \refthm{thm4.5.1} would imply existence of symplectic isotopy if we 
could show that for a generic path $\{J_t \}_{t \in [0,1]}$ the moduli space $\scrm
_{J_t}$ is non-empty. An obstruction to this is the fact that the projection $\pi
_\scrj:\scrm \to \scrj$ is not proper. This means that we must understand the 
structure of the total moduli space $\scrm$ ``at infinity''. In 
\refsubsection{5.2c} we have constructed a completion $\barm$ of $\scrm$ and 
eqquiped it with a natural stratification such that every stratum is a smooth
Banach manifold. In particular, the transversality technique developed in 
\refsection{2} can be applied to every such stratum. 

The next idea in the proof of \refthm{thm6.1.1} is to construct a path $\ti\gamma
_t \deff(C_t, J_t) \in \barr\scrm$ which goes piecewisely along some strata and 
which can be ``pushed'' into the ``main stratum'' $\scrm$ yielding the desired 
isotopy $(\Sigma_t, J_t)$. The main difficulty in realization this idea is to 
ensure that pushing $\ti\gamma_t$ into $\scrm$ we still remain in the same 
connected component of $\scrm$ so that the symplectic isotopy class is preserved.
This means that we are interested in describing possible connected components of 
$\scrm$ in a neighborhood of a given curve $(C^*, J^*) \in \barr\scrm$. Moreover, 
the positive solution of a symplectic isotopy problem would follow immediately from
the fact that locally exactly one such component exists. Indeed, it would be then
sufficient to construct {\sl any} path $\ti\gamma_t\deff(C_t, J_t) \in \barr\scrm$
connecting $\Sigma$ and $\Sigma'$. But existence of such a path follows easily
from \refthm{thm4.5.1} in the case $c_1(X)[C]>0$. 

\smallskip
The result of \refthm{thm6.1.1} is obtained via the proof of the local uniqueness
of such a component of $\scrm$ near a given $(C^*, J^*) \in \barr\scrm$ in the
special case when $C^*$ contains no multiple components. The restriction $d\le6$
in the theorem comes from the fact that in this case it is possible to avoid the 
appearance of multiple components in $C^*$. We are able to do so by demanding that 
the pseudoholomorphic curves $\Sigma_t$ in the isotopy path pass through fixed 
generic $3d-1$ points on $X=\cp^2$. Note that the number $3d-1$ is the maximal 
possible in \refthm{thm4.5.3}.

\newsubsection[6.2]{Local symplectic isotopy problem} As we have explained
in the previous paragraph, we are interested in the possible symplectic isotopy
classes of pseudoholomorphic curves $C$ in a neighborhood of a given singular curve
$C^*$ with no multiple components. The main difficulty in this case is, of course,
to understand the local behavior of curves $C$ near singular points of $C^*$. 
In this way we come to the following question.

\state The Local Symplectic Isotopy Problem. 
{\sl Let $B$ be the unit ball in $\rr^4$ equipped with the standard symplectic 
structure $\omega\st$, $J^*$ an $\omega\st$-tame almost complex structure, and 
$C^*\subset B$ a connected $J^*$-holomorphic curve in $B$ with a unique isolated 
singularity at $0\in B$ and without multiple components. Describe the possible 
symplectic isotopy classes of curves $C$ in $B$ which lie sufficiently close to 
$C^*$ \wrt the cycle topology and which have prescribed singularities, \eg
prescribed number of nodes and ordinary cusps.}

\medskip
We start with a construction of certain symplectic isotopy classes of nodal
pseudoholomorphic curves. For $C^*$ as above, let $C^*= \cup_i C^*_i$ be the 
decomposition
into irreducible components. Then there exist $J^*$-holomorphic parameterizations 
$u^*_i : S_i \to B$, $u^*_i(S_i) = C^*_i$. Shrinking $C^*_i$, if needed, we may 
assume that all $S_i$ are compact and smooth boundaries $\d S_i$, each 
consisting of finitely many circles. Note that the images of the boundary
circles are imbedded in $B$ and mutually disjoint. Further, we can also suppose 
that $u^*_i$ are $L^{1,p}$-smooth up to boundaries $\d S_i$. Set $S\deff \sqcup 
S_i$ and define $u^*: S \to B$ by $u^*\ogran_{S_i} \deff u^*_i$. Denote by $J^*
_S$ the complex structure on $S$ induced by $u^*: S \to B$ from $C^*$.

\newlemma{lem6.2.1} \sli The set $\scrp(S,B)_\nod$ of those $(u, J_S, J) \in 
\scrp(S,B)$ for which the map $u:S \to B$ is an immersion and the singularities 
of the image $C \deff u(S)$ are only nodal points is open and dense in $\scrp(S,B)$
and is connected;

\slii For $(u', J'_S, J')$ and $(u'', J''_S, J'') \in \scrp(S,B)_\nod$, 
sufficiently close to $(u^*, J^*_S, J^*) \in \scrp(S,B)$, the pseudoholomorphic 
curves $C' \deff u'(S)$ and $C' \deff u'(S)$ are symplectically isotopic;

\sliii For a fixed $J_S \in \scrj_S$ and $J\in \scrj(B)$, the subspace of nodal 
curves in each of the spaces $\scrp(S; B, J)$, $\scrp(S,J_S; B)$, and $\scrp(S,J_S; 
B, J)$ is open and dense in thecorresponding space.  
\end{lem}

\proof By results of \refsection{4}, the complement to $\scrp(S,B)_\nod$ in 
$\scrp(S,B)$ is closed and consists of submanifolds of real codimension at least
2. This shows \sli and implies \slii. Part \sliii is obtained similarly. 
\qed

\newdefi{def6.2.1} In the situation of \lemma{lem6.2.1}, we call $C_\nod=u(S)$ a 
{\sl maximal nodal deformation of $C^*$} and the number $\delta$ of nodes on $C
_\nod$ the {\sl nodal number of $C^*$} at the singular point $0\in C^*$. In other 
words, a maximal nodal deformation is a nodal pseudoholomorphic curve obtained from
$C^*= u^*(S)$ by a (sufficiently small) generic deformation of the parameterization
map $u^*: S \to X$, $C^*= u^*(S)$. 

Further, a {\sl canonical smoothing of $C^*$} is a $J^*$-holomorphic curve 
$C^\dag$ obtained from a maximal nodal nodal deformation $C_\nod$ by smoothing of 
all nodes. We use the notion of canonical smoothing for both the construction and 
the resulting curve. Further, we shall always assume that a canonical smoothing 
$C^\dag$ is sufficiently close to $C^*$ \wrt the cycle topology. 
\end{defi}

It follows immediately from \lemma{lem6.2.1} that the symplectic isotopy class
of a canonical smoothing of $C^*$ is well-defined.

\newprop{prop6.2.1a} Any two curves $C^\dag_1$ and $C^\dag_2$ obtained from $C^*$
by the construction of canonical smoothing are symplectically isotopic. Moreover,
such an isotopy can be be carried out sufficiently close to the identity map.
\end{prop}

\smallskip
Note that the number $\delta(C_\nod)$ of nodes on a maximal nodal deformation
$C_\nod$ of $C^*$ equals the nodal number $\delta(0, C^*)$ of $C^*$ at $0$. Observe
also that one can smooth {\sl some} number of nodes on $C_\nod$ producing further 
symplectic isotopy classes. It is easy to show that these new classes are 
determined by the set of the nodes on $C_\nod$ which are smoothed. We conjecture 
that these are all possible symplectic isotopy classes of nodal curves in a 
neighborhood of $C^*$ \wrt the cycle topology. 

\state Conjecture 2. {\sl (Local symplectic isotopy problem for nodal curves).
\it Let $J^*$ be a $C^2$-smooth $\omega\st$-tame almost
complex structure in $B \subset \rr^4$ and $C^* \subset B$ a $J^*$-holomorphic 
curve with a unique isolated singular point at $0\in B$ and without multiple 
components. Assume that $J$ is an almost complex structure in $B$ which is $C^{0,
\alpha}$-smooth for $\alpha>0$ and sufficiently close to $J^*$ \wrt the $C^{0,
\alpha}$-topology.

Then any {\it nodal} $J$-holomorphic curve $C$ sufficiently close to $C^*$ \wrt 
the cycle topology is symplectically isotopic to a $J^*$-holomorphic curve 
obtained from a maximal nodal deformation $C_\nod$ of $C^*$ by smoothing some 
number of nodes on $C_\nod$.}

\medskip
We give a proof the conjecture for the case of {\sl imbedded} curves. Observe 
that here we have only one candidate, namely the canonical smoothing.

\newthm{thm6.2.2} In the situation described in {\sl Conjecture 2}, let $C^\dag$ 
be $J^*$-holomorphic curve obtained by the canonical smoothing of $C^*$.

Let $J$ be an almost complex structure on $B$ sufficiently close to $J^*$ \wrt the
$C^{0,\alpha}$-topology and $C$ an imbedded $J$-holomorphic $C$ sufficiently close
to $C^*$ \wrt the cycle topology. Then there exist a homotopy $J_t$ which is
$C^0$-sufficiently close to $J^*$ and connects $J^*$ with $J$, and an isotopy
$C_t$ of $J_t$-holomorphic curves which connects $C^\dag$ with $C$ and is 
sufficiently close to $C^*$ \wrt the cycle topology.
\end{thm}

\smallskip
The proof will be given after some preparatory results. We shall always assume
that the hypotheses of the theorem are fulfilled. Denote by $S$ the real surface
parameterizing $C$. In other words $S$ is the curve $C$, considered as 
real oriented surface without complex structure.

\smallskip
Our first observation is that the theorem holds in the case when $C^*$ and the
approximating curve $C$ are holomorphic in the usual sense. The result is 
well-known, see \eg \cite{Mil}. Its proof is based on the main advantage of the 
holomorphic case: the fact that one can represent a holomorphic curve as the zero 
divisor of a holomorphic function. 

\newlemma{lem6.2.2a} Let $f^*$ be a holomorphic function in the ball $B$ in 
$\cc^2$ whose zero divisor is a holomorphic curve $C^*$ with a single singular 
point at $0 \in B$ and without multiple components. Assume that $f^*$ and $C^*$
are sufficiently smooth also at the boundary $\d B$. Then
\begin{itemize}
\item[\slip] a canonical smoothing $C$ is obtained as the zero divisor of
a sufficiently small perturbation $f$ of $f^*$;
\item[\sliip] for two generic sufficiently small perturbations $f_1$ and $f_2$
of $f^*$ their zero divisors $C_0$ and $C_1$ are non-singular and {\sl 
holomorphically isotopic}, \ie can be connected by a homotopy consisting
of holomorphic non-singular curves $C_t$.
\end{itemize}
\end{lem}

\smallskip
Denote by $\delta^*$ the nodal number of $C^*$ at $0\in C^*$. We may assume 
inductively that the claim of \refthm{thm6.2.2} holds for all curves $C'$ which 
satisfy the hypotheses of the theorem but have the nodal number $\delta(C')$ at
$0\in C'$ which is strictly less than $\delta^*$. Further, we assume that $\delta^*
\ge2$, since otherwise $\delta^*=1$ and $0\in C^*$ is a nodal point, 
the case covered by \refsubsection{5.3}.

\medskip
Recall that by the theorem of Micallef and White (see \lemma{lem1.2.1}) in a 
neighborhood of $0\in B$ there exists a $C^1$-diffeomorphism $\phi$ of $B\subset 
\rr^4$ such that $\phi(0)=0$, $\phi_*(J^*(0))=J\st$, the standard complex structure
in $\rr^4=\cc^2$, and such that $\phi(C^*)$ is a $J\st$-holomorphic curve. 
Obviously, we may also assume that $d\phi : T_0B \to T_0B$ is the identity map. 
This means that the form $\phi_*\omega\st$ coincides with $\omega\st$ at $0\in B$,
$\phi_* (\omega\st)\ogran_{T_0B}=\omega\st$, and similarly $\phi_*(J^*(0))=J^*(0)
=J\st$. Consequently, $\phi_*(J^*)$ is $\omega\st$-tame in a sufficiently small 
ball $B(r)$, $r\ll1$. Let us fix such a radius $r$. 

Moreover, since $C^*\subset B$ is imbedded outside $0$, we can additionally assume
that $\phi$ is {\sl smooth} outside $0\in B$.

Below, we translate the original situation by means of such $\phi$ and work with 
a {\sl holomorphic} curve $\phi(C^*)\cap B(r)$. This leads to the difficulty that 
$\phi_*(J^*)$ is apriori only {\sl continuous} at $0\in B(r)$. This requires an
additional control on the behavior of pairs $(C,J) \in \scrp(B)$ approximating
$(C^*, J^*)$.

\newlemma{lem6.2.3} Let $(u_n, J_{S,n}, J_n) \in \scrp(S, B)$ be a sequence such 
that $J_n$ converges to $J^*$ in the $C^{0,\alpha}$-topology with $0<\alpha <1$, 
and $C_n \deff u_n(S)$ converges to $C^*$ \wrt the cycle topology and \wrt the
$L^{1,p}$-topology near boundary $\d C_n = u_n(\d S)$. Further, let $\phi: B \to
B$ be the diffeomorphism introduced above. Then for all sufficiently big $n$

\sli $u_n: S \to B$ is an imbedding;

\slii there exists a sequence $J^*_n$ of $C^\ell$-smooth almost complex structures 
in $B$ such that
\begin{itemize}
\item $u_n$ are $(J_{S,n}, J^*_n)$-holomorphic;
\item $\phi_*(J^*_n)$ converges to $J\st$ in the $C^0$-topology in $B(r)$ and
in the $C^{0, \alpha}$-topology outside $0\in B(r)$.
\end{itemize}
\end{lem}

\proof The first part follows from \lemma{lem1.2.2a}, applied to a smaller
ball $B(\rho)$, $\rho<1$, and curves $C^* \cap B(\rho)$, $C_n \cap B(\rho)$.

\smallskip
Define $J^\sharp$ as the pull-back of $J\st$ \wrt $\phi$, $J^\sharp \deff \phi^*(
J\st)$. Then the second part is equivalent to the convergence $J^*_n \lrar 
J^\sharp$ in the appropriate topology. 

Fix some sufficiently small $\epsi>0$. Since $J^*(0)=J\st=J^\sharp(0)$, there 
exists a positive radius $\rho\ll r$ such that $\norm{J^* -J^\sharp}_{C^0(B(\rho))} 
< \epsi$. This implies that $\norm{J_n -J^\sharp}_{C^0(B(\rho))} < \epsi$ for all 
sufficiently big $n$.

Now observe that in $B\bs B(\rho)$ we have the $C^{1,\alpha}$ convergence $C_n 
\lrar C^*$. In particular, in $B\bs B(\rho)$ we have $C^{0,\alpha}$-convergence of 
tangent bundles $TC_n \lrar TC^*$. This implies that for $n \gg1$
we can extend every $J_n$ from $B(\rho)$ to $B(r)$ as a $C^\ell$-smooth structure 
$J^*_n$ which is defined {\sl along} $C_n$ and obeys the estimate
\begin{subequations}
\eqqno(6.2.1)
\begin{align}
&\norm{J^*_n -J^\sharp}_{C^0(C_n\cap B(r))} < \epsi,
\\
&\norm{J^*_n -J^\sharp}_{C^{0,\alpha}(C_n \cap ( B(r) \bs B(2\rho) ))} 
< \epsi.
\end{align}
\end{subequations}
Finally we extend the constructed $J^*_n$ from $C_n \cup B(\rho)$ to the whole
ball $B$ preserving the estimates \eqqref(6.2.1). \qed

\state Remark. In fact, below we shall merely make use of the weaker $C^0
$-convergence $\phi_*(J^*_n) \to J\st$. The H\"older $C^{0,\alpha} $-convergence 
$J_n \to J^*$ was used only to provide the $C^0$-convergence of 
tangent bundles $TC_n \to TC^*$ outside $0\in C^*$. In particular, it would be 
sufficient to have only $C^0$-convergence $J_n \to J^*$ in $B$ and the $C^{0,
\alpha}$-convergence outside $0\in B$. On the other hand, in the case when the 
convergence $J_n \to J^*$ is better, say in the $C^\ell$-topology with 
non-integer $\ell >1$, we could achieve just as as well the $C^\ell$-convergence in 
$B(r)$ outside $0$.

\medskip
\lemma{lem6.2.3} insures that we can reduce the problem to the case when $C^*$ is
holomorphic in the usual sense, \ie \wrt the structure $J\st$. Further, observe
that for the proof of \refthm{thm6.2.2} it is sufficient to show that for any 
sequence $(u_n, J_{S,n}, J_n)$ satisfying the hypotheses of \lemma{lem6.2.3} the 
curves $C_n \deff u_n(S)$ are symplectically isotopic to $C^\dag$ for $n\gg1$. An 
equivalent problem is to show that $\phi(C_n)$ are symplectically isotopic to 
$\phi(C^\dag)$ in $B(r)$. Thus we can replace our initial objects by their 
$\phi$-images in $B(r)$. For the sake of simplicity we maintain the original 
notations for these new objects, \eg $B$ for $B(r)$, $C^*$ and $C_n$ for 
respectively $\phi(C^*) \cap B(r)$ and $\phi(C_n) \cap B(r)$, $J^*$ and $J_n$ for 
respectively $\phi_*(J^*)\ogran_{B(r)}$ and $\phi_* (J_n)\ogran_{B(r)}$, and so 
on. On the other hand, $J\st$ and $\omega\st$ remain the standard structures in 
$B$. Observe that now we have the weaker $C^0$-convergence $J_n \lrar J^*$. 

\smallskip
Imbed $B$ in $\cp^2$ in the standard way so that $J\st$ becomes the standard
integrable structure, still denoted by $J\st$. Then we can extend $\omega$ to a 
global symplectic form on $\cp^2$ taming $J\st$. We maintain the notation $\omega$
for this extension.

We claim that $C^*$ also extends to $\cp^2$ as a compact closed pseudoholomorphic 
curve. Moreover, we claim that there exists an extension $\ti C^*$ with the 
following properties
\begin{itemize}
\item all irreducible components of $\ti C^*$ are rational, \ie 
parameterized by the sphere $S^2$;
\item except for the original singularity at $0\in \ti C^*$, all new singularities 
are only nodal points.
\end{itemize}

Indeed, every irreducible component of $C^*\subset B$ is
parameterized by a holomorphic map $u_i=u_i(z): \Delta \to B$ with $u_i(0)=0$. For
every $u_i(z)$ we take the Taylor polynomials $u_i^{(d)}(z)$ of degree $d$ chosen 
sufficiently high to satisfy the following conditions: 
\begin{itemize}
\item every $u_i^{(d)}(z)$ is non-multiple;
\item the images $u_i^{(d)}(\Delta)$ are pairwise distinct holomorphic discs.
\end{itemize}

Then every $u_i^{(d)}(z)$ can be considered as an algebraic map $f_i$ from $\cp^1=
S^2$ to $\cp^2$. Making a generic perturbation of $f_i$ outside $B$, we obtain 
desired curve $\ti C^* \subset \cp^2$ as the union of the images $\ti f_i(S^2)$ of
the perturbed maps. Observe that $d$ appears as the degree of every component 
$\ti C^*_i \deff \ti f_i(S^2)$. 

\newlemma{lem6.2.4} There exist an almost complex structure $\ti J^*$ and points 
$x_\alpha$ on $\ti C^*$ satisfying the following conditions:
\begin{itemize}
\item[{\sl (a)}] the points $x_\alpha$ are pairwise distinct, and there 
are exactly $3d-1$ of them on every component $\ti C^*_i$;
\item[{\sl (b)}] $\ti J^*$ is $C^\ell$-smooth and $\omega$-tame, $\ti C^*$ is 
$\ti J^*$-holomorphic, and $\ti J^*$ coincides with $J^*$ on $B$;
\item[{\sl (c)}] any $\ti J^*$-holomorphic curve $C'$ which
  \begin{itemize}
  \item passes through the fixed points $x_\alpha$;
  \item is sufficiently close to $\ti C^*$ \wrt the cycle topology;
  \item has the same number of singular points as $\ti C^*$;
  \item[$*$] has a singular point $x'\in C'$ with the nodal number $\delta^*$
 at $x'$
  \end{itemize} 
  must coincide with $\ti C^*$.
\end{itemize}
\end{lem}

The last property asserts that every pseudoholomorphic curve $C'\not=\ti C^*$ with 
the properties {\sl(c)} except $(*)$ has simpler singularities than $\ti C^*$. So 
the induction assumption can be applied to such a $C'$. 

\proof We use the results of {\sl Sections \ref{sec:2}} and {\sl\ref{sec:4}}. 
Fix non-singular points $x_\alpha$ on $\ti C^*$ such that condition {\sl(a)} 
is fulfilled. Let $\mbfx_i$ be the $(3d-1)$-tuple of the points lying on the
component $\ti C^*_i$. Denote by $\scrm'$ the space of pairs $(C',J')$, 
where $J'\in \scrj$ and $C'$ is $J'$-holomorphic curve $C'$ satisfying properties 
{\sl(c)} except $(*)$. Then by the genus formula \eqqref(1.2.1) any such curve 
$C'$ has only rational irreducible components $C'_i$, the number of which is the
same as for $\ti C^*$, and the degree of 
every component $C'_i$ is $d$. This means that $\scrm'$ is the fiber product
of the spaces $\scrm(S^2, \cp^2, d, \mbfx_i)$ of rational pseudoholomorphic curves
of degree $d$ in $\cp^2$ passing through $\mbfx_i$. The product is taken
over the space $\scrj$ of almost complex structure in $\cp^2$. By the transversality
technique of \refsection{2}, the space $\scrm'$ a Banach manifold. To compute
the Fredholm index of the natural projection $\pi'_\scrj: \scrm' \to \scrj$ observe
that the expected dimension of rational $J$-holomorphic curves in $\cp^2$ of 
degree $d$ passing through $3d-1$ fixed distinct points is $0$. This implies that
the index of the projection $\pi'_\scrj: \scrm' \to \scrj$ is also $0$.

Further, by results of \refsection{4} the condition $(*)$ defines a proper 
$C^\ell$-smooth submanifold $\scrm^*$ in $\scrm'$ of finite codimension, say $m$. 
Consequently, the index of the corresponding projection $\pi^*_\scrj\deff \pi'
_\scrj\ogran_{\scrm^*}: \scrm^* \to \scrj$ is negative. Using the transversality 
technique of \refsection{2} we can construct a $C^\ell$-smooth submanifold $Y 
\subset \scrm^*$ of dimension $m$ such that 
\begin{itemize}
\item $(\ti C^*, \ti J^*)\in Y$ for some $\ti J^*$ obeying the condition 
{\sl(b)} of the lemma;
\item $Y$ is transversal to $\scrm^*$;
\item the restricted projection $\pi'_\scrj\ogran_Y: Y \to \scrj$ is an 
imbedding.
\end{itemize}
Then $(\ti C^*, \ti J^*)$ is an isolated point of the intersection $Y \cap 
\scrm^*$. But this means that $\ti J^*$ has the desired properties. \qed

\medskip
Below we shall need a property which is a bit sharper than {\sl(c)} in 
\lemma{lem6.2.4}. Roughly speaking, it claims that one can recover a 
pseudoholomorphic curve $C$ in $B$ knowing its part $\bigl(\barr B \bs B(\half)
\bigr) \cap C$. 

\newdefi{def6.2.2} Denote by $A$ the spherical annulus $\barr B \bs B(\half)$.
It is a closed subset of the closed unit ball $\barr B \subset \cp^2$. For
closed subsets $Y_1, Y_2 \subset \cp^2$ we denote by $\dist_A(Y_1, Y_2)$ the
Hausdorff distance between $Y_1 \cap A$ and $Y_2 \cap A$,
$$
\dist_A(Y_1, Y_2) \deff \dist(Y_1\cap A, Y_2\cap A),
$$
if both $Y_1 \cap A$ and $Y_2 \cap A$ are non-empty. The standard distance 
function in $\cp^2$ is used as the base. If exactly one of the set $Y_i \cap A$ 
is empty, we set $\dist_A(Y_1, Y_2) \deff \diam(\cp^2)$. If $Y_1 \cap A= Y_2 
\cap A =\emptyset$, we define $\dist_A(Y_1, Y_2) \deff 0$. We call $\dist_A$
the {\sl $A$-distance}. 
\end{defi}

It is easy to see that $\dist_A$ is only a pseudo-distance function, \ie it is 
non-negative, symmetric, and has the triangle inequality property, but does not 
distinguish all closed subsets $Y_1 \not= Y_2 \subset \cp^2$ in general. It turns 
out that it induces the cycle topology on the set of pseudoholomorphic curves lying
in a sufficiently small  $\dist_A$-neighborhood of $\ti C^*$ provided only $C^1
$-smooth almost complex structures $J$ are used. More precise statement is given in

\newlemma{lem6.2.4a} There exists an $\epsi>0$ with the following property.

Let $J \in \scrj$ be a $C^1$-smooth almost complex structure which satisfies the
 condition $\norm{J-J^*}_{C^0(\cp^2)} \le \epsi$ and $C$ a $J$-holomorphic curve 
which is homologous to $\ti C^*$ and satisfies the condition $\dist_A(C, \ti C^*)
\le \epsi$. Then for any sequence $J_n$ of continuous almost complex 
structures $J_n$ converging to $J$ in the $C^0$-topology, $\norm{J_n-J}_{C^0(\cp
^2)} \lrar 0$, and any sequence of $J_n$-holomorphic curves $C_n$ the condition 
$\dist_A(C_n, C) \lrar 0$ implies that $C_n$ converges to $C$ in the cycle 
topology.
\end{lem}

\proof Consider a sequence of almost complex structures $J_n$ in $\cp^2$ which
converges to $J^*$ in the $C^0$-topology, and a sequence $C_n$ of closed $J_n
$-holomorphic curves homologous to $C^*$, for which $\lim \dist_A(C_n, \ti C^*)=
0$. Then $J_n$ are $\omega\st$-tame for all $n\gg1$. Hence we can apply the Gromov
compactness theorem (see \refthm{thm5.2.1}). This means that some subsequence,
still denoted $C_n$, converges to a $J^*$-holomorphic curve $C^+$ \wrt the cycle
topology. The condition $\lim \dist_A(C_n, \ti C^*)= 0$ implies that $\dist_A(C^+,
\ti C^*)=0$, which means that $\ti C^* \cap A = C^+ \cap A$.

Observe now that by the construction of $\ti C^*$ every irreducible component $\ti 
C^*_i$ of $\ti C^*$ meets the interior $\sf{Int}(A)$ of $A$. By the unique
continuation property of pseudoholomorphic curves proven in \lemma{lem1.2.4} 
\sliip, every component $\ti C^*_i$ is contained in $C^+$. Thus $\ti C^* \subset 
C^+$. Since $C^+$ is homologous to $\ti C^*$, we must have equality $\ti C^*= C^+$.
This means that $C_n$ converges to $\ti C^*$ in the cycle topology. In particular,
for every sufficiently big $n$ every irreducible component of $C_n$ meets the 
interior $\sf{Int}(A)$ of $A$. 

The latter property shows that the same argumentation can be used if we replace
$\ti C^*$ by any $C_n$ with $n\gg1$ and the lemma follows. \qed

\medskip
Now we are ready to complete the

\statep Proof. of \refthm{thm6.2.2}. It follows from the construction
of the extension $\ti C^*$ that the sequence $(C_n, J_n)$ can be extended 
to a sequence $(\ti C_n, \ti J_n)$ such that $\ti J_n$ is a sequence of
$\omega$-tamed almost complex structures in $\cp^2$ converging to $\ti J^*$
and $\ti C_n$ is a sequence of compact (\ie closed) $\ti J_n$-holomorphic curves 
converging to $\ti C^*$. Moreover, we may additionally assume that the curves
$\ti C_n$ pass through the marked points $\mbfx$ for all sufficiently big $n$.

Observe that all $\ti C_n$ are symplectically isotopic. We denote by $\ti S$ 
the closed oriented real surface parameterizing $\ti C_n$. It can be obtained
from the surface $S$ parameterizing $C_n$ by gluing in discs to fill out the holes 
in $S$.

Fix a sequence of homotopies $\{\ti J_{n,t}\}_{t\in [0,1]}$ of almost complex 
structures with the following properties:
\begin{itemize}
\item all $\ti J_{n,t}$ are $C^\ell$-smooth and depend $C^\ell$-smoothly on $t$;
\item every initial structure $J_{n,0}$ is $\ti J_n$;
\item for some small $\epsi_0>0$ the structures $\ti J_{n,t}$ are integrable in 
$B$ for all $t\in [1-\epsi_0, 1]$;
\item as $n$ goes to infinity, the structures $\ti J_{n,t}$ converge to $\ti J^*$
in the $C^0$-topology uniformly in $t\in [0,1]$, \ie 
$$
\lim_{n\lrar\infty}\; \mathop{\sup}\limits_{t\in [0,1]}\;
\norm{\ti J_{n,t} - \ti J^*}_{C^0(\cp^2)} =0;
$$
\item the homotopy $\{\ti J_{n,t}\}_{t \in [0,1]}$ is generic for every $n$.
\end{itemize} 

\smallskip
Now let us try to deform continuously every $\ti C_n$ inside a family $\ti J_{n,
t}$-holomorphic curves preserving the isotopy class. Since we want to control also
the local isotopy class we must impose the condition that the curves in the family
lie sufficiently close to $\ti C^*$. Apriori, it can occur that such a curve does 
not exist for all $t\in [0,1]$. Nevertheless, we can find the maximal subinterval 
where such a family of curves exists. Moreover, we allow that under the deformation
some nodal points appear. Let us formalize this observation.

\newprop{prop6.2.6} Fix a sufficiently small $\epsi>0$. Then for every $n\gg1$ 
there exists a $t^+_n \in (0,1]$ which is {\sl maximal} \wrt the following 
condition: 

For any $t <t^+_n$ there exists a $\ti J_{n,t}$-holomorphic curve $\ti C
_{n,t}$ such that
\begin{itemize}
\item $\ti C_{n,t}$ passes through the fixed points $\mbfx$ on $\cp^2$;
\item the curve $\ti C'_{n,t}$, obtained from $\ti C_{n,t}$ by smoothing of all 
singular points contained in $B$, is symplectically isotopic to $\ti C_n$;
\item $\dist_A(\ti C_{n,t}, C^*) < \epsi$.
\end{itemize}
\end{prop}

Recall that for $t$ sufficiently close to $1$ the structures $\ti J_{n,t}$ are 
integrable in $B$. So if $t^+_n =1$ for some $n$, then for some $t$ close to $1$ 
we obtain a {\sl holomorphic} curve $C_{n,t} \deff \ti C_{n,t} \cap B$ whose
smoothing is symplectically isotopic to the original curve $C_n$. In this case 
\refthm{thm6.2.2} follows from \lemma{lem6.2.2a}. We claim that it is always 
the case for $n\gg1$.

\smallskip
To show this, let us analyze the possible reasons which could cause the strict
inequality $t^+_n<1$ for a given $n\gg1$. 
Consider
an increasing sequence of parameters $t_\nu$ approaching to $t^+_n$. Then there
exists a sequence of $\ti J_{n, t_\nu}$-holomorphic curves $\ti C_{n, t_\nu}
$ with the properties from \propo{prop6.2.6}. In particular, all $\ti C_{n, 
t_\nu}$ are homologous to $\ti C^*$. Taking a subsequence, we may assume that
$\ti C_{n, t_\nu}$ converges to a $\ti J_{n, t^+_n}$-holomorphic curve $\ti C^+_n$
in the cycle topology. Note $\dist_A(\ti C^+_n, \ti C^*)\le \epsi$ by our 
construction. By \lemma{lem6.2.4a}, $\ti C^+_n$ is sufficiently close to 
$\ti C^*$ also \wrt the cycle topology. Consequently, near every nodal point of 
$\ti C^*$ there is exactly one nodal point of $\ti C^+_n$. 

Observe that $\ti C^+_n$ has no singular point $x^+_n \in \ti C^+_n$ with the 
nodal number $\ge \delta^*$ at $x^+_n$. Indeed, otherwise we can repeat the 
argumentation from the proof of \lemma{lem6.2.4} and show that $\ti C^+_n$ must 
consist of rational components the number of which is the same as that for 
$\ti C^*$. But the expected 
dimension of such curves in the space $\scrm'_{\ti J_{n, t^+_n}, \mbfx}$ with a 
singular point of this type is negative and less then $-1$. So the existence of 
$x^+_n\in \ti C^+_n$ with $\delta(x^+_n, \ti C^+_n) \ge \delta^*$ contradicts the 
genericity of the path $\ti J_{n, t}$. Thus all singularities of $\ti C^+_n$ are 
simpler than those of $\ti C^*$. By the induction assumption, the curve $\ti C'_n$ 
obtained as the canonical smoothing of all singular points of $\ti C^+_n$ 
contained in $B$ is symplectically isotopic to $\ti C_n$.

Let $u^+_n: S^+_n \to \cp^2$ be a {\sl normal} parameterization of $\ti C^+_n$
(see \refdefi{def5.2c.3}).
Consider the relative moduli space $\scrm_{h_n, \mbfx}(S^+_n, \cp^2)$ of $h_n(t)= 
\ti J_{n,t}$-holomorphic curves which are parameterized by $S^+_n$, are in the 
homology class $[\ti C^*]$, and pass through the fixed points $\mbfx$. This space
is non-empty since it contains $(\ti C^+_n, t^+_n)$. \refthm{thm4.5.3} provides 
that for some interval $t\in [t^+_n, t^{++}_n]$ with $t^{++}_n >t^+_n$ we can 
construct a path of $\ti J_{n,t}$-holomorphic curves $\ti C_{n, t}$ which lies in 
$\scrm_{h_n, \mbfx}(S^+_n, \cp^2)$ and starts at $\ti C^+_n$. Then the curves 
obtained from such $\ti C_{n, t}$ by smoothing of all singular points contained 
in $B$ will be symplectically isotopic to $\ti C_n$. Moreover,
if we would additionally have the strict inequality $\dist_A(\ti C^+_n, \ti C^*)< 
\epsi$, then $\dist_A(\ti C_{n, t}, \ti C^*)< \epsi$ for some $t\in ]t^+_n, 
t^{++}_n[$, and this would contradict the maximality of $t^+_n$.

\smallskip
Thus we may assume that $\dist_A(\ti C^+_n, \ti C^*)= \epsi$ for every $n$. Then 
for every $n\gg1$ we can fix $t^-_n \in [0, t^+_n]$ and a $\ti J_{n,t^-_n}
$-holomorphic curve $\ti C_{n,t^-_n}$ which has properties from \propo{prop6.2.6} 
and satisfies the additional condition 
$$
\msmall{\epsi\over2} \le \dist_A(\ti C_{n, t^-_n}, \ti C^*) \le \epsi.
$$ 
Taking a subsequence, we may assume that $\ti C_{n,t^-_n}$ converges to a $\ti J
$-holomorphic curve $\ti C^+$ in the cycle topology. Then 
$$
\msmall{\epsi\over2} \le \dist_A(\ti C^+, \ti C^*) 
\le \epsi.
$$
By {\sl Lemmas \ref{lem6.2.4}} and {\sl\ref{lem6.2.4a}}, $\ti C^+$ must have 
simpler singularities than $\ti C^*$ provided the constant $\epsi$ was chosen 
small enough. By the induction assumption, the curve $\ti C'$ obtained by canonical 
smoothing of all singular points of $\ti C^+$ lying in $B$ is symplectically 
isotopic to every $\ti C_n$, as also to every $\ti C_{n,t^-_n}$. On the other 
hand, $\ti J^*$ coincide in $B$ with the standard structure $J\st$. Thus $C' 
\deff \ti C' \cap B$ is a canonical smoothing of $C^*$ by \lemma{lem6.2.2a}.

\newsubsection[6.3]{Global symplectic isotopy in $\cp^2$} In this paragraph give

\nobreak
\statep Proof. of \refthm{thm6.1.1}. We proceed by making appropriate 
modifications of the argumentation used in the proof of \refthm{thm6.2.2}. Let 
$\Sigma$ be an imbedded surface in $\cp^2$ of degree $d\le 6$, such that $\omega\st
\ogran_{\Sigma}$ is non-degenerate. By \propo{prop6.1.2}, to prove the theorem it 
is sufficient to show that $\Sigma$ is symplectically isotopic to a non-singular 
algebraic curve of degree $d$.

\smallskip 
Find an $\omega\st$-tame almost complex structure $J_0$ making $\Sigma$ a $J_0
$-holomorphic curve, denoted by $C_0$. Fix $3d-1$ distinct points $\mbfx=(x_1,
\ldots, x_{3d-1})$ on $C_0$. Perturbing $C_0$ and the points, we may assume that 
$x_1,\ldots, x_{3d-1}$ are in generic position \wrt the standard structure $J\st$ 
in the following sense. For any positive degree $d'\le d$ and any closed oriented 
surface $S$, not necessary connected, the moduli space $\scrm_{J\st, \mbfx}(S, 
\cp^2; d')$ of $J\st$-holomorphic (and hence {\sl algebraic}) curves of degree 
$d'$ with normalization $S$ passing through $\mbfx$ is a (possibly empty)
complex space of the expected dimension.

Fix a generic path $h(t)$ of $\omega\st$-tame almost complex structures $J_t \deff 
h(t)$ connecting $J_0$ with $J\st=J_1$. Without loss of generality we may assume
that all $J_t$ are $C^\ell$-smooth and depend $C^\ell$-smoothly on $t$. 

\newprop{prop6.3.1} There exists a $t^+ \in (0,1]$ which is {\sl maximal} \wrt 
the following condition: 

For any $t <t^+$ there exists a $J_t$-holomorphic curve $C_t$ such that 
\begin{itemize}
\item[\slip] $C_t$ passes through the fixed points $\mbfx$ on $\cp^2$;
\item[\sliip] $C_t$ is non-multiple, but not necessarily irreducible;
\item[\sliiip] the curve $C'_t$, obtained from $C_t$ by smoothing of all singular
points, is symplectically isotopic to $C_0$.
\end{itemize}
\end{prop}

To prove the theorem, we must show that $t^+=1$ and that there exist a $J\st
$-holomorphic curve $C_1$ with the properties \slip--\slii in the proposition. 

\smallskip
Let $t_n$ be an increasing sequence converging to $t^+$. Fix $J_{t_n}$-holomorphic
curves $C_n$ with these properties. Property \slii implies that the $C_n$ have
the same degree $d$. Going to a subsequence we may assume that they converge
to a $J_{t^+}$-holomorphic curve $C^+$ in the cycle topology. 

We claim that $C^+$ contains no multiple components. To show this, it is sufficient
to consider only the case when $C^+$ has only two components $C^+_1$ and $C^+_2$ 
with multiplicities $m_1=1$ and $m_2=2$ respectively. Let $d_i$ be the degree of 
$C^+_i$, so that $d_1+2d_2=d$. Then the geometric genus $g_i$ of every $C_i$ is at
most $g_i \le{(d_i-1)(d_i-2) \over 2}\cdot\;$ By the genericity of the path $h(t)=
J_t$, each $C_i$ can contain at most $k_i \le 3d_i -1 + g_i \le{d_i(d_i+3) \over 
2}$ of the fixed points $\mbfx$, see \refsubsection{2.4}. Thus $C^+$ can contain
at most $\le{d_1(d_1+3) + d_2(d_2+3) \over 2}$ points. It is easy to show that for
$d\le 6$ this number is strictly less then $3d-1$. For example, in the worst case 
with $d=6$, $d_1=4$, and $d_2=1$ we would have on $C^+$ at most ${4\cdot(4+3) 
\over 2} + {1\cdot (1+3)\over 2}= 14 +2 =16$ the marked points $\mbfx$ instead 
of the necessary $3\cdot 6-1 =17$. Observe, that this argument remains valid also 
in the case $t^+=1$ and $J_1=J\st$.

Now the results of \refsubsection{6.2} show that the curve $C'$ obtained 
from $C^+$ by the 
canonical smoothing of all singular points is symplectically isotopic to $C_0$.
This implies the theorem in the case $t^+=1$. Indeed, in this case $C^+$ is the 
zero divisor of a homogeneous polynomial $F^+$ of degree $d$. Making a generic
perturbation of coefficients of $F^+$ we obtain a polynomial $F'$ whose
zero divisor is an algebraic (and hence $J\st$-holomorphic) curve $C'$ which
is symplectically isotopic to $C_0$.

In the case $t^+ <1$ we must show that for some $t^{++}>t^+$ there exists a 
$J_{t^{++}}$-holomorphic curve $C^{++}$ with the properties given in 
\propo{prop6.3.1}. To do this we fix a normal parameterization $u^+: S^+ \to 
C^+ \subset \cp^2$ (see \refdefi{def5.2c.3}) and consider the relative moduli 
space $\scrm_{h, \mbfx}(S^+, \cp^2, d)$ of $J_t=h(t)$-holomorphic curves which 
are parameterized by $S^+$, 
pass through $\mbfx$ and have the degree $d$. This space is non-empty because it 
contains $C^+$. The results of \refsubsection{4.5} imply that for some interval 
$t\in [t^+_n, t^{++}_n]$ with $t^{++}_n >t^+_n$ we can construct a path of 
$J_t$-holomorphic curves $C^+_t$ which lies in $\scrm_{h, \mbfx}(S^+, \cp^2, d)$ 
and starts at $C^+$. By \refsubsection{6.2}, the $C^+_t$ have the the properties 
\slip--\sliii from \propo{prop6.3.1}. This contradicts the maximality of $t^+$ 
and implies the statement of \refthm{thm6.1.1}. \qed

\medskip
\state Remark. In fact, the real homotopy $C_t$ from $C_0=\Sigma$ to an algebraic
curve $C_1$ has the property described in \refsubsection{6.1}. Namely, after 
fixing a generic homotopy $h(t)=J_t$, one tries to construct {\sl any} path $C_t$
of {\sl imbedded} $J_t$-holomorphic curves $C_t$. Such a path exists for some
interval $t\in [0, t')$. The saddle point property proven in
\refsubsection{4.5} removes the main difficulty in the construction of $C_t$:
the presence of local maxima in the corresponding moduli space $\scrm$. This means
that at end of this interval, when $t \lrar t'$, the curves $C_t$ go to infinity in
$\scrm$. By Gromov compactness, going along some sequence $t'_n \lrar 
t'$, we approach a $J_{t'}$-holomorphic curve $C'$ lying on some infinity 
stratum of $\barr \scrm$ parameterized by a new moduli space $\scrm'$. As we have 
shown in the proof, one can avoid the strata $\scrm'$ corresponding to curves with
multiple components. Now we continue to deform $C'$ as a path $C'_t$ along 
$\scrm'$, having in mind that the canonical (in fact, any) smoothing of singular 
points of $C'_t$ gives curves symplectically isotopic to $C_0$. The new path $C'_t$
continues until we come to the next infinity stratum $\scrm''$ of $\barr \scrm$,
and so on.

\medskip
We finish the paper with a remark on {\sl Conjecture 2} about the local symplectic
isotopy problem for nodal curves. The proof of this result would follow from the 
corresponding result for {\sl holomorphic} curves, which is essentially a local 
version of the Severi problem, see \propo{prop6.1.2}. Indeed, the proof of 
\refthm{thm6.2.2} could be applied after appropriate modification. 

\smallskip
\state Conjecture 3. {\sl (Local Severi-Harris problem). \it Let $C^*$ be 
a holomorphic curve in the ball $B \subset \cc^2$ with a unique isolated singular 
point at $0\in B$ and without multiple components. 

Then any {\it nodal} holomorphic curve $C\subset B$ sufficiently close to $C^*$ 
\wrt the cycle topology is holomorphically isotopic to a holomorphic curve 
$C^\dag$ obtained from a maximal nodal deformation $C_\nod$ of $C^*$ by smoothing 
some number of nodes on $C_\nod$.}

\smallskip
In view of the main results of this paper, the validity of {\sl Conjecture 3}, 
and hence of {\sl Conjecture 2}, seems quite plausible.

\vskip1cm plus 1cm

\newpage

\ifx\undefined\bysame
\newcommand{\bysame}{\leavevmode\hbox to3em{\hrulefill}\,}
\fi


\newpage
\tableofcontents

\end{document}